\documentclass[11pt]{amsart}

\usepackage{amsmath}
\usepackage{amssymb}
\usepackage{amsfonts}
\usepackage{geometry}
\usepackage{url}
\usepackage{color}
\usepackage{blkarray}
\usepackage{comment}
\usepackage{blindtext,rotating}

\usepackage{amsmath}
\usepackage{amssymb}
\usepackage{amsfonts}
\usepackage{geometry}
\usepackage{url}
\usepackage{setspace}
\usepackage{amsthm}
\usepackage[all,cmtip]{xy}
\usepackage{amsmath,amscd}
\usepackage{fourier-orns}

\usepackage{amsthm}
\theoremstyle{definition}
\newtheorem{theorem}{Theorem}[section]
\newtheorem{lemma}[theorem]{Lemma}
\newtheorem{proposition}[theorem]{Proposition}
\newtheorem{corollary}[theorem]{Corollary}

\newtheorem{conjecture}[theorem]{Conjecture}
\newtheorem{remark}[theorem]{Remark}

\newtheorem{assumption}[theorem]{Assumption}

\newcommand{\hstar}{\mathfrak{h}^*}
\newcommand{\oh}{\mathcal{O}}
\newcommand{\el}{\mathrm{L}}
\newcommand{\m}{\mathrm{M}}
\newcommand{\mbf}{\mathbf}
\newcommand{\s}{\mathfrak{sl}_2}
\newcommand{\B}{\mathcal{B}}
\newcommand{\G}{\mathrm{G}}

\newcommand{\Ind}{\operatorname{Ind}}

\newcommand{\Hom}{\operatorname{Hom}}

\newcommand{\Res}{\operatorname{Res}}
\newcommand{\Supp}{\operatorname{Supp}}
\newcommand{\Irr}{\operatorname{Irr}}
\newcommand{\rk}{\operatorname{rk}}
\newcommand{\Ext}{\operatorname{Ext}}
\newcommand{\Rad}{\operatorname{Rad}}
\newcommand{\Ann}{\operatorname{Ann}}
\newcommand{\Spec}{\operatorname{Spec}}
\newcommand{\KZ}{\operatorname{KZ}}

\begin{document}

\title{Irreducible representations of rational Cherednik algebras for exceptional Coxeter groups, Part I}
\maketitle

\begin{abstract} This paper and its sequel describe the irreducible representations of the rational Cherednik algebra $H_c(W)$ for a finite Coxeter group $W$ of type $H_4$, $F_4$ with equal parameters, $E_6$, $E_7$, and $E_8$, when $c$ is not a half-integer. Herein appear the decomposition matrices of Category $\oh_c(W)$ for $W=H_4$, $E_6$, and $E_7$, as well as a classification of finite-dimensional representations of $H_c(W)$ for these groups and for $H_c(F_4)$ with equal parameters.

%decomposition matrices of the blocks of Category $\oh_c(W)$ for $W$ of type $H_4$, $E_6$, and $E_7$. The character of any irreducible representation of $H_c(W)$ , as well as the classification of irreducible representations of any given dimension of support. Included is a classification of all finite-dimensional representations of $H_c(W)$. Additionally, a calculation of the finite-dimensional representations of $H_c(F_4)$ with equal parameters appears in this paper.
\end{abstract}

\section{ Introduction}

Rational Cherednik algebras exist for any complex reflection group $W$ but only the foundations of their representation theory are in place outside of the cyclotomic types ($W=G(\ell,1,n)$, of which types $A_n$ and $B_n$ are examples). Important elementary questions remain open in most types, for instance: 
\begin{itemize}
\item Which representations are finite-dimensional and what are their dimensions? 
\item What are the characters of the irreducible representations? 
\item What are the multiplicities of simple modules in standard (``Verma") modules? 
\end{itemize}
This paper and its sequel will provide the answers to these questions for $W$ one of the following exceptional Coxeter groups: $H_4$, $F_4$ with equal parameters, $E_6$, $E_7$, and $E_8$. 

The main result of this work is the computation of numerical data (decomposition matrices) associated to the rational Cherednik algebra $H_c(W)$ for all parameters $c=1/d$ with $d$ a divisor of one of the fundamental degrees of $W$ and $d\neq2$: calculating the decomposition matrices of the blocks of the Category $\oh$ of $H_c(W)$ gives answers to all three of the questions above. Up to an equivalence of categories this covers all cases where the Category $\oh$ of $H_c(W)$ is not semisimple, except for when $c$ is a half-integer. Finding decomposition numbers should be of interest in its own right. However, it turns out the decomposition matrices reveal more than expected. The decomposition matrix of $\oh_c(W)$ encodes the multiplicities $[\m(\tau):\el(\sigma)]$ of simples $\el(\sigma)$ appearing in the composition series of Vermas $\m(\tau)$ for a fixed parameter $c=1/d$. These multiplicities, which are nonnegative integers, turn out to be closely related to the decomposition numbers of the unipotent blocks of a finite group of Lie type over a field of characteristic $\ell$. If $G(q)$ is a finite group of Lie type with Weyl group $W$, then its representation theory in characteristic $\ell$ when $\ell$ divides the order of $G(q)$ and $\ell\neq p$ is, like the representation theory of $H_c(W)$, sensitive to the divisors $d$ of the fundamental degrees of $W$: what matters essentially is which cyclotomic polynomial $\Phi_d(q)$ is divisible by $\ell$, provided that $\ell$ divides exactly one $\Phi_d(q)$ and $\ell$ is large enough. It turns out that in all the examples in this paper where the decomposition matrices of the corresponding unipotent $G(q)$ blocks have been computed, the decomposition matrix of $H_c(W)$ embeds as a square submatrix in the decomposition matrix of $G(q)$, where $\ell$ divides $\Phi_d(q)$ and $c=1/d$. In particular, this is true (obviously) for all blocks of ``defect 1" but in particular for all the blocks for $E_6$, $E_7$, and $E_8$ at $c=1/4$ corresponding to those at $d=4$ computed recently by Dudas and Malle \cite{DM}.

The decomposition matrices thus make it clear that the representation theory of the rational Cherednik algebra for Weyl groups $W$ is closely related to the modular representation theory (in non-defining characteristic) of a finite group of Lie type with Weyl group $W$. The $E_n$ examples for $d=4$, along with the fact that in types $A_n$ and $B_n$, thanks to the $q$-Schur algebra which controls the numerics of both objects \cite{R}, the decomposition matrices for rational Cherednik algebras embed in those of the appropriate finite groups of Lie type, led the author to hope that the decomposition matrix of the rational Cherednik algebra should always embed in that of the associated finite group of Lie type for the appropriate $d$. Unfortunately this is not the case. In the sequel to this paper, an example in type $D_n$ will be explored where the decomposition matrix of the rational Cherednik algebra at $c=1/4$ is nearly identical to that of the finite special orthogonal group over a field of characteristic $\ell$ dividing $\Phi_4(q)$, but where there is one entry which is different between the two matrices. This could be explained by the fact that the generalized $q$-Schur algebra for type $D_n$ defined by Dipper-James is not a $1$-faithful quasi-hereditary cover. 

Based on the data in this paper and its sequel, I would like to make the following conjecture. %, which, in the spirit of the etymology of the word, is just an idea I'd like to throw out there. 
Let $d$ divide one of the fundamental degrees of $W$ and let $\ell | \Phi_d(q)$ be a prime. Let $\mathcal{B}$ be a block of $\oh_c(W)$ and take any $\tau\in\Irr W$ such that $\el(\tau)\in\mathcal{B}$. Let $G(q)$ be a finite group of Lie type with Weyl group $W$ and let $\mathcal{\tilde{B}}'$ be the unipotent block of $G(q)$ over a field of characteristic $\ell>>0$ such that $\mathcal{\tilde{B}}'$ contains the representation labeled by $\tau\otimes\mathrm{sgn}$. 

\conjecture For $\tau,\;\sigma$ lowest weights of indecomposable objects in $\mathcal{B}$, the rational Cherednik algebra decomposition number $[\m(\tau):\el(\sigma)]$ is a lower bound for the $G(q)$ decomposition number in row: $\tau\otimes\mathrm{sgn}$ and column: $\sigma\otimes\mathrm{sgn}$ of the unipotent block $\mathcal{\tilde{B}}'$.\\

%The decomposition matrix of the rational Cherednik algebra is upper unitriangular by convention; the decomposition matrix of the unipotent block of $G(q)$ is lower unitriangular by convention; the conjecture says that if the $G(q)$ decomposition matrix is rotated $180$ degrees so that it is upper unitriangular, the columns and rows corresponding to cuspidals are removed, and the labels are tensored by sign, then one gets the decomposition matrix of the rational Cherednik algebra either verbatim or after stripping away some extra multiplicities (mostly in the region of the upper left hand corner). 

Note that the square submatrix of $G(q)$ which is concerned is the one labeled by unipotent characters lying in the principal series. Also, note that the conjecture holds for those columns corresponding to simple representations of the Hecke algebra at a $d$'th root of unity, since this rectangular matrix embeds in both the rational Cherednik algebra decomposition matrix by \cite{CGG}, \cite{Chl} and the $G(q)$ decomposition matrix for $\ell>>0$ by James' Conjecture. One might imagine that what happens exactly is that the columns of the $G(q)$ decomposition matrix are nonnegative linear combinations of the columns of the $H_c(W)$ decomposition matrix.

Beyond the decomposition numbers themselves, which are an important invariant of the rational Cherednik algebra, the decomposition matrix also gives the characters of all irreducible representations of $H_c(W)$. The inverse of the decomposition matrix of multiplicities $[\m(\tau):\el(\sigma)]$ is the matrix of multiplicities $[\el(\tau):\m(\sigma)]$ which give the ``Verma-flag" of each simple representation. The character of $\el(\tau)$ may be read off the row labeled by $\tau$ in this matrix. It is a rational function in $t$ the order of whose pole at $t=1$ tells the dimension of support of $\el(\tau)$. The irreducible representations with support of dimension less than the rank of $W$ are of particular interest; these are the representations killed by the KZ functor. Among them, the representations whose support has dimension $0$ consist of the finite-dimensional representations. According to the philosophy of Miyachi, these should correspond to cuspidal representations for $G(q)$ (indeed, more generally, the minimal support representations in a block for the rational Cherednik algebra should correspond tos cuspidals for $G(q)$). Evaluating the character of a finite-dimensional representation at $t=1$ gives the dimension of the representation. Both the characters of finite-dimensional representations, which are polynomials in $\mathbb{Z}[t^{-1},t]$ symmetric under interchanging $t$ and $t^{-1}$, and the dimensions of finite-dimensional representations have meanings in other contexts. The characters should be Markov traces of braids of type $W$ \cite{GORS}, as Vivek Shende has explained to the author; and at least in some cases, the dimensions of finite-dimensional representations are dimensions of geometrically defined subspaces in the cohomology of the affine Springer fiber \cite{OY}, \cite{VV}.

The study of rational Cherednik algebras for exceptional Coxeter groups was begun by Chmutova for dihedral groups \cite{Chm} and Balagovic-Puranik for $H_3$ \cite{BP}. This paper treats $H_4$, $E_6$, $E_7$, and the finite-dimensional representations of $F_4$ with equal parameters. The sequel, which is mostly complete but still in progress, deals with $E_8$, decomposition matrices of $F_4$ with equal parameters, and some type $D_n$ examples for the purpose of comparing with the results of \cite{DM}. The outline of this paper is as follows: a summary of the finite-dimensional representations and their dimensions is given in the table below. The next section provides background and lemmas which will be needed for the computations. The remaining sections compute the decomposition matrices at all relevant parameters $c$ except for $c=1/2$.

\textbf{Acknowledgements.} I would like to thank Olivier Dudas, Pavel Etingof, and Alexei Oblomkov for enlightening conversations during different stages of this project. Thanks to Alexei Oblomkov for suggesting the problem to study the rational Cherednik algebras of $H_4$ and $F_4$.

%It is sweet to look up at the light.

\section{Finite-dimensional irreducible representations and their dimensions}

Let the parameter $c$ for $H_c(W)$ be $1/d$ for some positive integer $d$ dividing a fundamental degree of $W$. Then the following is a complete list of irreducible finite-dimensional representations  of $H_c(W)$ and their dimensions, for $W$ equal to $F_4$ (with equal parameters), $E_6$, $E_7$, and $H_4$, with two exceptions: (1) $d=2$, results are conjectural. These are marked with an asterisk. (2) In the case of $E_7$, it is known that the spherical representation $\el(1_a)$ is finite-dimensional when $d=2$ by \cite{VV}. This case is not studied in this paper and the finite-dimensional representations of $H_{\frac{1}{2}}(E_7)$ remain unknown.

The reason for only listing dimensions for $c=1/d$ is that for $c=r/d$, $r$ a positive integer coprime to $d$, there is the following dimension formula due to Rouquier \cite{R}: $$\dim\el_{\frac{r}{d}}(\tau)=r^n\dim\el_{\frac{1}{d}}(\tau)$$

In the case of $H_4$ where the characters of the group are defined not over $\mathbb{Q}$ but over $\mathbb{Q}(\sqrt{5})$, replacing $1/d$ with $r/d$ may have the effect that some lowest weights $\tau$ will have to be replaced with their Galois-conjugate. For finite-dimensional representations, this means that for some congruences of $r$ modulo $d$, $\mbf{3}$ will be interchanged with $\mbf{5}$, or $\mbf{11}$ with $\mbf{13}$, as lowest weights of the finite-dimensional representations. However, the dimensions themselves will not change.

\[
\begin{array}{ccccc} \begin{array}{c|c|c}
\underline{\mbf{F_4}}\\
\\
d&\el(\tau)&\dim\el(\tau)\\
\hline
12&\el(1_1)&1\\
\hline
8&\el(1_1)&6\\
\hline
6&\el(1_1)&20\\
&\el(2_1)&2\\
&\el(2_3)&2\\
\hline
4&\el(1_1)&96\\
&\el(4_1)&15\\
\hline
3&\el(1_1)&256\\
&\el(4_1)&64\\
\hline
2&\el(1_1)&1620\\
&\el(2_1)*&78\\
&\el(2_3)*&78\\
&\el(9_1)*&84\\
\end{array}

&

\hspace{7mm}

&

\begin{array}{c|c|c}
\underline{\mbf{E_6}}\\
\\
d&\el(\tau)&\dim\el(\tau)\\
\hline
12&\el(1_p)&1\\
\hline
9&\el(1_p)&8\\
\hline
6&\el(1_p)&92\\
&\el(6_p)&28\\
\hline
3&\el(1_p)&4152\\
&\el(6_p)&1680\\
&\el(15_p)&56\\
\end{array}
&

\hspace{7mm}

&
\begin{array}{c|c|c}
\underline{\mbf{E_7}}\\
\\
d&\el(\tau)&\dim\el(\tau)\\
\hline
18&\el(1_a)&1\\
\hline
14&\el(1_a)&9\\
\hline
10&\el(7_a')&36\\
\hline
6&\el(1_a)&3894\\
&\el(7_a')&1806\\
&\el(21_a)&84\\
&\el(15_a')&15\\
\hline
2&\el(1_a)&?\\
&?&?\\
&...&...\\
&?&?\\
\end{array}
\end{array}
\]

\vspace*{5mm}

\[
\begin{array}{ccc}
\begin{array}{c|c|c}
\underline{\mbf{H_4}}\\
\\
d&\el(\tau)&\dim\el(\tau)\\
\hline
30&\el(\mbf{1})&1\\
\hline
20&\el(\mbf{1})&6\\
\hline
15&\el(\mbf{1})&20\\
&\el(\mbf{5})&4\\
\hline
12&\el(\mbf{1})&50\\
\hline
10&\el(\mbf{1})&105\\
&\el(\mbf{5})&24\\
&\el(\mbf{13})&9\\
&\el(\mbf{3})&15\\
\hline
6&\el(\mbf{1})&800\\
&\el(\mbf{3})&175\\
&\el(\mbf{5})&175\\
\hline
5&\el(\mbf{1})&1620\\
&\el(\mbf{11})&84\\
&\el(\mbf{3})&384\\
&\el(\mbf{20})&60\\
\end{array}

&
\hspace*{7mm}
&

\begin{array}{c|c|c}
d&\el(\tau)&\dim\el(\tau)\\
\hline
4&\el(\mbf{1})&3450\\
&\el(\mbf{11})&300\\
&\el(\mbf{13})&300\\
\hline
3&\el(\mbf{1})&12800\\
&\el(\mbf{18})&300\\
&\el(\mbf{3})&2500\\
&\el(\mbf{5})&2500\\
\hline
2&\el(\mbf{1})&65625*\\
&\el(\mbf{3})*&8550\\
&\el(\mbf{5})*&8550\\
&\el(\mbf{11})*&5625\\
&\el(\mbf{13})*&5625\\
&\el(\mbf{27})*&825\\
\end{array}
\end{array}
\]

\section{Preliminaries}

\subsection{Background and notation}

Rational Cherednik algebras were defined by Etingof and Ginzburg in \cite{EG}.

The rational Cherednik algebra may be referred to in the text as an ``RCA" and irreducible representations as ``irreps."

Let $W$ be a finite Coxeter group. Let $\tau\in\Irr W.$ Let $\tau':=\tau\otimes\mathrm{sign}$ for any $\tau\in\Irr W$. Let $S\subset W$ be the set of all reflections in $W$ (not to be confused with the subset of simple reflections -- so for example, if $W=A_{n-1}=S_n$ then $S$ consists of all transpositions, not just adjacent transpositions). $\hstar$ denotes the standard representation, or reflection representation, of $W$, and it has dimension equal to the minimal number of generators of $W$, called the rank of $W$.

To $W$ one may associate a quadratic algebra called the rational Cherednik algebra, denoted $H_c(W)$, whose degree $0$ part is a semisimple algebra over $\mathbb{C}$. Fix $c\in\mathbb{C}$. As a $\mathbb{C}$-vector space, $H_c(W)\cong\mathbb{C}[\hstar]\otimes \mathbb{C}[W]\otimes\mathbb{C}[\mathfrak{h}]$ (this statement is known as the PBW theorem, by analogy with the similar statement of the triangular decomposition of the universal enveloping algebra of a complex semisimple Lie algebra). $H_c(W)$ is generated by $\mathfrak{h}$, $W$, and $\hstar$.  Elements of $\hstar$ will be denoted by $x$'s and elements of $\mathfrak{h}$ by $y$'s. As an algebra, $H_c(W)$ has the following description: all $x$'s commute with each other; all $y$'s commute with each other; there are relations between elements of $W$ and the two polynomial algebras in the usual way: moving an element of $W$ past $x$'s or $y$'s comes at the cost of the action of $W$ on the $x$'s and $y$'s. Thus $\mathbb{C}[\hstar]\rtimes W$ and $\mathbb{C}[\mathfrak{h}]\rtimes W$ are both subalgebras of $H_c(W)$ (note that $\hstar\subset\mathbb{C}[\mathfrak{h}]$ and $\mathfrak{h}\subset\mathbb{C}[\hstar]$). Finally, the most interesting relation is the quadratic relation between $x$'s and $y$'s which lands in the center of $\mathbb{C}[W]$, and this is where the parameter $c$ comes into play:
$$[y,x]=(y,x)-c\sum_{s\in S} (\alpha_s,x)(y,\alpha^\vee_s) s$$
Here, $(-,-)$ is the natural pairing between $\mathfrak{h}$ and $\hstar$, and $\alpha_s,$ $\alpha_s^\vee$ are eigenvectors for $s$ satisfying $(\alpha_s,\alpha_s^\vee)=2$.

$\oh_c(W)$ is the Category $\oh$ of representations of $H_c(W)$ which are finitely generated over $\mathbb{C}[\mathfrak{h}]$ and locally nilpotent for the action of $\mathfrak{h}$. It is a highest weight category in the sense of \cite{CPS} with irreducible objects labeled by $\Irr W$. However, instead of highest weights, the convention here will be of lowest weights. There is a poset structure on $\Irr W$ defined as follows. Let $n=\rk W$ and let $\{x_i\}$ and $\{y_i\}$, $i=1,...,n$, be dual orthonomal bases of $\hstar$ and $\mathfrak{h}$ respectively. $H_c(W)$ contains a canonical $\s$-triple $(E,H,F)$ with $E=\frac{1}{2}\sum\limits_{i=1}^n x_i^2$, $H=\sum\limits_{i=1}^n x_iy_i+y_ix_i$, and $F=-\frac{1}{2}\sum\limits_{i=1}^n y_i^2$. This $\s$-triple commutes with the action of $W$ on any representation in $\oh_c(W)$ and $H$ acts on the lowest weight $\tau$ of any indecomposable object in $\oh_c(W)$ by a scalar $h_c(\tau)$. The formula for $h_c(\tau)$ is:
$$h_c(\tau)=\frac{\dim\hstar}{2}-c\frac{\sum_{s\in S}\chi_{\tau}(s)}{\dim\tau}$$
In the case that $W$ has only one root length, this formula becomes 
$$h_c(\tau)=\frac{\dim\hstar}{2}-c|S|\frac{\chi_\tau(s)}{\dim\tau}$$
where $s$ is any reflection. One may then use $h_c(\tau)$ to put a partial order on $\Irr W$ via \cite{GGOR} :
$$\sigma>_c\tau\textrm{ if and only if }h_c(\sigma)-h_c(\tau)\in\mathbb{Z}_{>0}$$
In addition to bestowing a partial order on $\Irr W$, the $\s$-triple induces a natural $\mathbb{Z}$-grading on $H_c(W)$ and on the modules in $\oh_c(W)$ by assigning $W$ to degree $0$, $\hstar$ to degree $1$, and $\mathfrak{h}$ to degree $-1$. Imitating the style of \cite{BP}, we will visualize the values of $h_c(\tau)$ on $\tau\in\Irr W$ using ``weight lines," putting $\tau$ at the point $h_c(\tau)$ on the real number line; reps $\sigma$ and $\tau$ will be graphed on the same line only if $h_c(\tau)$ and $h_c(\sigma)$ are integer distances apart. This provides a helpful tool for seeing which Vermas can appear in the decomposition of a given simple and vice versa, and for quickly calculating the characters once the Verma-decompositions of simples are found.

$\m(\tau)$ denotes the Verma module with ``lowest weight" $\tau$:
$$\m(\tau):=H_c(W)\otimes_{\mathbb{C}[\hstar]\rtimes W}\tau$$
In this construction, $\tau$ is considered as a representation of $\mathbb{C}[\hstar]\rtimes W$ by letting $W$ act as it acts on $\tau$, and extending this to an action of $\mathbb{C}[\hstar]\rtimes W$ by having $\mathfrak{h}$ act by $0$. This representation $\tau$ is then induced up to the whole algebra $H_c(W)$. $\m(\tau)$ is $\mathbb{N}$-graded with a shift and generated in degree $h_c(\tau)$, and as a $W$-representation it is the direct sum of its graded pieces which are the tensor product of $\tau$ with symmetric powers of the standard representation:
$$\m(\tau)=\sum_{k=0}^\infty S^k\hstar\otimes\tau$$

$\el(\tau)$ denotes the unique simple quotient of $\m(\tau)$: $\el(\tau)=\m(\tau)/J(\tau)$ where $J(\tau)=\Rad\m(\tau)$. $J(\tau)$ has an interpretation as the kernel of the Shapovalov form, a $W$-invariant bilinear form on $\m(\tau)$, but this perspective will not enter into this paper. $\el(\tau)$ inherits a $\mathbb{N}$-grading with a shift from $\m(\tau)$ and is likewise has its lowest degree graded piece, equal to $\tau$, in degree $h_c(\tau)$.

The category $\oh_c(W)$ is finite-length, and any $H_c(W)$-module $M\in\oh_c(W)$ has a filtration whose successive quotients are simple modules. Then the \textit{decomposition number} $[\m(\tau):\el(\sigma)]$ denotes the multiplicity of $\el(\sigma)$ as a composition factor of $\m(\tau)$, a nonnegative integer. For any $\sigma,\;\tau\in\Irr W$, $[\m(\tau):\el(\tau)]=1$, and $[\m(\tau):\el(\sigma)]\neq0$ for $\sigma\neq\tau$ implies that $h_c(\sigma)-h_c(\tau)\in\mathbb{Z}_{>0}$. The \textit{decomposition matrix} is the matrix of multiplicities $[\m(\tau):\el(\sigma)]$, $\tau,\;\sigma\in\Irr W$. It is block diagonal so the different blocks are considered separately; these are the ``blocks" of the category $\oh_c(W)$. It is upper unitriangular if the rows and columns are ordered so that $h_c$ is non-decreasing down and to the right. Next to some of the rows is a star or a number in parentheses in typewriter font. A star next to row $\tau$ indicates that $\el(\tau)$ is finite-dimensional. A label of $\mathtt{(n)}$ next to row $\tau$ indicates that $\dim\Supp\el(\tau)=n$. In the course of working out the decomposition matrices, some versions of the matrix in different stages of completeness may appear in the text with bullets next to some rows -- the bullet by a row indicates that the irrep with lowest weight the character labeling that row has less than full support. These correspond exactly to the columns deleted by KZ functor.

A Grothendieck-group expression for the simple representation $\el(\tau)$ in terms of Verma modules can expressed by starting from $\m(\tau)$ and subtracting off some Verma modules, and then subtracting off some Verma modules from those, and so on, until one has subtracted off the submodule $J(\tau)$ expressed as a $\mathbb{Z}$-linear combination of Verma modules. This Verma-decomposition of $\el(\tau)$ is given by row $\tau$ of the inverse of the decomposition matrix, and the $\tau,\sigma$'th entry of the inverse of the decomposition matrix is $[\el(\tau):\m(\sigma)]$, the ``multiplicity" of $\m(\sigma)$ in $\el(\tau)$, an integer (which may be negative).

It suffices to consider $\oh_c(W)$ for only those parameters $c=1/d$ where $d$ divides one of the fundamental degrees of $W$ (i.e. the degrees of the generators of $\mathbb{C}[\hstar]^W$): first, $\oh_c(W)$ is semisimple unless $c=r/d$ with $d$ as above and $r$ coprime to $d$; second, equivalences of categories allow one to reduce to the case that the numerator is $1$ \cite{GGOR}.

As in any highest-weight category, there are no self-extensions between simple modules \cite{CPS}, \cite{BEG}:
$$\Ext(\el(\tau),\el(\tau))=0$$
and BGG reciprocity holds:
$$[\mathrm{P}(\sigma):\m(\tau)]=[\m(\tau):\el(\sigma)]$$

The graded character of a module $M$ of lowest weight $\tau$ is a Laurent series in $t$ times some possibly fractional power of $t$ (corresponding to the shift in the grading given by $h_c(\tau)$): $$\chi_M(t)=t^{h_c(\tau)}\sum_{k=0}^\infty \dim M[k+h_c(\tau)]t^k$$
For $M=\m(\tau)$ a Verma module, the formula reads $$\chi_{\m(\tau)}(t)=t^{h_c(\tau)}\dim\tau\sum_{k=0}^{\infty}{n-1+k\choose n-1}t^k$$ where $n=\rk W$. For $\el(\tau)$ a simple module, its character is given by the formula:
$$\chi_{\el(\tau)}(t)=\frac{\sum\limits_{\sigma\in\Irr W}\left([\el(\tau):\m(\sigma)]\dim\sigma\right) t^{h_c(\sigma)}}{(1-t)^n}$$
An important geometric object associated to a simple representation is its support, an algebraic subvariety of $\hstar$. The support of $\el(\tau)$ is defined as $$\Supp\el(\tau):=\Spec\left(\mathbb{C}[\hstar]/\Ann\el(\tau)\right)$$
On the one extreme, $\dim\Supp\el(\tau)=n$ if and only if $\Supp\el(\tau)=\hstar$, in which case we say that $\el(\tau)$ has ``full support;" on the other extreme, $\dim\Supp\el(\tau)=0$ if and only if $\Supp\el(\tau)=\{0\}$, in which case $\el(\tau)$ is finite-dimensional. $\dim\Supp\el(\tau)$ coincides with the power of $(1-t)$ appearing in the denominator of $\chi_{\el(\tau)}(t)$ when the character is written in lowest terms.

The modules whose support is neither the origin nor all of $\hstar$ have an interpretation in terms of finite-dimensional modules thanks to \cite{BE} -- they all arise via induction of finite-dimensional modules in the rational Cherednik algebras of parabolics $W/\subset W$ at the same parameter $c$. Induction and restriction between $\oh_c(W')$ and $\oh_c(W)$ was defined by \cite{BE}. Induction from a maximal parabolic $W'\subset W$ raises the dimension of support of an RCA module by $1$ while restriction to a maximal parabolic lowers dimension of support by \textit{at least 1}. In particular, a module of support bigger than the origin may be sent to $0$ upon restriction to some maximal parabolic; however, the finite-dimensional modules are exactly those which are sent to $0$ by restriction to every maximal parabolic. 

On the level of Grothendieck groups, the rules for induction and restriction are the same as those for the underlying groups $W$, $W'$ (Proposition 3.14, \cite{BE}). This was put to good use by \cite{BP} when they found decompositions of simples and Vermas in $\oh_{\frac{1}{2}}(H_3)$, and it will be much used in this paper in a similar spirit; in particular, the fact that induction and restriction take modules to modules. Namely, by finding an irreducible representation $\el(\tau)$ in $\oh_c(W')$ for $W'\subset W$ one may then produce the Grothendieck-group expression of a representation $\Ind_{W'}^W\el(\tau)$ and use this to get information about the Verma-decompositions of the irreducible representations appearing in it or conversely about the simple factors which are in its composition series. Often this will give a lower bound on certain decomposition numbers. Conversely, restriction also contains important information, since if a candidate expression in terms of Vermas for an irrep produces a virtual module upon restriction to a parabolic, then it cannot be correct. This will often be used to get upper bounds on decomposition numbers. Another trick is to induce and then restrict back down, or restrict then induce -- for example to get information about the number of copies of particular simples in the composition series of an induced module, in the first case; or to produce something in the parabolic that can then be induced back up if we don't know its rational Cherednik algebra well enough, in the second case. Playing around with induction and restriction also gives important information about the dimensions of supports of the various irreps in the block of $\oh_c(W)$ under consideration.

Crucial to the strategy of this paper is the structure of $\s$-module on any irreducible representation via the $\s$-triple $(E,H,F)$. This observation was important for the computations in \cite{BP} and will be critical in the analysis of modules in this paper as well. Since a finite-dimensional $\s$-representation must have integer weights with lowest weight in $\mathbb{Z}_{\leq0}$, the same is true for a finite-dimensional $H_c(W)$-representation: if $\el(\tau)$ is finite-dimensional then $h_c(\tau)\in\mathbb{Z}_{\leq0}$.

Also of critical importance is the faithful, exact functor $\KZ:\oh_c(W)\rightarrow\mathcal{H}_q(W)-mod$, which takes modules in Category $\oh_c(W)$ to modules over the Hecke algebra at a $q$th root of $1$ \cite{GGOR}. Here $q$ should be $e^{-2\pi i c}$. By the Double Centralizer Theorem  \cite{GGOR}, the characters labeling a block of the rational Cherednik algebra correspond to the characters labeling a block of the Hecke algebra at $q=e^{2\pi i c}$ after tensoring them with the sign character. But more than that is true: in fact, the decomposition matrix of a block of the Hecke algebra at $q=e^{2\pi ic}$ embeds into the block of the rational Cherednik algebra containing those characters tensored by sign. The effect of $\KZ$ on the decomposition matrix of $\B\subset\oh_c(W)$ is to delete the columns corresponding to the irreps of less than full support, and the result is the decomposition matrix of $\mathcal{H}_q(W)$ with respect to the canonical basic sets (\cite{Chl} Prop. 5.7 and Prop. 5.12; \cite{CGG}). These decomposition matrices of $\mathcal{H}_q$ have been calculated for all exceptional Coxeter groups and are collected into a series of tables at the end of \cite{GJ}. Thus the starting point in this paper for finding decomposition matrices for $H_c(W)$ will be to start from the decomposition matrices of $\mathcal{H}_q(W)$ at $q=e^{2\pi i c}$ as printed in \cite{GJ}, tensor the characters labeling rows and columns by the sign rep, embed this rectangular matrix into the square matrix of the corresponding block of $H_c(W)$, and then use the lemmas in the next section together with induction and restriction to recover the entries deleted by $\KZ$.

The ``defect" of a block is a nonnegative integer measuring, approximately, the complexity of the block (as the author learned from asking a question about this on MathOverflow). Let's say a block $\B$ of the rational Cherednik algebra has defect $n$ if $\KZ(\B)$ has defect $n$ as a block of the Hecke algebra. Blocks of defect $0$ are singletons -- a Verma which is simple and which has no nonzero homomorphism to any other Verma. Blocks of defect $1$ are those in which there is a unique irreducible of less than full support, and it is positioned at the leftmost end of the weight line for its block, and the decomposition matrix of the block has $1$'s on the diagonal and $1$'s just above the diagonal, and $0$'s everywhere else. The $h_c$-weight line coincides with a straight line Brauer tree in this case and embeds in the Brauer tree for the corresponding finite group of Lie type at the appropriate $d$ ($c=1/d$). The Verma-decomposition of any irreducible in such a block is simply the alternating sum of the Vermas to its right on the weight line, taken in order of increasing $h_c$-weight. This formula was proved by Rouquier in \cite{R}. One may obtain the list of all defect $n$ blocks of the Hecke algebra from the tables in \cite{GP}, and from there, read off the decomposition of the irreducibles in the blocks of defect $1$. The blocks of defect $1$ will still be listed in this paper, for the sake of completeness -- and the dimension of support of the irreducible of less than full support in each such block will be calculated and printed.

It is known by work of Varagnolo and Vasserot that the ``spherical representation" $\el(\mathrm{Triv})$ with lowest weight the trivial rep is finite-dimensional if and only if the denominator $d$ of the parameter $c=r/d$ is an elliptic number of $W$ \cite{VV}.

\subsection{Lemmas} Lemmas 3.2, 3.4, 3.5, and 3.6 will be referred to by their nicknames (E), (Symm), (dim Hom), and (RR) rather than by their numbers throughout the text. 

The following lemma was used by \cite{BP} in their study of rational Cherednik algebras for $H_3$, and will occasionally be used in a similar way in this study:
\lemma \cite{ES} If $h_c(\sigma)-h_c(\tau)=1$ and $\sigma\subset\hstar\otimes\tau$, then $\sigma$ generates an $H_c$-subrep of $\m(\tau)$.\\

Let $E=\frac{1}{2}\sum x_i^2$, part of the $\s$-triple $(E,H,F)\subset H_c(W)$.

Denote by $\el[n]$ the graded degree $n$ piece of an irreducible $H_c(W)$-module $\el=\el(\tau)$ and write $\el=\oplus_{k\geq0}\el[h_c(\tau)+k]$. 

\lemma \textbf{(E)} Let $\el$ be an irreducible $H_c(W)$-representation. Suppose there exists $v\in\el$ such that $E\cdot v=0$. Then $\el$ is finite-dimensional.

\begin{proof}
Say $v$ lies in $\el[n]$. Suppose $\el$ is not finite-dimensional -- then there are $x_1$ and $x_2$ in $\hstar$ such that $0\neq x_2x_1v\in\el[n+2]$. 

Indeed, suppose there is a vector $w\in\el$ with $\hstar\cdot w=0$. Say $w$ lies in graded degree $m$. Since $\el$ is irreducible, $w$ is a cyclic vector for $H_c(W)$. Then any vector in $\el[m+1]$ can be written as $z\cdot w$ for some $z$ in $H_c(W)$. Such a $z$ has graded degree $1$. By PBW theorem, it can be written as $xz'$ where $z'$ has graded degree $0$. Such an element $z'$ can be written, again by PBW theorem, as $z'=\sum_{j=1}^Ny^{a_j}(\sum_wc_{j,w}w)x^{b_j}$. Since $z'$ has degree $0$, $a_j=b_j$ for all $j=1,...,N$. Then the only terms of $z'$ that don't kill $v$ are those with $a_j=b_j=0$, i.e. an element of the group algebra of $W$. But if any $x\in\hstar$ applied to $v$ is $0$, then $(w^{-1}\cdot x)\cdot v=0$ for any $x\in\hstar$ and any $w\in W$, and so $x\cdot (w\cdot v)=(xw)\cdot v=(w(w^{-1}\cdot x))\cdot v=0$. So $z\cdot w=(xz')\cdot w=0$ for any $z$ of degree $1$, and therefore $\el[m+1]=0$.

So if $\el$ is infinite-dimensional we can find $x_2$ and $x_1$ as above. Then $E\cdot x_2x_1v=x_2x_1(E\cdot v)=0$, since $E$, $x_2$, $x_1$ all commute (they belong to the subalgebra $\mathbb{C}[\mathfrak{h}]$ which is a polynomial algebra in rank $W$ variables). So every application of $E$ starting from $\el[n]$ lowers the dimension of the graded pieces: $\dim E^k(\el[n])<\dim E^{k-1}(\el[n])$, and as the graded pieces are all finite-dimensional, eventually  for some $K\in\mathbb{N}$, $E^K(\el[n])=0$. So $\el$ contains a finite-dimensional $\s$-representation.

Now let $\el_{\mathrm{fin}}:=\{v\in\el\;|\;E^Nv=0\;\textrm{for all }N>>0\}$ be the finite-dimensional part of $\el$ as an $\s$-representation, and likewise let $\el_\infty:=\{v'\in\el\;|\;E^Nv'\neq0\textrm{ for all }N>0\}$ be the infinite-dimensional part of $\el$ as an $\s$-representation. Then $\el\cong\el_{\mathrm{fin}}\oplus\el_\infty$ as $\s$-reps. I claim that $\el_{\mathrm{fin}}$ is in fact an $H_c(W)$-subrepresentation of $\el$. Say $v'\in\el_\infty$. Suppose $v'=\prod_{j=1}^r x_j v$ for some $v\in\el_{\mathrm{fin}}$ and some elements $x_1,...,x_r\in\hstar$. Then $E^N v'=E^N\prod x_i v=\prod x_i E^N v=0$ for all $N>>0$, but then $v'\in\el_{\mathrm{fin}}$. Therefore $\hstar\cdot\el_{\mathrm{fin}}\subset\el_{\mathrm{fin}}$. Secondly, $\el_{\mathrm{fin}}$ is preserved under the action of $W$:  $E$ commutes with $W$, so $E^N(w\cdot v)=wE^Nv=0$ if $v\in\el_{\mathrm{fin}}$ and $N>>0$. Finally, take $y\in\mathfrak{h}$ and $v\in\el_{\mathrm{fin}}$, and let $N$ be large enough that $E^Nv=0$. By the quadratic relations defining $H_c(W)$,
\begin{align*}
Eyv&=\sum_{i=1}^nx_i^2yv=\sum_{i=1}^nx_i(x_iy)v=\sum_{i=1}^nx_iyx_iv+x_iC_i(w)v\\&=\sum_{i=1}^nyx_i^2v+(C_i(w)x_i+x_iC_i(w))v\\&=yEv+\sum_{i=1}^n(C_i(w)x_i+x_iC_i(w))v\\
\end{align*}
where $C_i(w)$ are the elements of $\mathbb{C}[W]$ given by the commutation relation between $y$ and $x_i$ (and now $x_i$ is the orthonormal basis of $\hstar$ as in the definition of the $\s$-triple $(E,h,F)$ introduced earlier). Therefore 
\begin{align*}
E^Kyv&=E^{K-1}\left(yEv+\sum_{i=1}^n(C_i(w)x_i+x_iC_i(w))v\right)\\&=E^{K-1}yEv+\sum_{i=1}^n(C_i(w)x_i+x_iC_i(w))E^{K-1}v\\&=E^{K-1}yEv
\end{align*}
 for any $K>N$. Consequently, 
 $$E^{2N}yv=E^NyE^Nv=0$$
 and so $yv\in\el_{\mathrm{fin}}$. To conclude: $\hstar$, $W$, and $\mathfrak{h}$ generate $H_c(W)$ and they all preserve $\el_{\mathrm{fin}}$, and therefore $\el_{\mathrm{fin}}$ is a nonzero $H_c$-subrepresentation of $\el$. But $\el$ is irreducible, so $\el_{\mathrm{fin}}=\el$, and $\el$ is finite-dimensional.

\end{proof}

\corollary Let $\el$ be an irreducible $H_c(W)$-representation. If $\dim\el[2]<\dim\el[0]$ then $\el$ is finite-dimensional.

Finite-dimensional representations, and more generally, minimal support representations in a block, possess symmetry in their composition factors with respect to tensoring lowest weights with sign:

\lemma \textbf{(Symm)} If $\el(\tau)$ is finite-dimensional, and 
$$\el(\tau)=\sum_{\sigma\in\Irr W}n_{\tau\sigma}M(\sigma)$$
is its decomposition into Vermas, then $$\sum_{\sigma\in\Irr W}n_{\tau\sigma}M(\sigma')=\pm\el(\tau)$$

More generally, the same lemma holds for $\el(\tau)$ an irreducible representation such that $\dim\Supp\el(\tau)$ is minimal over all dimensions of support of irreducibles $\el(\sigma)$ belonging to the same block as $\el(\tau)$; as Pavel Etingof explained to the author, this follows from the Cohen-Macaulay property of minimal support irreps proved in \cite{EGL}. This gives a kind of duality between minimal support irreps in ``dual blocks." A block might contain $\tau'$ if and only if it contains $\tau$, in which case the block is self-dual. If for any $\tau$ labeling objects in a block $\B$ of $\oh_c(W)$, objects with lowest weight $\tau'$ belong to another block $\B'$, the block $\B'$ will be called the dual block of $\B$. In either case, a minimal support irrep $\el(\tau)$ in $\B$ has a dual minimal support irrep $\el(\tau)^\vee$ whose lowest weight is $\sigma'$ where $\sigma$ is the lowest weight of maximal $h_c$-weight of the Vermas appearing in the Verma-decomposition of $\el(\tau)$. That is, if $\mu<_c\sigma$ for all $\mu$ such that $[\el(\tau):\m(\mu)]\neq0$ and $[\el(\tau):\m(\sigma)]=\pm1$ then $\el(\tau)^\vee=\el(\sigma')$, and the Verma-decomposition of $\el(\sigma')$ is given by $$\el(\tau)^\vee=\el(\sigma')=\pm\sum_{\mu\geq_c\tau}[\el(\tau):\m(\mu)]\m(\mu')$$

The following lemma equates dimensions of Hom spaces between two Vermas with lowest weights $\tau$, $\sigma$ and the two Vermas with those lowest weights tensored by the sign rep. The proof was communicated to the author by Pavel Etingof.

\lemma \textbf{(dim Hom)} Suppose $c$ is not a half-integer. Then $\dim\Hom(\m(\sigma),\m(\tau))=\dim\Hom(\m(\tau'),\m(\sigma'))$.
\begin{proof}
There is an isomorphism $H_c(W)\rightarrow H_{-c}(W)$ which sends $x$ to $x$, $y$ to $y$, and $s$ to $-s$ for reflections $s\in S$. The induced map $\oh_c(W)$ to $\oh_{-c}(W)$ sends $\m(\tau)$ to $\m(\tau')$. If $c$ is not a half-integer then the KZ functor is $0$-faithful by Theorem 5.3 of \cite{R}, and thus $\dim\Hom(\m(\sigma),\m(\tau))=\dim\Hom(\KZ\m(\sigma),\KZ\m(\tau))$. Let $\mathrm{S}(\tau)$ be the standard object labeled by $\tau$ in the category $\mathcal{H}_q$-mod that is the target of $\KZ$ functor. By Theorem 6.8 and Remark 6.9 of \cite{GGOR}, $\KZ\m(\tau)=\mathrm{S}(\tau)$ if $c>0$ and $\KZ\m(\tau)=\mathrm{S}(\tau)^*$ if $c<0$.
\end{proof}

In general, $[\m(\tau):\el(\sigma)]\geq\dim\Hom(\m(\sigma),\m(\tau))$. The following lemma gives a condition when equality holds and will be used in conjunction with (dim Hom) to recover decomposition numbers near the diagonal of the decomposition matrix.

\lemma \textbf{(RR)} Let $\tau<_c\sigma$. Suppose $\m(\sigma)$ is not linked to any $\m(\nu_1)$ that is linked to any $\m(\nu_2)$ that is linked to... that is linked to any $\m(\nu_k)$ that is linked to $\m(\tau)$, $\tau<\nu_k<_c...<_c\nu_2<_c\nu_1<_c\sigma$. Then $[\m(\tau):\el(\sigma)]=\dim\Hom(\m(\sigma),\m(\tau)).$

\begin{proof}
Let $\nu$ be linked by a chain to $\tau$ as above. The assumption implies that for all $i$, $\Ext^i(\m(\sigma),\m(\nu))=0$ .  

Let $\G(\nu)$ be any quotient of $\m(\nu)$. The first thing to show is that $\Ext^i(\m(\sigma),\G(\nu))=0$ also. By induction on $h_c(\sigma)-h_c(\nu)$: if $h_c(\sigma)-h_c(\nu)=0$ then $\sigma$ is not greater than $\nu$ in the partial order and so $\Ext^i(\m(\sigma),\G(\nu))=0$. Now assume $\Ext^i(\m(\sigma),\G(\nu'))=0$ for any $\nu'$ between $\tau$ and $\sigma$ such that $\nu'$ is linked by a chain to $\tau$ and such that $h_c(\sigma)-h_c(\nu')\leq n$. Consider $\tau<_c\nu<_c\sigma$ such that $h_c(\sigma)-h_c(\nu)>n$ and for any $\nu<_c\nu'<_c\sigma$, $h_c(\sigma)-h_c(\nu')\leq n$. Suppose $\m(\nu)$ is not simple, otherwise there is nothing to show. Then there is a short exact sequence $$0\longrightarrow \mathrm{N}\longrightarrow\m(\nu)\longrightarrow\G(\nu)\longrightarrow0$$
for some lowest-weight module $N$ of weight $\nu'$. Applying $\Hom(\m(\sigma),-)$ yields a long exact sequence of $\Ext$:
\begin{align*}
&0\longrightarrow \Hom(\m(\sigma),N)\longrightarrow \Hom(\m(\sigma),\m(\nu))\longrightarrow \Hom(\m(\sigma),\G(\nu))\\&\longrightarrow \Ext^1(\m(\sigma),N)\longrightarrow \Ext^1(\m(\sigma),\m(\nu))\longrightarrow\Ext^1(\m(\sigma),\G(\nu)) \longrightarrow ....
\end{align*}
By assumption, $\Ext^i(\m(\sigma),\m(\nu))=0$ vanishes for all $i$. Since $N$ has lowest weight $\nu'$ with $h_c(\sigma)-h_c(\nu')\leq n$, by induction $\Ext^i(\m(\sigma),N)$ vanishes for all $i$. Therefore $\Ext^i(\m(\sigma),\G(\nu))$ vanishes for all $i$ as well.

Now suppose that for some $\tau<_c\nu<_c\sigma$, there is a nonzero homomorphism\\ $\phi_\nu:\m(\nu)\rightarrow\m(\tau)$. Let $\G(\nu)$ be the image of $\phi_\nu$ in $\m(\tau)$. Then there is a short exact sequence 
$$0\longrightarrow\G(\nu)\longrightarrow\m(\tau)\longrightarrow\m(\tau)/\G(\nu)\longrightarrow0$$
giving rise to the long exact sequence of Ext with $\m(\sigma)$ fixed in the first argument:
\begin{align*}
&0\longrightarrow\Hom(\m(\sigma),\G(\nu))\longrightarrow\Hom(\m(\sigma),\m(\tau))\longrightarrow\Hom(\m(\sigma),\m(\tau)/\G(\nu))\\
&\longrightarrow\Ext^1(\m(\sigma),\G(\nu))\longrightarrow...
\end{align*}
By the previous paragraph, $\Hom(\m(\sigma),\G(\nu))=\Ext^1(\m(\sigma),\G(\nu))=0$. So there is an isomorphism $\Hom(\m(\sigma),\m(\tau))\cong\Hom(\m(\sigma),\m(\tau)/\G(\nu))$. This means that taking quotients of $\m(\tau)$ by any $\G(\nu)$ for $\nu$ between $\tau$ and $\sigma$ does not increase the multiplicity of $\sigma$ consisting of singular vectors in degree $h_c(\sigma)$. This proves the lemma.
\end{proof}

The conditions of the previous lemma can be phrased in terms of a game of ``Ricochet Robots" on the decomposition matrix, with a robot starting on the diagonal at $(\tau,\tau)$ and ending at $(\sigma,\sigma)$. The rules of moving a robot in the game Ricochet Robots is that the robot moves in a straight horizontal or vertical line until it hits a wall, at which point it starts another path in a straight line; this continues until the robot reaches its destination. In our case, a ``wall" is either (i) a nonzero decomposition number in the same row as the robot, or (ii) a $1$ on the diagonal in the same column as the robot. Additionally, we only allow the robot to move to the right and down. The Ricochet Robot interpretation of the conditions of (RR) is that there is no path starting from $(\tau,\tau)$ and ending at $(\sigma,\sigma)$ that does not go through the square $(\tau,\sigma)$ for a robot making moves to the right and down with the ``walls" specified above.

\textit{The general strategy for finding the decomposition matrices of the blocks of Category $\oh_c(W)$}: start by embedding the Hecke algebra decomposition matrix in the RCA decomposition matrix; use (dim Hom) together with (RR) to recover as many decomposition numbers as possible -- these will be mostly near the diagonal; then exploit induction and restriction from $H_c(W')$ for maximal parabolics $W'\subset W$.

\section{Simple representations of $H_c(H_4)$}

Every regular number of $H_4$ is an elliptic number, and so $\el(1)$ is always finite-dimensional when the Category O is not semisimple. There are $10$ regular numbers of $H_4$: $$30,\;20,\;15,\;12,\;10,\;6,\;5,\;4,\;3,\;2$$ When $c=\frac{1}{n}$ with $n\in\{3,5,15\}$, then $\oh_c$ has two blocks containing finite-dimensional modules; otherwise all finite-dimensional modules appear in a single block.

The labels for the thirty-four characters of $H_4$ used throughout this section are those used by Grove in his paper where the character table of $H_4$ first appeared \cite{Gr}. $\mbf{3}$ is the standard rep and $\mbf{1}$ is the trivial rep.

\subsection{\Large{$\mbf{c=\frac{1}{30}}$, $\mbf{\frac{1}{20}}$, $\mbf{\frac{1}{15}}$} \normalsize{and} \Large$\mbf{\frac{1}{12}}$}

 For these parameters with denominators ``cuspidal numbers" of $H_4$, all Vermas are simple unless they belong to the composition series of a finite-dimensional module \cite{BE}. In the cases of $c=\frac{1}{30},\;\frac{1}{20},$ and $\frac{1}{12}$, there's a unique nontrivial block, the one containing $\el(1)$. For $c=\frac{1}{15}$ there are two nontrivial blocks, one containing $\el(1)$ and the other containing the finite-dimensional rep $\el(5)$ (the lowest weight $5$ is the Galois-twist of the standard rep). In all cases the decomposition matrices have the form:
$$
\mathcal{L}=\begin{pmatrix}
1&1&0&0&0\\0&1&1&0&0\\0&0&1&1&0\\0&0&0&1&1\\0&0&0&0&1
\end{pmatrix}\qquad\qquad
\mathcal{M}=\begin{pmatrix}
1&-1&1&-1&1\\0&1&-1&1&-1\\0&0&1&-1&1\\0&0&0&1&-1\\0&0&0&0&1
\end{pmatrix}
$$
The appropriate reps labeling the rows and columns for a particular parameter can be read off of the table in Geck-Pfeiffer \cite{GP}, as these blocks map via KZ functor to the blocks of defect one in the corresponding Hecke algebras at $q=e^{2pic}$. The decompositions can also be determined without any ambiguity from tensor product decompositions $S^k\hstar\otimes\tau$: for each rep $\tau$, one can compute that for each $k\geq1$ no representation $\sigma>_c\tau$ with $h_c(\sigma)-h_c(\tau)=k$ appears in the decomposition of $S^k\hstar\otimes\tau$, except possibly when $k$ equals the distance from $h_c(\tau)$ to the next weight to its right on the weight line. For a given weight line there is one such potentially non-simple Verma for each weight on the line and these form a chain. In the case of a block containing $\el(1)$, the fact that $\el(1)$ is finite-dimensional forces the first such rep to the immediate right of $\el(1)$ to generate a subrep; this produces a domino effect on the reps to its right and it follows that the decomposition matrices have the form above. In the case of the block containing $\el(5)$ when $c=\frac{1}{15}$, the weights are spaced $1$ apart, in which case $\sigma\subset \hstar\otimes\tau$ implies that $\sigma$ consists of singular vectors in $\m(\tau)$. 

The decompositions, characters, and dimensions of the finite-dimensional reps at these parameters are:\\

{\Large$\mathbf{c=\frac{1}{30}}:$}\\

\begin{center}

\begin{picture}(200,50)
\put(0,30){\line(1,0){200}}
\put(0,30){\circle*{5}}
\put(50,30){\circle*{5}}
\put(100,30){\circle*{5}}
\put(150,30){\circle*{5}}
\put(200,30){\circle*{5}}
\put(0,40){$0$}
\put(48,40){$1$}
\put(98,40){$2$}
\put(148,40){$3$}
\put(198,40){$4$}
\put(-2,15){$\mathbf{1}$}
\put(48,15){$\mathbf{3}$}
\put(98,15){$\mathbf{7}$}
\put(148,15){$\mathbf{4}$}
\put(198,15){$\mathbf{2}$}
\end{picture}
\end{center}
 \begin{align*}
 %&\el(\mathbf{1})=\m(\mathbf{1})-\m(\mathbf{3})+\m(\mathbf{7})-\m(\mathbf{4})+\m(\mathbf{2})\\
 &\chi_{\el(\mathbf{1})}(w,t)=\chi_{\mathbf{1}},\;\;\;\chi_{\el(\mathbf{1})}(t)=1,\;\;\;\;\mathrm{dim}\;\el(\mathbf{1})=1\\
 \end{align*}
 \\
 
 {\Large$ \mathbf{c=\frac{1}{20}:}$} \\
 
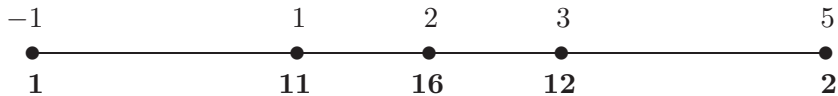
\begin{center}

\begin{picture}(300,50)
\put(0,30){\line(1,0){300}}
\put(0,30){\circle*{5}}
\put(100,30){\circle*{5}}
\put(150,30){\circle*{5}}
\put(200,30){\circle*{5}}
\put(300,30){\circle*{5}}
\put(-10,40){$-1$}
\put(98,40){$1$}
\put(148,40){$2$}
\put(198,40){$3$}
\put(298,40){$5$}
\put(-2,15){$\mathbf{1}$}
\put(93,15){$\mathbf{11}$}
\put(143,15){$\mathbf{16}$}
\put(193,15){$\mathbf{12}$}
\put(298,15){$\mathbf{2}$}
\end{picture}
\end{center}

 \begin{align*}
% &\el(\mathbf{1})=\m(\mathbf{1})-\m(\mathbf{11})+\m(\mathbf{16})-\m(\mathbf{12})+\m(\mathbf{2})\\
&\chi_{\el(\mathbf{1})}(w,t)=\chi_{\mathbf{1}}(t^{-1}+t)+\chi_{\mathbf{3}},\qquad\chi_{\el(\mathbf{1})}(t)=t^{-1}+4+t,\qquad\mathrm{dim}\;\el(\mathbf{1})=6\\
\end{align*}
\\

{\Large$\mathbf{c=\frac{1}{15}:}$}\\

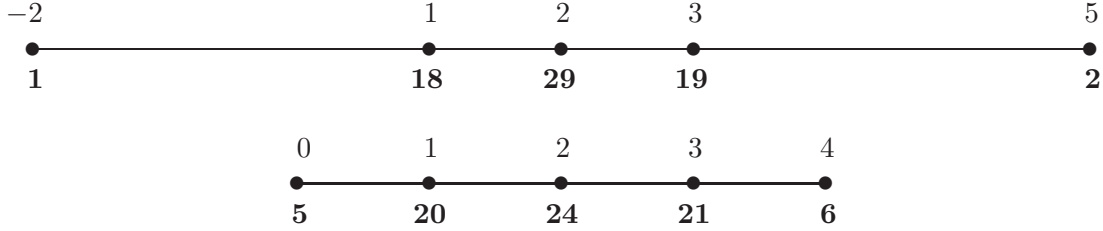
\begin{center}
\begin{picture}(400,50)
\put(0,30){\line(1,0){400}}
\put(0,30){\circle*{5}}
\put(150,30){\circle*{5}}
\put(200,30){\circle*{5}}
\put(250,30){\circle*{5}}
\put(400,30){\circle*{5}}
\put(-10,40){$-2$}
\put(148,40){$1$}
\put(198,40){$2$}
\put(248,40){$3$}
\put(398,40){$5$}
\put(-2,15){$\mathbf{1}$}
\put(143,15){$\mathbf{18}$}
\put(193,15){$\mathbf{29}$}
\put(243,15){$\mathbf{19}$}
\put(398,15){$\mathbf{2}$}
\end{picture}

\begin{picture}(200,50)
\put(0,30){\line(1,0){200}}
\put(0,30){\circle*{5}}
\put(50,30){\circle*{5}}
\put(100,30){\circle*{5}}
\put(150,30){\circle*{5}}
\put(200,30){\circle*{5}}
\put(0,40){$0$}
\put(48,40){$1$}
\put(98,40){$2$}
\put(148,40){$3$}
\put(198,40){$4$}
\put(-2,15){$\mathbf{5}$}
\put(44,15){$\mathbf{20}$}
\put(94,15){$\mathbf{24}$}
\put(144,15){$\mathbf{21}$}
\put(198,15){$\mathbf{6}$}
\end{picture}
\end{center}

\begin{align*}
%&\el(\mathbf{1})=\m(\mathbf{1})-\m(\mathbf{18})+\m(\mathbf{29})-\m(\mathbf{19})+\m(\mathbf{2})\\
&\chi_{\el(\mathbf{1})}(w,t)=\chi_{\mathbf{1}}(t^{-2}+t^2)+\chi_{\mathbf{3}}(t^{-1}+t)+\chi_{\mathbf{1}}+\chi_{\mathbf{11}}+,\qquad\chi_{\el(\mathbf{1})}(t)=t^{-2}+t^2+4(t^{-1}+t)+10\\&\mathrm{dim}\;\el(\mathbf{1})=20\\
\\
%&\el(\mathbf{5})=\m(\mathbf{5})-\m(\mathbf{20})+\m(\mathbf{24})-\m(\mathbf{21})+\m(\mathbf{6})\\
&\chi_{\el(\mathbf{5})}(w,t)=\chi_{\mathbf{5}},\qquad\chi_{\el(\mathbf{5})}(t)=4,\qquad\mathrm{dim}\;\el(\mathbf{5})=4
\\
\end{align*}

{\Large$\mathbf{c=\frac{1}{12}:}$}\\

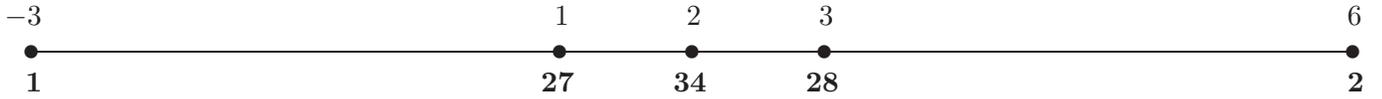
\begin{center}
\begin{picture}(500,50)
\put(0,30){\line(1,0){500}}
\put(0,30){\circle*{5}}
\put(200,30){\circle*{5}}
\put(250,30){\circle*{5}}
\put(300,30){\circle*{5}}
\put(500,30){\circle*{5}}
\put(-10,40){$-3$}
\put(198,40){$1$}
\put(248,40){$2$}
\put(298,40){$3$}
\put(498,40){$6$}
\put(-2,15){$\mathbf{1}$}
\put(193,15){$\mathbf{27}$}
\put(243,15){$\mathbf{34}$}
\put(293,15){$\mathbf{28}$}
\put(498,15){$\mathbf{2}$}
\end{picture}
\end{center}

\begin{align*}
%& \el(\mathbf{1})=\m(\mathbf{1})-\m(\mathbf{27})+\m(\mathbf{34})-\m(\mathbf{28})+\m(\mathbf{2})\\
&\chi_{\el(\mathbf{1})}(w,t)=\chi_{\mathbf{1}}(t^{-3}+t^3)+\chi_{\mathbf{3}}(t^{-2}+t^2)+(\chi_{\mathbf{1}}+\chi_{\mathbf{11}})(t^{-1}+t)+\chi_{\mathbf{3}}+\chi_{\mathbf{18}}\\&\chi_{\el(\mathbf{1})}(t)=t^{-3}+t^3+4(t^{-2}+t^2)+10(t^{-1}+t)+20\\&\mathrm{dim}\;\el(\mathbf{1})=50\\
\end{align*}

\subsection{\Large$\mbf{c=\frac{1}{10}}$} There are two blocks containing finite-dimensional irreps -- the principal block, and a block of defect $1$ in which lies the irrep with lowest weight the standard rep of $H_4$, which we'll nickname ``the standard block."

\subsubsection{The principal block}

\begin{center}
\begin{picture}(400,80)
\put(20,60){\line(1,0){360}}
\put(20,60){\circle*{5}}
\put(110,60){\circle*{5}}
\put(140,60){\circle*{5}}
\put(170,60){\circle*{5}}
\put(200,60){\circle*{5}}
\put(230,60){\circle*{5}}
\put(260,60){\circle*{5}}
\put(290,60){\circle*{5}}
\put(380,60){\circle*{5}}
\put(10,70){$-4$}
\put(100,70){$-1$}
\put(138,70){$0$}
\put(168,70){$1$}
\put(198,70){$2$}
\put(228,70){$3$}
\put(258,70){$4$}
\put(288,70){$5$}
\put(378,70){$8$}
\put(17,45){$\mathbf{1}$}
\put(107,45){$\mathbf{5}$}
\put(133,45){$\mathbf{13}$}
\put(163,45){$\mathbf{31}$}
\put(193,45){$\mathbf{26}$}
\put(193,30){$\mathbf{30}$}
\put(193,15){$\mathbf{33}$}
\put(223,45){$\mathbf{32}$}
\put(253,45){$\mathbf{14}$}
\put(287,45){$\mathbf{6}$}
\put(377,45){$\mathbf{2}$}
\end{picture}
\end{center}

\[
\begin{blockarray}{ccccccccccccc}
\begin{block}{cc(ccccccccccc)}
\star& \mbf{1}&1&\cdot&\cdot&1&\cdot&\cdot&1&\cdot&\cdot&\cdot&\cdot\\
\star&\mbf{5}&\cdot&1&\cdot&1&\cdot&1&\cdot&\cdot&\cdot&\cdot&\cdot\\
\star&\mbf{13}&\cdot&\cdot&1&1&1&\cdot&\cdot&\cdot&\cdot&\cdot&\cdot\\
\mathtt{(1)}&\mbf{31}&\cdot&\cdot&\cdot&1&1&1&1&1&\cdot&\cdot&\cdot\\
&\mbf{33}&\cdot&\cdot&\cdot&\cdot&1&\cdot&\cdot&1&1&\cdot&\cdot\\
&\mbf{30}&\cdot&\cdot&\cdot&\cdot&\cdot&1&\cdot&1&\cdot&1&\cdot\\
&\mbf{26}&\cdot&\cdot&\cdot&\cdot&\cdot&\cdot&1&1&\cdot&\cdot&1\\
&\mbf{32}&\cdot&\cdot&\cdot&\cdot&\cdot&\cdot&\cdot&1&1&1&1\\
&\mbf{14}&\cdot&\cdot&\cdot&\cdot&\cdot&\cdot&\cdot&\cdot&1&\cdot&\cdot\\
&\mbf{6}&\cdot&\cdot&\cdot&\cdot&\cdot&\cdot&\cdot&\cdot&\cdot&1&\cdot\\
&\mbf{2}&\cdot&\cdot&\cdot&\cdot&\cdot&\cdot&\cdot&\cdot&\cdot&\cdot&1\\
\end{block}
\\
\begin{block}{cc(ccccccccccc)}
\star&\mbf{1}& 1&  \cdot&  \cdot& -1&  1&  1&  \cdot& -1&  \cdot&  \cdot&  1\\
\star&\mbf{5}& \cdot & 1 & \cdot &-1 & 1&  \cdot&  1 &-1 & 1 & \cdot & \cdot\\
\star&\mbf{13}& \cdot & \cdot & 1 &-1 & \cdot & 1  &1 &-1  &\cdot & 1&  \cdot\\
\mathtt{(1)}&\mbf{31}& \cdot&  \cdot & \cdot & 1 &-1 &-1& -1&  2& -1& -1& -1\\
&\mbf{33}& \cdot & \cdot & \cdot & \cdot & 1 & \cdot & \cdot &-1 & 1&  \cdot & 1\\
&\mbf{30}& \cdot & \cdot & \cdot & \cdot  &\cdot & 1 & \cdot &-1 & \cdot & 1 & 1\\
&\mbf{26}& \cdot & \cdot & \cdot&  \cdot & \cdot & \cdot & 1 &-1  &1 & 1 & \cdot\\
&\mbf{32}& \cdot & \cdot & \cdot & \cdot & \cdot & \cdot & \cdot & 1 &-1& -1& -1\\
&\mbf{14}& \cdot&  \cdot & \cdot & \cdot & \cdot & \cdot  &\cdot  &\cdot  &1  &\cdot & \cdot\\
&\mbf{6}& \cdot & \cdot & \cdot & \cdot & \cdot & \cdot & \cdot & \cdot & \cdot & 1 & \cdot\\
&\mbf{2}& \cdot & \cdot & \cdot & \cdot & \cdot & \cdot & \cdot & \cdot & \cdot & \cdot & 1\\
\end{block}
\end{blockarray}
\]

The entries in the first four columns of $\left([\m(\tau):\el(\sigma)]\right)$ are obtained from the rest of the matrix (the Hecke algebra decomposition matrix for $q=e^{\frac{2\pi i}{10}}$) by applying (dim Hom) and (RR).

\remark It is possible to compute this entire decomposition matrix by hand without reference to the Hecke algebra. The author originally took this approach. To do this one must play with tensor products of the lowest weights $\tau$ with symmetric powers of the standard rep, use induction from the rational Cherednik algebra of $H_3$ at the same parameter, apply (E), apply (dim Hom), make whatever ad hoc arguments about dimensions of graded pieces of a representation to rule out this or that possible decomposition of a simple representation, and so on. Eventually one can find the whole decomposition matrix honestly and laboriously in this way and it coincides with that obtained by the slicker method here. 

\begin{align*}
&\chi_{\el(\mbf{1})}(w,t)=\chi_{\mathbf{1}}(t^{-4}+t^4)+\chi_{\mathbf{3}}(t^{-3}+t^3)+(\chi_{\mathbf{1}}+\chi_{\mathbf{11}})(t^{-2}+t^2)\\&\qquad\qquad\qquad+(\chi_{\mathbf{3}}+\chi_{\mathbf{18}})(t^{-1}+t)+\chi_{\mathbf{1}}+\chi_{\mathbf{11}}+\chi_{\mathbf{27}}\\
&\chi_{\el(\mbf{1})}(t)=t^{-4}+t^4+4(t^{-3}+t^3)+10(t^{-2}+t^2)+20(t^{-1}+t)+35\\
&\dim\el(\mbf{1})=105\\
\\
&\chi_{\el(\mathbf{5})}(w,t)=\chi_\mathbf{5}(t^{-1}+t)+\chi_\mathbf{20}\\
&\chi_{\el(\mbf{5})}(t)=4(t^{-1}+t)+16\\
&\dim\el(\mbf{5})=24\\
\\
&\chi_{\el(\mathbf{13})}(w,t)=\chi_{\mathbf{13}}\\
&\chi_{\el(\mbf{13})}(t)=9\\
&\dim\el(\mbf{13})=9\\
\\
&\chi_{\el(\mbf{31})}(t)=\frac{t^4 + 3t^3 + 6t^2 + 14t + 36}{1-t}\\
&\dim\Supp\el(\mbf{31})=1\\
\end{align*}

\subsubsection{The standard block} This block is of defect $1$.

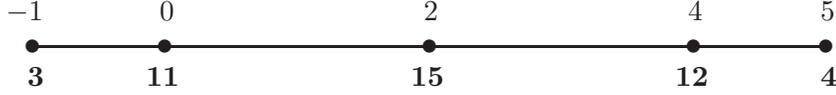
\begin{center}

\begin{picture}(300,50)
\put(0,20){\line(1,0){300}}
\put(0,20){\circle*{5}}
\put(50,20){\circle*{5}}
\put(150,20){\circle*{5}}
\put(250,20){\circle*{5}}
\put(300,20){\circle*{5}}
\put(-10,30){$-1$}
\put(48,30){$0$}
\put(148,30){$2$}
\put(248,30){$4$}
\put(298,30){$5$}
\put(-2,5){$\mathbf{3}$}
\put(43,5){$\mathbf{11}$}
\put(143,5){$\mathbf{15}$}
\put(243,5){$\mathbf{12}$}
\put(298,5){$\mathbf{4}$}
\end{picture}
\end{center}

\begin{align*}
&\chi_{\el(\mbf{3})}(w,t)=\chi_{\mbf{3}}(t^{-1}+t)+\chi_{\mbf{1}}+\chi_{\mbf{11}}\\
&\chi_{\el(\mbf{3})}(t)=4(t^{-1}+t)+7\\
&\dim\el(\mbf{3})=15\\
\end{align*}

\subsection{\Large$\mbf{c=\frac{1}{6}}$}
The decomposition matrix for $c=\frac{1}{6}$ is identical to that of the principal block for $c=\frac{1}{10}$ up to a relabeling of the $W$-irreps and swapping the order of a couple of the sign-invariant reps in the middle. Not only that, but the dimension of supports of the modules match up in the two blocks.

\remark This entire decomposition matrix may also be computed by hand as in the case of $c=\frac{1}{10}$.

\begin{center}
\begin{picture}(400,80)
\put(0,60){\line(1,0){400}}
\put(0,60){\circle*{5}}
\put(100,60){\circle*{5}}
\put(160,60){\circle*{5}}
\put(200,60){\circle*{5}}
\put(240,60){\circle*{5}}
\put(300,60){\circle*{5}}
\put(400,60){\circle*{5}}
\put(-10,70){$-8$}
\put(90,70){$-3$}
\put(158,70){$0$}
\put(198,70){$2$}
\put(238,70){$4$}
\put(298,70){$7$}
\put(395,70){$12$}
\put(-3,45){$\mathbf{1}$}
\put(97,45){$\mathbf{3}$}
\put(97,30){$\mathbf{5}$}
\put(153,45){$\mathbf{27}$}
\put(193,45){$\mathbf{22}$}
\put(193,30){$\mathbf{25}$}
\put(193,15){$\mathbf{26}$}
\put(233,45){$\mathbf{28}$}
\put(297,45){$\mathbf{4}$}
\put(297,30){$\mathbf{6}$}
\put(397,45){$\mathbf{2}$}
\end{picture}
\end{center}

\[
\begin{blockarray}{ccccccccccccc}
\begin{block}{cc(ccccccccccc)}
\star & \mbf{1}&1&\cdot&\cdot&1&\cdot&\cdot&1&\cdot&\cdot&\cdot&\cdot\\
\star&\mbf{5}&\cdot&1&\cdot&1&\cdot&1&\cdot&\cdot&\cdot&\cdot&\cdot\\
\star&\mbf{3}&\cdot&\cdot&1&1&1&\cdot&\cdot&\cdot&\cdot&\cdot&\cdot\\
\mathtt{(1)}&\mbf{27}&\cdot&\cdot&\cdot&1&1&1&1&1&\cdot&\cdot&\cdot\\
&\mbf{26}&\cdot&\cdot&\cdot&\cdot&1&\cdot&\cdot&1&\cdot&1&\cdot\\
&\mbf{25}&\cdot&\cdot&\cdot&\cdot&\cdot&1&\cdot&1&1&\cdot&\cdot\\
&\mbf{22}&\cdot&\cdot&\cdot&\cdot&\cdot&\cdot&1&1&\cdot&\cdot&1\\
&\mbf{28}&\cdot&\cdot&\cdot&\cdot&\cdot&\cdot&\cdot&1&1&1&1\\
&\mbf{6}&\cdot&\cdot&\cdot&\cdot&\cdot&\cdot&\cdot&\cdot&1&\cdot&\cdot\\
&\mbf{4}&\cdot&\cdot&\cdot&\cdot&\cdot&\cdot&\cdot&\cdot&\cdot&1&\cdot\\
&\mbf{2}&\cdot&\cdot&\cdot&\cdot&\cdot&\cdot&\cdot&\cdot&\cdot&\cdot&1\\
\end{block}
\\
\begin{block}{cc(ccccccccccc)}
\star&\mbf{1}& 1&  \cdot&  \cdot& -1&  1&  1&  \cdot& -1&  \cdot&  \cdot&  1\\
\star&\mbf{5}& \cdot & 1 & \cdot &-1 & 1&  \cdot&  1 &-1 & 1 & \cdot & \cdot\\
\star&\mbf{3}& \cdot & \cdot & 1 &-1 & \cdot & 1  &1 &-1  &\cdot & 1&  \cdot\\
\mathtt{(1)}&\mbf{27}& \cdot&  \cdot & \cdot & 1 &-1 &-1& -1&  2& -1& -1& -1\\
&\mbf{26}& \cdot & \cdot & \cdot & \cdot & 1 & \cdot & \cdot &-1 & 1&  \cdot & 1\\
&\mbf{25}& \cdot & \cdot & \cdot & \cdot  &\cdot & 1 & \cdot &-1 & \cdot & 1 & 1\\
&\mbf{22}& \cdot & \cdot & \cdot&  \cdot & \cdot & \cdot & 1 &-1  &1 & 1 & \cdot\\
&\mbf{28}& \cdot & \cdot & \cdot & \cdot & \cdot & \cdot & \cdot & 1 &-1& -1& -1\\
&\mbf{6}& \cdot&  \cdot & \cdot & \cdot & \cdot & \cdot  &\cdot  &\cdot  &1  &\cdot & \cdot\\
&\mbf{4}& \cdot & \cdot & \cdot & \cdot & \cdot & \cdot & \cdot & \cdot & \cdot & 1 & \cdot\\
&\mbf{2}& \cdot & \cdot & \cdot & \cdot & \cdot & \cdot & \cdot & \cdot & \cdot & \cdot & 1\\
\end{block}
\end{blockarray}
\]

\begin{align*}
&\chi_{\el(\mathbf{1})}(w,t)=\chi_\mathbf{1}(t^{-8}+t^8)+\chi_\mathbf{3}(t^{-7}+t^7)+(\chi_\mathbf{1}+\chi_\mathbf{11})(t^{-6}+t^6)+(\chi_\mathbf{3}+\chi_\mathbf{18})(t^{-5}+t^5)\\&\qquad\qquad\qquad+(\chi_\mathbf{1}+\chi_\mathbf{11}+\chi_\mathbf{27})(t^{-4}+t^4)+(\chi_\mathbf{3}+\chi_\mathbf{18}+\chi_\mathbf{31})(t^{-3}+t^3)\\&\qquad\qquad\qquad+(\chi_\mathbf{1}+\chi_\mathbf{11}+\chi_\mathbf{13}+\chi_\mathbf{20}+\chi_\mathbf{26}+\chi_\mathbf{27})(t^{-2}+t^2)\\
&\qquad\qquad+(\chi_\mathbf{3}+\chi_\mathbf{5}+\chi_\mathbf{18}+\chi_\mathbf{24}+2\chi_\mathbf{31})(t^{-1}+t)+\chi_\mathbf{1}+\chi_\mathbf{11}+\chi_\mathbf{13}+2\chi_\mathbf{20}+\chi_\mathbf{26}+\chi_\mathbf{27}+\chi_\mathbf{33}\\
&\chi_{\el(\mbf{1})}(t)=t^{-8}+t^8 + 4(t^{-7}+t^7) + 10(t^{-6}+t^6) + 20(t^{-5}+t^5) + 35(t^{-4}+t^4)\\
&\qquad\qquad\qquad + 56(t^{-3}+t^3) + 84(t^{-2}+t^2) + 120(t^{-1}+t) + 140\\
&\dim\el(\mbf{1})=800\\
\\
&\chi_{\el(\mathbf{3})}(w,t)=\chi_\mathbf{3}(t^{-3}+t^3)+(\chi_\mathbf{1}+\chi_\mathbf{7}+\chi_\mathbf{11})(t^{-2}+t^2)+(2\chi_\mathbf{3}+\chi_\mathbf{16}+\chi_\mathbf{18})(t^{-1}+t)\\&\qquad\qquad\qquad+\chi_\mathbf{1}+\chi_\mathbf{7}+2\chi_\mathbf{11}+\chi_\mathbf{29}\\
&\chi_{\el(\mbf{3})}(t)=4(t^{-3}+t^3) + 16(t^{-2}+t^2) + 40(t^{-1}+t) + 55\\
&\dim\el(\mbf{3})=175\\
\\
&\chi_{L(\mathbf{5})}(w,t)=\chi_\mathbf{5}(t^{-3}+t^3)+\chi_\mathbf{20}(t^{-2}+t^2)+(\chi_\mathbf{5}+\chi_\mathbf{31})(t^{-1}+t)+\chi_\mathbf{13}+\chi_\mathbf{20}+\chi_\mathbf{30}\\
&\chi_{\el(\mbf{5})}(t)=4(t^{-3}+t^3) + 16(t^{-2}+t^2) + 40(t^{-1}+t) + 55\\
&\dim\el(\mbf{5})=175\\
\\
&\chi_{\el(\mbf{27})}(t)=\frac{t^9 + 3t^8 + 6t^7 + 10t^6 + 15t^5 + 29t^4 + 52t^3 + 84t^2 + 75t + 25}{1-t}\\
&\dim\Supp\el(\mbf{27})=1\\
\end{align*}

\subsection{\Large{$\mbf{c=\frac{1}{5}}$} \normalsize{and} \Large{$\mbf{\frac{1}{3}}$}}

The rational Cherednik algebra at both of these parameters has two blocks containing finite-dimensional irreps, and all four blocks in question have the same decomposition matrix $([\m(\tau):\el(\sigma)])$ and the same dimensions of supports in each row:

\[
\begin{blockarray}{cccccccccc}
\begin{block}{c(ccccccccc)}
\star & 1 & \cdot & 1 & 1 & \cdot & \cdot & \cdot & \cdot & \cdot \\
\star & \cdot & 1 & 1 & \cdot & 1 & \cdot & \cdot & \cdot & \cdot \\
\mathtt{(2)} & \cdot & \cdot & 1 & 1 & 1 & 1 & 1 & \cdot & \cdot \\
& \cdot & \cdot & \cdot & 1 & \cdot & \cdot & 1 & \cdot & 1 \\
& \cdot & \cdot & \cdot & \cdot & 1 & \cdot & 1 & 1 & \cdot\\
\mathtt{(2)} & \cdot & \cdot & \cdot & \cdot & \cdot & 1 & 1 & \cdot & \cdot \\
& \cdot & \cdot & \cdot & \cdot & \cdot & \cdot & 1 & 1 & 1 \\
& \cdot & \cdot & \cdot & \cdot & \cdot & \cdot & \cdot & 1 & \cdot \\
& \cdot & \cdot & \cdot & \cdot & \cdot & \cdot & \cdot & \cdot & 1\\
\end{block}
\end{blockarray}
\]
\\
The inverse matrix $([\el(\tau):\m(\sigma)])$ is then
\[
\begin{blockarray}{cccccccccc}
\begin{block}{c(ccccccccc)}
\star & 1 & \cdot & -1 & \cdot & 1 & 1 & -1 & \cdot & 1 \\
\star & \cdot & 1 & -1 & 1 & \cdot & 1 & -1 & 1 & \cdot \\
\mathtt{(2)} & \cdot & \cdot & 1 & -1 & -1 & -1 & 2 & -1 & -1 \\
& \cdot & \cdot & \cdot & 1 & \cdot & \cdot & -1 & 1 & \cdot \\
& \cdot & \cdot & \cdot & \cdot & 1 & \cdot & -1 & \cdot & 1\\
\mathtt{(2)} & \cdot & \cdot & \cdot & \cdot & \cdot & 1 & -1 & 1 & 1 \\
& \cdot & \cdot & \cdot & \cdot & \cdot & \cdot & 1 & -1 & -1 \\
& \cdot & \cdot & \cdot & \cdot & \cdot & \cdot & \cdot & 1 & \cdot \\
& \cdot & \cdot & \cdot & \cdot & \cdot & \cdot & \cdot & \cdot & 1\\
\end{block}
\end{blockarray}
\]
\\
The labelings of rows and columns for these four blocks and their $h_c$-weight lines are as follows:
\begin{align*}
&c=\frac{1}{5}\qquad\{\mbf{1},\;\mbf{11},\;\mbf{18},\;\mbf{9},\;\mbf{26},\;\mbf{8},\;\mbf{19},\;\mbf{14},\;\mbf{2}\},\qquad\{\mbf{3},\;\mbf{20},\;\mbf{31},\;\mbf{24},\;\mbf{34},\;\mbf{17},\;\mbf{32},\;\mbf{21},\;\mbf{4}\}\\
\end{align*}

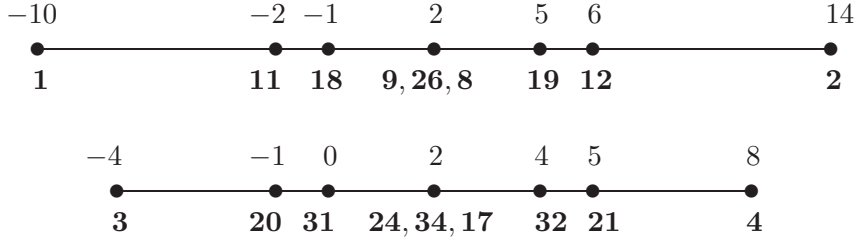
\begin{center}
\begin{picture}(300,50)
\put(0,30){\line(1,0){300}}
\put(0,30){\circle*{5}}
\put(90,30){\circle*{5}}
\put(110,30){\circle*{5}}
\put(150,30){\circle*{5}}
\put(190,30){\circle*{5}}
\put(210,30){\circle*{5}}
\put(300,30){\circle*{5}}
\put(-12,40){$-10$}
\put(80,40){$-2$}
\put(100,40){$-1$}
\put(148,40){$2$}
\put(188,40){$5$}
\put(208,40){$6$}
\put(298,40){$14$}
\put(-2,15){$\mbf{1}$}
\put(80,15){$\mbf{11}$}
\put(103,15){$\mbf{18}$}
\put(130,15){$\mbf{9},\mbf{26},\mbf{8}$}
\put(185,15){$\mbf{19}$}
\put(205,15){$\mbf{12}$}
\put(298,15){$\mbf{2}$}
\end{picture}
\end{center}

\begin{center}
\begin{picture}(240,50)
\put(0,30){\line(1,0){240}}
\put(0,30){\circle*{5}}
\put(60,30){\circle*{5}}
\put(80,30){\circle*{5}}
\put(120,30){\circle*{5}}
\put(160,30){\circle*{5}}
\put(180,30){\circle*{5}}
\put(240,30){\circle*{5}}
\put(-12,40){$-4$}
\put(50,40){$-1$}
\put(78,40){$0$}
\put(118,40){$2$}
\put(158,40){$4$}
\put(178,40){$5$}
\put(238,40){$8$}
\put(-2,15){$\mbf{3}$}
\put(50,15){$\mbf{20}$}
\put(70,15){$\mbf{31}$}
\put(95,15){$\mbf{24},\mbf{34},\mbf{17}$}
\put(158,15){$\mbf{32}$}
\put(178,15){$\mbf{21}$}
\put(238,15){$\mbf{4}$}
\end{picture}
\end{center}

\begin{align*}
&c=\frac{1}{3}\qquad\{\mbf{1},\;\mbf{18},\;\mbf{27},\;\mbf{15},\;\mbf{33},\;\mbf{9},\;\mbf{28},\;\mbf{19},\;\mbf{2}\},\qquad\{\mbf{3},\;\mbf{5},\;\mbf{20},\;\mbf{16},\;\mbf{17},\;\mbf{10},\;\mbf{21},\;\mbf{6},\;\mbf{4}\}  \\
\end{align*}

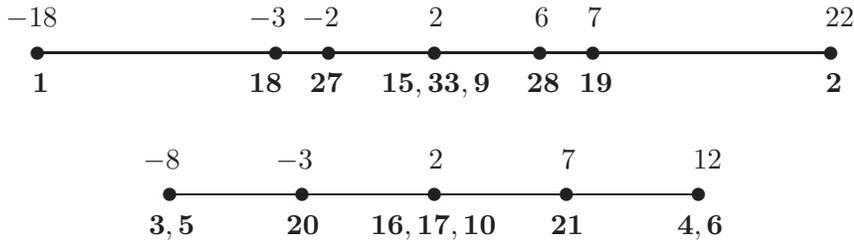
\begin{center}
\begin{picture}(300,50)
\put(0,30){\line(1,0){300}}
\put(0,30){\circle*{5}}
\put(90,30){\circle*{5}}
\put(110,30){\circle*{5}}
\put(150,30){\circle*{5}}
\put(190,30){\circle*{5}}
\put(210,30){\circle*{5}}
\put(300,30){\circle*{5}}
\put(-12,40){$-18$}
\put(80,40){$-3$}
\put(100,40){$-2$}
\put(148,40){$2$}
\put(188,40){$6$}
\put(208,40){$7$}
\put(298,40){$22$}
\put(-2,15){$\mbf{1}$}
\put(80,15){$\mbf{18}$}
\put(103,15){$\mbf{27}$}
\put(130,15){$\mbf{15},\mbf{33},\mbf{9}$}
\put(185,15){$\mbf{28}$}
\put(205,15){$\mbf{19}$}
\put(298,15){$\mbf{2}$}
\end{picture}
\end{center}

\begin{center}
\begin{picture}(200,50)
\put(0,30){\line(1,0){200}}
\put(0,30){\circle*{5}}
\put(50,30){\circle*{5}}
\put(100,30){\circle*{5}}
\put(150,30){\circle*{5}}
\put(200,30){\circle*{5}}
\put(-10,40){$-8$}
\put(40,40){$-3$}
\put(98,40){$2$}
\put(148,40){$7$}
\put(198,40){$12$}
\put(-8,15){$\mbf{3},\mbf{5}$}
\put(44,15){$\mbf{20}$}
\put(76,15){$\mbf{16},\mbf{17},\mbf{10}$}
\put(144,15){$\mbf{21}$}
\put(192,15){$\mbf{4},\mbf{6}$}
\end{picture}
\end{center}

The argument for how to complete the matrix $([\el(\tau):\m(\sigma)])$ above from the corresponding decomposition matrix of the Hecke algebra, which contains all but the first, second, third, and sixth columns, is the same in all four cases up to relabeling the rows and columns. We give the argument in the abstract in a way that applies to all four of the cases at hand. Label the reps $\tau_1,...,\tau_9$. Then $\tau_k=\tau_{9+1-k}'$ for $k=1,2,3$, while $\tau_4,\;\tau_5,\;$ and $\tau_6$ are sign-invariant. We have $\tau_n\geq_c\tau_{n-1}$ for each $n=2,..,9$, and so it follows from $[\m(\tau_8):\el(\tau_9)]=0$ that $\dim\Hom(\m(\tau_9),\m(\tau_8))=0$. Since $\el(\tau_7)=\m(\tau_7)-\m(\tau_8)-\m(\tau_9)$ but $\tau_8$ and $\tau_9$ are not linked, it must be that $\dim\Hom(\m(\tau_8),\m(\tau_7))=\dim\Hom(\m(\tau_9),\m(\tau_7))=1$. Then by (dim Hom), $\dim\Hom(\m(\tau_2),\m(\tau_1))=0$ and $\dim\Hom(\m(\tau_3),\m(\tau_1)=\dim\Hom(\m(\tau_3),\m(\tau_2))=1$. This gives the first three columns of the decomposition matrix, since for the first step of a composition series the multiplicity $[\m(\tau):\el(\sigma)]$ and the dimension $\dim\Hom(\m(\sigma),\m(\tau))$ are the same. As for the column labeled $\tau_6$: $[\m(\tau_4):\el(\tau_6)]=[\m(\tau_5):\el(\tau_6)]=0$ since $h_c(\tau_k)=2$ for $k=4,5,6$ and reps with the same $h_c$-weight can't be linked. From the decomposition $\el(\tau_6)=\m(\tau_6)-\m(\tau_7)+\m(\tau_8)+\m(\tau_9),$ we see that $1=\dim\Hom(\m(\tau_7),\m(\tau_6))=\dim\Hom(\m(\tau_6),\m(\tau_3))$. Since there are no reps with $h_c$-weight intermediate between $\tau_3$ and $\tau_6$, it must be that $[\m(\tau_3):\el(\tau_6)]=1$. For the top two entries of column $\tau_6$, we look instead to the decompositions of $\el(\tau_1)$ and $\el(\tau_2)$ into Vermas. We know so far that $\el(\tau_1)=\m(\tau_1)-\m(\tau_3)+\m(\tau_5)...$ and $\el(\tau_2)=\m(\tau_2)-\m(\tau_3)+\m(\tau_4)$. In all four cases of blocks and parameters, it's necessary that $[\el(\tau_i):\m(\tau_6]\geq1$ for $i=1,2$ in order to have $\dim\el(\tau_i)[2]\geq\dim\el(\tau_i)[-2]$, which must be the case since $\el(\tau_i)$ is an $\s$-representation. On the other hand, $\dim\Hom(\m(\tau_6),\m(\tau_3))=1$ implies that $[\el(\tau_i):\m(\tau_6]\leq1$, since $c_6\m(\tau_6)$ occurring next in the decomposition of $\el(\tau_i)=\m(\tau_i)-\m(\tau_3)....$ with $c_6>0$ would have to correspond to $c_6$ copies of $\tau_6$ generating $H_c$-subreps inside $\m(\tau_3)$. Thus $[\el(\tau_i):\m(\tau_6]=1$ for $i=1,2$, and it follows that $[\m(\tau_i):\el(\tau_6)]=0$. This completes the decomposition matrix $([\m(\tau):\el(\sigma)])$.

\underline{$c=\frac{1}{5}$} There are four finite-dimensional irreps with the following decompositions, characters, and dimensions:
\begin{align*}
\el(\mbf{1})&=\m(\mbf{1})-\m(\mbf{18})+\m(\mbf{26})+\m(\mbf{8})-\m(\mbf{19})+\m(\mbf{2})\\
\chi_{\el(\mbf{1})}(w,t)&=\chi_\mbf{1}(t^{-10}+t^{10})+\chi_\mbf{3}(t^{-9}+t^9)+(\chi_\mbf{1}+\chi_\mbf{11})(t^{-8}+t^8)+(\chi_\mbf{3}+\chi_\mbf{18})(t^{-7}+t^7)\\&\qquad\qquad+(\chi_\mbf{1}+\chi_\mbf{11}+\chi_\mbf{27})(t^{-6}+t^6)+(\chi_\mbf{3}+\chi_\mbf{18}+\chi_\mbf{31})(t^{-5}+t^5)\\
&\qquad\qquad+(\chi_\mbf{1}+\chi_\mbf{11}+\chi_\mbf{13}+\chi_\mbf{20}+\chi_\mbf{26}+\chi_\mbf{27})(t^{-4}+t^4)\\
&\qquad\qquad+(\chi_\mbf{3}+\chi_\mbf{5}+\chi_\mbf{18}+\chi_\mbf{24}+2\chi_\mbf{31})(t^{-3}+t^3)\\
&\qquad\qquad+(\chi_\mbf{1}+\chi_\mbf{11}+\chi_\mbf{13}+2\chi_\mbf{20}+\chi_\mbf{26}+2\chi_\mbf{27}+\chi_\mbf{33})(t^{-2}+t^2)\\
&\qquad\qquad+(\chi_\mbf{3}+\chi_\mbf{5}+\chi_\mbf{18}+\chi_\mbf{24}+3\chi_\mbf{31}+\chi_\mbf{34})(t^{-1}+t)\\
&\qquad\qquad+\chi_\mbf{1}+\chi_\mbf{11}+2\chi_\mbf{13}+2\chi_\mbf{20}+\chi_\mbf{22}+\chi_\mbf{26}+2\chi_\mbf{27}+\chi_\mbf{30}+\chi_\mbf{33}\\
\chi_{\el(\mbf{1})}(t)&=t^{-10}+t^{10}+4(t^{-9}+t^9)+10(t^{-8}+t^8)+20(t^{-7}+t^7)+35(t^{-6}+t^6)+56(t^{-5}+t^5)
\\&\qquad+84(t^{-4}+t^4)+120(t^{-3}+t^3)+165(t^{-2}+t^2)+204(t^{-1}+t)+222\\
\dim\el(\mbf{1})&=1620\\
\\
\el(\mathbf{11})&=\m(\mathbf{11})-\m(\mathbf{18})+\m(\mathbf{8})+\m(\mathbf{9})-\m(\mathbf{19})+\m(\mathbf{12})\\
\chi_{\el(\mathbf{11})}&=\chi_{\mathbf{11}}(t^{-2}+t^2)+(\chi_{\mathbf{3}}+\chi_{\mathbf{16}})(t^{-1}+t)+\chi_{\mathbf{1}}+\chi_{\mathbf{7}}+\chi_{\mathbf{11}}+\chi_{\mathbf{15}}\\
\chi_{\el(\mathbf{11})}(t)&=9t^{-2}+20t^{-1}+26+20t+9t^2\\
\mathrm{dim}\;\el(\mathbf{11})&=84\\
\\
\el(\mbf{3})&=\m(\mbf{3})-\m(\mbf{31})+\m(\mbf{34})+\m(\mbf{17})-\m(\mbf{32})+\m(\mbf{4})\\
\chi_{\el(\mbf{3})}(w,t)&=\chi_\mbf{3}(t^{-4}+t^4)+(\chi_\mbf{1}+\chi_\mbf{7}+\chi_\mbf{11})(t^{-3}+t^3)+(2\chi_\mbf{3}+\chi_\mbf{16}+\chi_\mbf{18})(t^{-2}+t^2)\\
&\qquad+(\chi_\mbf{1}+\chi_\mbf{7}+2\chi_\mbf{11}+\chi_\mbf{27}+\chi_\mbf{29})(t^{-1}+t)+2\chi_\mbf{3}+\chi_\mbf{16}+2\chi_\mbf{18}+\chi_\mbf{34}\\
\chi_{\el(\mbf{3})}(t)&=4(t^{-4}+t^4)+16(t^{-3}+t^3)+40(t^{-2}+t^2)+80(t^{-1}+t)+104\\
\dim\el(\mbf{3})&=384\\
\\
\el(\mathbf{20})&=\m(\mathbf{20})-\m(\mathbf{31})+\m(\mathbf{17})+\m(\mathbf{24})-\m(\mathbf{32})+\m(\mathbf{21})\\
\chi_{\el(\mathbf{20})}(w,t)&=\chi_{\mathbf{20}}(t^{-1}+t)+\chi_\mathbf{5}+\chi_\mathbf{24}\\
\chi_{\el(\mathbf{20})}(t)&=16t^{-1}+28+16t\\
\mathrm{dim}\;\el(\mathbf{20})&=60
\end{align*}

The reps with bulleted rows in the decomposition matrix are the reps which don't have full support but aren't finite-dimensional. These are $\el(\mbf{18}),\;\el(\mbf{8}),\;\el(\mbf{31}),\;\el(\mbf{17})$. Computing their graded characters we find:
$$\dim\Supp\el(\tau)=2,\qquad\tau\in\{\mbf{18},\;\mbf{8},\;\mbf{31},\;\mbf{17}\}$$

Additionally, in the case of $c=1/5$ there are two more self-dual blocks of defect $1$ \cite{GP}:
\begin{align*}
&\el(\mbf{13})=\m(\mbf{13})-\m(\mbf{22})+\m(\mbf{14})\\
&\chi_{\el(\mbf{13})}(t)=9t^{-2}\frac{t^6 + 2t^5 + 3t^4 + 4t^3 + 3t^2 + 2t + 1}{(1-t)^2}\\
&\dim\Supp\el(\mbf{13})=2\\
\\
&\el(\mbf{5})=\m(\mbf{5})-\m(\mbf{10})+\m(\mbf{6})\\
&\chi_{\el(\mbf{5})}(t)=4t^{-4}\frac{t^{10} + 2t^9 + 3t^8 + 4t^7 + 5t^6 + 6t^5 + 5t^4 + 4t^3 + 3t^2 + 2t + 1}{(1-t)^2}\\
&\dim\Supp\el(\mbf{5})=2\\
\end{align*}

\newgeometry{top=2cm,bottom=2cm}
\underline{$c=\frac{1}{3}$} There are four finite-dimensional irreps with the following decompositions, characters, and dimensions:
\begin{align*}
\el(\mbf{1})&=\m(\mbf{1})-\m(\mbf{27})+\m(\mbf{33})+\m(\mbf{9})-\m(\mbf{28})+\m(\mbf{2})\\
%\chi_{\el(\mbf{1})}(w,t)&=\\
\chi_{\el(\mbf{1})}(t)&=t^{-18}+t^{18} + 4(t^{-17}+t^{17}) + 10(t^{-16}+t^{16}) + 20(t^{-15}+t^{15}) + 35(t^{-14}+t^{14})\\
&\qquad + 56(t^{-13}+t^{13}) + 84(t^{-12}+t^{12}) + 120(t^{-11}+t^{11}) + 165(t^{-10}+t^{10})\\
&\qquad + 220(t^{-9}+t^9) + 286(t^{-8}+t^8) + 364(t^{-7}+t^7) + 455(t^{-6}+t^6) + 560(t^{-5}+t^5) \\
&\qquad+ 680(t^{-4}+t^4) + 816(t^{-3}+t^3) + 944(t^{-2}+t^2) + 1040(t^{-1}+t) + 1080\\
\dim\el(\mbf{1})&=12,800\\
\\
\el(\mathbf{18})&=\m(\mbf{18})-\m(\mbf{27})+\m(\mbf{15})+\m(\mbf{9})-\m(\mbf{28})+\m(\mbf{19})\\
\chi_{\el(\mathbf{18})}(w,t)&=\chi_\mathbf{18}(t^{-3}+t^3)+(\chi_\mathbf{11}+\chi_\mathbf{29})(t^{-2}+t^2)+(\chi_\mathbf{3}+\chi_\mathbf{16}+\chi_\mathbf{18}+\chi_\mathbf{23})(t^{-1}+t)\\&\qquad+\chi_\mathbf{1}+\chi_\mathbf{7}+\chi_\mathbf{8}+\chi_\mathbf{9}+\chi_\mathbf{11}+\chi_\mathbf{15}+\chi_\mathbf{29}\\
\chi_{\el(\mathbf{18}}(t)&=16(t^{-3}+t^3)+39(t^{-2}+t^2)+60(t^{-1}+t)+70\\
\dim\el(\mathbf{18})&=300\\
\\
\el(\mbf{3})&=\m(\mbf{3})-\m(\mbf{20})+\m(\mbf{10})+\m(\mbf{17})-\m(\mbf{21})+\m(\mbf{4})\\
\chi_{\el(\mbf{3})}(w,t)&=\chi_\mbf{3}(t^{-8}+t^8)+(\chi_\mbf{1}+\chi_\mbf{7}+\chi_\mbf{11})(t^{-7}+t^7)+(2\chi_\mbf{3}+\chi_\mbf{16}+\chi_\mbf{18})(t^{-6}+t^6)\\&\qquad+(\chi_\mbf{1}+\chi_\mbf{7}+2\chi_\mbf{11}+\chi_\mbf{27}+\chi_\mbf{29})(t^{-5}+t^5)+(2\chi_\mbf{3}+\chi_\mbf{16}+2\chi_\mbf{18}+\chi_\mbf{31}+\chi_\mbf{34})(t^{-4}+t^4)\\
&\qquad+(\chi_\mbf{1}+\chi_\mbf{7}+2\chi_\mbf{11}+\chi_\mbf{13}+\chi_\mbf{26}+2\chi_\mbf{27}+\chi_\mbf{29}+\chi_\mbf{30}+\chi_\mbf{33})(t^{-3}+t^3)\\
&\qquad+(2\chi_\mbf{3}+\chi_\mbf{16}+2\chi_\mbf{18}+\chi_\mbf{24}+3\chi_\mbf{31}+\chi_\mbf{32}+\chi_\mbf{34})(t^{-2}+t^2)\\
&\qquad+(\chi_\mbf{1}+\chi_\mbf{7}+2\chi_\mbf{11}+\chi_\mbf{13}+\chi_\mbf{20}+\chi_\mbf{21}+2\chi_\mbf{26}+2\chi_\mbf{27}+\chi_\mbf{29}+\chi_\mbf{30}+2\chi_\mbf{33})(t^{-1}+t)\\
&\qquad+2\chi_\mbf{3}+\chi_\mbf{5}+\chi_\mbf{16}+\chi_\mbf{17}+2\chi_\mbf{18}+3\chi_\mbf{24}+3\chi_\mbf{31}+\chi_\mbf{32}+\chi_\mbf{34}\\
\chi_{\el(\mbf{3})}(t)&=4(t^{-8}+t^8) + 16(t^{-7}+t^7) + 40(t^{-6}+t^6) + 80(t^{-5}+t^5) + 140(t^{-4}+t^4) \\&\qquad+ 208(t^{-3}+t^3) + 272(t^{-2}+t^2) + 320(t^{-1}+t) + 340\\
\dim\el(\mbf{3})&=2500\\
\\
\el(\mbf{5})&=\m(\mbf{5})-\m(\mbf{20})+\m(\mbf{10})+\m(\mbf{16})-\m(\mbf{21})+\m(\mbf{6})\\
\chi_{\el(\mbf{5})}(w,t)&=\chi_\mbf{5}(t^{-8}+t^8)+\chi_\mbf{20}(t^{-7}+t^7)+(\chi_\mbf{5}+\chi_\mbf{31})(t^{-6}+t^6)+(\chi_\mbf{13}+\chi_\mbf{20}+\chi_\mbf{27}+\chi_\mbf{30})(t^{-5}+t^5)\\
&\qquad+(\chi_\mbf{5}+\chi_\mbf{18}+2\chi_\mbf{31}+\chi_\mbf{34})(t^{-4}+t^4)\\
&\qquad+(\chi_\mbf{11}+\chi_\mbf{13}+\chi_\mbf{20}+\chi_\mbf{25}+2\chi_\mbf{27}+\chi_\mbf{29}+\chi_\mbf{30}+\chi_\mbf{33})(t^{-3}+t^3)\\
&\qquad+(\chi_\mbf{3}+\chi_\mbf{5}+\chi_\mbf{10}+\chi_\mbf{16}+\chi_\mbf{17}+2\chi_\mbf{18}+\chi_\mbf{23}+2\chi_\mbf{31}+2\chi_\mbf{34})(t^{-2}+t^2)\\
&\qquad+(\chi_\mbf{1}+\chi_\mbf{7}+\chi_\mbf{8}+\chi_\mbf{9}+2\chi_\mbf{11}+\chi_\mbf{13}+\chi_\mbf{15}+\chi_\mbf{20}+\chi_\mbf{22}+2\chi_\mbf{25}\\&\qquad\qquad\qquad+2\chi_\mbf{27}+2\chi_\mbf{29}+\chi_\mbf{30}+\chi_\mbf{33})(t^{-1}+t)\\
&\qquad+2\chi_\mbf{3}+\chi_\mbf{5}+\chi_\mbf{10}+2\chi_\mbf{16}+\chi_\mbf{17}+2\chi_\mbf{18}+3\chi_\mbf{23}+2\chi_\mbf{31}+2\chi_\mbf{34}\\
\chi_{\el(\mbf{5})}(t)&=4(t^{-8}+t^8) + 16(t^{-7}+t^7) + 40(t^{-6}+t^6) + 80(t^{-5}+t^5) + 140(t^{-4}+t^4) \\&\qquad+ 208(t^{-3}+t^3) + 272(t^{-2}+t^2) + 320(t^{-1}+t) + 340\\
\dim\el(\mbf{5})&=2500\\
\end{align*}
\restoregeometry

The reps with bulleted rows in the decomposition matrix are the reps which don't have full support but aren't finite-dimensional. These are $\el(\mbf{27}),\;\el(\mbf{9}),\;\el(\mbf{20}),\;\el(\mbf{10})$. Computing their graded characters we find:
$$\dim\Supp\el(\tau)=2,\qquad\tau\in\{\mbf{27},\;\mbf{9},\;\mbf{20},\;\mbf{10}\}$$

\subsection{\Large{$\mbf{c=\frac{1}{4}}$}}

\begin{center}
\begin{picture}(400,50)
\put(0,30){\line(1,0){400}}
\put(0,30){\circle*{5}}
\put(120,30){\circle*{5}}
\put(152,30){\circle*{5}}
\put(200,30){\circle*{5}}
\put(248,30){\circle*{5}}
\put(280,30){\circle*{5}}
\put(400,30){\circle*{5}}
\put(-10,40){$-13$}
\put(110,40){$-3$}
\put(150,40){$-1$}
\put(198,40){$2$}
\put(246,40){$5$}
\put(278,40){$7$}
\put(394,40){$17$}
\put(-2,15){$\mbf{1}$}
\put(103,15){$\mbf{11},\mbf{13}$}
\put(146,15){$\mbf{27}$}
\put(176,15){$\mbf{23},\mbf{24},\mbf{10}$}
\put(240,15){$\mbf{28}$}
\put(265,15){$\mbf{14},\mbf{12}$}
\put(398,15){$\mbf{2}$}
\end{picture}
\end{center}

\[
\begin{blockarray}{ccccccccccccc}
&& \mbf{1} & \mbf{11} & \mbf{13} & \mbf{27} & \mbf{23} & \mbf{24} & \mbf{10} & \mbf{28} & \mbf{14} & \mbf{12} & \mbf{2}\\
\begin{block}{cc(ccccccccccc)}
\star & \mbf{1} & 1 & 1 & 1 & 1 & \cdot & \cdot & 1 & \cdot & \cdot & \cdot &\cdot \\
\star & \mbf{11} & \cdot & 1 & \cdot & 1 & \cdot & 1 & \cdot & \cdot & \cdot & \cdot &\cdot\\
\star & \mbf{13} & \cdot & \cdot & 1 & 1 & 1 & \cdot & \cdot & \cdot & \cdot & \cdot & \cdot\\
& \mbf{27} & \cdot & \cdot & \cdot & 1 & 1 & 1 & 1 & 1 & \cdot & \cdot & \cdot \\
& \mbf{23} &  \cdot & \cdot & \cdot & \cdot & 1 & \cdot & \cdot & 1 & 1 & \cdot & \cdot\\
& \mbf{24} & \cdot & \cdot & \cdot & \cdot & \cdot & 1 & \cdot & 1 & \cdot & 1 & \cdot \\
\mathtt{(1)} & \mbf{10} & \cdot & \cdot & \cdot & \cdot & \cdot & \cdot & 1 & 1 & \cdot & \cdot & \cdot\\
& \mbf{28} & \cdot & \cdot & \cdot & \cdot & \cdot & \cdot & \cdot & 1 & 1 & 1 & 1 \\
& \mbf{14} & \cdot & \cdot & \cdot & \cdot & \cdot & \cdot & \cdot & \cdot & 1 & \cdot & 1 \\
& \mbf{12} & \cdot & \cdot & \cdot & \cdot & \cdot & \cdot & \cdot & \cdot & \cdot & 1 & 1\\
& \mbf{2} & \cdot & \cdot & \cdot & \cdot & \cdot & \cdot & \cdot & \cdot & \cdot & \cdot & 1\\
\end{block}
\end{blockarray}
\]
\\

\[
\begin{blockarray}{ccccccccccccc}
&& \mbf{1} & \mbf{11} & \mbf{13} & \mbf{27} & \mbf{23} & \mbf{24} & \mbf{10} & \mbf{28} & \mbf{14} & \mbf{12} & \mbf{2}\\
\begin{block}{cc(ccccccccccc)}
\star & \mbf{1} & 1 & -1 & -1 & 1 & \cdot & \cdot & -2 & 1 & -1 & -1 & 1 \\
\star & \mbf{11} & \cdot & 1 & \cdot & -1 & 1 & \cdot & 1 & -1 & \cdot & 1 &\cdot\\
\star & \mbf{13} & \cdot & \cdot & 1 & -1 & \cdot & 1 & 1 & -1 & 1 & \cdot & \cdot\\
& \mbf{27} & \cdot & \cdot & \cdot & 1 & -1 & -1 & -1 & 2 & -1 & -1 & \cdot \\
& \mbf{23} &  \cdot & \cdot & \cdot & \cdot & 1 & \cdot & \cdot & -1 & \cdot & 1 & \cdot\\
& \mbf{24} & \cdot & \cdot & \cdot & \cdot & \cdot & 1 & \cdot & -1 & 1 & \cdot & \cdot \\
\mathtt{(1)} & \mbf{10} & \cdot & \cdot & \cdot & \cdot & \cdot & \cdot & 1 & -1 & 1 & 1 & -1\\
& \mbf{28} & \cdot & \cdot & \cdot & \cdot & \cdot & \cdot & \cdot & 1 & -1 & -1 & 1 \\
& \mbf{14} & \cdot & \cdot & \cdot & \cdot & \cdot & \cdot & \cdot & \cdot & 1 & \cdot & -1 \\
& \mbf{12} & \cdot & \cdot & \cdot & \cdot & \cdot & \cdot & \cdot & \cdot & \cdot & 1 & -1\\
& \mbf{2} & \cdot & \cdot & \cdot & \cdot & \cdot & \cdot & \cdot & \cdot & \cdot & \cdot & 1\\
\end{block}
\end{blockarray}
\]
\\

All but columns $\mbf{1},\mbf{11},\mbf{13}$, and $\mbf{10}$ of the matrix $([\m(\tau):\el(\sigma)])$ appear in the Hecke algebra decomposition matrix. Column $\mbf{1}$ is immediate, and columns $\mbf{11},\;\mbf{13}$ follow from Lemmas (RR) and (dim Hom). As for column $\mbf{10}$: $[\m(\mbf{23}):\el(\mbf{10})]=[\m(\mbf{24}):\el(\mbf{10})]=0$ since all three of these reps have the same $h_c$-weight. By (RR) and (dim Hom), $[\m(\mbf{27}):\el(\mbf{10})]=1$. To determine the remaining entries of this column (in the first three rows), first observe that $\el(\mbf{1})$ must be finite-dimensional since $4$ is an elliptic number of $H_4$. Since we know the decomposition into Vermas up through $h_c$-weight $0$, we can the graded character, a polynomial, by calculating the dimensions of the graded pieces; multiplying this polynomial by $\frac{(1-t)^4}{(1-t)^4}$ we have:
$$\chi_{\el(\mbf{1})}(t)=\frac{t^{-13}-18t^{-3}+25t^{-1}-16t^2+25t^5-18t^7+t^{17}}{(1-t)^4}$$
and so $-2$ must be the coefficient of $\mbf{10}$ in the decomposition of $\el(\mbf{1})$ as neither $\m(\mbf{23})$ nor $\m(\mbf{24})$ occurs.  Likewise we may determine the decomposition of $\el(\mbf{11})$ and $\el(\mbf{13})$ into Vermas by looking at the decomposition up to $h_c$-weight $2$ and applying (E) to conclude that $\el(\mbf{11})$ and $\el(\mbf{13})$ are finite-dimensional; a look at the graded character of these modules then shows that the coefficient of $\m(\mbf{10})$ must be $1$ in both $\el(\mbf{11})$ and $\el(\mbf{13})$. Now by taking the dot product of these top three rows of $([\el(\tau):\m(\sigma)])$ with the column $\mbf{10}$ of $([\m(\tau):\el(\sigma)]$ to be $0$, we find that $[\m(\mbf{13}):\el(\mbf{10})]=[\m(\mbf{11}):\el(\mbf{10})]=0$ and $[\m(\mbf{1}):\el(\mbf{10})]=1$. This completes the decomposition matrix.

Decompositions, characters, and dimensions of the finite-dimensional irreps:
\begin{align*}
&L(\mathbf{1})=M(\mathbf{1})-M(\mathbf{11})-M(\mathbf{13})+M(\mathbf{27})-2M(\mathbf{10})+M(\mathbf{28})-M(\mathbf{14})-M(\mathbf{12})+M(\mathbf{2})\\
&\chi_{L(\mathbf{1})}(t)=t^{-13}+t^{13}+4(t^{-12}+t^{12})+10(t^{-11}+t^{11})+20(t^{-10}+t^{10})+35(t^{-9}+t^9)\\
&\qquad\qquad\qquad+56(t^{-8}+t^8)+84(t^{-7}+t^7)+120(t^{-6}+t^6)+165(t^{-5}+t^5)\\
&\qquad\qquad\qquad+220(t^{-4}+t^4)+268(t^{-3}+t^3)+292(t^{-2}+t^2)+300(t^{-1}+t)+300\\
&\mathrm{dim}\;L(\mathbf{1})=3450\\
\\
&L(\mathbf{11})=M(\mathbf{11})-M(\mathbf{27})+M(\mathbf{10})+M(\mathbf{23})-M(\mathbf{28})+M(\mathbf{12})\\
&L(\mathbf{13})=M(\mathbf{13})-M(\mathbf{27})+M(\mathbf{10})+M(\mathbf{24})-M(\mathbf{28})+M(\mathbf{14})\\
&\chi_{L(\mathbf{11})}(t)=\chi_{L(\mathbf{13})}(t)=9t^{-3}+36t^{-2}+65t^{-1}+80+65t+36t^2+9t^3\\
&\mathrm{dim}\;\chi_{L(\mathbf{11})}=\mathrm{dim}\;\chi_{L(\mathbf{13})}=300
\end{align*}

In addition, the irrep $\el(\mbf{10})$ has less than full support, and the dimension of its support is $1$:
\begin{align*}
\chi_{\el(\mbf{10})}(t)&=\frac{8t^2-25t^5+18t^7-t^{17}}{(1-t)^4}\\
&=\frac{8t^2+24t^3+48t^4+\sum_{n=0}^9{11-n\choose2}t^{5+n}}{(1-t)}\\
\end{align*}
 
\subsection{\Large{$\mathbf{c=\frac{1}{2}}$}} The $h_c$-weight line of the principal block is as follows.

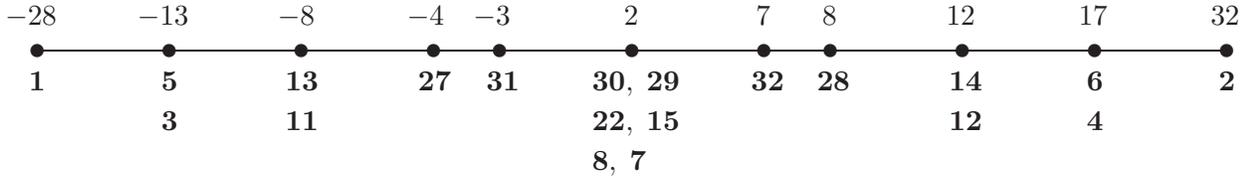
\begin{center}
\begin{picture}(450,90)
\put(0,70){\line(1,0){450}}
\put(0,70){\circle*{5}}
\put(50,70){\circle*{5}}
\put(100,70){\circle*{5}}
\put(150,70){\circle*{5}}
\put(175,70){\circle*{5}}
\put(225,70){\circle*{5}}
\put(275,70){\circle*{5}}
\put(300,70){\circle*{5}}
\put(350,70){\circle*{5}}
\put(400,70){\circle*{5}}
\put(450,70){\circle*{5}}
\put(-12,80){$-28$}
\put(38,80){$-13$}
\put(91,80){$-8$}
\put(140,80){$-4$}
\put(165,80){$-3$}
\put(222,80){$2$}
\put(272,80){$7$}
\put(297,80){$8$}
\put(344,80){$12$}
\put(394,80){$17$}
\put(444,80){$32$}
\put(-3,55){$\mbf{1}$}
\put(47,55){$\mbf{5}$}
\put(47,40){$\mbf{3}$}
\put(94,55){$\mbf{13}$}
\put(94,40){$\mbf{11}$}
\put(144,55){$\mbf{27}$}
\put(170,55){$\mbf{31}$}
\put(210,55){$\mbf{30},\;\mbf{29}$}
\put(210,40){$\mbf{22},\;\mbf{15}$}
\put(210,25){$\mbf{8},\;\mbf{7}$}
\put(270,55){$\mbf{32}$}
\put(295,55){$\mbf{28}$}
\put(345,55){$\mbf{14}$}
\put(345,40){$\mbf{12}$}
\put(397,55){$\mbf{6}$}
\put(397,40){$\mbf{4}$}
\put(447,55){$\mbf{2}$}
\end{picture}
\end{center}

It is basically possible to deduce that $\el(\bf{27})$, $\el(\bf{3})$, and $\el(\bf{5})$ are finite-dimensional by direct arguments. Finite-dimensionality of $\el(\bf{27})$ follows from (E) and Lemma 3.1 but one should check that $E$ does in fact kill something in degree $0$. 
\begin{align*}
&\el(\mbf{27})=\m(\mbf{27})-\m(\mbf{31})+\m(\mbf{7})+\m(\mbf{8})+\m(\mbf{15})-\m(\mbf{32})+\m(\mbf{28})\\
&\chi_{\el(\mbf{27})}(t)=25(t^{-4}+t^4)+64(t^{-3}+t^3)+106(t^{-2}+t^2)+140(t^{-1}+t)+155\\
&\dim\el(\mbf{27})=825\\
\end{align*}

$\el(\bf{3})$ and $\el(\bf{5})$ require some work: for $\el(\bf{3})$, induce $\el(\bf{3_+})$ from $\oh_{\frac{1}{2}}(H_3)$:
$$\Ind\el(\bf{3_-})=\m(\bf{3})-\m(\bf{13})-\m(\bf{27})+\m(\bf{31})-\m(\mbf{32})+\m(\mbf{28})+\m(\mbf{14})-\m(\mbf{4})$$
Thus $\dim\Hom(\m(\mbf{13}),\m(\mbf{3}))\geq1$. The decomposition of $S^5\hstar\otimes\mbf{3}$ into irreducibles implies that $\dim\Hom(\m(\mbf{13}),\m(\mbf{3}))\leq1$. Since $\mbf{11}^{\oplus2}\subset S^5\hstar\otimes\mbf{3}$ it is not immediately obvious that $\dim\Hom(\m(\mbf{11}),\m(\mbf{3}))=0$. However, considering the decomposition of $S^7\hstar\otimes\mbf{3}$ into irreducibles we see that it does not contain $\mbf{15}$ as a subrep, while $S^2\hstar\otimes\mbf{11}$ does -- and as $$\el(\mbf{3})[-6]=S^7\hstar\otimes\mbf{3}-S^2\hstar\otimes\mbf{13}-c_{11}S^2\hstar\otimes\mbf{11}$$ can't contain $\mbf{15}$ with negative multiplicity, $-\m(\mbf{11})$ must not appear in the composition series. This now implies that $\dim\Hom(\m(\mathbf{27}),\m(\mathbf{3}))\geq1$. Again we must argue that it equals $1$; this can be done by considering the decomposition of $$\el(\mbf{3})[-2]=S^{11}\hstar\otimes\mbf{3}-S^6\hstar\otimes\mbf{13}-c_{27}S^2\hstar\otimes\mbf{27}$$ into $H_4$-irreps. The irrep $\mbf{22}$ appears with multiplicity $3$ in $S^{11}\hstar\otimes\mbf{3}$ and with multiplicity $2$ in $S^6\hstar\otimes\mbf{13}$, implying $c_{27}\leq 1$. Thus $\el(\bf{3})=\m(\bf{3})-\m(\bf{13})-\m(\bf{27})+\m(\bf{31})...$ where the rest of the decomposition consists of Vermas $\m$ with $h_c(\m)\geq2$. And so from this beginning of the decomposition we may calculate $\dim\el(\mbf{3})[-1]=\dim\el(\mbf{3})[1]$. So long as the operator $\mathrm{E}$ kills something in $\el(\mbf{3})[0]$, then by (E), $\el(\bf{3})$ is finite-dimensional. To check this, one should decompose $\el(\mbf{3})[0]$ as a $H_4$-representation; hopefully it would turn out that there is a $6$-dimensional representation which is not $\mbf{7}$ which will not appear in $\el(\mbf{3})[2]$.

To determine the full decomposition of $\el(\mbf{3})$ into Vermas, write 
$\el(\mbf{3})=\m(\mbf{3})-\m(\mbf{13})-\m(\mbf{27})+\m(\mbf{31})+M+\m(\mbf{32})-\m(\mbf{28})-\m(\mbf{14})+\m(\mbf{4})$ where $M$ is some linear combination of Vermas whose lowest weight $H_4$-irreps have $h_c$-weight $2$. Now restrict this expression to $\oh_{\frac{1}{2}}(H_3)$ and note that restriction of a finite-dimensional representation is always $0$. For ease on the eyes we drop the notation of Vermas and simply write the restriction of their lowest weight $H_4$-reps (it is the same rule on the level of the Grothendieck group anyways): 

\begin{align*}
0&=\Res^{H_4}_{H_3}\el(\mbf{3})\\
&=\Res^{H_4}_{H_3}\mbf{3}-\Res^{H_4}_{H_3}\mbf{13}-\Res^{H_4}_{H_3}\mbf{27}+\Res^{H_4}_{H_3}\mbf{31}+\Res^{H_4}_{H_3}M\\&\qquad\qquad+\Res^{H_4}_{H_3}\mbf{32}-\Res^{H_4}_{H_3}\mbf{28}-\Res^{H_4}_{H_3}\mbf{14}+\Res^{H_4}_{H_3}\mbf{4}\\
&=(\mbf{1_+}\oplus\mbf{3_-})-(\mbf{1_+}\oplus\mbf{\tilde{3}_-}\oplus\mbf{5_+})-(\mbf{1_+}\oplus\mbf{3_-}\oplus\mbf{\tilde{3}_-}\oplus\mbf{4_+}\oplus\mbf{4_-}\oplus\mbf{5_+}^{\oplus2})\\
&\qquad\qquad+(\mbf{1_+}\oplus\mbf{3_-}^{\oplus2}\oplus\mbf{\tilde{3}_-}^{\oplus2}\oplus\mbf{4_+}\oplus\mbf{4_-}\oplus\mbf{5_+}^{\oplus2}\oplus\mbf{5_-})+\Res^{H_4}_{H_3}M\\
&\qquad\qquad+(\mbf{1_-}\oplus\mbf{3_+}^{\oplus2}\oplus\mbf{\tilde{3}_+}^{\oplus2}\oplus\mbf{4_+}\oplus\mbf{4_-}\oplus\mbf{5_+}\oplus\mbf{5_-}^{\oplus2})\\
&\qquad\qquad-(\mbf{1_-}\oplus\mbf{3_+}\oplus\mbf{\tilde{3}_+}\oplus\mbf{4_+}\oplus\mbf{4_-}\oplus\mbf{5_-}^{\oplus2})-(\mbf{1_-}\oplus\mbf{\tilde{3}_+}\oplus\mbf{5_-})+(\mbf{1_-}\oplus\mbf{3_+})\\
&=\mbf{3_+}^{\oplus2}\oplus\mbf{3_-}^{\oplus2}+\Res^{H_4}_{H_3}M
\end{align*}

Therefore (after examining the restrictions of the possible constituents of $M$):
\begin{align*}
&\el(\mbf{3})=\m(\mbf{3})-\m(\mbf{13})-\m(\mbf{27})+\m(\mbf{31})-2\m(\mbf{7})+\m(\mbf{32})-\m(\mbf{28})-\m(\mbf{14})+\m(\mbf{4})\\
&\chi_{\el(\mbf{3})}(t)=4(t^{-13}+t^{13}) + 16(t^{-12}+t^{12}) + 40(t^{-11}+t^{11}) + 80(t^{-10}+t^{10}) + 140(t^{-9}+t^9) \\&\qquad\qquad + 215(t^{-8}+t^8) + 300(t^{-7}+t^7) + 390(t^{-6}+t^6)+ 480(t^{-5}+t^5) + 540(t^{-4}+t^4) \\&\qquad\qquad+ 576(t^{-3}+t^3) + 594(t^{-2}+t^2) + 600(t^{-1}+t^1) + 600\\
&\dim\el(\mbf{3})=8550\\
\end{align*}

Using similar arguments starting with $\Ind\el(\mbf{\tilde{3}_-})$ in place of $\Ind\el(\bf{3_-})$, one can show that $\el(\bf{5})$ is finite-dimensional as well, with the following decomposition into Vermas and the same graded character and dimension as $\el(\mbf{3})$:
\begin{align*}
&\el(\mbf{5})=\m(\mbf{5})-\m(\mbf{11})-\m(\mbf{27})+\m(\mbf{31})-2\m(\mbf{8})+\m(\mbf{32})-\m(\mbf{28})-\m(\mbf{12})+\m(\mbf{6})\\
&\chi_{\el(\mbf{5})}(t)=\chi_{\el(\mbf{3})}(t)\\
&\dim\el(\mbf{5})=8550\\
\end{align*}

When $c=\frac{1}{2}$ (dim Hom) can't be applied; nonetheless the statement may still hold. If we assume it does, we can determine the finite-dimensional irreps and their decompositions.
\assumption $\dim\Hom (\m_{\frac{1}{2}}(\sigma),\m_{\frac{1}{2}}(\tau))=\dim\Hom(\m_{\frac{1}{2}}(\tau'),\m_{\frac{1}{2}}(\sigma')).$\\

Conjecturally, if things in fact work the same way for $c=\frac{1}{2}$ as they do otherwise, then the decomposition matrix of $H_{\frac{1}{2}}(H_4)$ should be as follows. This is only conjectural, as the proof of (dim Hom), and probably the relevant statements of [CGG] as well, uses faithfulness of KZ functor which need not hold when $c=\frac{1}{2}$:

\[
\begin{blockarray}{cccccccccccccccccccccc}
&&\mbf{1}&\mbf{5}&\mbf{3}&\mbf{13}&\mbf{11}&\mbf{27}&\mbf{31}&\mbf{30}&\mbf{29}&\mbf{22}&\mbf{15}&\mbf{8}&\mbf{7}&\mbf{32}&\mbf{28}&\mbf{14}&\mbf{12}&\mbf{6}&\mbf{4}&\mbf{2}\\
\begin{block}{cc(cccccccccccccccccccc)}
\star&\mbf{1}&1 &\cdot &\cdot& \cdot& \cdot& \cdot& 1 &\cdot &\cdot &1 &\cdot &1 &1 &1& \cdot& \cdot& \cdot& \cdot& \cdot &1\\
\star&\mbf{5}&\cdot &1 &\cdot &\cdot &1 &1 &1 &\cdot &1 &\cdot &\cdot &1 &\cdot &\cdot &\cdot& \cdot& \cdot& 1 &\cdot &\cdot\\
\star&\mbf{3}&\cdot& \cdot& 1& 1 &\cdot& 1& 1& 1& \cdot& \cdot &\cdot &\cdot &1 &\cdot& \cdot& \cdot &\cdot &\cdot &1 &\cdot\\
\star&\mbf{13}&\cdot &\cdot &\cdot &1 &\cdot &\cdot &1 &1 &\cdot &\cdot &1 &1 &\cdot &1 &\cdot &1& \cdot& \cdot& 1& \cdot\\
\star&\mbf{11}&\cdot &\cdot &\cdot &\cdot &1 &\cdot &1 &\cdot &1 &\cdot &1& \cdot& 1 &1 &\cdot &\cdot &1& 1 &\cdot &\cdot\\
\star&\mbf{27}&\cdot& \cdot& \cdot& \cdot &\cdot &1 &1& 1& 1& 1& \cdot& \cdot& \cdot& 1& 1& \cdot& \cdot& 1 &1& \cdot\\
\mathtt{(1)}&\mbf{31}&\cdot &\cdot &\cdot &\cdot &\cdot& \cdot& 1 &1 &1 &1& 1& 1& 1& 3 &1 &1& 1& 1 &1 &1\\
\mathtt{(3)}&\mbf{30}&\cdot &\cdot &\cdot &\cdot& \cdot& \cdot& \cdot& 1& \cdot& \cdot& \cdot& \cdot& \cdot& 1& 1& 1& \cdot& \cdot &2& \cdot\\
\mathtt{(3)}&\mbf{29}&\cdot& \cdot& \cdot& \cdot &\cdot &\cdot &\cdot &\cdot &1 &\cdot &\cdot &\cdot &\cdot &1 &1 &\cdot &1 &2 &\cdot &\cdot\\
\mathtt{(2)}&\mbf{22}&\cdot &\cdot &\cdot& \cdot& \cdot &\cdot &\cdot& \cdot& \cdot& 1& \cdot& \cdot& \cdot& 1& 1& \cdot& \cdot& \cdot &\cdot &1\\
\mathtt{(2)}&\mbf{15}&\cdot &\cdot &\cdot &\cdot &\cdot &\cdot &\cdot &\cdot& \cdot& \cdot& 1& \cdot& \cdot& 1& \cdot& 1& 1& \cdot& \cdot &\cdot\\
\mathtt{(1)}&\mbf{8}&\cdot& \cdot& \cdot& \cdot& \cdot &\cdot& \cdot& \cdot& \cdot& \cdot& \cdot &1& \cdot& 1 &\cdot &1 &\cdot &\cdot &\cdot &1\\
\mathtt{(1)}&\mbf{7}&\cdot& \cdot& \cdot &\cdot& \cdot& \cdot& \cdot& \cdot &\cdot& \cdot &\cdot &\cdot &1 &1& \cdot &\cdot &1& \cdot& \cdot& 1\\
\mathtt{(3)}&\mbf{32}&\cdot &\cdot& \cdot& \cdot &\cdot &\cdot &\cdot &\cdot &\cdot &\cdot &\cdot &\cdot &\cdot &1& 1 &1& 1& 1& 1& 1\\
&\mbf{28}&\cdot &\cdot &\cdot& \cdot& \cdot& \cdot& \cdot &\cdot& \cdot& \cdot& \cdot &\cdot &\cdot &\cdot &1 &\cdot &\cdot &1 &1 &\cdot\\
&\mbf{14}&\cdot &\cdot& \cdot& \cdot &\cdot &\cdot& \cdot& \cdot &\cdot &\cdot &\cdot &\cdot &\cdot &\cdot &\cdot &1 &\cdot& \cdot& 1 &\cdot\\
&\mbf{12}&\cdot &\cdot &\cdot& \cdot& \cdot& \cdot& \cdot& \cdot& \cdot& \cdot& \cdot& \cdot &\cdot &\cdot &\cdot &\cdot &1 &1 &\cdot &\cdot\\
&\mbf{6}&\cdot &\cdot &\cdot &\cdot& \cdot& \cdot& \cdot& \cdot &\cdot &\cdot &\cdot& \cdot &\cdot& \cdot &\cdot& \cdot& \cdot& 1& \cdot& \cdot\\
&\mbf{4}&\cdot &\cdot &\cdot &\cdot &\cdot &\cdot &\cdot &\cdot &\cdot &\cdot& \cdot& \cdot& \cdot& \cdot& \cdot &\cdot &\cdot &\cdot& 1 &\cdot\\
&\mbf{2}&\cdot &\cdot& \cdot& \cdot& \cdot& \cdot &\cdot &\cdot &\cdot &\cdot &\cdot& \cdot& \cdot& \cdot &\cdot &\cdot& \cdot& \cdot& \cdot& 1\\
\end{block}
\\
\begin{block}{cc(cccccccccccccccccccc)}
\star&\mbf{1}& 1&  \cdot&  \cdot & \cdot & \cdot & \cdot &-1 & 1 & 1 & \cdot & 1 & \cdot & \cdot &-1 & \cdot & \cdot & \cdot & \cdot & \cdot & 1\\
\star&\mbf{5}& \cdot & 1 & \cdot & \cdot &-1 &-1 & 1 & \cdot & \cdot & \cdot & \cdot &-2 & \cdot & 1& -1 & \cdot &-1 & 1 & \cdot& \cdot\\
\star&\mbf{3}& \cdot & \cdot & 1 &-1 & \cdot& -1 & 1 & \cdot & \cdot & \cdot & \cdot & \cdot& -2 & 1 &-1 &-1 & \cdot & \cdot & 1 & \cdot\\
\star&\mbf{13}& \cdot & \cdot & \cdot & 1 & \cdot & \cdot &-1 & \cdot & 1 & 1 & \cdot & \cdot & 1 &-1 & \cdot & 1 & \cdot & \cdot & \cdot & \cdot\\
\star&\mbf{11}& \cdot & \cdot & \cdot & \cdot & 1 & \cdot &-1 & 1 & \cdot & 1 & \cdot & 1 & \cdot &-1 & \cdot & \cdot & 1 & \cdot & \cdot & \cdot\\
\star&\mbf{27}& \cdot & \cdot & \cdot & \cdot & \cdot & 1 &-1 & \cdot & \cdot & \cdot & 1 & 1 & 1 &-1 & 1 & \cdot & \cdot & \cdot & \cdot & \cdot\\
\mathtt{(1)}&\mbf{31}& \cdot & \cdot & \cdot & \cdot & \cdot & \cdot  &1 &-1& -1& -1 &-1 &-1& -1 & 3 &-1 &-1 &-1 & \cdot & \cdot &-1\\
\mathtt{(3)}&\mbf{30}& \cdot & \cdot & \cdot & \cdot & \cdot & \cdot & \cdot & 1 & \cdot & \cdot & \cdot & \cdot & \cdot &-1 & \cdot & \cdot & 1 & \cdot &-1 & 1\\
\mathtt{(3)}&\mbf{29}& \cdot & \cdot & \cdot & \cdot  &\cdot & \cdot & \cdot & \cdot & 1 & \cdot & \cdot & \cdot & \cdot &-1 & \cdot & 1 & \cdot& -1 & \cdot & 1\\
\mathtt{(2)}&\mbf{22}& \cdot & \cdot & \cdot & \cdot & \cdot & \cdot & \cdot & \cdot & \cdot & 1 & \cdot & \cdot & \cdot &-1 & \cdot & 1 & 1 & \cdot & \cdot & \cdot\\
\mathtt{(2)}&\mbf{15}& \cdot & \cdot & \cdot & \cdot & \cdot  &\cdot & \cdot & \cdot & \cdot & \cdot & 1 & \cdot & \cdot &-1 & 1 & \cdot & \cdot & \cdot & \cdot & 1\\
\mathtt{(1)}&\mbf{8}& \cdot & \cdot & \cdot & \cdot & \cdot & \cdot & \cdot & \cdot & \cdot & \cdot & \cdot & 1 & \cdot &-1 & 1 & \cdot & 1& -1 & \cdot & \cdot\\
\mathtt{(1)}&\mbf{7}& \cdot & \cdot & \cdot & \cdot & \cdot & \cdot & \cdot & \cdot & \cdot & \cdot & \cdot & \cdot & 1 &-1 & 1 & 1 & \cdot & \cdot& -1 & \cdot\\
\mathtt{(3)}&\mbf{32}& \cdot & \cdot & \cdot & \cdot & \cdot & \cdot  &\cdot & \cdot & \cdot & \cdot & \cdot & \cdot & \cdot & 1 &-1& -1& -1 & 1 & 1& -1\\
&\mbf{28}& \cdot & \cdot & \cdot & \cdot & \cdot & \cdot & \cdot & \cdot & \cdot & \cdot & \cdot & \cdot & \cdot & \cdot & 1 & \cdot & \cdot &-1& -1 & \cdot\\
&\mbf{14}& \cdot & \cdot & \cdot & \cdot & \cdot & \cdot & \cdot & \cdot & \cdot & \cdot & \cdot & \cdot & \cdot & \cdot & \cdot & 1 & \cdot & \cdot &-1 & \cdot\\
&\mbf{12}& \cdot & \cdot & \cdot & \cdot & \cdot & \cdot & \cdot & \cdot & \cdot & \cdot & \cdot & \cdot & \cdot & \cdot & \cdot & \cdot & 1 &-1 & \cdot & \cdot\\
&\mbf{6}& \cdot & \cdot & \cdot & \cdot & \cdot & \cdot  &\cdot & \cdot & \cdot & \cdot & \cdot & \cdot & \cdot & \cdot & \cdot & \cdot & \cdot & 1 & \cdot & \cdot\\
&\mbf{4}& \cdot & \cdot & \cdot & \cdot & \cdot  &\cdot  &\cdot & \cdot & \cdot  &\cdot  &\cdot  &\cdot & \cdot&  \cdot&  \cdot & \cdot & \cdot  &\cdot & 1 & \cdot\\
&\mbf{2}& \cdot & \cdot & \cdot & \cdot & \cdot & \cdot & \cdot & \cdot & \cdot&  \cdot & \cdot&  \cdot & \cdot & \cdot & \cdot&  \cdot & \cdot & \cdot & \cdot & 1\\
\end{block}
\end{blockarray}
\]

The characters of the finite-dimensional reps $\el(\mbf{1})$, $\el(\mbf{13})$, and $\el(\mbf{11})$ would then be:
\begin{align*}
&\chi_{\el(\mbf{1})}(t)=3775 + 3700(t^{-1}+t) + 3510(t^{-2}+t^2) + 3240(t^{-3}+t^3) + 2925(t^{-4}+t^4) + 2600(t^{-5}+t^5)\\
&\qquad\qquad + 2300(t^{-6}+t^6) + 2024(t^{-7}+t^7) + 1771(t^{-8}+t^8) + 1540(t^{-9}+t^9) + 1330(t^{-10}+t^{10})\\
&\qquad\qquad + 1140(t^{-11}+t^{11}) + 969(t^{-12}+t^{12}) + 816(t^{-13}+t^{13}) + 680(t^{-14}+t^{14}) + 560(t^{-15}+t^{15})\\
&\qquad\qquad + 455(t^{-16}+t^{16}) + 364(t^{-17}+t^{17}) + 286(t^{-18}+t^{18}) + 220(t^{-19}+t^{19}) + 165(t^{-20}+t^{20})\\
&\qquad\qquad + 120(t^{-21}+t^{21}) + 84(t^{-22}+t^{22}) + 56(t^{-23}+t^{23}) + 35(t^{-24}+t^{24}) + 20(t^{-25}+t^{25}) \\
&\qquad\qquad+ 10(t^{-26}+t^{26}) + 4(t^{-27}+t^{27}) + t^{-28}+t^{28}\\
&\dim\el(\mbf{1})=65625\\
\\
&\chi_{\el(\mbf{11})}(t)=\chi_{\el(\mbf{13})}(t)=765 + 720(t^{-1}+t) + 612(t^{-2}+t^2) + 468(t^{-3}+t^3) + 315(t^{-4}+t^4)\\
&\qquad\qquad\qquad\qquad\qquad + 180(t^{-5}+t^5) + 90(t^{-6}+t^6) + 36(t^{-7}+t^7) + 9(t^{-8}+t^8)\\
&\dim\el(\mbf{11})=\dim\el(\mbf{13})=5625\\
\end{align*}

Additionally, according to \cite{GP} there are also two dual blocks of defect $1$, each containing only two irreps. Their characters come out to be:

\begin{align*}
&\el(\mbf{18})=\m(\mbf{18})-\m(\mbf{21})\\
&\el(\mbf{20})=\m(\mbf{20})-\m(\mbf{19})\\
&\chi_{\el(\mbf{18})}(t)=\chi_{\el(\mbf{20})}(t)=16t^{-5.5}\frac{\sum\limits_{i=0}^9t^i}{(1-t)^3}\\
&\dim\Supp\el(\mbf{18})=\dim\Supp\el(\mbf{20})=3\\
\end{align*}

\section{Finite-dimensional irreducible representations of $H_c(F_4)$}

The dimensions of $\el_c(\mathbf{1})$ for $H_c(F_4)$ with equal parameters were calculated in \cite{OY}. The finite-dimensional representations can be identified and their dimensions calculated by hand for $F_4$ with equal parameters by methods similar to those in \cite{BP}. The results of such calculations are given below.

The labels for the $F_4$-irreps are those given by Sage in the Sage for iPad app. The numbers labeling the irreps are not totally random; they are increasing with the dimensions of the irreps. $\mbf{9}$ denotes the standard rep, $\mbf{1}$ the trivial rep, $\mbf{6}$ and $\mbf{8}$ are $2$-dimensional reps, and $\mbf{23}$ is $9$-dimensional. 

\proposition
The following is a complete list of the finite-dimensional irreps for $H_{\frac{1}{d}}(F_4)$ with equal parameters, except in the case of $c=\frac{1}{2}$ where the finite-dimensionality of the three non-spherical reps listed is shown to be likely. If they are finite-dimensional then their dimensions and characters are as claimed.
\begin{align*}
&\mbf{c=\frac{1}{12}:}\qquad\chi_{\el(\mbf{1})}(w,t)=\chi_\mbf{1}\\
& \qquad\qquad\dim \el(\mbf{1})=1\\
\\
&\mbf{c=\frac{1}{8}:}\qquad\chi_{\el(\mbf{1})}(w,t)=\chi_\mbf{1}(t^{-1}+t)+\chi_{\hstar}\\
&\qquad\qquad\chi_{\el(\mbf{1})}(t)=t^{-1}+4+t\\
&\qquad\qquad\dim\el(\mbf{1})=6\\
\\
&\mbf{c=\frac{1}{6}:}\qquad\chi_{\el(\mbf{1})}(w,t)=\chi_\mbf{1}(t^{-2}+t^2)+\chi_{\hstar}(t^{-1}+t)+\chi_{S^2\hstar}\\
&\qquad\qquad\chi_{\el(\mbf{1})}(t)=t^{-2}+t^2+4(t^{-1}+t)+10\\
&\qquad\qquad\dim\el(\mbf{1})=20\\
\\
&\qquad\qquad\chi_\el(\mbf{6})(w,t)=\chi_\mbf{6}\\
&\qquad\qquad\dim\el(\mbf{6})=2\\
\\
&\qquad\qquad\chi_\el(\mbf{8})(w,t)=\chi_\mbf{8}\\
&\qquad\qquad\dim\el(\mbf{8})=2\\
\\
&\mbf{c=\frac{1}{4}:}\qquad\chi_{\el(\mbf{1})}=\chi_\mbf{1}(t^{-4}+t^4)+\chi_{\hstar}(t^{-3}+t^3)+\chi_{S^2\hstar}(t^{-2}+t^2)+\chi_{S^3\hstar}(t^{-1}+t)\\
&\qquad\qquad\qquad\qquad\qquad+\chi_{S^4\hstar}-\chi_\mbf{23}\\
&\qquad\qquad\chi_{\el(\mbf{9})}(t)=t^{-4}+t^4+4(t^{-3}+t^3)+10(t^{-2}+t^2)+20(t^{-1}+t)+26\\
&\qquad\qquad\dim\el(\mbf{1})=96\\
\\
&\qquad\qquad\chi_{\el(\mbf{9})}(w,t)=\chi_\mbf{9}(t^{-1}+t)+\chi_{\hstar}-\chi_\mbf{23}\\
&\qquad\qquad\chi_{\el(\mbf{9})}(t)=4(t^{-1}+t)+7\\
&\qquad\qquad\dim\el(\mbf{9})=15\\
\\
&\mbf{c=\frac{1}{3}:}\qquad_{\el(\mbf{1})}(w,t)=\chi_\mbf{1}(t^{-6}+t^6)+\chi_{\hstar}(t^{-5}+t^5)+\chi_{S^2\hstar}(t^{-4}+t^4)+\chi_{S^3\hstar}(t^{-3}+t^3)\\
&\qquad\qquad\qquad\qquad\qquad+(\chi_{S^4\hstar}-\chi_{\mbf{6}\oplus\mbf{8}})(t^{-2}+t^2)+(\chi_{S^5\hstar}-\chi_{\hstar\otimes(\chi_\mbf{6}\oplus\chi_\mbf{8})})(t^{-1}+t)\\
&\qquad\qquad\qquad\qquad\qquad+\chi_{S^6\hstar}-\chi_{S^2\hstar\otimes(\chi_\mbf{6}\oplus\chi_\mbf{8})}\\
&\qquad\qquad\chi_{\el(\mbf{1})}(t)=t^{-6}+t^6+4(t^{-5}+t^5)+10(t^{-4}+t^4)+20(t^{-3}+t^3)+31(t^{-2}+t^2)\\
&\qquad\qquad\qquad\qquad\qquad+40(t^{-1}+t)+44\\
&\qquad\qquad\dim\el(\mbf{1})=256\\
\\
&\quad\quad\qquad\chi_{\el(\mbf{9})}(w,t)=\chi_\mbf{9}(t^{-2}+t^2)+\chi_{\hstar\otimes\mbf{9}}(t^{-1}+t)+\chi_{S^2\otimes\hstar}-\chi_{\mbf{16}\oplus\mbf{18}}\\
&\qquad\qquad\chi_{\el(\mbf{9})}(t)=4(t^{-1}+t^2)+16(t^{-1}+t)+24\\
&\qquad\qquad\dim\el(\mbf{9})=64
\end{align*}
\begin{align*}
&\mbf{c=\frac{1}{2}:}\qquad\chi_{\el(\mbf{1})}(w,t)=\left(\sum_{k=0}^8\chi_{S^k\hstar}(t^{k-10}+t^{10-k})\right)+(\chi_{S^9\hstar}-\chi_{\mbf{16}\oplus\mbf{18}})(t^{-1}+t)\\
&\qquad\qquad\qquad\qquad\qquad\qquad+\chi_{S^{10}\hstar}-\chi_{\hstar\otimes(\mbf{16}\oplus\mbf{18})}\\
&\qquad\qquad\chi_{\el(\mbf{1})}(t)=\sum_{k=0}^8{3+k\choose3}(t^{k-10}+t^{10-k})+204(t^{-1}+t)+222\\
&\qquad\qquad\dim\el(\mbf{1})=1620\\
\\
%&\qquad\qquad\chi_{\el(\mbf{6})}(w,t)=\chi_{\mbf{6}}(t^{-4}+t^4)+\chi_{\mbf{18}}(t^{-3}+t^3)+\chi_{\mbf{6}\oplus\mbf{22}}(t^{-2}+t^2)+\chi_{\mbf{9}\oplus\mbf{10}\oplus\mbf{18}}(t^{-1}+t)\\
%&\qquad\qquad\qquad\qquad\qquad+\chi_{S^4\hstar\otimes\mbf{6}}-\chi_{S^2\hstar\otimes\mbf{23}}+\chi_{\hstar\otimes\mbf{16}}\\
&\qquad\qquad\chi_{\el(\mbf{6})}(w,t)=\chi_\mbf{6}(t^{-4}+t^4)+\chi_{\hstar\otimes\mbf{6}}(t^{-3}+t^3)+(\chi_{S^2\hstar\otimes\mbf{6}}-\chi_\mbf{23})(t^{-2}+t^2)\\
&\qquad\qquad\qquad\qquad\qquad+(\chi_{S^3\hstar\otimes\mbf{6}}-\chi_{\hstar\otimes\mbf{23}}+\chi_\mbf{16})(t^{-1}+t)+\chi_{S^4\hstar\otimes\mbf{6}}-\chi_{S^2\hstar\otimes\mbf{23}}+\chi_{\hstar\otimes\mbf{16}}\\
&\qquad\qquad\chi_{\el(\mbf{6})}(t)=2(t^{-4}+t^4)+8(t^{-3}+t^3)+11(t^{-2}+t^2)+12(t^{-1}+t)+12\\
&\qquad\qquad\dim\el(\mbf{6})=78\\
\\
&\qquad\qquad\chi_{\el(\mbf{8})}(w,t)=\chi_\mbf{8}(t^{-4}+t^4)+\chi_{\hstar\otimes\mbf{8}}(t^{-3}+t^3)+(\chi_{S^2\hstar\otimes\mbf{8}}-\chi_\mbf{23})(t^{-2}+t^2)\\
&\qquad\qquad\qquad\qquad\qquad+(\chi_{S^3\hstar\otimes\mbf{8}}-\chi_{\hstar\otimes\mbf{23}}+\chi_\mbf{18})(t^{-1}+t)+\chi_{S^4\hstar\otimes\mbf{8}}-\chi_{S^2\hstar\otimes\mbf{23}}+\chi_{\hstar\otimes\mbf{18}}\\
&\qquad\qquad\chi_{\el(\mbf{8})}(t)=\chi_{\el(\mbf{6})}(t)\\
&\qquad\qquad\dim\el(\mbf{8})=78\\
\\
&\qquad\qquad\chi_{\el(\mbf{23})}(w,t)=\chi_{\mbf{23}}(t^{-2}+t^2)+(\chi_{\hstar\otimes\mbf{23}}-\chi_{\mbf{16}\otimes\mbf{18}})(t^{-1}+t)+\chi_{S^2\hstar\otimes\mbf{23}}-\chi_{\hstar\otimes(\mbf{16}\oplus\mbf{18})}\\
&\qquad\qquad\chi_{\el(\mbf{23})}(t)=9(t^{-2}+t^2)+20(t^{-1}+t)+26\\
&\qquad\qquad\dim\el(\mbf{23})=84
\end{align*}

\begin{proof}
The character of $\el_{\frac{1}{12}}(\mbf{1})=1$ because $12$ is the Coxeter number of $F_4$, and this always holds $c=1/h$ with $h$ the Coxeter number. It's also immediate from the weight line for $c=\frac{1}{12}$.

\vspace*{.2cm}

\large
\begin{center}
\aldine\\
\end{center}
\normalsize

\vspace*{.3cm}

When $c=\frac{1}{8}$ only $\tau=\mbf{1}$ satisfies $h_c(\tau)\in\mathbb{Z}_{\leq0}$, so only $\el(\mbf{1})$ can be finite-dimensional. As $8$ is an elliptic number of $F_4$, $\el(\mbf{1})$ is in fact finite-dimensional. There is no $\tau$ with $h_c(\tau)=0$ so $\el(\mbf{1})[0]=\chi_{\hstar}$.

\vspace*{.2cm}

\large
\begin{center}
\aldine\\
\end{center}
\normalsize

\vspace*{.3cm}

$\el_{\frac{1}{6}}\mbf{1}$ is finite-dimensional because $6$ is an elliptic number of $F_4$ \cite{VV}. $S^2\hstar=\mbf{23}\oplus\mbf{1}$ so $\mbf{6},\;\mbf{8}$ don't appear in the Verma-decomposition of $\el_{\frac{1}{6}}\mbf{1}$. In order to have $\dim\el(\mbf{1})[-1]=\dim\el(\mbf{1})[1],$ it must be that $\el(\mbf{1})=\m(\mbf{1})-c_{16}\m(\mbf{16})-c_{18}\m(\mbf{18})...$ with $c_{16}+c_{18}=2$. Since $S^3\hstar=\mbf{9}\oplus\mbf{16}\oplus\mbf{18}$, $c_{16}=c_{18}=1$.

$\el_{\frac{1}{6}}(\mbf{6})$ and $\el_{\frac{1}{6}}(\mbf{8})$ are easily seen to be finite-dimensional: their lowest weights occur in degree $0$, and $\hstar\otimes\mbf{6}=\mbf{18}$, $\hstar\otimes\mbf{8}=\mbf{16}$. Then by Lemma 3.1, $\mbf{18}$ generates a subrep of $\mbf{6}$ and $\mbf{16}$ generates a subrep of $\mbf{8}$, so the dimension of the graded piece in degree $1$ of both $\el_{\frac{1}{6}}(\mbf{6})$ and $\el_{\frac{1}{6}}(\mbf{8})$ is $0$.

\vspace*{.2cm}

\large
\begin{center}
\aldine\\
\end{center}
\normalsize

\vspace*{.3cm}

$F_4$ has two finite-dimensional representations at $c=\frac{1}{4}$, and their lowest weights are the trivial rep and the standard rep.

When $c=\frac{1}{4}$ there are five $\tau\in\Irr F_4$ satisfying $h_c(\tau)\in\mathbb{Z}_\leq0$: $\tau\in\{\mbf{1},\mbf{9},\;\mbf{23},\;\mbf{6},\;\mbf{8}\}$. $h_c(\mbf{23})=0$ but there is no $\tau$ with $h_c(\tau)=1$, so $0=\dim\el(\mbf{23})[-1]<\dim\el(\mbf{23})[1]$ and $\el(\mbf{23})$ cannot be finite-dimensional. The weight line containing $\mbf{6}$ and $\mbf{8}$ has no $\tau$ with $h_c(\mbf{\tau})=0$ or $1$ but $h_c(\mbf{6})=h_c(\mbf{8})=-1$, so $2=\dim\el(\mbf{6})[-1]=\dim\el(\mbf{8})[-1]<\dim\el(\mbf{6})[1]=\dim\el(\mbf{8})[1]$ and $\el(\mbf{6})$ and $\el(\mbf{8})$ cannot be finite-dimensional. 

$\el_{\frac{1}{4}}(\mbf{1})$ is finite-dimensional because $4$ is an elliptic number of $F_4$. $S^3\hstar=\mbf{9}\oplus\mbf{16}\oplus\mbf{18}$ so $\dim\Hom(\m(\mbf{9}),\m(\mbf{1}))\leq1$. If $\dim\Hom(\m(\mbf{9}),\m(\mbf{1}))=1$ then the equality $\dim\el(\mbf{1})[-1]=\dim\el(\mbf{1})[1]$ is satisified and the dimension of spherical rep is $81$. If $\dim\Hom(\m(\mbf{9}),\m(\mbf{1}))=0$ then the only way for $\dim\el(\mbf{1})[-1]=\dim\el(\mbf{1})[1]$ to hold is if $\dim\Hom(\m(\mbf{23}),\m(\mbf{1}))=1$, in which case the dimension of the spherical rep is $96$. By the calculation of $\dim\el(\mbf{1})$ in \cite{OY}, the correct option is the second one, and the character of $\el(\mbf{1})$ follows.

$\el_{\frac{1}{4}}(\mbf{9})$ is also finite-dimensional: $h_c(\mbf{23})-h_c(\mbf{9})=1$ and $\mbf{23}\subset\hstar\otimes\mbf{9}=\mbf{1}\oplus\mbf{14}\oplus\mbf{23}$, so by Lemma 3.1, $\dim\Hom(\m(\mbf{23}),\m(\mbf{9}))=1$ and $\mbf{23}$ generates a subrep of $\m(\mbf{9})$. Then $\dim\el(\mbf{9})[-1]=4={5\choose3}\cdot4-{4\choose3}\cdot9=\dim\el(\mbf{9})[1]$, and the operator $E$ kills $\mbf{1}\subset\el(\mbf{1})[0]$. By (E) then, $\el(\mbf{9})$ is finite-dimensional and its character is as claimed.

\vspace*{.2cm}

\large
\begin{center}
\aldine\\
\end{center}
\normalsize

\vspace*{.3cm}

$\el_{\frac{1}{3}}(\mbf{1})$ is finite-dimensional since $3$ is an elliptic number of $F_4$ \cite{VV} . Writing $\el(\mbf{1})=\m(\mbf{1})-c_6\m(\mbf{6})-c_8\m(\mbf{8})$... and setting $\dim\el(\mbf{1})[-1]=\dim\el(\mbf{1})[1]$ shows that $c_6+c_8=2$. The graded character of $\el(\mbf{1})$ follows. To determine the correct multiplicities $c_6$ and $c_8$, consider parabolic induction of the spherical rep from the two copies of $A2\times A_1$ in $F_4$ -- inducing $\el(\mathrm{Triv})$ from one copy produces the expression $\m(\mbf{1})-\m(\mbf{8})...$ while inducing from the other gives $\m(\mbf{1})-\m(\mbf{6})...$. Thus $c_6=c_8=1$.

Let $\mbf{1}_{+-}$ be the one-dimensional irrep of $A_2\times A_1$ on which $A_2$ acts trivially and $A_1$ acts by sign. Inducing $\el(\mbf{1}_{+-})$ from the Grothendieck group of the rational Cherednik algebra at $c=\frac{1}{3}$ for one copy of $A_2\times A_1$ gives rise to a module containing the expression $\m(\mbf{9})-\m(\mbf{18})$ while inducing the same representation from $A_2\times A_2\times A_1$ embedded in the other way gives something containing $\m(\mbf{9})-\m(\mbf{16})$. The combined multiplicities of $\m(\mbf{16})$ and $\m(\mbf{18})$ in $\el(\mbf{9})$ cannot be less than $-2$ or else the result will be that $\dim\el(\mbf{9})[1]<\dim\el(\mbf{9})[-1]$ in violation of $\s$-theory. With $\el(\mbf{9})=\m(\mbf{9})-\m(\mbf{16})-\m(\mbf{18})...$, $\dim\el(\mbf{9})[-1]=\dim\el(\mbf{9})[1]$. The character in degrees $-2$, $-1$, and $0$, $1$ follow, while finite-dimensionality follows from (E), and so $\el(\mbf{9})$ has the stated character and dimension.

\vspace*{.2cm}

\large
\begin{center}
\aldine\\
\end{center}
\normalsize

\vspace*{.3cm}

The rational Cherednik algebra of $B_3$ with equal parameters $c=\frac{1}{2}$ has two finite-dimensional representations: $\el(\mathrm{Triv})$ and a one-dimensional representation spanned by the one-dimensional rep $\tau\in\Irr B_3$ such that $h_c(\tau)=0$. Its decomposition into Vermas is easy to compute, and inducing $\el(\tau)$ from $H_{\frac{1}{2}}(B_3)$ to $H_{\frac{1}{2}}(F_4)$ gives an expression whose lowest-$h_c$-weight terms are $\m(\mbf{6})-\m(\mbf{23})+\m(\mbf{16})...$ Thus the decomposition of $\el(\mbf{6})$ up through degree $1$ is $\m(\mbf{6})-\m(\mbf{23})+c_{18}\m(\mbf{18})+c_{16}\m(\mbf{16})$ where $c_{18}\leq0$ and $c_{16}$ is $0$ or $1$. The inequality $\dim\el(\mbf{6})[-1]\leq\dim\el(\mbf{6})[1]$ only holds if $c_{18}+c_{16}\geq1$, so $c_{18}=0$ and $c_{16}=1$. Then $\dim\el(\mbf{6})[-1]=\dim\el(\mbf{6})[1]$, and furthermore, \textit{if} $E$ kills a $1$-dimensional rep in graded degree $[0]$ then by (E) $\el(\mbf{6})$ is finite-dimensional. By switching long and short roots, $\el(\mbf{8})$ would be finite-dimensional as well with the character given in the theorem.

Finite-dimensionality of $\el_{\frac{1}{2}}(\mbf{23})$ follows from Lemma 3.1 and (E).
$$\hstar\otimes\mbf{23}=\mbf{9}\oplus\mbf{16}\oplus\mbf{18}\oplus\mbf{25}$$
and $h_c(\mbf{16})-h_c(\mbf{23})=h_c(\mbf{18})-h_c(\mbf{23})$=1, so $\mbf{16}$ and $\mbf{18}$ each generate a subrep of $\m(\mbf{23})$ by Lemma 3.1. Then $\dim\el(\mbf{23})[1]=20\cdot9-10\cdot16=20=4\cdot9-16=\dim\el(\mbf{23})[1]$, and \textit{if} $E$ kills something in $\el(\mbf{23})[0]$ then by (E) $\el(\mbf{23})$ would be finite-dimensional.

The dimension of $\el(\mbf{1})$ is known thanks to \cite{OY}: $\dim\el(\mbf{1})=1620$. This is enough to pin down the character of $\el(\mbf{1})$. Counting dimensions of weight spaces,
\begin{align*}
&\dim\el(\mbf{1})[-1]={12\choose3}+2(c_6+c_8){6\choose3}+9c_{23}{4\choose3}+8(c_{16}+c_{18})\\
&\dim\el(\mbf{1})[1]={14\choose3}+2(c_6+c_8){8\choose3}+9c_{23}{6\choose3}+8(c_{16}+c_{18}){5\choose3}\\
\end{align*}
where $c_\tau$ is the multiplicity of $\m(\tau)$ in $\el(\mbf{1})$. Since $\dim\el(\mbf{1})[-1]=\dim\el(\mbf{1})[1]$, it follows that $-2(c_{23}+1)=c_6+c_8+c_{16}+c_{18}$. $c_6$ and $c_8$ are either $0$ or negative. Then it's enough to check that the only way to get $1620$ as the dimension is if $c_{16}+c_{18}=-2$ and $c_{23}=c_6=c_8=0$. It's then necessary to have $c_{16}=c_{18}$ in order that the restriction to all parabolics will be $0$.

\end{proof}

\section{Simple representations of $H_c(E_6)$}
In this section will be found the decompositions, characters, dimensions of support, and dimensions of all irreducible representations of $H_c(E_6)$ for $c=1/d$ when $d>2$ divides one of the fundamental degrees of $E_6$. The $d$ that need to be considered are $d\in\{3,\;4,\;5,\;6,\;8,\;9,\;12\}$. The result is a description of Category $\oh_c(E_6)$ when $\oh_c(E_6)$ is not semisimple for all $c$ except half-integers. For each $c=1/d$ as above, the upper-unitriangular decomposition matrix of $\oh_c(E_6)$ recording the multiplicities $[\m(\tau):\el(\sigma)]$ has only $0$'s and $1$'s for entries.

\subsection{$c=\frac{1}{12},\;\frac{1}{9},$ or $\frac{1}{8}$.} In each case, there is a single block of defect $1$. All other blocks have defect $0$, i.e. consist of a single Verma which is simple. Category $\oh_c$ contains a unique irrep of less than full support, which in all three cases is $\el(1_p)$. \\

\begin{align*}
&\mbf{c=\frac{1}{12}:}\qquad
\begin{picture}(240,10)
\put(0,5){\line(1,0){240}}
\put(0,5){\circle*{5}}
\put(40,5){\circle*{5}}
\put(80,5){\circle*{5}}
\put(120,5){\circle*{5}}
\put(160,5){\circle*{5}}
\put(200,5){\circle*{5}}
\put(240,5){\circle*{5}}
\put(-2,20){$0$}
\put(38,20){$1$}
\put(78,20){$2$}
\put(118,20){$3$}
\put(158,20){$4$}
\put(198,20){$5$}
\put(238,20){$6$}
\put(-2,-10){$1_p$}
\put(38,-10){$6_p$}
\put(75,-10){$15_p$}
\put(115,-10){$20_p$}
\put(154,-10){$15_p'$}
\put(196,-10){$6_p'$}
\put(236,-10){$1_p'$}
\end{picture}\\
\\
%&\qquad\el(1_p)=\m(1_p)-\m(6_p)+\m(15_p)-\m(20_p)+\m(15_p')-\m(6_p')+\m(1_p')\\
&\qquad\qquad\qquad\chi_{\el(1_p)}(t)=1, \qquad\dim\el(1_p)=1\\
\\
\\
&\mbf{c=\frac{1}{9}:}\qquad
\begin{picture}(320,20)
\put(0,5){\line(1,0){320}}
\put(0,5){\circle*{5}}
\put(80,5){\circle*{5}}
\put(120,5){\circle*{5}}
\put(160,5){\circle*{5}}
\put(200,5){\circle*{5}}
\put(240,5){\circle*{5}}
\put(320,5){\circle*{5}}
\put(-5,20){$-1$}
\put(78,20){$1$}
\put(118,20){$2$}
\put(158,20){$3$}
\put(198,20){$4$}
\put(238,20){$5$}
\put(318,20){$7$}
\put(-2,-10){$1_p$}
\put(75,-10){$20_p$}
\put(115,-10){$64_p$}
\put(155,-10){$90_s$}
\put(195,-10){$64_p'$}
\put(235,-10){$20_p'$}
\put(316,-10){$1_p'$}
\end{picture}\\
\\
%&\qquad\el(1_p)=\m(1_p)-\m(20_p)+\m(64_p)-\m(90_s)+\m(64_p')-\m(20_p')+\m(1_p')\\
&\qquad\qquad\qquad\chi_{\el(1_p)}(t)=t^{-1}+6+t, \qquad\dim\el(1_p)=8\\
\\
\\
&\mbf{c=\frac{1}{8}:}\qquad
\begin{picture}(360,20)
\put(0,5){\line(1,0){360}}
\put(0,5){\circle*{5}}
\put(120,5){\circle*{5}}
\put(160,5){\circle*{5}}
\put(200,5){\circle*{5}}
\put(240,5){\circle*{5}}
\put(360,5){\circle*{5}}
\put(-8,20){$-1.5$}
\put(115,20){$1.5$}
\put(155,20){$2.5$}
\put(195,20){$3.5$}
\put(235,20){$4.5$}
\put(355,20){$7.5$}
\put(-2,-10){$1_p$}
\put(115,-10){$30_p$}
\put(155,-10){$81_p$}
\put(195,-10){$81_p'$}
\put(235,-10){$30_p'$}
\put(356,-10){$1_p'$}
\end{picture}\\
\\
%&\qquad\el(1_p)=\m(1_p)-\m(30_p)+\m(81_p)-\m(81_p')+\m(30_p')-\m(1_p')\\
%&\qquad\qquad\qquad\chi_{\el(1_p)}(t)=t^{-1.5}\frac{1+5t+15t^2+5t^3+t^4}{1-t}\\
&\qquad\qquad\qquad\dim\Supp\el(1_p)=1\\
\end{align*}

\subsection{{\Large$\mbf{c=\frac{1}{6}}$}}

 There is a single nontrivial block, the block containing the spherical rep. It has five irreps of less than full support, of which two (with lowest weights the trivial rep and the standard rep) are finite-dimensional, as given by the following theorem: 
 
 \theorem The decomposition matrices $\mathcal{M}=\left([\m(\tau):\el(\sigma)]\right)$ and $\mathcal{L}=\left([\el(\tau):\m(\sigma)]\right)$ are:
\[
\begin{blockarray}{cccccccccccccccccc}
&& 1_p & 6_p & 20_p & 30_p & 15_q & 24_p & 60_p & 80_s & 60_s & 24_p' & 60_p' & 30_p' & 15_q' & 20_p' & 6_p' & 1_p'\\
\begin{block}{cc(cccccccccccccccc)}
\star&1_p & 1 & \cdot & 1 & \cdot & 1 & \cdot & 1 & \cdot & \cdot & \cdot & \cdot & \cdot & \cdot & \cdot & \cdot & \cdot \\
\star&6_p & \cdot & 1 & 1 & 1 & \cdot & \cdot & \cdot & \cdot & \cdot & \cdot & \cdot & \cdot & \cdot & \cdot & \cdot & \cdot \\
\mathtt{(2)}&20_p & \cdot & \cdot & 1 & 1 & \cdot & 1 & 1 & 1 & \cdot & \cdot & \cdot & \cdot & \cdot & \cdot & \cdot & \cdot \\
&30_p & \cdot & \cdot & \cdot & 1 & \cdot & \cdot & \cdot & 1 & \cdot & 1 & \cdot & \cdot & \cdot & \cdot & \cdot & \cdot \\
\mathtt{(1)}&15_q & \cdot & \cdot & \cdot & \cdot & 1 & \cdot & 1 & \cdot & 1 & \cdot & \cdot & \cdot & \cdot & \cdot & \cdot & \cdot \\
\mathtt{(2)}&24_p & \cdot & \cdot & \cdot & \cdot & \cdot & 1 & \cdot & 1 & \cdot & \cdot & \cdot & 1 & \cdot & \cdot & \cdot & \cdot \\
&60_p & \cdot & \cdot & \cdot & \cdot & \cdot & \cdot & 1 & 1 & 1 & \cdot & 1 & \cdot & \cdot & \cdot & \cdot & \cdot \\
&80_s & \cdot & \cdot & \cdot & \cdot & \cdot & \cdot & \cdot & 1 & \cdot & 1 & 1 & 1 & \cdot & 1 & \cdot & \cdot \\
&60_s & \cdot & \cdot & \cdot & \cdot & \cdot & \cdot & \cdot & \cdot & 1 & \cdot & 1 & \cdot & 1 & \cdot & \cdot & \cdot \\
&24_p' & \cdot & \cdot & \cdot & \cdot & \cdot & \cdot & \cdot & \cdot & \cdot & 1 & \cdot & \cdot & \cdot & 1& \cdot & \cdot \\
&60_p' & \cdot & \cdot & \cdot & \cdot & \cdot & \cdot & \cdot & \cdot & \cdot & \cdot & 1 & \cdot & 1 & 1 & \cdot & 1\\
&30_p' & \cdot & \cdot & \cdot & \cdot & \cdot & \cdot & \cdot & \cdot & \cdot & \cdot & \cdot & 1 & \cdot & 1 & 1 & \cdot\\
&15_q' & \cdot & \cdot & \cdot & \cdot & \cdot & \cdot & \cdot & \cdot & \cdot & \cdot & \cdot & \cdot & 1& \cdot & \cdot & 1\\
&20_p' & \cdot & \cdot & \cdot & \cdot & \cdot & \cdot & \cdot & \cdot & \cdot & \cdot & \cdot & \cdot & \cdot & 1 & 1 & 1\\
&6_p' & \cdot & \cdot & \cdot & \cdot & \cdot & \cdot & \cdot & \cdot & \cdot & \cdot & \cdot & \cdot & \cdot & \cdot & 1 & \cdot\\
&1_p' & \cdot & \cdot & \cdot & \cdot & \cdot & \cdot & \cdot & \cdot & \cdot & \cdot & \cdot & \cdot & \cdot & \cdot & \cdot & 1\\
\end{block}
\end{blockarray}
\]
\\
\[
\begin{blockarray}{cccccccccccccccccc}
&& 1_p & 6_p & 20_p & 30_p & 15_q & 24_p & 60_p & 80_s & 60_s & 24_p' & 60_p' & 30_p' & 15_q' & 20_p' & 6_p' & 1_p'\\
\begin{block}{cc(cccccccccccccccc)}
\star & 1_p & 1 & \cdot & -1 & 1 & -1 & 1 & 1 & -2 & \cdot & 1 & 1 & 1 & -1 & -1 & \cdot & 1 \\
\star & 6_p & \cdot & 1 & -1 & \cdot & \cdot & 1 & 1 & -1 & -1 & 1 & 1 & \cdot & \cdot & -1 & 1 & \cdot \\
\mathtt{(2)} & 20_p & \cdot & \cdot & 1 & -1 & \cdot & -1 & -1 & 2 & 1 & -1 & -2 & -1 & 1 & 2 & -1 & -1 \\
& 30_p & \cdot & \cdot & \cdot & 1 & \cdot & \cdot & \cdot & -1 & \cdot & \cdot & 1& 1 & -1 & -1 & \cdot & 1 \\
\mathtt{(1)} & 15_q & \cdot & \cdot & \cdot & \cdot & 1 & \cdot & -1 & 1 & \cdot & -1 & \cdot & -1 & \cdot & 1 & \cdot & -1 \\
\mathtt{(2)} & 24_p & \cdot & \cdot & \cdot & \cdot & \cdot & 1 & \cdot & -1 & \cdot & 1 & 1 & \cdot & -1 & -1 & 1 & 1 \\
& 60_p & \cdot & \cdot & \cdot & \cdot & \cdot & \cdot & 1 & -1 & -1 & 1 & 1 & 1 & \cdot & -2 & 1 & 1 \\
& 80_s & \cdot & \cdot & \cdot & \cdot & \cdot & \cdot & \cdot & 1 & \cdot & -1 & -1 & -1 & 1 & 2 & -1 & -2 \\
& 60_s & \cdot & \cdot & \cdot & \cdot & \cdot & \cdot & \cdot & \cdot & 1 & \cdot & -1 & \cdot & \cdot & 1 & -1 & \cdot \\
& 24_p' & \cdot & \cdot & \cdot & \cdot & \cdot & \cdot & \cdot & \cdot & \cdot & 1 & \cdot & \cdot & \cdot & -1 & 1 & 1 \\
& 60_p' & \cdot & \cdot & \cdot & \cdot & \cdot & \cdot & \cdot & \cdot & \cdot & \cdot & 1 & \cdot & -1 & -1 & 1 & 1\\
& 30_p' & \cdot & \cdot & \cdot & \cdot & \cdot & \cdot & \cdot & \cdot & \cdot & \cdot & \cdot & 1 & \cdot & -1 & \cdot & 1\\
& 15_q' & \cdot & \cdot & \cdot & \cdot & \cdot & \cdot & \cdot & \cdot & \cdot & \cdot & \cdot & \cdot & 1& \cdot & \cdot & -1\\
& 20_p' & \cdot & \cdot & \cdot & \cdot & \cdot & \cdot & \cdot & \cdot & \cdot & \cdot & \cdot & \cdot & \cdot & 1 & -1 & -1\\
& 6_p' & \cdot & \cdot & \cdot & \cdot & \cdot & \cdot & \cdot & \cdot & \cdot & \cdot & \cdot & \cdot & \cdot & \cdot & 1 & \cdot\\
& 1_p' & \cdot & \cdot & \cdot & \cdot & \cdot & \cdot & \cdot & \cdot & \cdot & \cdot & \cdot & \cdot & \cdot & \cdot & \cdot & 1\\
\end{block}
\end{blockarray}
\]\\

\begin{proof}
The weight line for this block has integer weights:

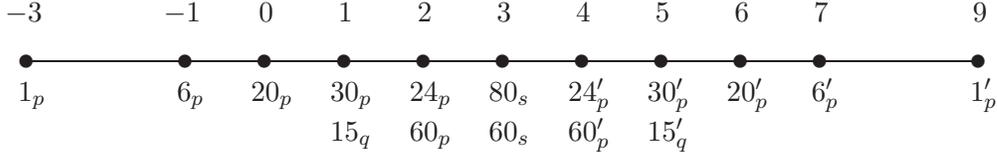
\begin{center}
\begin{picture}(360,65)
\put(0,40){\line(1,0){360}}
\put(0,40){\circle*{5}}
\put(60,40){\circle*{5}}
\put(90,40){\circle*{5}}
\put(120,40){\circle*{5}}
\put(150,40){\circle*{5}}
\put(180,40){\circle*{5}}
\put(210,40){\circle*{5}}
\put(240,40){\circle*{5}}
\put(270,40){\circle*{5}}
\put(300,40){\circle*{5}}
\put(360,40){\circle*{5}}
\put(-8,55){$-3$}
\put(52,55){$-1$}
\put(88,55){$0$}
\put(118,55){$1$}
\put(148,55){$2$}
\put(178,55){$3$}
\put(208,55){$4$}
\put(238,55){$5$}
\put(268,55){$6$}
\put(298,55){$7$}
\put(358,55){$9$}
\put(-3,25){$1_p$}
\put(57,25){$6_p$}
\put(85,25){$20_p$}
\put(115,25){$30_p$}
\put(115,10){$15_q$}
\put(145,25){$24_p$}
\put(145,10){$60_p$}
\put(175,25){$80_s$}
\put(175,10){$60_s$}
\put(205,25){$24_p'$}
\put(205,10){$60_p'$}
\put(235,25){$30_p'$}
\put(235,10){$15_q'$}
\put(265,25){$20_p'$}
\put(297,25){$6_p'$}
\put(357,25){$1_p'$}
\end{picture}\\
\end{center}

All of matrix $\mathcal{M}$ save the columns labeled $1_p,\;6_p,\;20_p,\;15_q,\;24_p$ have been reprinted verbatim from the decomposition matrix of the Hecke algebra at $e=6$. The missing columns of $\mathcal{M}$ killed by KZ functor, corresponding to the columns labeled by the bulleted and starred rows, can be reconstructed from the Hecke submatrix as follows. First, by (RR) the row $24_p'$ of $\mathcal{M}$, $\dim\Hom(\m(20_p'),\m(24_p'))=1$ while $60_p',\;30_p',\;15_q'$ don't consist of singular vectors in $\m(20_p')$. Then by (dim H,om) $\dim\Hom(\m(24_p),\m(15_q))=0$ and since $24_p$ is the next weight to the right of $\m(15_q)$, then $[\m(15_q),\el(24_p)]=0$. Likewise $[\m(30_p),\el(24_p)]=0$. Also, $[\m(30_p),\el(15_q)]=0$ since $h_c(30_p)=h_c(15_q)$. That completes the rows $15_q$ and $30_p$ of $\mathcal{M}$ and also implies that $[\m(20_p):\el(24_p)]=1$.

From row $6_p'$ of $\mathcal{M}$, $1_p'$ doesn't consist of singular vectors in $\m(6_p')$. Then by (dim Hom), neither does $6_p$ in $\m(1_p')$. Since $6_p$ is the rep immediately to the right of $1_p$ on the weight line, the decomposition of $\el(1_p)$ into Vermas skips $6_p$ and $[\m(1_p),\el(6_p)]=0$. From row $20_p'$ it then follows that $\dim\Hom(\m(20_p),\m(1_p))=\dim\Hom(\m(20_p),\m(6_p))=1$ and $[\m(1_p):\el(20_p)]=[\m(6_p):\el(20_p)]=1.$ 

Now it only remains to find the remaining five entries of $\mathcal{M}$ in columns $15_q,\;24_p$ and rows $1_p,\;6_p,\;20_p$. Row $15_q'$ shows that 
\begin{align*}&\dim\Hom(\m(6_p'),\m(15_q'))=0=\dim\Hom(\m(15_q),\m(6_p))\\&\dim\Hom(\m(20_p'),\m(15_q'))=0=\dim\Hom(\m(15_q),\m(20_p))\\&\dim\Hom(\m(1_p'),\m(15_q'))=1=\dim\Hom(\m(15_q),\m(1_p))\end{align*}
So $[\m(20_p):\el(15_q)]=0$, and in the decomposition $\el(6_p)=\m(6_p)-\m(20_p)+c_{15}\m(15_q)...$, the only possibility is $c_{15}=0$, and thus $[\el(6_p):\m(15_q)]=0$. It then follows that $[\m(6_p):\el(15_q)]=0$ since the dot product of row $6_p$ of $\mathcal{M}$ with column $15_q$ of $\mathcal{L}$ must be $0$. Since $\el(1)=\m(1_p)-\m(20_p)+\m(30_p)...$ and $15_q$ consists of singular vectors in $\m(1_p)$ but not in $\m(20_p)$, the decomposition continues: $\el(1)=\m(1_p)-\m(20_p)+\m(30_p)-\m(15_q)...$, we have $[\el(1_p):\m(15_q)]=-1$ and $[\m(1_p):\el(15_q)]=1.$ 

So the only question left is what are $[\m(1_p):\el(24_p)]$ and $[\m(6_p):\el(24_p)]$? Since we know the decomposition of $\el(1_p)$ up through Vermas of $h_c$-weight $0$, we may calculate its graded character which is symmetric in $t$ and $t^{-1}$:
\begin{align*}
\chi_{\el(1_p)}(t)&=t^{-3}+t^3+6(t^{-2}+t^2)+21(t^{-1}+t)+36\\=&\frac{t^{-3}-20+15t+84t^2-160t^3+84t^4+15t^5-20t^6+t^9}{(1-t)^6}
\end{align*}
Seeing that the coefficient of $t^2$ is $84$ and $[\el(1_p):\m(60_p)]=1$, it must be that \\$[\el(1_p):\m(24_p)]=1$. From the condition that row $1_p$ of $\mathcal{M}$ dotted with column $24_p$ of $\mathcal{L}$ must be $0$, we get $[\el(1_p):\m(24_p)]=0$. Next, $[\el(6_p):\m(24_p)]$ must be at least $1$ or else $\dim\el(6)[2]<0$. On the other hand, $24_p$ only occurs once as a subspace of singular vectors in $\m(20_p)$, so $[\el(6_p):\m(24_p)]$ is at most $1$. So $[\el(6_p):\m(24_p)]=1$, and so $[\m(6_p):\el(24_p)]=0$.

This completes the calculation of the decomposition matrix $\mathcal{M}$. Inverting $\mathcal{M}$ we obtain $\mathcal{L}$ and the characters of all irreps, from which we find the supports of all $\el(\tau)$ such that KZ$(\el(\tau))=0$, and in particular, the characters and dimensions of the finite-dim irreps:

\begin{align*}
%&\el(1_p)=\m(1_p)-\m(20_p)+\m(30_p)-\m(15_q)+\m(24_p)+\m(60_p)-2\m(80_s)\\&\qquad\qquad\qquad+\m(24_p')+\m(60_p')+\m(30_p')-\m(15_q')-\m(20_p')+\m(1_p')\\
&\chi_{\el(1_p)}(t)=t^{-3}+t^3+6(t^{-2}+t^2)+21(t^{-1}+t)+36\\
&\dim\el(1_p)=92\\
\\
%&\el(6_p)=\m(6_p)-\m(20_p)+\m(24_p)+\m(60_p)-\m(80_s)-\m(60_s)+\m(60_p')+\m(24_p')\\
%&\qquad\qquad-\m(20_p')+\m(6_p')\\
&\chi_{\el(6_p)}(t)=6(t^{-1}+t)+16\\
&\dim\el(6_p)=28\\
\end{align*}

\begin{itemize}
\item$\dim\Supp\el(20_p)=2$
\item$\dim\Supp\el(15_q)=1$
\item$\dim\Supp\el(24_p)=2$
\end{itemize}

\begin{comment} 
&\el(20_p)=\m(20_p)-\m(30_p)-\m(24_p)-\m(60_p)+2\m(80_s)+\m(60_s)-\m(24_p')-2\m(60_p')\\
&\qquad\qquad-\m(30_p')+\m(15_q')+2\m(20_p')-\m(6_p')-\m(1_p')\\
&\chi_{\el(20_p)}(t)=t^{\color{red}h_c(20_p)}\frac{-t^5 - 4t^4 - 16t^3 - 4t^2 + 50t + 20}{(1-t)^2}\\
\\
&\el(15_q)=\m(15_q)-\m(60_p)+\m(80_s)-\m(24_p')-\m(30_p')+\m(20_p')-\m(1_p')\\
&\chi_{\el(15_q)}(t)=t^{\color{red}h_c(15_q)}\frac{t^3 + 5t^2 + 15t + 15}{1-t}\\
\\
&\el(24_p)=\m(24_p)-\m(80_s)+\m(24_p')+\m(60_p')-\m(15_q')-\m(20_p')+\m(6_p')+\m(1_p')\\
&\chi_{\el(24_p)}(t)=t^{\color{red}h_c(24_p)}\frac{t^3 + 4t^2 + 16t + 24}{(1-t)^2}\\
\end{align*}
\end{comment}

\end{proof}

\subsection{{\Large$\mbf{c=\frac{1}{5}}$}}

There are two dual blocks which KZ maps to blocks of defect $1$ in the Hecke algebra. Thus there are two reps of less than full support, and their lowest weights are the trivial rep and the standard rep:
\begin{align*}
&\begin{picture}(360,30)
\put(0,0){\line(1,0){360}}
\put(0,0){\circle*{5}}
\put(180,0){\circle*{5}}
\put(240,0){\circle*{5}}
\put(270,0){\circle*{5}}
\put(360,0){\circle*{5}}
\put(-12,15){$-4.2$}
\put(175,15){$1.8$}
\put(235,15){$3.8$}
\put(265,15){$4.8$}
\put(355,15){$7.8$}
\put(-2,-15){$1_p$}
\put(175,-15){$24_p$}
\put(235,-15){$81_p'$}
\put(265,-15){$64_p'$}
\put(356,-15){$6_p'$}
\end{picture}\\
%\\
%&\el(1_p)=\m(1_p)-\m(24_p)+\m(81_p')-\m(64_p')+\m(6_p')\\
%&\chi_{\el(1_p)}(t)=t^{-4.2}\frac{1+4t+10t^2+20t^3+35t^4+56t^5+60t^6+24t^7+6t^8}{(1-t)^2}\\
%&\dim\Supp\el(1_p)=2\\
\\
&\begin{picture}(360,30)
\put(0,0){\line(1,0){360}}
\put(0,0){\circle*{5}}
\put(90,0){\circle*{5}}
\put(120,0){\circle*{5}}
\put(180,0){\circle*{5}}
\put(360,0){\circle*{5}}
\put(-12,15){$-1.8$}
\put(85,15){$1.2$}
\put(115,15){$2.2$}
\put(175,15){$4.2$}
\put(351,15){$10.2$}
\put(-2,-15){$6_p$}
\put(85,-15){$64_p$}
\put(115,-15){$81_p$}
\put(175,-15){$24_p'$}
\put(356,-15){$1_p'$}
\end{picture}\\
\\
\\
%&\el(6_p)=\m(6_p)-\m(64_p)+\m(81_p)-\m(24_p')+\m(1_p')\\
%&\chi_{\el(6_p)}(t)=t^{-1.8}\frac{6+24t+60t^2+56t^3+35t^4+20t^5+10t^6+4t^7+t^8}{(1-t)^2}\\
&\bullet\;\dim\Supp\el(1_p)=2\\
&\bullet\;\dim\Supp\el(6_p)=2\\
\end{align*}

\subsection{{\Large$\mbf{c=\frac{1}{4}}$}} There are two blocks of nonzero defect: the principal block, and a block of defect $1$ anchored by the minimal support irrep $\el(20_p)$. 

\theorem The decomposition matrix $\big([M(\tau):L(\sigma)]\big)$ of the principal block is

\begin{center}
\[
\begin{blockarray}{ccccccccccccccc}
\begin{block}{cc(ccccccccccccc)}
\mathtt{(2)} & 1_p &1 & \cdot & -1 & \cdot & \cdot & \cdot & 1 & \cdot & -1 & \cdot & 1 & \cdot & \cdot \\
\mathtt{(2)} & 6_p & \cdot & 1 & -1 & \cdot & -1 & 1 & 1 & 1 & -1 & -1 & \cdot & 1 & \cdot\\
\mathtt{(3)} & 15_q & \cdot & \cdot & 1 & \cdot & \cdot & \cdot & -1 & -1 & 1 & 1 & -1 & -1 & \cdot\\
\mathtt{(2)} & 15_p & \cdot & \cdot & \cdot & 1 & -1 & \cdot & 1 & \cdot & \cdot & -1 & \cdot & \cdot & 1\\
\mathtt{(3)} & 81_p & \cdot & \cdot & \cdot & \cdot & 1 & -1 & -1 & \cdot & 1 & 1 & \cdot & -1 & -1\\
& 90_s & \cdot & \cdot & \cdot & \cdot & \cdot & 1 & \cdot & \cdot & -1 & \cdot & \cdot & 1 & \cdot\\
& 80_s & \cdot & \cdot & \cdot & \cdot & \cdot & \cdot & 1 & \cdot & -1 & -1 & 1 & 1 & 1\\
& 10_s & \cdot & \cdot & \cdot & \cdot & \cdot & \cdot & \cdot & 1 & \cdot & -1 & \cdot & 1 & \cdot\\
& 81_p' & \cdot & \cdot & \cdot & \cdot & \cdot & \cdot & \cdot & \cdot & 1 & \cdot & -1 & -1 & \cdot\\
& 15_q' & \cdot & \cdot & \cdot & \cdot & \cdot & \cdot & \cdot & \cdot & \cdot & 1 & \cdot & -1 & -1\\
& 15_p' & \cdot & \cdot & \cdot & \cdot & \cdot & \cdot & \cdot & \cdot & \cdot & \cdot & 1 & \cdot & \cdot\\
& 6_p' & \cdot & \cdot & \cdot & \cdot & \cdot & \cdot & \cdot & \cdot & \cdot & \cdot & \cdot & 1 & \cdot\\
& 1_p'& \cdot & \cdot & \cdot & \cdot & \cdot & \cdot & \cdot & \cdot & \cdot & \cdot & \cdot & \cdot & 1\\
\end{block}
\end{blockarray}
\]
\\
\[
\begin{blockarray}{ccccccccccccccc}
\begin{block}{cc(ccccccccccccc)}
\mathtt{(2)} & 1_p &1 & \cdot & 1 & \cdot & \cdot & \cdot & \cdot & 1 & \cdot & \cdot & \cdot & \cdot & \cdot\\
\mathtt{(2)} & 6_p & \cdot & 1 & 1 & \cdot & 1 & \cdot & 1 & \cdot & \cdot & \cdot & \cdot & \cdot & \cdot\\
\mathtt{(3)} & 15_q & \cdot & \cdot & 1 & \cdot & \cdot & \cdot & 1 & 1 & \cdot & 1 & \cdot & \cdot & \cdot\\
\mathtt{(2)} & 15_p & \cdot & \cdot & \cdot & 1 & 1 & 1 & \cdot & \cdot & \cdot & \cdot & \cdot & \cdot & \cdot\\
\mathtt{(3)} & 81_p & \cdot & \cdot & \cdot & \cdot & 1 & 1 & 1 & \cdot & 1 & \cdot & \cdot & \cdot & \cdot\\
& 90_s & \cdot & \cdot & \cdot & \cdot & \cdot & 1 & \cdot & \cdot & 1 & \cdot & 1 & \cdot & \cdot\\
& 80_s & \cdot & \cdot & \cdot & \cdot & \cdot & \cdot & 1 & \cdot & 1 & 1 & \cdot & 1 & \cdot\\
& 10_s & \cdot & \cdot & \cdot & \cdot & \cdot & \cdot & \cdot & 1 & \cdot & 1 & \cdot & \cdot & 1\\
& 81_p' & \cdot & \cdot & \cdot & \cdot & \cdot & \cdot & \cdot & \cdot & 1 & \cdot & 1 & 1 & \cdot\\
& 15_q' & \cdot & \cdot & \cdot & \cdot & \cdot & \cdot & \cdot & \cdot & \cdot & 1 & \cdot & 1 & 1\\
& 15_p' & \cdot & \cdot & \cdot & \cdot & \cdot & \cdot & \cdot & \cdot & \cdot & \cdot & 1 & \cdot & \cdot\\
& 6_p' & \cdot & \cdot & \cdot & \cdot & \cdot & \cdot & \cdot & \cdot & \cdot & \cdot & \cdot & 1 & \cdot\\
& 1_p'& \cdot & \cdot & \cdot & \cdot & \cdot & \cdot & \cdot & \cdot & \cdot & \cdot & \cdot & \cdot & 1\\
\end{block}
\end{blockarray}
\]
\end{center}

\begin{proof}

All but the five leftmost columns come from the decomposition matrix of the Hecke algebra of $E_6$ with $e=4$. Observe that $\m(15_p')$ and $\m(6_p')$ are simple; by (dim Hom), $6_p$ and $15_p$ can't consist of singular vectors in $\m(1_p)$, nor can $15_p$ in $\m(6_p)$. Since $\el(15_q')=\m(15_q')-\m(6_p')-\m(1_p')$ and all Vermas to the right of $15_q$ are simple, $\dim\Hom(\m(15_q),\m(6_p))=\dim\Hom(\m(6_p'),\m(15_q'))=1$ and likewise  \\$\dim\Hom(\m(15_q),\m(1_p))=1$. This gives column $15_q$. Since $15_q'$ doesn't generate a subrep of $\m(81_p')$, and the remaining Vermas to the right are simple, column $81_p$ follows by the same reasoning.

The  weight line for this block together with the inverse of the decomposition matrix then gives the characters of the simples in this block:
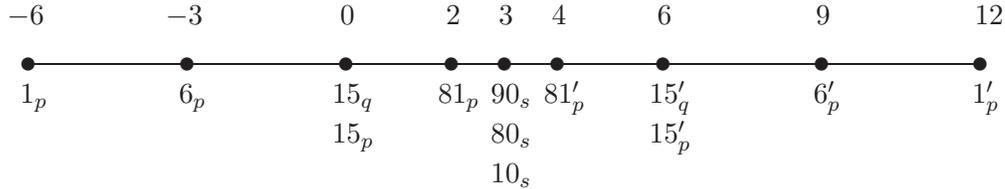
\begin{center}
\begin{picture}(360,75)
\put(0,50){\line(1,0){360}}
\put(0,50){\circle*{5}}
\put(60,50){\circle*{5}}
\put(120,50){\circle*{5}}
\put(160,50){\circle*{5}}
\put(180,50){\circle*{5}}
\put(200,50){\circle*{5}}
\put(240,50){\circle*{5}}
\put(300,50){\circle*{5}}
\put(360,50){\circle*{5}}
\put(-8,65){$-6$}
\put(52,65){$-3$}
\put(118,65){$0$}
\put(158,65){$2$}
\put(178,65){$3$}
\put(198,65){$4$}
\put(238,65){$6$}
\put(298,65){$9$}
\put(358,65){$12$}
\put(-3,35){$1_p$}
\put(57,35){$6_p$}
\put(115,35){$15_q$}
\put(115,20){$15_p$}
\put(155,35){$81_p$}
\put(175,35){$90_s$}
\put(175,20){$80_s$}
\put(175,5){$10_s$}
\put(195,35){$81_p'$}
\put(235,35){$15_q'$}
\put(235,20){$15_p'$}
\put(297,35){$6_p'$}
\put(357,35){$1_p'$}
\end{picture}\\
\end{center}

After calculating characters, we conclude that the dimensions of support of the simples $\el(\tau)$ such that $\mathrm{KZ}(\el(\tau))=0$ are all $2$ or $3$, as printed next to the corresponding rows of the matrices $\mathcal{L}$ and $\mathcal{M}$. 
\end{proof}

In addition to the principal block, there's a single block of defect $1$:

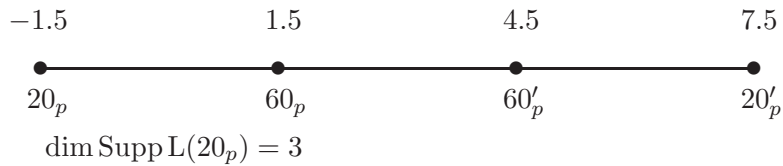
\begin{align*}
&\begin{picture}(270,30)
\put(0,0){\line(1,0){270}}
\put(0,0){\circle*{5}}
\put(90,0){\circle*{5}}
\put(180,0){\circle*{5}}
\put(270,0){\circle*{5}}
\put(-12,15){$-1.5$}
\put(85,15){$1.5$}
\put(175,15){$4.5$}
\put(265,15){$7.5$}
\put(-5,-15){$20_p$}
\put(85,-15){$60_p$}
\put(175,-15){$60_p'$}
\put(265,-15){$20_p'$}
\end{picture}\\
\\
%&\el(20_p)=\m(20_p)-\m(60_p)+\m(60_p')-\m(20_p')\\
%&\chi_{\el(20_p)}(t)=t^{-1.5}\frac{20+60t+120t^2+140t^3+60t^4+20t^6}{(1-t)^3}\\
&\dim\Supp\el(20_p)=3\\
\end{align*}

\subsection{{\Large$\mbf{c=\frac{1}{3}}$}} There is a single big block containing all irreducible representations of $H_{\frac{1}{3}}(E_6)$, save for $\el(81_p)$ and $\el(81_p')$ each of which sits in a defect zero block by itself. 

\theorem The matrices below are the decomposition matrix of the principal block of $H_{\frac{1}{3}}(E_6)$ and its inverse, decorated with the dimensions of support of the irreps killed by the KZ functor. The dimensions of support occurring in this block are $0$, $2$, $4$, and $6$.

\begin{center}
\[
\begin{blockarray}{cccccccccccccccccccccccc}
\begin{block}{cc(cccccccccccccccccccccc)}
\star&1_p&1&1&1&\cdot&1&\cdot&\cdot&\cdot&\cdot&\cdot&\cdot&\cdot&\cdot&\cdot&\cdot&\cdot&\cdot&1&\cdot&\cdot&\cdot&\cdot\\
\star&6_p&\cdot&1&1&1&\cdot&1&1&\cdot&\cdot&\cdot&\cdot&\cdot&1&\cdot&1&\cdot&\cdot&1&\cdot&\cdot&\cdot&\cdot\\
\mathtt{(4)}&20_p&\cdot&\cdot&1&\cdot&1&\cdot&1&1&1&\cdot&1&\cdot&1&1&1&\cdot&\cdot&1&\cdot&\cdot&\cdot&\cdot\\
\mathtt{(2)}&30_p&\cdot&\cdot&\cdot&1&\cdot&\cdot&1&\cdot&1&1&\cdot&\cdot&1&\cdot&1&\cdot&\cdot&\cdot&\cdot&\cdot&\cdot&\cdot\\
\mathtt{(2)}&15_q&\cdot&\cdot&\cdot&\cdot&1&\cdot&\cdot&\cdot&1&\cdot&1&\cdot&\cdot&1&\cdot&\cdot&\cdot&\cdot&\cdot&\cdot&\cdot&1\\
\star&15_p&\cdot&\cdot&\cdot&\cdot&\cdot&1&1&\cdot&\cdot&\cdot&\cdot&1&\cdot&1&1&1&\cdot&\cdot&\cdot&\cdot&\cdot&\cdot\\
\mathtt{(4)}&64_p&\cdot&\cdot&\cdot&\cdot&\cdot&\cdot&1&1&1&1&1&1&\cdot&1&1&1&\cdot&\cdot&1&\cdot&\cdot&\cdot\\
\mathtt{(4)}&24_p&\cdot&\cdot&\cdot&\cdot&\cdot&\cdot&\cdot&1&\cdot&\cdot&1&\cdot&\cdot&\cdot&\cdot&\cdot&\cdot&\cdot&1&\cdot&\cdot&\cdot\\
\mathtt{(4)}&60_p&\cdot&\cdot&\cdot&\cdot&\cdot&\cdot&\cdot&\cdot&1&1&1&\cdot&1&\cdot&1&\cdot&\cdot&\cdot&1&\cdot&1&1\\
&80_s&\cdot&\cdot&\cdot&\cdot&\cdot&\cdot&\cdot&\cdot&\cdot&1&\cdot&\cdot&\cdot&\cdot&1&1&1&\cdot&1&\cdot&1&\cdot\\
&60_s&\cdot&\cdot&\cdot&\cdot&\cdot&\cdot&\cdot&\cdot&\cdot&\cdot&1&\cdot&\cdot&1&1&1&\cdot&1&1&1&1&1\\
\mathtt{(4)}&20_s&\cdot&\cdot&\cdot&\cdot&\cdot&\cdot&\cdot&\cdot&\cdot&\cdot&\cdot&1&\cdot&\cdot&\cdot&1&\cdot&\cdot&1&\cdot&\cdot&\cdot\\
\mathtt{(2)}&10_s&\cdot&\cdot&\cdot&\cdot&\cdot&\cdot&\cdot&\cdot&\cdot&\cdot&\cdot&\cdot&1&\cdot&1&\cdot&\cdot&\cdot&\cdot&\cdot&1&\cdot\\
\mathtt{(2)}&24_p'&\cdot&\cdot&\cdot&\cdot&\cdot&\cdot&\cdot&\cdot&\cdot&\cdot&\cdot&\cdot&\cdot&1&\cdot&1&\cdot&\cdot&\cdot&1&\cdot&\cdot\\
&60_p'&\cdot&\cdot&\cdot&\cdot&\cdot&\cdot&\cdot&\cdot&\cdot&\cdot&\cdot&\cdot&\cdot&\cdot&1&1&1&1&\cdot&1&1&\cdot\\
&64_p'&\cdot&\cdot&\cdot&\cdot&\cdot&\cdot&\cdot&\cdot&\cdot&\cdot&\cdot&\cdot&\cdot&\cdot&\cdot&1&1&\cdot&1&1&1&\cdot\\
&30_p'&\cdot&\cdot&\cdot&\cdot&\cdot&\cdot&\cdot&\cdot&\cdot&\cdot&\cdot&\cdot&\cdot&\cdot&\cdot&\cdot&1&\cdot&\cdot&\cdot&1&\cdot\\
&15_q'&\cdot&\cdot&\cdot&\cdot&\cdot&\cdot&\cdot&\cdot&\cdot&\cdot&\cdot&\cdot&\cdot&\cdot&\cdot&\cdot&\cdot&1&\cdot&1&\cdot&\cdot\\
&15_p'&\cdot&\cdot&\cdot&\cdot&\cdot&\cdot&\cdot&\cdot&\cdot&\cdot&\cdot&\cdot&\cdot&\cdot&\cdot&\cdot&\cdot&\cdot&1&\cdot&1&\cdot\\
&20_p'&\cdot&\cdot&\cdot&\cdot&\cdot&\cdot&\cdot&\cdot&\cdot&\cdot&\cdot&\cdot&\cdot&\cdot&\cdot&\cdot&\cdot&\cdot&\cdot&1&1&1\\
&6_p'&\cdot&\cdot&\cdot&\cdot&\cdot&\cdot&\cdot&\cdot&\cdot&\cdot&\cdot&\cdot&\cdot&\cdot&\cdot&\cdot&\cdot&\cdot&\cdot&\cdot&1&1\\
&1_p'&\cdot&\cdot&\cdot&\cdot&\cdot&\cdot&\cdot&\cdot&\cdot&\cdot&\cdot&\cdot&\cdot&\cdot&\cdot&\cdot&\cdot&\cdot&\cdot&\cdot&\cdot&1\\
\end{block}
\end{blockarray}
\]
\end{center}

\newgeometry{top=2cm, bottom=3cm, left=4cm, right=4cm}
\begin{turn}{270}
\begin{minipage}{\linewidth}
\[
\begin{blockarray}{cccccccccccccccccccccccc}
&&1_p&6_p&20_p&30_p&15_q&15_p&64_p&24_p&60_p&80_s&60_s&20_s&10_s&24_p'&60_p'&64_p'&30_p'&15_q'&15_p'&20_p'&6_p'&1_p'\\
\begin{block}{cc(cccccccccccccccccccccc)}
\star&1_p&1&-1&\cdot&1&-1&1&-1&1&1&-1&\cdot&\cdot&-1&1&1&-1&1&-1&1&\cdot&-1&1\\
\star&6_p&\cdot&1&-1&-1&1&-1&2&-1&-1&\cdot&\cdot&-1&2&-1&-1&2&-1&1&-1&-1&1&\cdot\\
\mathtt{(4)}&20_p&\cdot&\cdot&1&\cdot&-1&\cdot&-1&\cdot&1&\cdot&\cdot&1&-2&1&1&-2&1&-2&1&2&-2&\cdot\\
\mathtt{(2)}&30_p&\cdot&\cdot&\cdot&1&\cdot&\cdot&-1&1&\cdot&\cdot&\cdot&1&-1&1&1&-2&1&-1&1&1&-1&\cdot\\
\mathtt{(2)}&15_q&\cdot&\cdot&\cdot&\cdot&1&\cdot&\cdot&\cdot&-1&1&\cdot&\cdot&1&-1&-1&1&-1&1&-1&\cdot&1&-1\\
\star&15_p&\cdot&\cdot&\cdot&\cdot&\cdot&1&-1&1&1&\cdot&-1&\cdot&-1&1&1&-1&\cdot&\cdot&1&\cdot&\cdot&\cdot\\
\mathtt{(4)}&64_p&\cdot&\cdot&\cdot&\cdot&\cdot&\cdot&1&-1&-1&\cdot&1&-1&1&-2&-2&3&-1&1&-2&-1&2&-1\\
\mathtt{(4)}&24_p&\cdot&\cdot&\cdot&\cdot&\cdot&\cdot&\cdot&1&\cdot&\cdot&-1&\cdot&\cdot&1&1&-1&\cdot&\cdot&1&\cdot&\cdot&1\\
\mathtt{(4)}&60_p&\cdot&\cdot&\cdot&\cdot&\cdot&\cdot&\cdot&\cdot&1&-1&-1&\cdot&-1&1&2&-1&\cdot&-1&2&\cdot&-1&1\\
&80_s&\cdot&\cdot&\cdot&\cdot&\cdot&\cdot&\cdot&\cdot&\cdot&1&\cdot&\cdot&\cdot&\cdot&-1&\cdot&\cdot&1&-1&\cdot&1&-1\\
&60_s&\cdot&\cdot&\cdot&\cdot&\cdot&\cdot&\cdot&\cdot&\cdot&\cdot&1&\cdot&\cdot&-1&-1&1&\cdot&\cdot&-2&\cdot&1&-2\\
\mathtt{(4)}&20_s&\cdot&\cdot&\cdot&\cdot&\cdot&\cdot&\cdot&\cdot&\cdot&\cdot&\cdot&1&\cdot&\cdot&\cdot&-1&1&\cdot&\cdot&1&-1&\cdot\\
\mathtt{(2)}&10_s&\cdot&\cdot&\cdot&\cdot&\cdot&\cdot&\cdot&\cdot&\cdot&\cdot&\cdot&\cdot&1&\cdot&-1&1&\cdot&1&-1&-1&1&\cdot\\
\mathtt{(2)}&24_p'&\cdot&\cdot&\cdot&\cdot&\cdot&\cdot&\cdot&\cdot&\cdot&\cdot&\cdot&\cdot&\cdot&1&\cdot&-1&1&\cdot&1&\cdot&-1&1\\
&60_p'&\cdot&\cdot&\cdot&\cdot&\cdot&\cdot&\cdot&\cdot&\cdot&\cdot&\cdot&\cdot&\cdot&\cdot&1&-1&\cdot&-1&1&1&-2&1\\
&64_p'&\cdot&\cdot&\cdot&\cdot&\cdot&\cdot&\cdot&\cdot&\cdot&\cdot&\cdot&\cdot&\cdot&\cdot&\cdot&1&-1&\cdot&-1&-1&2&-1\\
&30_p'&\cdot&\cdot&\cdot&\cdot&\cdot&\cdot&\cdot&\cdot&\cdot&\cdot&\cdot&\cdot&\cdot&\cdot&\cdot&\cdot&1&\cdot&\cdot&\cdot&-1&1\\
&15_q'&\cdot&\cdot&\cdot&\cdot&\cdot&\cdot&\cdot&\cdot&\cdot&\cdot&\cdot&\cdot&\cdot&\cdot&\cdot&\cdot&\cdot&1&\cdot&-1&1&\cdot\\
&15_p'&\cdot&\cdot&\cdot&\cdot&\cdot&\cdot&\cdot&\cdot&\cdot&\cdot&\cdot&\cdot&\cdot&\cdot&\cdot&\cdot&\cdot&\cdot&1&\cdot&-1&1\\
&20_p'&\cdot&\cdot&\cdot&\cdot&\cdot&\cdot&\cdot&\cdot&\cdot&\cdot&\cdot&\cdot&\cdot&\cdot&\cdot&\cdot&\cdot&\cdot&\cdot&1&-1&\cdot\\
&6_p'&\cdot&\cdot&\cdot&\cdot&\cdot&\cdot&\cdot&\cdot&\cdot&\cdot&\cdot&\cdot&\cdot&\cdot&\cdot&\cdot&\cdot&\cdot&\cdot&\cdot&1&-1\\
&1_p'&\cdot&\cdot&\cdot&\cdot&\cdot&\cdot&\cdot&\cdot&\cdot&\cdot&\cdot&\cdot&\cdot&\cdot&\cdot&\cdot&\cdot&\cdot&\cdot&\cdot&\cdot&1\\
\end{block}
\end{blockarray}
\]
\end{minipage}
\end{turn}

\restoregeometry

\begin{proof}

For reference throughout the calculations, here is the weight line for this block:

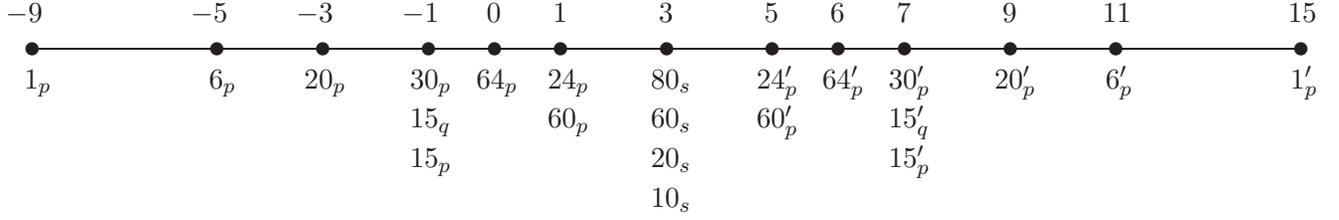
\begin{center}
\begin{picture}(500,90)
\put(0,70){\line(1,0){480}}
\put(0,70){\circle*{5}}
\put(70,70){\circle*{5}}
\put(110,70){\circle*{5}}
\put(150,70){\circle*{5}}
\put(175,70){\circle*{5}}
\put(200,70){\circle*{5}}
\put(240,70){\circle*{5}}
\put(280,70){\circle*{5}}
\put(305,70){\circle*{5}}
\put(330,70){\circle*{5}}
\put(370,70){\circle*{5}}
\put(410,70){\circle*{5}}
\put(480,70){\circle*{5}}
\put(-10,80){$-9$}
\put(60,80){$-5$}
\put(100,80){$-3$}
\put(140,80){$-1$}
\put(172,80){$0$}
\put(197,80){$1$}
\put(237,80){$3$}
\put(277,80){$5$}
\put(302,80){$6$}
\put(327,80){$7$}
\put(367,80){$9$}
\put(405,80){$11$}
\put(475,80){$15$}
\put(-3,55){$1_p$}
\put(67,55){$6_p$}
\put(103,55){$20_p$}
\put(143,55){$30_p$}
\put(143,40){$15_q$}
\put(143,25){$15_p$}
\put(168,55){$64_p$}
\put(195,55){$24_p$}
\put(195,40){$60_p$}
\put(234,55){$80_s$}
\put(234,40){$60_s$}
\put(234,25){$20_s$}
\put(234,10){$10_s$}
\put(274,55){$24_p'$}
\put(274,40){$60_p'$}
\put(299,55){$64_p'$}
\put(324,55){$30_p'$}
\put(324,40){$15_q'$}
\put(324,25){$15_p'$}
\put(364,55){$20_p'$}
\put(406,55){$6_p'$}
\put(476,55){$1_p'$}
\end{picture}
\end{center}

Importing the columns from the Hecke algebra decomposition matrix of $E_6$ for $e=3$ in \cite{GJ} and then applying (dim Hom) brings the decomposition matrix of $\oh_{\frac{1}{3}}(E_6)$ to the following state:
\footnotesize
\[
\begin{blockarray}{cccccccccccccccccccccccc}
&&1_p&6_p&20_p&30_p&15_q&15_p&64_p&24_p&60_p&80_s&60_s&20_s&10_s&24_p'&60_p'&64_p'&30_p'&15_q'&15_p'&20_p'&6_p'&1_p'\\
\begin{block}{cc(cccccccccccccccccccccc)}
\bullet&1_p&1&1&{\color{red}?}&{\color{red}?}&{\color{red}?}&{\color{red}?}&{\color{red}?}&{\color{red}?}&{\color{red}?}&\cdot&\cdot&{\color{red}?}&{\color{red}?}&{\color{red}?}&\cdot&\cdot&\cdot&1&\cdot&\cdot&\cdot&\cdot\\
\bullet&6_p&\cdot&1&1&1&{\color{red}?}&1&{\color{red}?}&{\color{red}?}&{\color{red}?}&\cdot&\cdot&{\color{red}?}&{\color{red}?}&{\color{red}?}&1&\cdot&\cdot&1&\cdot&\cdot&\cdot&\cdot\\
\bullet&20_p&\cdot&\cdot&1&\cdot&1&\cdot&1&{\color{red}?}&{\color{red}?}&\cdot&1&{\color{red}?}&{\color{red}?}&{\color{red}?}&1&\cdot&\cdot&1&\cdot&\cdot&\cdot&\cdot\\
\bullet&30_p&\cdot&\cdot&\cdot&1&\cdot&\cdot&1&{\color{red}?}&{\color{red}?}&1&\cdot&{\color{red}?}&{\color{red}?}&{\color{red}?}&1&\cdot&\cdot&\cdot&\cdot&\cdot&\cdot&\cdot\\
\bullet&15_q&\cdot&\cdot&\cdot&\cdot&1&\cdot&\cdot&\cdot&1&\cdot&1&\cdot&{\color{red}?}&{\color{red}?}&\cdot&\cdot&\cdot&\cdot&\cdot&\cdot&\cdot&1\\
\bullet&15_p&\cdot&\cdot&\cdot&\cdot&\cdot&1&1&{\color{red}?}&{\color{red}?}&\cdot&\cdot&{\color{red}?}&{\color{red}?}&{\color{red}?}&1&1&\cdot&\cdot&\cdot&\cdot&\cdot&\cdot\\
\bullet&64_p&\cdot&\cdot&\cdot&\cdot&\cdot&\cdot&1&1&1&1&1&{\color{red}?}&{\color{red}?}&{\color{red}?}&1&1&\cdot&\cdot&1&\cdot&\cdot&\cdot\\
\bullet&24_p&\cdot&\cdot&\cdot&\cdot&\cdot&\cdot&\cdot&1&\cdot&\cdot&1&{\color{red}?}&{\color{red}?}&{\color{red}?}&\cdot&\cdot&\cdot&\cdot&1&\cdot&\cdot&\cdot\\
\bullet&60_p&\cdot&\cdot&\cdot&\cdot&\cdot&\cdot&\cdot&\cdot&1&1&1&\cdot&1&{\color{red}?}&1&\cdot&\cdot&\cdot&1&\cdot&1&1\\
&80_s&\cdot&\cdot&\cdot&\cdot&\cdot&\cdot&\cdot&\cdot&\cdot&1&\cdot&\cdot&\cdot&\cdot&1&1&1&\cdot&1&\cdot&1&\cdot\\
&60_s&\cdot&\cdot&\cdot&\cdot&\cdot&\cdot&\cdot&\cdot&\cdot&\cdot&1&\cdot&\cdot&1&1&1&\cdot&1&1&1&1&1\\
\bullet&20_s&\cdot&\cdot&\cdot&\cdot&\cdot&\cdot&\cdot&\cdot&\cdot&\cdot&\cdot&1&\cdot&{\color{red}?}&\cdot&1&\cdot&\cdot&1&\cdot&\cdot&\cdot\\
\bullet&10_s&\cdot&\cdot&\cdot&\cdot&\cdot&\cdot&\cdot&\cdot&\cdot&\cdot&\cdot&\cdot&1&{\color{red}?}&1&\cdot&\cdot&\cdot&\cdot&\cdot&1&\cdot\\
\bullet&24_p'&\cdot&\cdot&\cdot&\cdot&\cdot&\cdot&\cdot&\cdot&\cdot&\cdot&\cdot&\cdot&\cdot&1&\cdot&1&\cdot&\cdot&\cdot&1&\cdot&\cdot\\
&60_p'&\cdot&\cdot&\cdot&\cdot&\cdot&\cdot&\cdot&\cdot&\cdot&\cdot&\cdot&\cdot&\cdot&\cdot&1&1&1&1&\cdot&1&1&\cdot\\
&64_p'&\cdot&\cdot&\cdot&\cdot&\cdot&\cdot&\cdot&\cdot&\cdot&\cdot&\cdot&\cdot&\cdot&\cdot&\cdot&1&1&\cdot&1&1&1&\cdot\\
&30_p'&\cdot&\cdot&\cdot&\cdot&\cdot&\cdot&\cdot&\cdot&\cdot&\cdot&\cdot&\cdot&\cdot&\cdot&\cdot&\cdot&1&\cdot&\cdot&\cdot&1&\cdot\\
&15_q'&\cdot&\cdot&\cdot&\cdot&\cdot&\cdot&\cdot&\cdot&\cdot&\cdot&\cdot&\cdot&\cdot&\cdot&\cdot&\cdot&\cdot&1&\cdot&1&\cdot&\cdot\\
&15_p'&\cdot&\cdot&\cdot&\cdot&\cdot&\cdot&\cdot&\cdot&\cdot&\cdot&\cdot&\cdot&\cdot&\cdot&\cdot&\cdot&\cdot&\cdot&1&\cdot&1&\cdot\\
&20_p'&\cdot&\cdot&\cdot&\cdot&\cdot&\cdot&\cdot&\cdot&\cdot&\cdot&\cdot&\cdot&\cdot&\cdot&\cdot&\cdot&\cdot&\cdot&\cdot&1&1&1\\
&6_p'&\cdot&\cdot&\cdot&\cdot&\cdot&\cdot&\cdot&\cdot&\cdot&\cdot&\cdot&\cdot&\cdot&\cdot&\cdot&\cdot&\cdot&\cdot&\cdot&\cdot&1&1\\
&1_p'&\cdot&\cdot&\cdot&\cdot&\cdot&\cdot&\cdot&\cdot&\cdot&\cdot&\cdot&\cdot&\cdot&\cdot&\cdot&\cdot&\cdot&\cdot&\cdot&\cdot&\cdot&1\\
\end{block}
\end{blockarray}
\]
\normalsize

The way to acquire all the missing decomposition numbers systematically and efficiently is to weave in and out of picking up Verma-decompositions of simples and dropping off the new entries they make available into the decomposition matrix above, the overall route going from the bottom of the matrix to the top.

\large
\begin{center}
\decofourleft\decosix\decofourleft\\
\end{center}
\normalsize

First, an easy acquisition: $\el(15_p)$ is finite-dimensional. If $\alpha=[\m(15_p):\el(24_p)]$ and $\beta=[\m(15_p):\el(60_p)]$, then $\el(15_p)=\m(15_p)-\m(64_p)+(1-\alpha)\m(24_p)+(1-\beta)\m(60_p)...$ Counting dimensions of graded pieces,
\begin{align*}
\dim\el(15_p)[-1]&=15\\
\dim\el(15_p)[0]&={6\choose5}15-64=26\\
\dim\el(15_p)[1]&={7\choose5}15-{6\choose5}64+(1-\alpha)24+(1-\beta)60\\
&=-69+(1-\alpha)24+(1-\beta)60\\
\end{align*}
The dimension cannot be negative, so $\alpha=\beta=0$. Then $$\dim\el(15_p)[2]={8\choose5}15-{7\choose5}64+{6\choose5}84=0$$
Therefore $\el(15_p)$ is finite-dimensional, and $$\dim\el(15_p)=15+26+15=56$$
The full decomposition of $\el(15_p)$ into Vermas follows by applying (Symm) and observing that the dot product of the Verma vector of $\el(15_p)$ with appropriate columns of the decomposition matrix above must be $0$ (columns $6_p'$, $30_p'$, and $64_p'$ do the job) :
\begin{align*}
\el(15_p)&=\m(15_p)-\m(64_p)+\m(24_p)+\m(60_p)-\m(60_s)-\m(10_s)\\&\qquad\qquad+\m(60_p')+\m(24_p')-\m(64_p')+\m(15_p')
\end{align*}

\vspace*{.4cm}

\large
\begin{center}
\decofourleft\decofourleft
\end{center}
\normalsize

\vspace*{.5cm}

Now begin building up the missing entries of the decomposition matrix from the bottom up. The starting point is the Verma-decomposition of $\el(24_p')$, given by inverting the square submatrix from $24_p'$ to its right and below and taking the top row of that inverse of the submatrix. The decomposition numbers in this submatrix survived the KZ functor and have been given by the Hecke algebra's decomposition matrix.
$$\el(24_p')=\m(24_p')-\m(64_p')+\m(30_p')+\m(15_p')-\m(6_p')+\m(1_p')$$
The graded character of $\el(24_p')$ has a pole of order $2$ at $t=1$, so
$$\dim\Supp\el(24_p')=2$$

\vspace*{.4cm}

\large
\begin{center}
\decofourleft\decofourleft
\end{center}
\normalsize

\vspace*{.5cm}

Now restrict $\el(24_p')$ to $A_5$. $\Res^{E_6}_{A_5}\el(24_p')$ is a nonzero $H_{\frac{1}{3}}(A_5)$-module with lowest weight the partition $(3^2)$. By \cite{Wi}, the minimal support irreps of $H_{\frac{1}{3}}(A_5)$ have lowest weights $(6)$ and $(3^2)$. Since $(6)<_c(3^2)$ it follows that $\Res\el(24_p')=\el(3^2)$, and the dimension of support of this module is $1$.

Write $\alpha=[\m(10_s):\el(24_p')]$. Then 
\begin{align*}
\el(10_s)&=\left(\m(10_s)-\m(60_p')+\m(64_p')+\m(15_q')-\m(15_p')-\m(20_p')+\m(6_p')\right)-\alpha\el(24_p')\\
&=:(*)-\alpha\el(24_p')\\
\end{align*}
$\Res^{E_6}_{A_5}(*)=\el(3^2)$ so $\alpha\in\{0,1\}$. Suppose $\alpha=1$, i.e. $\el(10_s)=(*)-\el(24_p')$, then 
$$\Ind^{E_6}_{A_2\times A_2\times A_1}\Res^{E_6}_{A_2\times A_2\times A_1}\el(10_s)=\m(1_p)-2\m(6_p)+...$$
where all elided $\m(\tau)$ are such that $\tau>_c6_p$. Subtracting the factor $\el(1_p)=\m(1_p)-\m(6_p)...$ leaves an expression with negative leading term $-\m(6_p)+...$, contradiction. So $\alpha=0$ and $(*)=\el(10_s)$. The graded character gives $\dim\Supp\el(10_s)=2$.

\vspace*{.4cm}

\large
\begin{center}
\decofourleft\decofourleft
\end{center}
\normalsize

\vspace*{.5cm}

With the same easy argument, $\el(20_p)$: let $\alpha=[\m(20_s):\el(24_p')]$ and write $$\el(20_s)=\m(20_s)-\m(64_p')+\m(30_p')+\m(20_p')-\m(6_p')-\alpha\el(24_p')$$ Then compute that $$\Ind^{E_6}_{A_5}\Res^{E_6}_{A_5}\el(20_s)=-\alpha\m(15_q)+...$$ where $...$ denotes terms $\m(\tau)$ with $\tau\geq_c15_q$. Therefore $\alpha=0$. Its graded character gives that $\dim\Supp\el(20_s)=4$.

\vspace*{.4cm}

\large
\begin{center}
\decofourleft\decofourleft
\end{center}
\normalsize

\vspace*{.5cm}

Next on the list: $\el(60_p)$. Set $\alpha=[\m(60_p):\el(24_p')]$ and write 
\begin{align*}
\el(60_p)&=\m(60_p)-\m(80_s)-\m(60_s)-\m(10_s)+\m(24_p')+2\m(60_p')-\m(64_p')-\m(15_q')\\
&\qquad+2\m(15_p')-\m(6_p')+\m(1_p')-\alpha\el(24_p')
\end{align*}
Consider $\Ind^{E_6}_{A_2\times A_2\times A_1}\Res^{E_6}_{A_2\times A_2\times A_1}\el(60_p)$: it starts with a $-$ sign unless $\alpha<3$. If $\alpha=2$, it's $4\m(20_p)-4\m(30_p)-4\m(15_q)...$, but $[\el(20_p):\m(30_p)]=0$ so this is only a virtual module. If $\alpha=1$, it's $\m(6_p)+3\m(20_p)-3\m(30_p)$, but $[\el(6_p):\m(30_p)]=-1$ so after killing the factor of $\el(6_p)$ one is again left merely with a virtual module. Therefore $\alpha=0$. The character of $\el(60_p)$ has a pole of order $4$ at $t=1$, so $\dim\Supp\el(60_p)=4$.

\vspace*{.4cm}

\large
\begin{center}
\decofourleft\decofourleft
\end{center}
\normalsize

\vspace*{.5cm}
By (dim Hom) and (RR), $[\m(24_p):\el(10_s)]=[\m(24_p):\el(20_s)]=0$. The value of $\alpha:=[\m(24_p):\el(24_p')]$ will tell the rest of the Verma-decomposition of $\el(24_p)$:
\begin{align*}
\el(24_p)=\m(24_p)-\m(60_s)+\m(24_p')+\m(60_p')-\m(64_p')+\m(15_p')+\m(1_p')-\alpha\el(24_p')\\
\end{align*}
It turns out that $\alpha$ must be $0$, otherwise $\Ind\circ\Res$ with respect to the parabolic \linebreak$A_2\times A_2\times A_1$ produces out of the supposed $\el(24_p)$ an expression with negative leading term. Computing the character of $\el(24_p)$ gives $\dim\Supp\el(24_p)=4$.

\vspace*{.4cm}

\large
\begin{center}
\decofourleft\decofourleft
\end{center}
\normalsize

\vspace*{.5cm}

There are three missing entries in the decomposition matrix for row $64_p$. The first comes for free from (dim Hom), (RR) now that the decomposition numbers in row $20_s$ are known: $[\m(64_p):\el(20_s)]=1$. Playing around with $\Ind$ and $\Res,$ 
\begin{align*}
M:&=\frac{1}{2}\left(\Ind^{E_6}_{A_2\times A_2\times A_1}\Res^{E_6}_{A_2\times A_2\times A_1}\el(24p)\right)-2\el(24_p)-2\el(20_s)\\
&=\m(64_p)-\m(24_p)-\m(60_p)+\m(60_s)-\m(20_s)+\m(10_s)-2\m(24_p')-2\m(60_p')\\
&\qquad+3\m(64_p')-\m(30_p')+\m(15_q')-2\m(15_p')-\m(20_p')+2\m(6_p')-\m(1_p')\\
\end{align*}
 coincides with $\el(64_p)$ up through $\m(20_s)$. So for some $a,\;b\geq0$, $$\el(64_p)=M-a\el(10_s)-b\el(24_p')$$
Then  $$\Res^{E_6}_{A_2\times A_2\times A_1}\left(M-a\el(10_s)\right)=-a\el(\mathrm{Triv})+...$$ implies $a=0$, and
 $$\Ind^{E_6}_{A_2\times A_2\times A_1}\Res^{E_6}_{A_2\times A_2\times A_1}\left(M-b\el(24_p')\right)=-b\el(6_p)+...$$ implies $b=0$. Therefore $M=\el(64_p)$, and computing its character, $\dim\Supp\el(64_p)=4$. It follows that $[\m(64_p):\el(10_s)]=0$ and $[\m(64_p):\el(24_p')]=1$.

\vspace*{.4cm}

\large
\begin{center}
\decofourleft\decofourleft
\end{center}
\normalsize

\vspace*{.5cm}

Since $\el(15_p)$ has already been determined and all rows of the decomposition matrix below row $15_p$ are now known, $[\m(15_p):\el(20_s)]=1$ and $[\m(15_p):\el(10_s)]=[\m(15_p):\el(24_p')]=0$.

\vspace*{.4cm}

\large
\begin{center}
\decofourleft\decofourleft
\end{center}
\normalsize

\vspace*{.5cm}

Row $30_p$ of the decomposition matrix has five missing entries. The first three, $[\m(30_p):\el(24_p)]$, $[\m(30_p):\el(60_p)]$, and $[\m(30_p):\el(20_s)]$ are easily recovered by inducing $\el(3^2)$ from $H_{\frac{1}{3}}(A_5)$ to $E_6$.
\begin{align*}
\Ind^{E_6}_{A_5}\el(3^2)&=\m(30_p)+\m(15_q)-\m(64_p)+\m(24_p)-\m(60_p)+\m(80_s)+\m(20_s)+2\m(10_s)\\
&\qquad+2\m(24_p')-2\m(60_p')-\m(64_p')+2\m(30_p')+2\m(15_q')-\m(20_p')+\m(1_p')\\
\end{align*}
This module contains $\el(30_p)$ and $\el(15_q)$ in its composition series. The first few terms of $\el(15_q)$ (up through the coefficient of $\m(20_s)$) are already known: $\el(15_q)=\m(15_q)-\m(60_p)+\m(80_s)...$, and subtracting off this expression for $\el(15_q)$ reveals that the Verma-decomposition of $\el(30_p)$ through the coefficient of $\m(20_s)$ is $\el(30_p)=\m(30_p)-\m(64_p)+\m(24_p)+\m(20_s)...$ There cannot be any other composition factors of $\Ind\el(3^2)$ until $10_s$ since there are no more simple modules of $2$-dimensional support until then, and $\dim\Supp\Ind\el(3^2)=2$. So $[\m(30_p):\el(24_p)]=0$, $[\m(30_p):\el(60_p)]=1$, and $[\m(30_p):\el(20_s)]=0$. \\

\vspace*{.4cm}

\large
\begin{center}
\decofourleft\decofourleft
\end{center}
\normalsize

\vspace*{.5cm}

Row $15_q$ has two missing entries: $\alpha=[\m(15_q):\el(10_s)]$ and $\beta=[\m(15_q):\el(24_p')]$. Set
$$\alpha=[\m(15_q):\el(10_s)],\qquad\beta=[\m(15_q):\el(24_p')]$$
and write:
$$\el(15_q)=\m(15_q)-\m(60_p)+\m(80_s)+\m(10_s)-\m(60_p')+\m(15_q')-\alpha\el(10_s)-\beta\el(24_p')$$
Likewise, set $\gamma=[\m(30_p):\el(10_s)]$ and $\delta=[\m(30_p):\el(24_p')]$ so that 
$$\el(30_p)=\m(30_p)-\m(64_p)+\m(24_p)+\m(20_s)-\m(64_p')+\m(30_p')-\gamma\el(10_s)-\delta\el(24_p')$$
Then 
\begin{align*}
\Res^{E_6}_{A_5}\el(15_q)&=\el(6)+\el(3^2)-(\alpha+\beta)\el(3^2)\\
\Res^{E_6}_{A_5}\el(30_p)&=\el(6)+\el(3^2)-(\gamma+\delta)\el(3^2)
\end{align*}
and so  $\alpha+\beta\leq1$ and $\gamma+\delta\leq1$.

The parabolic $A_2\times A_2\times A_1$ contains some information about $\el(15_q)$ and $\el(30_p)$. The restrictions of $\el(10_s)$ and $\el(24_p')$ do not contain any Vermas in common; they are identical in the first two components $A_2\times A_2$ but whereas $\Res\el(10_s)$ only contains terms labeled with $(\lambda,\mu,(2))$, $\Res\el(24_p')$ only contains terms labeled with $(\lambda,\mu,(1^2))$. For the rest of this section, let
$$\mathbf{1}=\Res^{E_6}_{A_2\times A_2\times A_1}\el(10_s),\qquad\mathbf{2}=\Res^{E_6}_{A_2\times A_2\times A_1}\el(24_p')$$
$\Res^{E_6}_{A_2\times A_2\times A_1}\el(30_p)=(1-\gamma)\mathbf{1}+(1-\delta)\mathbf{2}$ and $\Res^{E_6}_{A_2\times A_2\times A_1}\el(15_q)=(1-\alpha)\mathbf{1}+(1-\beta)\mathbf{2}$. 

Inducing $\mathbf{2}$ up to $E_6$ and then restricting it back down will force one of $\Res\el(15_q)$, $\Res\el(30_p)$ to be equal to $\mathbf{2}$.
\begin{align*}
\Ind^{E_6}_{A_2\times A_2\times A_1}\mathbf{2}&=\m(6_p)-\m(20_p)+\m(30_p)+2\m(15_q)-\m(15_p)+\m(24_p)-2\m(60_p)\\
&\qquad+\m(80_s)+\m(20_s)+2\m(10_s)+2\m(24_p')-\m(60_p')-2\m(64_p')\\
&\qquad+2\m(30_p')+\m(15_q')+\m(15_p')-\m(6_p')+\m(1_p')\\
\end{align*}
This module contains $\el(6_p)$ in its composition series, but $\el(6_p)$ begins with $\el(6_p)=\m(6_p)-\m(20_p)-\m(30_p)...$, so also $2\el(30_p)$ belongs to the composition series of $\Ind\mathbf{2}$. If we write $\epsilon=[\m(6_p):\el(15_q)]$ then $[\el(6_p):\m(15_q)]=1-\epsilon$. But $2\m(15_q)$ appears above in $\Ind\mathbf{2}$, and $[\el(30_p):\m(15_q)]=0$, so it must happen that $\el(15_q)$ is also a composition factor (possibly with multiplicity) of $\Ind\mathbf{2}$. Restricting $\Ind\mathbf{2}$ back down to $A_2\times A_2\times A_1$ there are only two copies of $\mathbf{1}$
$$\Res^{E_6}_{A_2\times A_2\times A_1}\Ind^{E_6}_{A_2\times A_2\times A_1}\mathbf{2}=2\cdot\mathbf{1}+4\cdot\mathbf{2}$$
Consequently, at least one of $\Res\el(15_q),\;\Res\el(30_p)$ is $\mathbf{2}$ rather than $\mathbf{1}$ or $\mathbf{1}+\mathbf{2}$. So at least one of $\alpha,\;\gamma$ equals $1$. 

The problem of determining $\alpha,\;\beta,\;\gamma,\;\delta$ will be bracketed for a moment. But note that since $\el(15_q)$ and $\el(30_p)$ are composition factors of $\Ind\mathbf{2}$, a module of support $2$-dimensional support, and moreover, they both restrict to modules of $1$-dimensional support in a maximal parabolic, $\dim\Supp\el(15_q)=\dim\Supp\el(30_p)=2$.

\vspace*{.4cm}

\large
\begin{center}
\decofourleft\decofourleft
\end{center}
\normalsize

\vspace*{.5cm}

It is now possible to prove that $\el(6_p)$ is finite-dimensional and find its dimension and Verma-decomposition. Subtract off the expression for $2\el(30_p)+\el(15_q)$ up through the coefficient of $\m(20_s)$ from $\Ind^{E_6}_{A_2\times A_2\times A_1}\mathbf{2}$, and call it $M$:
\begin{align*}
M&=\m(6_p)-\m(20_p)-\m(30_p)+\m(15_q)-\m(15_p)+2\m(64_p)-\m(24_p)-\m(60_p)\\
&\qquad-\m(20_s)...
\end{align*}
Up through $\m(20_s)$, $M$ coincides with $\el(6_p)+\gamma\el(15_q)$ where $\gamma=[\m(6_p):\el(15_q)]$, by the same dimension of support and restriction arguments given in the previous paragraph. However, if $\gamma>0$ then the dimensions of the graded pieces of $\el(6_p)$ fail to satisfy the inequality $\dim\el(6_p)[-1]\leq\dim\el(6_p)[1]$ guaranteed by $\s$-theory. The only option is that $\gamma=0$, and $M$ is equal to $\el(6_p)$ up through the term $\m(20_s)$, which occurs in degree $3$. The dimensions of graded pieces of $\el(6_p)$ in degrees $-2$ through $2$ are:
\begin{align*}
\dim\el(6_p)[-2]=\dim\el(6_p)[2]&=216\\
\dim\el(6_p)[-1]=\dim\el(6_p)[1]&=306\\
\dim\el(6_p)[0]&=340\\
\end{align*}
By (E), $\el(6_p)$ is finite-dimensional. Computing the dimensions of graded pieces in degrees $[-5],\;[-4],$ and $\;[-3]$ as well, 
$$\dim\el(6_p)=2\left(6+36+106+216+306\right)+340=1680$$
and the character of $\el(6_p)$ is of course $$\chi_{\el(6_p)}(t)=6(t^{-5}+t^5)+36(t^{-4}+t^4)+106(t^{-3}+t^3)+216(t^{-2}+t^2)+306(t^{-1}+t)+340$$
When $\chi_{\el(6_p)}(t)$ is expressed as a rational function with denominator $(1-t)^6$, the coefficient of $t^3$ in the numerator is $0$; this forces $[\el(6_p):\m(10_s)]=-2$. Then by (Symm),
\begin{align*}
\el(6_p)&=\m(6_p)-\m(20_p)-\m(30_p)+\m(15_q)-\m(15_p)+2\m(64_p)-\m(24_p)-\m(60_p)\\
&\qquad-\m(20_s)+2\m(10_s)-\m(24_p')-\m(60_p')+2\m(64_p')-\m(30_p')+\m(15_q')\\
&\qquad-\m(15_p')-\m(20_p')+\m(6_p')\\
\end{align*}

The expression for $\el(6_p)$ in terms of Vermas also recovers the decomposition number  $[\m(6_p):\el(64_p)]=1$ (and $[\m(6_p):\el(15_q)]=0$ was recovered above). The rest of the missing numbers in this row depend on those in row $20_p$ which haven't been found yet.

\vspace*{.4cm}

\large
\begin{center}
\decofourleft\decofourleft
\end{center}
\normalsize

\vspace*{.5cm}

To understand $\el(20_p)$ it is necessary to understand the structure of $H_{\frac{1}{3}}(D_5)$. First of all, $\dim\Supp\el(20_p)\leq4$, for $\Ind^{E_6}_{D_5}\Res^{E_6}_{D_5}\el(24_p)=\m(20_p)-\m(15_q)+\m(64_p)+\m(24_p)-\m(60_p)-\m(60_s)+\m(20_p')+\m(6_p')+\m(1_p')$ is a module of dimension $4$ support containing $\el(20_p)$ in its composition series. So now dimensions of supports of all irreducible $H_{\frac{1}{3}}(E_6)$ modules are known -- there are four dimensions of support, $0$, $2$, $4$, and $6$. Therefore $H_{\frac{1}{3}}(D_5)$ has no modules of even-dimensional support, or else inducing such would give a module of odd-dimensional support in $H_{\frac{1}{3}}(E_6)$. On the other hand, it is easy to check that no parabolic of $D_5$ of rank greater than $2$ has a finite-dimensional irrep, mostly just by considering $h_{\frac{1}{3}}(\textrm{Triv})$ which in almost every case is positive. In the case of $A_3\times A_1$, $h_{\frac{1}{3}}(\textrm{Triv})=-\frac{1}{3}$, so there could be a finite-dimensional irrep if some $\tau\in\Irr A_3\times A_1$ satisfied $h_{\frac{1}{3}}(\tau)=0$; however, the dimension of $\tau$ would need to be a multiple of $7$, and there is no such $\tau$. It follows from Corollary 3.24 in \cite{BE} that $H_{\frac{1}{3}}(D_5)$ has no module of $1$-dimensional support. On the other hand, the parabolic $A_2$ of $D_5$ produces a finite-dimensional irrep for the rational Cherednik algebra at $c=\frac{1}{3}$ because $3$ is the Coxeter number of $A_2$. This means $H_{\frac{1}{3}}(D_5)$ contains irreps of $3$-dimensional support. To conclude, the only dimensions of support for $H_{\frac{1}{3}}(D_5)$-modules are $3$ and $5$.

The irreducible representations of minimal support (dim Supp $3$) for $H_{\frac{1}{3}}(D_5)$ all have length $3$ and each falls into its own block of defect $1$. They are given by restricting the five irreps of dim $4$ support from $H_{\frac{1}{3}}(E_6)$: 
\begin{align*}
&\Res\el(20_s)=\el(3:1^2)=\m(3:1^2)-\m(2,1:1^2)+\m(1^3:1^2)\\
&\Res\el(60_p)=\el(3:2)=\m(3:2)-\m(2,1:2)+\m(1^3:2)\\
&\Res\el(24_p)=\el(4,1:-)=\m(4,1:-)-\m(3,2:-)+\m(1^5:-)\\
&\Res\el(64_p)=\el(4:1)=\m(4:1)-\m(2^2:1)+\m(1^4:1)\\
&\Res\el(20_p)=\el(5:-)=\m(5:-)-\m(2^2,1:-)+\m(2,1^3:-)\\
\end{align*}
$\el(3:2)$ and $\el(3:1^2)$ are a dual pair, $\el(5:-)$ and $\el(1^5:-)$ are a dual pair, and $\el(4,1:-)$ is self-dual.

$\Ind\el(4,1:-)=\Ind\Res\el(24_p)$ as two paragraphs ago contains $\el(20_p)$ once in its composition series. Computing the restriction back down to $D_5,$
\begin{align*}
\Res\Ind\el(4,1:-)&=\el(5:-)+3\el(4,1:-)+2\el(4:1)+\el(3:1^2)\\
&=\Res\el(20_p)+3\Res\el(24_p)+2\Res\el(64_p)+\Res\el(20_s)\\
\end{align*}
Therefore
$$\el(20_p)=\Ind\el(4,1:-)-3\el(24_p)-2\el(64_p)-\el(20_s)-a\el(10_s)-b\el(24_p')$$
for some $a,b\geq0$ ( -- since $\el(10_s),\;\el(24_p')$ have $2$-dimensional support, their restriction to $\oh_{\frac{1}{3}}(D_5)$ is $0$. Also, from the known entries of row $20_p$ in the decomposition matrix, none of $\el(15_q)$, $\el(15_p)$, $\el(30_p)$ is a composition factor of $\Ind\el(4,1:-)$). To determine $a$ and $b$, take the right-hand-side of the equation above, restrict it to $A_5$, then induce the result back up to $E_6$ to get $\m(6_p)+2\m(20_p)-(a+b+1)\m(30_p)...$ Since $[\el(20_p):\m(30_p)]=0$ and $[\el(6_p):\m(30_p)]=1$, this will be merely a virtual module unless $a=b=0$. So $\el(20_p)$ is given by the equation above, and explicitly by the Verma-decomposition:
\begin{align*}
\el(20_p)&=\m(20_p)-\m(15_q)-\m(64_p)+\m(60_p)+\m(20_s)-2\el(10_s)+\m(24_p')+\m(60_p')\\
&\qquad-2\m(64_p')+\m(30_p')-2\m(15_q')+\m(15_p')+2\m(20_p')-2\m(6_p')
\end{align*}

\vspace*{.4cm}

\begin{center}
\large
\decofourleft\decofourleft
\end{center}
\normalsize

\vspace*{.5cm}

There is enough information now to determine the decomposition of $\el(1_p)$. As above, let $\mathbf{1}$ denote the module $\Res^{E_6}_{A_2\times A_2\times A_1}\el(10_s)$ and $\mathbf{2}$ the module $\Res^{E_6}_{A_2\times A_2\times A_1}\el(24_p')$. Consider
\begin{align*}
\Ind^{E_6}_{A_2\times A_2\times A_1}\mathbf{1}&=\m(1_p)-\m(6_p)+2\m(30_p)+\m(15_q)+\m(15_p)-2\m(64_p)+2\m(24_p)\\
&\qquad-\m(60_p)+\m(80_s)+\m(20_s)+2\m(10_s)+\m(24_p')-2\m(60_p')+\m(30_p')\\
&\qquad+2\m(15_q')-\m(15_p')-\m(20_p')+\m(6_p')\\
\end{align*}
The dimension of support of $\Ind^{E_6}_{A_2\times A_2\times A_1}\mathbf{1}$ is $2$. Its possible composition factors are the finite-dimensional irreps $\el(1_p),\;\el(6_p),$ and $\el(15_p),$ and the irreps with $2$-dimensional support: $\el(30_p),\;\el(15_q),\;\el(10_s)$, and $\el(24_p')$. Since it's already known that $[\m(1_p):\el(6_p)]=1$, $\el(6_p)$ does not occur. Furthermore, it follows that $[\el(1_p):\m(20_p)]=0$, for given the expression above it can be at most $0$, and if it were less than $0$ then $\el(20_p)$ would be a composition factor, which it is not. $\el(15_p)$ could occur; and $[\el(1_p):\m(15_p)]\leq1$.

The dimension of $\el(1_p)$ was calculated by Oblomkov and Yun in \cite{OY} using an algorithm deriving from the geometry of the affine Springer fiber; their calculation gives $\dim\el(1_p)=4152$. On the other hand, by counting up dimensions of graded pieces, all of which except the ones in degrees $-1$ and $0$ are known now, the dimension of $\el(1_p)$ is 
\begin{align*}
\dim\el(1_p)&=2\left(\sum_{n=h_{\frac{1}{3}}(1_p)}^{-1}\dim\el(1_p)[n]\right)+\dim\el(1_p)[0]\\
&=2\left(\sum_{n=0}^{8}{n+5\choose5}-6{n+1\choose5}\right)+2\left(c_{30_p}30+c_{15_q}15+c_{15_p}15\right)\\
&\qquad\qquad+{14\choose5}-6{10\choose5}+\left(c_{30_p}30+c_{15_q}15+c_{15_p}15\right){6\choose5}+c_{64}64
\end{align*}
where $c_{\tau}$ denotes $[\el(1_p):\m(\tau)]$. $\el(30_p),\;\el(15_q)$. and $\el(15_p)$ start in degree $-1$ and $\el(64_p)$ starts in degree $0$. From $\Ind\mathbf{1}$, $c_{30_p}\leq2,\;c_{15_q}\leq1,$ and $c_{15_p}\leq1$. $c_{64_p}$ depends on the other three, since $\m(64_p)$ isn't a composition factor of $\Ind\mathbf{1}$, and to find $\el(1_p)$ it's enough to subtract off the correct multiples of $\el(30_p),\; \el(15_q),\;\el(15_p)$ from $\Ind\mathbf{1}$. Subtracting off too many reps from among those three quickly produces things of too small a dimension. There is only one combination which produces $4152$: $$c_{30_p}=1,\qquad c_{15_q}=-1,\qquad c_{15_p}=1$$
So $\el(30_p)$ appears once in the composition series of $\Ind\mathbf{1}$, $\el(15_q)$ appears twice, and $\el(15_p)$ does not appear. Subtracting from $\Ind\mathbf{1}$ with these multiplicities what is known of the Verma-decompositions of $\el(30_p)$ and $\el(15_q)$ at this stage, so through the coefficient of $\m(20_s)$, it follows that the decomposition of $\el(1_p)$ through $\m(20_s)$ is $$\el(1_p)=\m(1_p)-\m(6_p)+\m(30_p)-\m(15_q)+\m(15_p)-\m(64_p)+\m(24_p)+\m(60_p)-\m(80_s)...$$
Calculating the graded character which is symmetric in $t$ and $t^{-1}$ since $\el(1_p)$ is finite-dimensional and so is determined by the dimensions of graded pieces up through degree $0$, it turns out that when written as a polynomial over $(1-t)^6$, the coefficient of $t^3$ in the numerator is $-90$. This implies that $c_{10_s}=-1$ since $c_{80_s}=-1,\;c_{60_s}=c_{20_s}=0$. Now the rest of the Verma-decomposition may be written by (Symm):
\begin{align*}
\el(1_p)&=\m(1_p)-\m(6_p)+\m(30_p)-\m(15_q)+\m(15_p)-\m(64_p)+\m(24_p)+\m(60_p)-\m(80_s)\\
&\;\;\;-\m(10_s)+\m(24_p')+\m(60_p')-\m(64_p')+\m(30_p')-\m(15_q')+\m(15_p')-\m(6_p')+\m(1_p')\\
\end{align*}
One may check that the restriction to any parabolic of the expression found for $\el(1_p)$ is $0$, as must be true for any finite-dimensional module. Using this Verma-decomposition, the missing decomposition numbers $[\m(1_p):\el(\tau)]$ may all be filled in now except for when $\tau=10_s$ or $24_p'$. They are all $0$ except for when $\tau=20_p$ or $15_q$, in which case they are $1$.

\vspace*{.4cm}

\huge
\begin{center}
\leafleft\\
\end{center}
\normalsize

\vspace*{.5cm}

The Verma-decomposition of $\el(1_p)$ makes it possible to finish the Verma-decompositions of $\el(15_q)$, $\el(30_p)$, and $\el(20_p)$, and consequently, to find the remaining ten decomposition numbers in the first five rows and the two columns $10_s$ and $24_p'$. $\el(15_q)$ appears twice in $\Ind\mathbf{1}$ while $\el(30_p)$ appears once. As one of $\tau\in\{30_p,\;15_q\}$ satisfies $\Res^{E_6}_{A_2\times A_2\times A_1}\el(\tau)=\mathbf{2}$, and $\Res\Ind\mathbf{1}=4\cdot\mathbf{1}+2\cdot\mathbf{2}$, it follows that the other of $\el(30_p),\;\el(15_q)$ must restrict to $\mathbf{1}$ or else the coefficient of $\mathbf{2}$ in $\Res\Ind\mathbf{1}$ would be at least $3$. This means that the two-by-two square where rows $30_p,\;15_q$ intersect columns $10_s,\;24_p'$ in the decomposition matrix is either $\begin{pmatrix}1&0\\0&1\end{pmatrix}$ or $\begin{pmatrix}0&1\\1&0\end{pmatrix}$. The Verma-decomposition $\el(20_p)=\m(20_p)-\m(15_q)-\m(64_p)+\m(60_p)+\m(20_s)-2\m(10_s)+\m(24_p')...$, and likewise the Verma-decompositions of $\el(1_p),\;\el(6_p)$ up through $\m(24_p')$, then determine what the remaining decomposition numbers would have to be in each case without any ambiguity. If the two-by-two square is $\begin{pmatrix}0&1\\1&0\end{pmatrix}$ then $[\m(1_p):\el(24_p')]$ turns out to be negative. This is impossible. So the correct option is $\begin{pmatrix}1&0\\0&1\end{pmatrix}$ and the top 5 rows of columns $10_s$ and $24_p'$ read:
$$\begin{pmatrix}
0&0\\
1&0\\
1&1\\
1&0\\
0&1\\
\end{pmatrix}$$

That completes the calculation of the decomposition matrix. Any remaining characters now follow from the character formula, the inverse of the decomposition matrix, and the $h_{\frac{1}{3}}$-weight line, and we get the dimensions of support for $\el(\tau)$ printed next to row $\tau$ of the matrices displayed in the theorem.

\end{proof}

The graded characters and the dimensions of the finite-dimensional irreps of $H_{\frac{1}{3}}(E_6)$ are:

\begin{align*}
&\chi_{\el(1_p)}(t)=t^{-9}+t^9+6(t^{-8}+t^8)+21(t^{-7}+t^7)+56(t^{-6}+t^6)+120(t^{-5}+t^5)\\
&\qquad\qquad\;\;\;+216(t^{-4}+t^4)+336(t^{-3}+t^3)+456(t^{-2}+t^2)+561(t^{-1}+t)+606\\
&\dim\el(1_p)=4152\\
\\
&\chi_{\el(6_p)}(t)=6(t^{-5}+t^5)+36(t^{-4}+t^4)+106(t^{-3}+t^3)+216(t^{-2}+t^2)+306(t^{-1}+t)+340\\
&\dim\el(6_p)=1680\\
\\
&\chi_{\el(15_p)}(t)=15(t^{-1}+t)+26\\
&\dim\el(15_p)=56\\
\end{align*}

\section{Simple representations of $H_c(E_7)$}

All decomposition numbers for $H_c(E_7)$, $c=1/d$ with $d>2$ dividing one of the fundamental degrees of $E_7$, are either $0$'s or $1$'s. The decomposition matrices of $\oh_c(W)$ for all such $d$ are found in this section. Note that $d=10$ is not an elliptic number of $E_7$, yet nonetheless there is a finite-dimensional representation when $c=1/10$; its lowest weight is the standard representation of $E_7$.

%REMARK!! the ``amplitude" of a "defect 1" block is always the same for a given parameter c, except in the case of E8, c=1/15, or H4, c=1/15.

\subsection{{\Large$\mbf{c=\frac{1}{18}}$}}

\begin{center}
\begin{picture}(280,50)
\put(0,30){\line(1,0){280}}
\put(0,30){\circle*{5}}
\put(40,30){\circle*{5}}
\put(80,30){\circle*{5}}
\put(120,30){\circle*{5}}
\put(160,30){\circle*{5}}
\put(200,30){\circle*{5}}
\put(240,30){\circle*{5}}
\put(280,30){\circle*{5}}
\put(-2,40){$0$}
\put(38,40){$1$}
\put(78,40){$2$}
\put(118,40){$3$}
\put(158,40){$4$}
\put(198,40){$5$}
\put(238,40){$6$}
\put(278,40){$7$}
\put(-3,15){$1_a$}
\put(37,15){$7_a'$}
\put(74,15){$21_a$}
\put(114,15){$35_a'$}
\put(154,15){$35_a$}
\put(194,15){$21_a'$}
\put(234,15){$7_a$}
\put(276,15){$1_a'$}
\end{picture}
\end{center}

\begin{align*}
%\el(1_a)&=\m(1_a)-\m(7_a')+\m(21_a)-\m(35_a')+\m(35_a)-\m(21_a')+\m(7_a)-\m(1_a')\\
\chi_{\el(1_a)}(t,w)&=\chi_{1_a}\\
\dim\el(1_a)&=1
\end{align*}

\subsection{{\Large$\mbf{c=\frac{1}{14}}$}}

\begin{center}
\begin{picture}(360,50)
\put(0,30){\line(1,0){360}}
\put(0,30){\circle*{5}}
\put(80,30){\circle*{5}}
\put(120,30){\circle*{5}}
\put(160,30){\circle*{5}}
\put(200,30){\circle*{5}}
\put(240,30){\circle*{5}}
\put(280,30){\circle*{5}}
\put(360,30){\circle*{5}}
\put(-8,40){$-1$}
\put(78,40){$1$}
\put(118,40){$2$}
\put(158,40){$3$}
\put(198,40){$4$}
\put(238,40){$5$}
\put(278,40){$6$}
\put(358,40){$8$}
\put(-3,15){$1_a$}
\put(74,15){$27_a$}
\put(114,15){$105_a'$}
\put(154,15){$189_a$}
\put(194,15){$189_a'$}
\put(232,15){$105_a$}
\put(273,15){$27_a'$}
\put(356,15){$1_a'$}
\end{picture}
\end{center}

\begin{align*}
%\el(1_a)&=\m(1_a)-\m(27_a)+\m(105_a')-\m(189_a)+\m(189_a')-\m(105_a)+\m(27_a')-\m(1_a')\\
\chi_{\el(1_a)}(t,w)&=\chi_{1_a}(t^{-1}+t)+\chi_{7_a'}\\
\dim\el(1_a)&=9\\
\end{align*}

\subsection{{\Large$\mbf{c=\frac{1}{12}}$}}

\begin{center}
\begin{picture}(320,50)
\put(0,30){\line(1,0){320}}
\put(0,30){\circle*{5}}
\put(120,30){\circle*{5}}
\put(160,30){\circle*{5}}
\put(200,30){\circle*{5}}
\put(240,30){\circle*{5}}
\put(280,30){\circle*{5}}
\put(320,30){\circle*{5}}
\put(-8,40){$-1.75$}
\put(118,40){$1.25$}
\put(158,40){$2.25$}
\put(198,40){$3.25$}
\put(238,40){$4.25$}
\put(278,40){$5.25$}
\put(318,40){$6.25$}
\put(-3,15){$1_a$}
\put(114,15){$56_a'$}
\put(154,15){$210_a$}
\put(194,15){$336_a'$}
\put(234,15){$280_a$}
\put(276,15){$120_a'$}
\put(316,15){$21_b$}
\end{picture}

\begin{picture}(320,50)
\put(0,30){\line(1,0){320}}
\put(0,30){\circle*{5}}
\put(40,30){\circle*{5}}
\put(80,30){\circle*{5}}
\put(120,30){\circle*{5}}
\put(160,30){\circle*{5}}
\put(200,30){\circle*{5}}
\put(320,30){\circle*{5}}
\put(-2,40){$.75$}
\put(38,40){$1.75$}
\put(78,40){$2.75$}
\put(118,40){$3.75$}
\put(158,40){$4.75$}
\put(198,40){$5.75$}
\put(318,40){$8.75$}
\put(-3,15){$21_b'$}
\put(37,15){$120_a$}
\put(74,15){$280_a'$}
\put(114,15){$336_a$}
\put(154,15){$210_a'$}
\put(194,15){$56_a$}
\put(316,15){$1_a'$}
\end{picture}
\end{center}

\begin{itemize}
\item$\dim\Supp\el(1_a)=1$
\item$\dim\Supp\el(21_b')=1$
\end{itemize}

\begin{comment}
\begin{align*}
\el(1_a)&=\m(1_a)-\m(56_a')+\m(210_a)-\m(336_a')+\m(280_a)-\m(120_a')+\m(21_b)\\
\chi_{\el(1_a)}(t)&=t^{-1.75}\frac{21t^2+6t+1}{1-t}\\
\\
\el(21_b')&=\m(21_b')-\m(120_a)+\m(280_a')-\m(336_a)+\m(210_a')-\m(56_a)+\m(1_a')\\
\chi_{\el(21_b')}(t)&=t^{.75}\frac{t^2+6t+21}{1-t}\\
\end{align*}

\end{comment}

\subsection{{\Large$\mbf{c=\frac{1}{10}}$}}

\begin{center}
\begin{picture}(360,50)
\put(0,30){\line(1,0){360}}
\put(0,30){\circle*{5}}
\put(120,30){\circle*{5}}
\put(200,30){\circle*{5}}
\put(240,30){\circle*{5}}
\put(280,30){\circle*{5}}
\put(320,30){\circle*{5}}
\put(360,30){\circle*{5}}
\put(-12,40){$-2.8$}
\put(117,40){$.2$}
\put(195,40){$2.2$}
\put(235,40){$3.2$}
\put(275,40){$4.2$}
\put(315,40){$5.2$}
\put(355,40){$6.2$}
\put(-3,15){$1_a$}
\put(114,15){$21_b'$}
\put(194,15){$189_c'$}
\put(232,15){$420_a$}
\put(273,15){$405_a'$}
\put(316,15){$189_b$}
\put(356,15){$35_b'$}
\end{picture}

\begin{picture}(360,50)
\put(0,30){\line(1,0){360}}
\put(0,30){\circle*{5}}
\put(40,30){\circle*{5}}
\put(80,30){\circle*{5}}
\put(120,30){\circle*{5}}
\put(160,30){\circle*{5}}
\put(240,30){\circle*{5}}
\put(360,30){\circle*{5}}
\put(-3,40){$.8$}
\put(35,40){$1.8$}
\put(75,40){$2.8$}
\put(115,40){$3.8$}
\put(155,40){$4.8$}
\put(235,40){$6.8$}
\put(355,40){$9.8$}
\put(-6,15){$35_b$}
\put(32,15){$189_b'$}
\put(72,15){$405_a$}
\put(112,15){$420_a'$}
\put(152,15){$189_c$}
\put(235,15){$21_b$}
\put(356,15){$1_a'$}
\end{picture}

\begin{picture}(360,50)
\put(0,30){\line(1,0){360}}
\put(0,30){\circle*{5}}
\put(40,30){\circle*{5}}
\put(120,30){\circle*{5}}
\put(160,30){\circle*{5}}
\put(200,30){\circle*{5}}
\put(240,30){\circle*{5}}
\put(320,30){\circle*{5}}
\put(360,30){\circle*{5}}
\put(-8,40){$-1$}
\put(38,40){$0$}
\put(118,40){$2$}
\put(158,40){$3$}
\put(198,40){$4$}
\put(238,40){$5$}
\put(318,40){$7$}
\put(358,40){$8$}
\put(-3,15){$7_a'$}
\put(34,15){$27_a$}
\put(114,15){$168_a$}
\put(154,15){$378_a'$}
\put(194,15){$378_a$}
\put(232,15){$168_a'$}
\put(313,15){$27_a'$}
\put(356,15){$7_a$}
\end{picture}
\end{center}

\begin{itemize}
\item$\dim\Supp\el(1_a)=1$
\item$\dim\Supp\el(35_b)=1$
\item$\dim\Supp\el(7_a')=0$\\
$\chi_{\el(7_a')}(t)=7(t^{-1}+t)+22$\\
$\dim\el(7_a')=36$
\end{itemize}

\begin{comment}
\begin{align*}
\el(1_a)&=\m(1_a)-\m(21_b')+\m(189_c')-\m(420_a)+\m(405_a')-\m(189_b)+\m(35_b')\\
\chi_{\el(1_a)}(t)&=t^{-2.8}\frac{35t^3+21t^2+6t+1}{1-t}\\
\\
\el(35_b)&=\m(35_b)-\m(189_b')+\m(405_a)-\m(420_a')+\m(189_c)-\m(21_b)+\m(1_a')\\
\chi_{\el(35_b)}(t)&=t^{.8}\frac{t^3+6t^2+21t+35}{1-t}\\
\\
\el(7_a')&=\m(7_a')-\m(27_a)+\m(168_a)-\m(378_a')+\m(378_a)-\m(168_a')+\m(27_a')-\m(7_a)\\
\end{align*}
\end{comment}

\newpage
\subsection{{\Large$\mbf{c=\frac{1}{9}}$}}

\begin{center}
\begin{picture}(360,50)
\put(0,30){\line(1,0){360}}
\put(0,30){\circle*{5}}
\put(120,30){\circle*{5}}
\put(180,30){\circle*{5}}
\put(210,30){\circle*{5}}
\put(240,30){\circle*{5}}
\put(300,30){\circle*{5}}
\put(360,30){\circle*{5}}
\put(-12,40){$-3.5$}
\put(116,40){$.5$}
\put(174,40){$2.5$}
\put(204,40){$3.5$}
\put(234,40){$4.5$}
\put(294,40){$6.5$}
\put(354,40){$8.5$}
\put(-3,15){$1_a$}
\put(115,15){$35_b$}
\put(173,15){$280_b$}
\put(203,15){$512_a'$}
\put(233,15){$315_a$}
\put(295,15){$56_a$}
\put(355,15){$7_a$}
\end{picture}

\begin{picture}(360,50)
\put(0,30){\line(1,0){360}}
\put(0,30){\circle*{5}}
\put(60,30){\circle*{5}}
\put(120,30){\circle*{5}}
\put(150,30){\circle*{5}}
\put(180,30){\circle*{5}}
\put(240,30){\circle*{5}}
\put(360,30){\circle*{5}}
\put(-12,40){$-1.5$}
\put(56,40){$.5$}
\put(114,40){$2.5$}
\put(144,40){$3.5$}
\put(174,40){$4.5$}
\put(234,40){$6.5$}
\put(354,40){$10.5$}
\put(-3,15){$7_a'$}
\put(55,15){$56_a'$}
\put(113,15){$315_a'$}
\put(143,15){$512_a$}
\put(173,15){$280_b'$}
\put(235,15){$35_b'$}
\put(355,15){$1_a'$}
\end{picture}
\end{center}

\begin{itemize}
\item$\dim\Supp\el(1_a)=1$
\item$\dim\Supp\el(7_a')=1$
\end{itemize}

\begin{comment}
\begin{align*}
\el(1_a)&=\m(1_a)-\m(35_b)+\m(280_b)-\m(512_a')+\m(315_a)-\m(56_a)+\m(7_a)\\
\chi_{\el(1)}(t)&=t^{-3.5}\frac{7t^6+42t^5+91t^4+56t^3+21t^2+6t+1}{1-t}\\
\\
\el(7_a')&=\m(7_a')-\m(56_a')+\m(315_a')-\m(512_a)+\m(280_b')-\m(35_b')+\m(1_a')\\
\chi_{\el(7_a')}(t)&=t^{-1.5}\frac{t^6+6t^5+21t^4+56t^3+91t^2+42t+7}{1-t}\\
\end{align*}
\end{comment}

\subsection{{\Large$\mbf{c=\frac{1}{8}}$}}

\begin{itemize}
\item$\dim\Supp\el(1_a)=2$
\item$\dim\Supp\el(21_b')=2$
\item$\dim\Supp\el(7_a')=2$
\item$\dim\Supp\el(27_a)=2$
\end{itemize}

\begin{comment}
\begin{align*}
\el(1_a)&=\m(1_a)-\m(105_b)+\m(216_a')-\m(280_b')+\m(189_b)-\m(21_b)\\
\chi_{\el(1_a)}(t)&=t^{-4.375}\frac{21t^7+105t^6+126t^5+70t^4+35t^3+15t^2+5t+1}{(1-t)^2}\\
\\
\el(21_b')&=\m(21_b')-\m(189_b')+\m(280_b)-\m(216_a)+\m(105_b')-\m(1_a')\\
\chi_{\el(21_b')}(t)&=t^{-.625}\frac{t^7+5t^6+15t^5+35t^4+70t^3+126t^2+105t+21}{(1-t)^2}\\
\\
\el(7_a')&=\m(7_a')-\m(120_a)+\m(189_c')-\m(105_c;)+\m(56_a)-\m(27_a')\\
\chi_{\el(7_a')}(t)&=t^{-2.125}\frac{27t^5+79t^4+125t^3+105t^2+35t+7}{(1-t)^2}\\
\\
\el(27_a)&=\m(27_a)-\m(56_a')+\m(105_c)-\m(189_c)+\m(120_a')-\m(7_a)\\
\chi_{\el(27_a)}(t)&=t^{-.875}\frac{7t^5+35t^4+105t^3+125t^2+79t+27}{(1-t)^2}\\
\end{align*}
\end{comment}

\subsection{{\Large$\mbf{c=\frac{1}{7}}$}}

\begin{itemize}
\item$\dim\Supp\el(1_a)=1$
\item$\dim\Supp\el(15_a')=1$
\end{itemize}

\begin{comment}
\begin{align*}
\el(1_a)&=\m(1_a)-\m(27_a)+\m(120_a)-\m(405_a)+\m(512_a')-\m(216_a)+\m(15_a)\\
\chi_{\el(1_a)}(t)&=t^{-5.5}\frac{15t^6+90t^5+99t^4+56t^3+21t^2+6t+1}{1-t}\\
\\
\el(15_a')&=\m(15_a')-\m(216_a')+\m(512_a)-\m(405_a')+\m(120_a')-\m(27_a')+\m(1_a')\\
\chi_{\el(15_a')}(t)&=t^{.5}\frac{t^6+6t^5+21t^4+56t^3+99t^2+90t+15}{1-t}\\
\end{align*}
\end{comment}

\subsection{{\Large$\mbf{c=\frac{1}{6}}$}}

\theorem The decomposition matrix and its inverse for the principal block of $H_{\frac{1}{6}}(E_7)$ are as follows. There are four finite-dimensional irreducible representations, three irreducible representations each of $1$- and $2$-dimensional support, and four irreducible representations of $3$-dimensional support; all other irreps have full support.

\vspace*{2in}
\newgeometry{right=2cm,top=1cm,bottom=1cm}
\begin{turn}{270}
\begin{minipage}{\linewidth}
\[
\begin{blockarray}{cccccccccccccccccccccccccccccccccccc}
%1_a & 7_a' & 21_b' & 35_b & 21_a & 105_a' & 15_a' & 210_a & 168_a & 105_b & 315_a' & 280_a' & 70_a' & 
%105_c & 420_a & 210_b & 84_a & 105_c' & 420_a' & 210_b' & 84_a' & 315_a & 280_a & 70_a & 210_a' & 168_a' & 105_b' & 105_a & 15_a & 35_b' & 21_a' & 21_b & 7_a & 1_a'\\
\begin{block}{cc(cccccccccccccccccccccccccccccccccc)}
\star&1_a&1&\cdot&1&1&\cdot&\cdot&1&\cdot&1&1&\cdot&\cdot&1&\cdot&\cdot&1&\cdot&\cdot&\cdot&\cdot&\cdot&\cdot&\cdot&\cdot&\cdot&\cdot&\cdot&\cdot&\cdot&\cdot&\cdot&\cdot&\cdot&\cdot\\
\star&7_a'&\cdot&1&1&1&\cdot&1&\cdot&1&1&1&1&\cdot&\cdot&1&\cdot&\cdot&\cdot&\cdot&\cdot&\cdot&\cdot&\cdot&\cdot&\cdot&\cdot&\cdot&\cdot&\cdot&\cdot&\cdot&\cdot&\cdot&\cdot&\cdot\\
\texttt{(1)}&21_b'&\cdot&\cdot&1&\cdot&\cdot&\cdot&\cdot&1&1&1&1&\cdot&\cdot&\cdot&\cdot&1&\cdot&1&\cdot&\cdot&\cdot&1&\cdot&\cdot&\cdot&\cdot&\cdot&\cdot&\cdot&\cdot&\cdot&\cdot&\cdot&\cdot\\
\texttt{(1)}&35_b&\cdot&\cdot&\cdot&1&\cdot&\cdot&\cdot&\cdot&1&1&1&\cdot&1&1&\cdot&1&1&\cdot&\cdot&1&\cdot&\cdot&\cdot&\cdot&\cdot&\cdot&\cdot&\cdot&\cdot&\cdot&\cdot&\cdot&\cdot&\cdot\\
\star&21_a&\cdot&\cdot&\cdot&\cdot&1&1&\cdot&1&\cdot&\cdot&\cdot&1&\cdot&\cdot&\cdot&\cdot&\cdot&\cdot&\cdot&\cdot&\cdot&\cdot&\cdot&\cdot&\cdot&\cdot&\cdot&\cdot&\cdot&\cdot&\cdot&\cdot&\cdot&\cdot\\
\texttt{(1)}&105_a'&\cdot&\cdot&\cdot&\cdot&\cdot&1&\cdot&1&1&\cdot&1&1&\cdot&\cdot&1&\cdot&\cdot&1&\cdot&\cdot&\cdot&\cdot&\cdot&\cdot&\cdot&\cdot&\cdot&\cdot&\cdot&\cdot&\cdot&\cdot&\cdot&\cdot\\
\star&15_a'&\cdot&\cdot&\cdot&\cdot&\cdot&\cdot&1&\cdot&\cdot&1&\cdot&\cdot&\cdot&\cdot&\cdot&1&\cdot&\cdot&\cdot&\cdot&1&\cdot&\cdot&\cdot&\cdot&\cdot&\cdot&\cdot&\cdot&\cdot&\cdot&\cdot&\cdot&\cdot\\
\texttt{(3)}&210_a&\cdot&\cdot&\cdot&\cdot&\cdot&\cdot&\cdot&1&\cdot&\cdot&1&1&\cdot&1&1&\cdot&\cdot&1&1&\cdot&\cdot&1&\cdot&\cdot&\cdot&\cdot&\cdot&\cdot&\cdot&\cdot&\cdot&\cdot&\cdot&\cdot\\
\texttt{(3)}&168_a&\cdot&\cdot&\cdot&\cdot&\cdot&\cdot&\cdot&\cdot&1&\cdot&1&\cdot&\cdot&\cdot&1&1&1&1&\cdot&1&\cdot&1&\cdot&\cdot&\cdot&\cdot&1&\cdot&\cdot&\cdot&\cdot&\cdot&\cdot&\cdot\\
\texttt{(2)}&105_b&\cdot&\cdot&\cdot&\cdot&\cdot&\cdot&\cdot&\cdot&\cdot&1&1&\cdot&\cdot&\cdot&\cdot&1&\cdot&\cdot&\cdot&1&1&1&\cdot&\cdot&\cdot&1&\cdot&\cdot&\cdot&\cdot&\cdot&\cdot&\cdot&\cdot\\
&315_a'&\cdot&\cdot&\cdot&\cdot&\cdot&\cdot&\cdot&\cdot&\cdot&\cdot&1&\cdot&\cdot&1&1&\cdot&\cdot&\cdot&1&1&\cdot&1&\cdot&\cdot&\cdot&1&1&\cdot&\cdot&\cdot&\cdot&1&\cdot&\cdot\\
&280_a'&\cdot&\cdot&\cdot&\cdot&\cdot&\cdot&\cdot&\cdot&\cdot&\cdot&\cdot&1&\cdot&\cdot&1&\cdot&\cdot&\cdot&1&\cdot&\cdot&\cdot&1&\cdot&\cdot&\cdot&\cdot&\cdot&\cdot&\cdot&\cdot&\cdot&\cdot&\cdot\\
\texttt{(2)}&70_a'&\cdot&\cdot&\cdot&\cdot&\cdot&\cdot&\cdot&\cdot&\cdot&\cdot&\cdot&\cdot&1&\cdot&\cdot&1&\cdot&\cdot&\cdot&1&\cdot&\cdot&\cdot&1&\cdot&\cdot&\cdot&\cdot&\cdot&\cdot&\cdot&\cdot&\cdot&\cdot\\
\texttt{(2)}&105_c&\cdot&\cdot&\cdot&\cdot&\cdot&\cdot&\cdot&\cdot&\cdot&\cdot&\cdot&\cdot&\cdot&1&\cdot&\cdot&\cdot&\cdot&1&\cdot&\cdot&\cdot&\cdot&\cdot&\cdot&\cdot&\cdot&\cdot&\cdot&\cdot&\cdot&1&\cdot&\cdot\\
&420_a&\cdot&\cdot&\cdot&\cdot&\cdot&\cdot&\cdot&\cdot&\cdot&\cdot&\cdot&\cdot&\cdot&\cdot&1&\cdot&\cdot&1&1&\cdot&\cdot&1&1&\cdot&1&\cdot&1&\cdot&\cdot&\cdot&\cdot&1&\cdot&\cdot\\
&210_b&\cdot&\cdot&\cdot&\cdot&\cdot&\cdot&\cdot&\cdot&\cdot&\cdot&\cdot&\cdot&\cdot&\cdot&\cdot&1&\cdot&\cdot&\cdot&1&1&1&\cdot&1&\cdot&1&1&\cdot&\cdot&1&\cdot&\cdot&\cdot&\cdot\\
\texttt{(3)}&84_a&\cdot&\cdot&\cdot&\cdot&\cdot&\cdot&\cdot&\cdot&\cdot&\cdot&\cdot&\cdot&\cdot&\cdot&\cdot&\cdot&1&\cdot&\cdot&1&\cdot&\cdot&\cdot&\cdot&\cdot&\cdot&1&\cdot&1&\cdot&\cdot&\cdot&\cdot&\cdot\\
\texttt{(3)}&105_c'&\cdot&\cdot&\cdot&\cdot&\cdot&\cdot&\cdot&\cdot&\cdot&\cdot&\cdot&\cdot&\cdot&\cdot&\cdot&\cdot&\cdot&1&\cdot&\cdot&\cdot&1&\cdot&\cdot&1&\cdot&\cdot&\cdot&\cdot&\cdot&\cdot&\cdot&\cdot&\cdot\\
&420_a'&\cdot&\cdot&\cdot&\cdot&\cdot&\cdot&\cdot&\cdot&\cdot&\cdot&\cdot&\cdot&\cdot&\cdot&\cdot&\cdot&\cdot&\cdot&1&\cdot&\cdot&1&1&\cdot&1&1&\cdot&1&\cdot&\cdot&\cdot&1&\cdot&\cdot\\
&210_b'&\cdot&\cdot&\cdot&\cdot&\cdot&\cdot&\cdot&\cdot&\cdot&\cdot&\cdot&\cdot&\cdot&\cdot&\cdot&\cdot&\cdot&\cdot&\cdot&1&\cdot&\cdot&\cdot&1&\cdot&1&1&\cdot&1&1&\cdot&1&\cdot&1\\
&84_a'&\cdot&\cdot&\cdot&\cdot&\cdot&\cdot&\cdot&\cdot&\cdot&\cdot&\cdot&\cdot&\cdot&\cdot&\cdot&\cdot&\cdot&\cdot&\cdot&\cdot&1&\cdot&\cdot&\cdot&\cdot&1&\cdot&\cdot&\cdot&1&\cdot&\cdot&\cdot&\cdot\\
&315_a&\cdot&\cdot&\cdot&\cdot&\cdot&\cdot&\cdot&\cdot&\cdot&\cdot&\cdot&\cdot&\cdot&\cdot&\cdot&\cdot&\cdot&\cdot&\cdot&\cdot&\cdot&1&\cdot&\cdot&1&1&1&1&\cdot&1&\cdot&1&1&\cdot\\
&280_a&\cdot&\cdot&\cdot&\cdot&\cdot&\cdot&\cdot&\cdot&\cdot&\cdot&\cdot&\cdot&\cdot&\cdot&\cdot&\cdot&\cdot&\cdot&\cdot&\cdot&\cdot&\cdot&1&\cdot&1&\cdot&\cdot&1&\cdot&\cdot&1&\cdot&\cdot&\cdot\\
&70_a&\cdot&\cdot&\cdot&\cdot&\cdot&\cdot&\cdot&\cdot&\cdot&\cdot&\cdot&\cdot&\cdot&\cdot&\cdot&\cdot&\cdot&\cdot&\cdot&\cdot&\cdot&\cdot&\cdot&1&\cdot&\cdot&\cdot&\cdot&\cdot&1&\cdot&\cdot&\cdot&1\\
&210_a'&\cdot&\cdot&\cdot&\cdot&\cdot&\cdot&\cdot&\cdot&\cdot&\cdot&\cdot&\cdot&\cdot&\cdot&\cdot&\cdot&\cdot&\cdot&\cdot&\cdot&\cdot&\cdot&\cdot&\cdot&1&\cdot&\cdot&1&\cdot&\cdot&1&1&1&\cdot\\
&168_a'&\cdot&\cdot&\cdot&\cdot&\cdot&\cdot&\cdot&\cdot&\cdot&\cdot&\cdot&\cdot&\cdot&\cdot&\cdot&\cdot&\cdot&\cdot&\cdot&\cdot&\cdot&\cdot&\cdot&\cdot&\cdot&1&\cdot&1&\cdot&1&\cdot&1&1&1\\
&105_b'&\cdot&\cdot&\cdot&\cdot&\cdot&\cdot&\cdot&\cdot&\cdot&\cdot&\cdot&\cdot&\cdot&\cdot&\cdot&\cdot&\cdot&\cdot&\cdot&\cdot&\cdot&\cdot&\cdot&\cdot&\cdot&\cdot&1&\cdot&1&1&\cdot&1&1&1\\
&105_a&\cdot&\cdot&\cdot&\cdot&\cdot&\cdot&\cdot&\cdot&\cdot&\cdot&\cdot&\cdot&\cdot&\cdot&\cdot&\cdot&\cdot&\cdot&\cdot&\cdot&\cdot&\cdot&\cdot&\cdot&\cdot&\cdot&\cdot&1&\cdot&\cdot&1&\cdot&1&\cdot\\
&15_a&\cdot&\cdot&\cdot&\cdot&\cdot&\cdot&\cdot&\cdot&\cdot&\cdot&\cdot&\cdot&\cdot&\cdot&\cdot&\cdot&\cdot&\cdot&\cdot&\cdot&\cdot&\cdot&\cdot&\cdot&\cdot&\cdot&\cdot&\cdot&1&\cdot&\cdot&\cdot&\cdot&1\\
&35_b'&\cdot&\cdot&\cdot&\cdot&\cdot&\cdot&\cdot&\cdot&\cdot&\cdot&\cdot&\cdot&\cdot&\cdot&\cdot&\cdot&\cdot&\cdot&\cdot&\cdot&\cdot&\cdot&\cdot&\cdot&\cdot&\cdot&\cdot&\cdot&\cdot&1&\cdot&\cdot&1&1\\
&21_a'&\cdot&\cdot&\cdot&\cdot&\cdot&\cdot&\cdot&\cdot&\cdot&\cdot&\cdot&\cdot&\cdot&\cdot&\cdot&\cdot&\cdot&\cdot&\cdot&\cdot&\cdot&\cdot&\cdot&\cdot&\cdot&\cdot&\cdot&\cdot&\cdot&\cdot&1&\cdot&\cdot&\cdot\\
&21_b&\cdot&\cdot&\cdot&\cdot&\cdot&\cdot&\cdot&\cdot&\cdot&\cdot&\cdot&\cdot&\cdot&\cdot&\cdot&\cdot&\cdot&\cdot&\cdot&\cdot&\cdot&\cdot&\cdot&\cdot&\cdot&\cdot&\cdot&\cdot&\cdot&\cdot&\cdot&1&1&1\\
&7_a&\cdot&\cdot&\cdot&\cdot&\cdot&\cdot&\cdot&\cdot&\cdot&\cdot&\cdot&\cdot&\cdot&\cdot&\cdot&\cdot&\cdot&\cdot&\cdot&\cdot&\cdot&\cdot&\cdot&\cdot&\cdot&\cdot&\cdot&\cdot&\cdot&\cdot&\cdot&\cdot&1&\cdot\\
&1_a'&\cdot&\cdot&\cdot&\cdot&\cdot&\cdot&\cdot&\cdot&\cdot&\cdot&\cdot&\cdot&\cdot&\cdot&\cdot&\cdot&\cdot&\cdot&\cdot&\cdot&\cdot&\cdot&\cdot&\cdot&\cdot&\cdot&\cdot&\cdot&\cdot&\cdot&\cdot&\cdot&\cdot&1\\
\end{block}
\end{blockarray}
\]
\end{minipage}
\end{turn}
%\end{sideways}
\restoregeometry
\newpage

\newgeometry{top=.5cm,bottom=0cm}
\footnotesize
\begin{turn}{270}
\begin{minipage}{\linewidth}
\[
\begin{blockarray}{ccccccccccccccccccccccccccccccccccc}
%1_a & 7_a' & 21_b' & 35_b & 21_a & 105_a' & 15_a' & 210_a & 168_a & 105_b & 315_a' & 280_a' & 70_a' & 
%105_c & 420_a & 210_b & 84_a & 105_c' & 420_a' & 210_b' & 84_a' & 315_a & 280_a & 70_a & 210_a' & 168_a' & 105_b' & 105_a & 15_a & 35_b' & 21_a' & 21_b & 7_a & 1_a'\\
\begin{block}{c(cccccccccccccccccccccccccccccccccc)}
1_a&1&\cdot&-1&-1&\cdot&\cdot&-1&1&1&2&-2&-1&\cdot&2&1&-1&\cdot&-2&-1&1&\cdot&2&1&\cdot&-1&-1&-2&\cdot&1&1&\cdot&1&\cdot&-1\\
7_a'&\cdot&1&-1&-1&\cdot&-1&\cdot&1&2&1&-2&\cdot&1&1&\cdot&-2&-1&-1&\cdot&2&1&2&\cdot&-1&-1&-2&-1&1&\cdot&1&\cdot&1&-1&\cdot\\
21_b'&\cdot&\cdot&1&\cdot&\cdot&\cdot&\cdot&-1&-1&-1&2&1&\cdot&-1&-1&1&1&2&\cdot&-2&\cdot&-2&\cdot&1&1&2&2&-1&-1&-2&\cdot&-1&1&1\\
35_b&\cdot&\cdot&\cdot&1&\cdot&\cdot&\cdot&\cdot&-1&-1&1&\cdot&-1&-2&\cdot&2&\cdot&1&1&-1&-1&-3&-1&\cdot&2&2&2&-1&-1&-1&\cdot&-2&1&1\\
21_a&\cdot&\cdot&\cdot&\cdot&1&-1&\cdot&\cdot&1&\cdot&\cdot&\cdot&\cdot&\cdot&\cdot&-1&-1&\cdot&\cdot&1&1&\cdot&\cdot&\cdot&\cdot&-1&\cdot&1&\cdot&\cdot&-1&\cdot&\cdot&\cdot\\
105_a'&\cdot&\cdot&\cdot&\cdot&\cdot&1&\cdot&-1&-1&\cdot&1&\cdot&\cdot&\cdot&\cdot&1&1&1&\cdot&-2&-1&-1&\cdot&1&\cdot&2&1&-1&\cdot&-1&1&-1&1&\cdot\\
15_a'&\cdot&\cdot&\cdot&\cdot&\cdot&\cdot&1&\cdot&\cdot&-1&1&\cdot&\cdot&-1&-1&\cdot&\cdot&1&1&\cdot&\cdot&-1&\cdot&\cdot&\cdot&\cdot&1&\cdot&-1&\cdot&\cdot&\cdot&\cdot&\cdot\\
210_a&\cdot&\cdot&\cdot&\cdot&\cdot&\cdot&\cdot&1&\cdot&\cdot&-1&-1&\cdot&\cdot&1&\cdot&\cdot&-2&\cdot&1&\cdot&1&\cdot&-1&\cdot&-1&-2&\cdot&1&2&\cdot&1&-1&-1\\
168_a&\cdot&\cdot&\cdot&\cdot&\cdot&\cdot&\cdot&\cdot&1&\cdot&-1&\cdot&\cdot&1&\cdot&-1&-1&-1&\cdot&2&1&2&\cdot&-1&-1&-3&-2&2&1&2&-1&2&-2&-1\\
105_b&\cdot&\cdot&\cdot&\cdot&\cdot&\cdot&\cdot&\cdot&\cdot&1&-1&\cdot&\cdot&1&1&-1&\cdot&-1&-1&1&\cdot&2&\cdot&\cdot&-1&-1&-2&1&1&1&\cdot&1&-1&-1\\
315_a'&\cdot&\cdot&\cdot&\cdot&\cdot&\cdot&\cdot&\cdot&\cdot&\cdot&1&\cdot&\cdot&-1&-1&\cdot&\cdot&1&1&-1&\cdot&-2&\cdot&1&1&1&3&-1&-2&-2&\cdot&-2&2&2\\
280_a'&\cdot&\cdot&\cdot&\cdot&\cdot&\cdot&\cdot&\cdot&\cdot&\cdot&\cdot&1&\cdot&\cdot&-1&\cdot&\cdot&1&\cdot&\cdot&\cdot&\cdot&\cdot&\cdot&\cdot&\cdot&1&\cdot&-1&-1&\cdot&\cdot&\cdot&1\\
70_a'&\cdot&\cdot&\cdot&\cdot&\cdot&\cdot&\cdot&\cdot&\cdot&\cdot&\cdot&\cdot&1&\cdot&\cdot&-1&\cdot&\cdot&\cdot&\cdot&1&1&\cdot&\cdot&-1&-1&\cdot&1&\cdot&\cdot&\cdot&1&-1&\cdot\\
105_c&\cdot&\cdot&\cdot&\cdot&\cdot&\cdot&\cdot&\cdot&\cdot&\cdot&\cdot&\cdot&\cdot&1&\cdot&\cdot&\cdot&\cdot&-1&\cdot&\cdot&1&1&\cdot&-1&\cdot&-1&\cdot&1&\cdot&\cdot&1&\cdot&-1\\
420_a&\cdot&\cdot&\cdot&\cdot&\cdot&\cdot&\cdot&\cdot&\cdot&\cdot&\cdot&\cdot&\cdot&\cdot&1&\cdot&\cdot&-1&-1&\cdot&\cdot&1&\cdot&\cdot&\cdot&\cdot&-2&\cdot&2&1&\cdot&1&-1&-2\\
210_b&\cdot&\cdot&\cdot&\cdot&\cdot&\cdot&\cdot&\cdot&\cdot&\cdot&\cdot&\cdot&\cdot&\cdot&\cdot&1&\cdot&\cdot&\cdot&-1&-1&-1&\cdot&\cdot&1&2&1&-2&\cdot&-1&1&-2&2&1\\
84_a&\cdot&\cdot&\cdot&\cdot&\cdot&\cdot&\cdot&\cdot&\cdot&\cdot&\cdot&\cdot&\cdot&\cdot&\cdot&\cdot&1&\cdot&\cdot&-1&\cdot&\cdot&\cdot&1&\cdot&1&\cdot&-1&\cdot&-1&1&\cdot&1&\cdot\\
105_c'&\cdot&\cdot&\cdot&\cdot&\cdot&\cdot&\cdot&\cdot&\cdot&\cdot&\cdot&\cdot&\cdot&\cdot&\cdot&\cdot&\cdot&1&\cdot&\cdot&\cdot&-1&\cdot&\cdot&\cdot&1&1&\cdot&-1&-1&\cdot&-1&1&1\\
420_a'&\cdot&\cdot&\cdot&\cdot&\cdot&\cdot&\cdot&\cdot&\cdot&\cdot&\cdot&\cdot&\cdot&\cdot&\cdot&\cdot&\cdot&\cdot&1&\cdot&\cdot&-1&-1&\cdot&1&\cdot&1&\cdot&-1&\cdot&\cdot&-2&1&2\\
210_b'&\cdot&\cdot&\cdot&\cdot&\cdot&\cdot&\cdot&\cdot&\cdot&\cdot&\cdot&\cdot&\cdot&\cdot&\cdot&\cdot&\cdot&\cdot&\cdot&1&\cdot&\cdot&\cdot&-1&\cdot&-1&-1&1&\cdot&2&-1&1&-2&-1\\
84_a'&\cdot&\cdot&\cdot&\cdot&\cdot&\cdot&\cdot&\cdot&\cdot&\cdot&\cdot&\cdot&\cdot&\cdot&\cdot&\cdot&\cdot&\cdot&\cdot&\cdot&1&\cdot&\cdot&\cdot&\cdot&-1&\cdot&1&\cdot&\cdot&-1&1&-1&\cdot\\
315_a&\cdot&\cdot&\cdot&\cdot&\cdot&\cdot&\cdot&\cdot&\cdot&\cdot&\cdot&\cdot&\cdot&\cdot&\cdot&\cdot&\cdot&\cdot&\cdot&\cdot&\cdot&1&\cdot&\cdot&-1&-1&-1&1&1&1&\cdot&2&-2&-2\\
280_a&\cdot&\cdot&\cdot&\cdot&\cdot&\cdot&\cdot&\cdot&\cdot&\cdot&\cdot&\cdot&\cdot&\cdot&\cdot&\cdot&\cdot&\cdot&\cdot&\cdot&\cdot&\cdot&1&\cdot&-1&\cdot&\cdot&\cdot&\cdot&\cdot&\cdot&1&\cdot&-1\\
70_a&\cdot&\cdot&\cdot&\cdot&\cdot&\cdot&\cdot&\cdot&\cdot&\cdot&\cdot&\cdot&\cdot&\cdot&\cdot&\cdot&\cdot&\cdot&\cdot&\cdot&\cdot&\cdot&\cdot&1&\cdot&\cdot&\cdot&\cdot&\cdot&-1&\cdot&\cdot&1&\cdot\\
210_a'&\cdot&\cdot&\cdot&\cdot&\cdot&\cdot&\cdot&\cdot&\cdot&\cdot&\cdot&\cdot&\cdot&\cdot&\cdot&\cdot&\cdot&\cdot&\cdot&\cdot&\cdot&\cdot&\cdot&\cdot&1&\cdot&\cdot&-1&\cdot&\cdot&\cdot&-1&1&1\\
168_a'&\cdot&\cdot&\cdot&\cdot&\cdot&\cdot&\cdot&\cdot&\cdot&\cdot&\cdot&\cdot&\cdot&\cdot&\cdot&\cdot&\cdot&\cdot&\cdot&\cdot&\cdot&\cdot&\cdot&\cdot&\cdot&1&\cdot&-1&\cdot&-1&1&-1&2&1\\
105_b'&\cdot&\cdot&\cdot&\cdot&\cdot&\cdot&\cdot&\cdot&\cdot&\cdot&\cdot&\cdot&\cdot&\cdot&\cdot&\cdot&\cdot&\cdot&\cdot&\cdot&\cdot&\cdot&\cdot&\cdot&\cdot&\cdot&1&\cdot&-1&-1&\cdot&-1&1&2\\
105_a&\cdot&\cdot&\cdot&\cdot&\cdot&\cdot&\cdot&\cdot&\cdot&\cdot&\cdot&\cdot&\cdot&\cdot&\cdot&\cdot&\cdot&\cdot&\cdot&\cdot&\cdot&\cdot&\cdot&\cdot&\cdot&\cdot&\cdot&1&\cdot&\cdot&-1&\cdot&-1&\cdot\\
15_a&\cdot&\cdot&\cdot&\cdot&\cdot&\cdot&\cdot&\cdot&\cdot&\cdot&\cdot&\cdot&\cdot&\cdot&\cdot&\cdot&\cdot&\cdot&\cdot&\cdot&\cdot&\cdot&\cdot&\cdot&\cdot&\cdot&\cdot&\cdot&1&\cdot&\cdot&\cdot&\cdot&-1\\
35_b'&\cdot&\cdot&\cdot&\cdot&\cdot&\cdot&\cdot&\cdot&\cdot&\cdot&\cdot&\cdot&\cdot&\cdot&\cdot&\cdot&\cdot&\cdot&\cdot&\cdot&\cdot&\cdot&\cdot&\cdot&\cdot&\cdot&\cdot&\cdot&\cdot&1&\cdot&\cdot&-1&-1\\
21_a'&\cdot&\cdot&\cdot&\cdot&\cdot&\cdot&\cdot&\cdot&\cdot&\cdot&\cdot&\cdot&\cdot&\cdot&\cdot&\cdot&\cdot&\cdot&\cdot&\cdot&\cdot&\cdot&\cdot&\cdot&\cdot&\cdot&\cdot&\cdot&\cdot&\cdot&1&\cdot&\cdot&\cdot\\
21_b&\cdot&\cdot&\cdot&\cdot&\cdot&\cdot&\cdot&\cdot&\cdot&\cdot&\cdot&\cdot&\cdot&\cdot&\cdot&\cdot&\cdot&\cdot&\cdot&\cdot&\cdot&\cdot&\cdot&\cdot&\cdot&\cdot&\cdot&\cdot&\cdot&\cdot&\cdot&1&-1&-1\\
7_a&\cdot&\cdot&\cdot&\cdot&\cdot&\cdot&\cdot&\cdot&\cdot&\cdot&\cdot&\cdot&\cdot&\cdot&\cdot&\cdot&\cdot&\cdot&\cdot&\cdot&\cdot&\cdot&\cdot&\cdot&\cdot&\cdot&\cdot&\cdot&\cdot&\cdot&\cdot&\cdot&1&\cdot\\
1_a'&\cdot&\cdot&\cdot&\cdot&\cdot&\cdot&\cdot&\cdot&\cdot&\cdot&\cdot&\cdot&\cdot&\cdot&\cdot&\cdot&\cdot&\cdot&\cdot&\cdot&\cdot&\cdot&\cdot&\cdot&\cdot&\cdot&\cdot&\cdot&\cdot&\cdot&\cdot&\cdot&\cdot&1\\
\end{block}
\end{blockarray}
\]
\end{minipage}
\end{turn}
%\end{sideways}
\restoregeometry
\normalsize
\newpage

\begin{proof}

The picture below shows the weights of the grading element $h_c$:

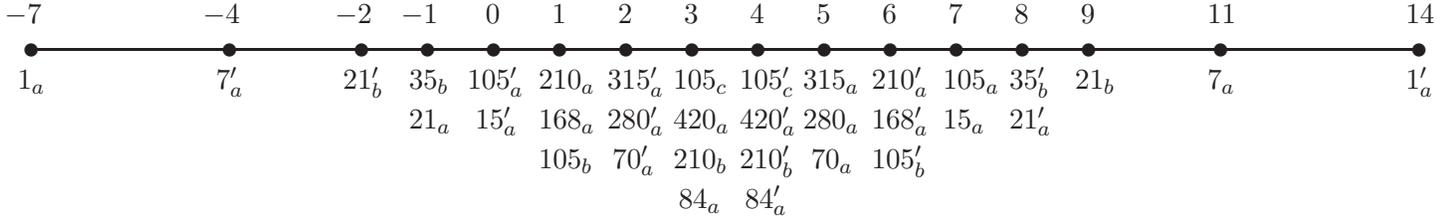
\begin{center}
\begin{picture}(525,100)
\put(0,80){\line(1,0){525}}
\put(0,80){\circle*{5}}
\put(75,80){\circle*{5}}
\put(125,80){\circle*{5}}
\put(150,80){\circle*{5}}
\put(175,80){\circle*{5}}
\put(200,80){\circle*{5}}
\put(225,80){\circle*{5}}
\put(250,80){\circle*{5}}
\put(275,80){\circle*{5}}
\put(300,80){\circle*{5}}
\put(325,80){\circle*{5}}
\put(350,80){\circle*{5}}
\put(375,80){\circle*{5}}
\put(400,80){\circle*{5}}
\put(450,80){\circle*{5}}
\put(525,80){\circle*{5}}
\put(-10,90){$-7$}
\put(65,90){$-4$}
\put(115,90){$-2$}
\put(140,90){$-1$}
\put(172,90){$0$}
\put(197,90){$1$}
\put(222,90){$2$}
\put(247,90){$3$}
\put(272,90){$4$}
\put(297,90){$5$}
\put(322,90){$6$}
\put(347,90){$7$}
\put(372,90){$8$}
\put(397,90){$9$}
\put(445,90){$11$}
\put(520,90){$14$}
\put(-5,65){$1_a$}
\put(70,65){$7_a'$}
\put(118,65){$21_b'$}
\put(143,65){$35_b$}
\put(143,50){$21_a$}
\put(165,65){$105_a'$}
\put(168,50){$15_a'$}
\put(192,65){$210_a$}
\put(192,50){$168_a$}
\put(192,35){$105_b$}
\put(218,65){$315_a'$}
\put(218,50){$280_a'$}
\put(220,35){$70_a'$}
\put(243,65){$105_c$}
\put(243,50){$420_a$}
\put(243,35){$210_b$}
\put(245,20){$84_a$}
\put(268,65){$105_c'$}
\put(268,50){$420_a'$}
\put(268,35){$210_b'$}
\put(270,20){$84_a'$}
\put(292,65){$315_a$}
\put(292,50){$280_a$}
\put(295,35){$70_a$}
\put(318,65){$210_a'$}
\put(318,50){$168_a'$}
\put(318,35){$105_b'$}
\put(345,65){$105_a$}
\put(345,50){$15_a$}
\put(370,65){$35_b'$}
\put(370,50){$21_a'$}
\put(395,65){$21_b$}
\put(445,65){$7_a$}
\put(520,65){$1_a'$}
\end{picture}
\end{center}

\normalsize

The upper left hand corner of the decomposition matrix looks as follows after (dim Hom), (RR) have been applied. Columns $315_a'$, $280_a'$, $420_a$, and $210_b$ are also columns of the Hecke algebra decomposition matrix:

\[
\small
\begin{blockarray}{cccccccccccccccccccc}
&& 1_a & 7_a' & 21_b' & 35_b & 21_a & 105_a' & 15_a' & 210_a & 168_a & 105_b & 315_a' & 280_a' & 70_a' & 
105_c & 420_a & 210_b & 84_a & 105_c' \\
\begin{block}{cc(cccccccccccccccccc}
\star & 1_a & 1 & \cdot & 1 & 1 & \cdot & \cdot & 1 & {\color{red}?} &  {\color{red}?}  &  {\color{red}?}  & \cdot & \cdot &  {\color{red}?}  &  {\color{red}?}  & \cdot & 1 &  {\color{red}?}  &  {\color{red}?} \\
\star? & 7_a' & \cdot & 1 & 1 & 1 & \cdot & 1 & \cdot &  {\color{red}?}  &  {\color{red}?}  & 1 & 1 & \cdot &  {\color{red}?}  &  {\color{red}?}  & \cdot & \cdot &  {\color{red}?}  &  {\color{red}?} \\
\bullet & 21_b' & \cdot & \cdot & 1 & \cdot & \cdot & \cdot & \cdot & 1 & 1 & 1 & 1 & \cdot & \cdot &  {\color{red}?}  & \cdot & 1 &  {\color{red}?}  &  {\color{red}?}  \\
\bullet & 35_b & \cdot & \cdot & \cdot & 1 & \cdot & \cdot & \cdot & \cdot & 1 & 1 & 1 & \cdot & 1 &  {\color{red}?}  & \cdot & 1 &  {\color{red}?}  &  {\color{red}?} \\
\star? & 21_a & \cdot & \cdot & \cdot & \cdot & 1 & 1 & \cdot &  {\color{red}?}  &  {\color{red}?}  & \cdot & \cdot & 1 & \cdot &  {\color{red}?}  & \cdot & \cdot &  {\color{red}?}  &  {\color{red}?} \\
\bullet & 105_a' & \cdot & \cdot & \cdot & \cdot & \cdot & 1 & \cdot & 1 & 1 & \cdot & 1 & 1 & \cdot &  {\color{red}?}  & 1 & \cdot &  {\color{red}?}  &  {\color{red}?} \\
\star & 15_a' & \cdot & \cdot & \cdot & \cdot & \cdot & \cdot & 1 & \cdot & \cdot & 1 & \cdot & \cdot & \cdot &  {\color{red}?}  & \cdot & 1 & \cdot &  {\color{red}?} \\
\bullet & 210_a & \cdot & \cdot & \cdot & \cdot & \cdot & \cdot & \cdot & 1 & \cdot & \cdot & 1 & 1 & \cdot &  {\color{red}?}  & 1 & \cdot & \cdot &  {\color{red}?} \\
\bullet & 168_a & \cdot & \cdot & \cdot & \cdot & \cdot & \cdot & \cdot & \cdot & 1 & \cdot & 1 & \cdot & \cdot &  {\color{red}?}  & 1 & 1 & 1 &  {\color{red}?} \\
\bullet & 105_b & \cdot & \cdot & \cdot & \cdot & \cdot & \cdot & \cdot & \cdot & \cdot & 1 & 1 & \cdot & \cdot &  {\color{red}?}  & \cdot & 1 & \cdot &  {\color{red}?} \\
& 315_a' & \cdot & \cdot & \cdot & \cdot & \cdot & \cdot & \cdot & \cdot & \cdot & \cdot & 1 & \cdot & \cdot & 1 & 1 & \cdot & \cdot &  {\color{red}?} \\
& 280_a' & \cdot & \cdot & \cdot & \cdot & \cdot & \cdot & \cdot & \cdot & \cdot & \cdot & \cdot & 1 & \cdot & \cdot & 1 & \cdot & \cdot &  {\color{red}?} \\
\bullet & 70_a' & \cdot & \cdot & \cdot & \cdot & \cdot & \cdot & \cdot & \cdot & \cdot & \cdot & \cdot & \cdot  & 1 & \cdot & \cdot & 1 & \cdot & \cdot\\
\bullet & 105_c & \cdot & \cdot & \cdot & \cdot & \cdot & \cdot & \cdot & \cdot & \cdot & \cdot & \cdot & \cdot & \cdot & 1 & \cdot & \cdot & \cdot &  {\color{red}?} \\
& 420_a & \cdot & \cdot & \cdot & \cdot & \cdot & \cdot & \cdot & \cdot & \cdot & \cdot & \cdot & \cdot & \cdot & \cdot & 1 & \cdot & \cdot & 1\\
& 210_b & \cdot & \cdot & \cdot & \cdot & \cdot & \cdot & \cdot & \cdot & \cdot & \cdot & \cdot & \cdot & \cdot & \cdot & \cdot & 1 & \cdot & \cdot\\
\bullet & 84_a & \cdot & \cdot & \cdot & \cdot & \cdot & \cdot & \cdot & \cdot & \cdot & \cdot & \cdot & \cdot & \cdot & \cdot & \cdot & \cdot & 1 & \cdot\\
\bullet & 105_c' & \cdot & \cdot & \cdot & \cdot & \cdot & \cdot & \cdot & \cdot & \cdot & \cdot & \cdot & \cdot & \cdot & \cdot & \cdot & \cdot & \cdot & 1\\
\end{block}
\end{blockarray}
\]
\normalsize

The remaining columns of the matrix are all preserved by KZ functor and can be seen in the table for $E_7$, $e=6$ in \cite{GJ}.

First, $[\m(315_a'):\el(105_c')],\;[\m(280_a'):\el(105_c')],\;[\el(105_c):\m(105_c')]$ are all $0$ as may be seen by examining restrictions to $H_{\frac{1}{6}}(E_6)$. Set $\alpha:=[\el(105_c):\m(105_c')]$; then $\el(105_c)=\m(105_c)-\alpha\m(105_c')-\m(420_a')...$ is the decomposition through $h_c$-weight $4$. Then 
\begin{align*}
\Res\el(105_c)&=\Res\m(105_c)-\alpha\Res\m(105_c')-\Res\m(420_a')...\\
&=\m(81_p)+\m(24_p')-\alpha(\m(81_p')+\m(24_p))-\m(64_p')-\m(81_p')-\m(24_p')\\
&\qquad-\m(81_p)-\m(90_s)-\m(80_s)\\
&=-\alpha(81_p'+24_p)-64_p'-81_p'-90_s-80_s...
\end{align*}
 can only be a module if $\alpha=0$, because there is no representation $\tau$ with $h_c(\tau)>h_c(\el(105_c'))=4$ such that $\Res(\tau)$ contains $24_p$ as a summand. Likewise with $[\m(280_a'):\el(105_c')]:=\alpha$, writing the decomposition of $\el(280_a')$ through $h_c$-weight $4$ reps as $\el(280_a')=\m(280_a')-\m(420_a)+(1-\alpha)\m(105_c')...$, 
$$\Res\el(280_a')=\m(15_p)+\m(30_p)-\m(24_p)-\m(81_p')-\m(80_s)+(1-\alpha)(\m(24_p)+\m(81_p'))$$
The Verma $\m(15_p)$ does not belong to the same block as $\m(24_p)$, while $\m(30_p)$ does not contain any submodule with lowest weight $\m(24_p)$. So the only way the expression above can be that of a module is if $\alpha=0$ so that the $-\m(24_p)$ coming from $-\Res\m(420_a)$ is canceled out. And now for $[\m(315_a'):\el(105_c')]:=\alpha$, write $\el(315_a')=\m(315_a')-\m(105_c)-\m(420_a)+(1-\alpha)\m(105_c')+\m(420_a')+\m(210_b')....$, the Verma-decomposition through $h_c$-weight $4$. Then 
$$\Res\el(315_a')=\m(30_p)+\m(60_p)-\m(81_p)-\m(24_p')-\m(24_p)-\m(81_p')-\m(90_s)...$$
are the terms in $\Res\el(315_a')$ coming from those $\tau$ with $h_c(\tau)\leq3$. Since neither $\el(30_p)$ nor $\el(60_p)$ contains $\m(24_p)$ in its composition series, and $105_c'$ is the only $\tau$ with $h_c(\tau)>3$ such that $24_p\subset\Res(\tau)$, it must be that $\alpha=0$ so that $+\Res\m(105_c')$ cancels out the term $-\m(24_p)$.

The Verma-decompositions and characters of the simples of less than full support which reside to the south of row $315_a'$ now follow from the calculations above:

\begin{align*}
&\el(70_a')=\m(70_a')-\m(210_b)+\m(84_a')+\m(315_a)-\m(210_a')-\m(168_a')+\m(105_a)+\m(21_b)-\m(7_a)\\
&\chi_{\el(70_a')}(t)=t^2\frac{7t^4 + 35t^3 + 84t^2 + 140t + 70}{(1-t)^2}\\
&\dim\Supp\el(70_a')=2\\
\\
&\el(105_c)=\m(105_c)-\m(420_a')+\m(315_a)+\m(280_a)-\m(210_a')-\m(105_b)+\m(15_a)+\m(21_b)-\m(1_a')\\
&\chi_{\el(105_c)}(t)=t^3\frac{t^6 + 5t^5 + 15t^4 + 35t^3 + 70t^2 + 105t + 105}{(1-t)^2}\\
&\dim\Supp\el(105_c)=2\\
\\
&\el(84_a)=\m(84_a)-\m(210_b')+\m(70_a)+\m(168_a')-\m(105_a)-\m(35_b')+\m(21_a')+\m(7_a)\\
&\chi_{\el(84_a)}(t)=t^3\frac{7t^4 + 28t^3 + 70t^2 + 126t + 84}{(1-t)^3}\\
&\dim\Supp\el(84_a)=3\\
\\
&\el(105_c')=\m(105_c')-\m(315_a)+\m(168_a')+\m(105_b')-\m(15_a)-\m(35_b')-\m(21_b)+\m(7_a)+\m(1_a')\\
&\chi_{\el(105_c')}(t)=t^4\frac{t^6 + 4t^5 + 10t^4 + 27t^3 + 63t^2 + 105t + 105}{(1-t)^3}\\
&\dim\Supp\el(105_c')=3\\
\end{align*}

Next, finite-dimensionality of two irreps comes for free. First, $\el(1_a)$ must be finite-dimensional since $6$ is an elliptic number of $E_7$ \cite{VV}. Secondly, $105_b$ generates a subrep of $\m(15_a')$. Then $\dim\el(15_a')[1]=15\cdot7-105=0$, so $\el(15_a')$ is finite-dimensional and of dimension $15.$ The graded character of $\el(15_a')$ is:
$$\chi_{\el(15_a')}(t)=15=\frac{15-105t+315t^2-525t^3+525t^4-315t^5+105t^6-15t^7}{(1-t)^7}$$
from which it is easy to deduce (using the columns of the decomposition matrix given by the Hecke algebra decomposition matrix and (Symm)) that:
\begin{align*}
\el(15_a')&=\m(15_a')-\m(105_b)+\m(315_a')-\m(105_c)-\m(420_a)\\&\qquad\qquad+\m(105_c')+\m(420_a')-\m(315_a)+\m(105_b')-\m(15_a)\\
\end{align*}

The spherical representation $\el(1_a)$ has Verma-decomposition beginning: $\el(1_a)=\m(1_a)-\m(21_b')-\m(35_b)-\m(15_a')...$ As this is the decomposition of $\el(1_a)$ up through $h_c$-weight $0$, it determines the graded character and dimension of $\el(1_a)$:
\begin{align*}
\chi_{\el(1_a)}(t)&=t^{-7}+t^7+7(t^{-6}+t^6)+28(t^{-5}+t^5)+84(t^{-4}+t^4)+210(t^{-3}+t^3)\\&\qquad+441(t^{-2}+t^2)+742(t^{-1}+t)+868\\
\\
\dim\el(1_a)&=\chi_{\el(1_a)}(1)=3894\\
\end{align*}
Multiplying $\chi_{\el(1_a)}(t)$ by $(1-t)^7$, the coefficient of $t^k$ then records the signed dimension with multiplicities of the lowest weights $\tau$ with $h_c(\tau)=k$ of the Vermas $\m(\tau)$ in the Verma-decomposition of $\el(1_a)$:
\begin{align*}
\chi_{\el(1_a)}(t)&=\frac{t^{-7}-21t^{-2}-35t^{-1}-15+588t-910t^2+420t^3-420t^4+910t^5-588t^6}{(1-t)^7}\\
&\qquad+\frac{15t^7+35t^8+21t^9-t^{14}}{(1-t)^7}\\
\end{align*}

Consider the restriction of $\el(1_a)$ to $H_{\frac{1}{6}}(E_6)$ for $\el(1_a)$ in its decomposition found so far. $$\Res 1_a-\Res 21_b'-\Res 35_b'-\Res 15_a'=-2(20_p)-2(15_q)$$
Since $\el(1_a)$ is finite-dimensional, $\Res\el(1_a)=0$.  Let $c_{210},\;c_{168},\;c_{105}$ be the signed multiplicities of $\m(210_a),\;\m(168_a),\;\m(105_a)$ in $\el(1_a)$. $105_b$ is the only rep to the right of $h_c$-weight $0$ in this block whose restriction to $E_6$ contains $15_q$. Therefore $c_{105_b}=2$. Both $\Res 210_a$ and $\Res 168_a$ contain one copy of $20_p$. Moreover, $c_{210}210+c_{168}168=378$, from the coefficient of $t$ in the character of $\el(1_a)$. So $c_{210}=c_{168}=1$. Next, for the Vermas at $h_c$-weight $2$: computing $\Res$ of the decomposition of $\el(1_a)$ so far gives $15_p+3(30_p)+2(64_p)+24_p+81_p+3(60_p)$, and as $280_a'$ is the only rep to the right of $h_c$-weight $1$ whose restriction contains $15_p$, it's necessary that $[\el(1_a):\m(280_a')]=-1$ in order for the restriction of $\el(1_a)$ to be $0$. $30_p$ is contained only in $\Res 280_a'$, $\Res 315_a'$ for reps to the right of $h_c$-weight $1$, so it must be that $[\el(1_a):\m(315_a')]=-2$. This contributes $-910t^2$ to the numerator of the graded character of $\el(1_a)$ and thus $[\el(1_a):\m(70_a')]=0$. Continuing to the $h_c$-weight $3$ irreps, the restriction to $E_6$ of the decomposition of $\el(1_a)$ is $24_p-64_p-2(81_p)+60_p-2(80_s)-90_s$. 
$420_a$ is the last chance to kill $64_p$ occuring in the restriction to $E_6$, while $210_b$ is the last chance to kill $60_p$ in the restriction, so $[\el(1_a):\m(420_a)]=1$ and $[\el(1_a):\m(210_b)]=-1$. Then as $420=420-210+c_{84_a}84+c_{105_c}105_c$
it must be that $[\el(1_a):\m(84_a)]=0$ and $[\el(1_a):\m(105_c)]=2$.
(Symm) then gives the second half of the Verma-decomposition of $\el(1_a)$:
\begin{align*}
\el(1_a)&=\m(1_a)-\m(21_b')-\m(35_b)-\m(15_a')+\m(210_a)+\m(168_a)+2\m(105_b)-2\m(315_a')\\
&-\m(280_a')+\m(420_a)-\m(210_b)+2\m(105_c)-2\m(105_c')+\m(210_b')-\m(420_a')+\m(280_a)\\
&+2\m(315_a)-2\m(105_b')-\m(168_a')-\m(210_a')+\m(15_a)+\m(35_b')+\m(21_b)-\m(1_a')
\end{align*}

To determine the decomposition of $\el(7_a'),$ induce the spherical representation of $H_{\frac{1}{6}}(E_6)$:
\begin{align*}
\Ind\el(1_p)&=\Ind\left(\m(1_p)-\m(20_p)+\m(30_p)-\m(15_q)+\m(24_p)+\m(60_p)-2\m(80_s)\right.\\&\left.\qquad\qquad+\m(24_p')+\m(60_p')+\m(30_p')-\m(15_q')-\m(20_p')+\m(1_p')\right)\\
&=\m(1_a)+\m(1_a')+\m(7_a)+\m(7_a')-2\m(35_b)-2\m(35_b')-\m(105_a')-\m(105_a)\\&-\m(15_a')-\m(15_a)+\m(168_a)+\m(168_a')+\m(84_a)+\m(84_a')-\m(420_a)\\&-\m(420_a')+\m(105_c')+\m(105_c)+\m(105_b)+\m(105_b')-\m(210_b)-\m(210_b')\\&+\m(70_a')+\m(70_a)+\m(280_a')+\m(280_a)
\end{align*}
Whatever this module is, $\el(1_a)$ appears in its composition series into simples. Subtracting the expression for $\el(1_a)$ from $\Ind\el(1_p)$ (and writing simply ``$\tau$" instead of ``$\m(\tau)$" everywhere) gives:
\begin{align*}
&\Ind\el(1_p)-\el(1_a)=7_a'+21_b'-35_b-105_a'-210_a-105_b+70_a'+2(315_a')+2(280_a')\\
&\qquad-2(420_a)-105_c+84_a+84_a'-2(210_b')+3(105_c')+70_a-2(315_a)+3(105_b')+2(168_a')\\
&\qquad+210_a'-2(15_a)-105_a-3(35_b')-21_b+7_a+2(1_a')\\
\end{align*}

The decomposition numbers found so far imply that $$\el(21_b')=\m(21_b')-\m(210_a)-\m(105_b)-\m(168_a)+2\m(315_a')+\m(280_a')....$$ is the decomposition of $\el(21_b')$ up through $h_c$-weight $2$.
Moreover, $\dim\Hom(\m(21_b'),\m(7_a'))=\dim\Hom(\m(7_a),\m(21_b))=1$, and so $\el(7_a')=\m(7_a')-\m(21_b')...$. Consequently, $2\el(21_b')$ occurs as a summand of $\Ind\el(1_p)-\el(1_a)$; subtracting the terms of its decomposition known so far gives:
\begin{align*}
\Ind\el(1_p)-\el(1_a)-2\el(21_b')&=\m(7_a')-\m(21_b')-\m(35_b)-\m(105_a')+\m(105_b)+2\m(168_a)\\&\qquad+\m(210_a)
-2\m(315_a')+\m(70_a')+...
\end{align*}
where the $...$ denotes the tail of reps with lowest $h_c$-weight bigger than $2$. On the one hand, this expression contains the beginning of $\el(7_a')$, and in particular $\dim\el(7_a')[-1]=406$. Then $\dim\el(7_a')[1]\geq406$ by $\s$-module theory. On the other hand, the expression in $\Ind\el(1_p)-\el(1_a)-2\el(21_b')$ implies $\dim\el(7_a')[1]\leq406$, since subtracting off additional composition factors could not raise the dimension. So $\dim\el(7_a')[1]=406=\dim\el(7_a')[-1]$ and the expression above coincides with the decomposition of $\el(7_a')$ up through its  Verma-factors whose lowest weights have $h_c$-weight $1$:
$$\el(7_a')=\m(7_a')-\m(21_b')-\m(35_b)-\m(105_a')+\m(105_b)+2\m(168_a)+\m(210_a)...$$ 
Next, consider what happens with the Vermas with $h_c$-weight $2$. The columns $315_a'$ and $280_a'$ of the decomposition matrix are given by the Hecke algebra decomposition matrix; taking dot product with the vector of Verma multiplicities in $\el(7_a')$ shows that $\m(315_a')$ must appear with coefficient $-2$ and $\m(280_a')$ with coefficient $0$. In view of the single occurrence of $\m(70_a')$ in $\Ind\el(1_p)-\el(1_a)-2\el(21_b')$, and also taking into account that all the preceding terms belong to $\el(7_a')$, it follows that $\m(70)$ appears with either coefficient $0$ or $1$ in $\el(7_a')$. But $\dim\el(7_a')[2]\geq\dim\el(7_a')[-2]=175$ and the only way $\dim\el(7_a')[2]\geq175$ is if $+\m(70_a')$ belongs to the composition series of $\el(7_a')$, in which case $\dim\el(7_a')[2]=175$. Thus up through the $h_c$-weight 2 irreps, the decomposition is:
$$\el(7_a')=\m(7_a')-\m(21_b')-\m(35_b)-\m(105_a')+\m(105_b)+2\m(168_a)+\m(210_a)-2\m(315_a')+\m(70_a')...$$ 
The dimension of $\el(7_a')$ in graded degree $0$ is $532$ and thus $\dim\el(7_a')[2]<\dim\el(7_a')[0]$; now by (E), $\el(7_a')$ is finite-dimensional.

The graded character and dimension of $\el(7_a')$ may then be calculated using the Verma-decomposition up through weight $0$ to calculate the dimension of each graded piece:
\begin{align*}
\chi_{\el(7_a')}(t)&=7(t^{-4}+t^4)+49(t^{-3}+t^3)+175(t^{-2}+t^2)+406(t^{-1}+t)+532\\
&=\frac{-399(t^3-t^4) - 560(t^2-t^5) + 651(t-t^6) - 105(1-t^7) - 35(t^{-1}-t^8)}{(1-t)^7} \\
&\qquad+\frac{- 21(t^{-2}-t^9) + 7(t^{-4}-t^{11})}{(1-t)^7}\\
\\
\dim\el(7_a')&=\chi_{\el(7_a')}(1)=1806\\
\end{align*}

Using the coefficients in the numerator to determine the total signed dimension of lowest weights of Vermas with multiplicities at a given $h_c$-weight, together with the fact that columns of the decomposition matrix %such as the columns $420_a$ and $210_b$ must 
have dot product $0$ with the vector of Verma-multiplicities of $\el(7_a')$, one arrives at the Verma-decomposition of $\el(7_a')$:
\begin{align*}
\el(7_a')&=\m(7_a')-\m(21_b')-\m(35_b)-\m(105_a')+\m(210_a)+2\m(168_a)+\m(105_b)\\&\qquad-2\m(315_a')+\m(70_a')+\m(105_c)-2\m(210_b)-\m(84_a)-\m(105_c')\\&\qquad+2\m(210_b')+\m(84_a')+2\m(315_a)-\m(70_a)-\m(210_a')-2\m(168_a')\\&\qquad-\m(105_b')+\m(105_a)+\m(35_b')+\m(21_b')-\m(7_a)\\
\end{align*}

Next, the decompositions of $\el(1_a)$ and $\el(7_a')$ provide information relevant to the Verma-decomposition of $\el(21_b')$: recall that $\Ind_{E_6}^{E_7}\el(1_p)=\el(1_a)+\el(7_a')+2\el(21_b')+?$ on the level of K-theory, and that $\dim\Supp\Ind_{E_6}^{E_7}\el(1_p)=1$. So, unless $\el(21_b')$ is finite-dimensional, $\dim\Supp\el(21_b')=1$. Indeed, the decomposition numbers already known give enough of the Verma-decomposition of $\el(21_b')$ to see that $\dim\el(21_b')[1]>\dim\el(21_b')[-1]$, so $\el(21_b')$ is not finite-dimensional and therefore $\dim\Supp\el(21_b')=1$. Additionally, $\el(21_b')$ belongs to the composition series of the following module whose support has dimension $1$:
\begin{align*}
&\frac{1}{2}\left(\Ind\el(1_p)-\el(1_a)-\el(7_a')\right)=21_b'-210_a-168_a-105_b+2(315_a')+280_a'\\
&\qquad-105_c-420_a+210_b+84_a+2(105_c')-2(210_b')-2(315_a)+70_a+210_a'\\
&\qquad+2(168_a')+2(105_b')-105_a-15_a-2(35_b')-21_b+7_a+1_a'\\
\end{align*}
This could be $\el(21_b')$ itself, or it could contain some other simples $\el(\tau)$ with $\dim\Supp\el(\tau)=0$ or $1$ in its composition series. %(It is easy to check that no other finite-dimensional simples appear). 
The Verma-decomposition of $\el(21_b')$ coincides with that of $\frac{1}{2}\left(\Ind\el(1_p)-\el(1_a)-\el(7_a')\right)$ up through $h_c$-weight $2$. Any additional simple factor $\el(\tau)$ of $\frac{1}{2}\left(\Ind\el(1_p)-\el(1_a)-\el(7_a')\right)$ would thus have to satisfy the condition $h_c(\tau)\geq3$. %(in particular no finite-dimensional representation could appear, since they are generated in graded degree at most $0$). 
There are only three such $\tau$ of less than full support: $\tau=105_c,$ $84_a,$ and $105_c'$. However, $\el(105_c),$ $\el(84_a),$ and $\el(105_c')$ all have at least $2$-dimensional support, so they cannot occur in the composition series of $\frac{1}{2}\left(\Ind\el(1_p)-\el(1_a)-\el(7_a')\right)$. Therefore $\el(21_b')=\frac{1}{2}\left(\Ind\el(1_p)-\el(1_a)-\el(7_a')\right)$ with decomposition as above. %And from this we learn that $[\m(21_b'):\el(84_a)]=0$.

There is a fourth finite-dimensional irrep in this block: $\el(21_a)$. It occurs in the socle of $\Ind_{E_6}^{E_7}\el(6_p)$ :
\begin{align*}
\Ind\el(6_p)&=\Ind\m(6_p)-\Ind\m(20_p)+\Ind\m(24_p)+\Ind\m(60_p)-\Ind\m(80_s)-\Ind\m(60_s)\\
&\qquad+\Ind\m(60_p')+\Ind\m(24_p')-\Ind\m(20_p')+\Ind\m(6_p')\\
&=\m(21_a)+\m(21_a')+\m(7_a)+\m(7_a')-\m(35_b)-\m(35_b')-\m(210_a)-\m(210_a')\\
&\qquad-\m(21_b)-\m(21_b')+\m(105_c)+\m(105_c')+\m(105_b)+\m(105_b')+\m(168_a)\\
&\qquad+\m(168_a')+\m(70_a)+\m(70_a')-\m(210_b)-\m(210_b')\\
\end{align*}
Out of all the Vermas $\m(\tau)$ appearing in this expression, $\m(7_a')$ has the minimal $h_c$-weight, so $\m(7_a')$ belongs to the socle. Subtracting $\el(7_a')$ from this expression, and just writing the Vermas with $h_c$-weight less than $3$:
$$\Ind\el(6_p)-\el(7_a')=\m(21_a)+\m(105_a')-\m(168_a)-2\m(210_a)+2\m(315_a')...$$
$\el(21_a)$ then belongs to the socle as well since $\m(7_a')$ and $\m(21_a)$ are not linked, and since $21_a$ has minimal $h_c$-weight of all the terms excepting $7_a'$. The decomposition numbers known at this point give that
$$\el(105_a')=\m(105_a')-\m(210_a)-\m(168_a)+\m(315_a')...$$
where ``$...$" consists of some linear combination of Vermas of $h_c$-weight at least $3$. Also, $\el(21_a)=\m(21_a)-\m(105_a')...$ where the rest of the Vermas appearing in $\el(21_a)$ have lowest weights $\tau$ satisfying $h_c(\tau)\geq1$ $h_c$-weight at least $1$.Thus $2\el(105_a')$ must be a summand of $\Ind\el(6_p)-\el(7_a')$. Subtracting it, we find that up through Vermas of $h_c$-weight $2$,
$$\Ind\el(6_p)-\el(7_a')-2\el(105_a')=\m(21_a)-\m(105_a')+\m(168_a)+....$$
Now either this expression is $\el(21_a)$ or $\el(21_a)+\el(168_a)$, up through the Vermas of weight $2$ in the decomposition into simple modules. Write
$$\el(21_a)=\m(21_a)-\m(105_a')+c_{168}\m(168_a)...$$
with $c_{168}=0$ or $1$, and calculate the dimension of the graded piece at weight $1$:
$$\el(21_a)[1]=21{8\choose6}-105{7\choose6}+c_{168}168=-147+c_{168}168$$
Since the dimension cannot be negative, $c_{168}=1$ and $\m(21_a)-\m(105_a')+\m(168_a)+....$ is the decomposition of $\el(21_a)$ up through $h_c$-weight $2$. Then 
$$\dim\el(21_a)[2]=21{9\choose6}-105{8\choose6}+168{7\choose6}=0$$
and therefore $\el(21_a)$ is finite-dimensional, and $$\dim\el(21_a)=2\cdot21+(7\cdot21-105)=84$$
The graded character of $\el(21_a)$ is 
\begin{align*}
\chi_{\el(21_a)}(t)&=21t^{-1}+42+21t\\
&=\frac{21t^{-1}-105+168t-294t^3+294t^4-168t^6+105t^7-21t^8}{(1-t)^7}\\
\end{align*}
The columns $420_a$ and $210_b$ of the decomposition matrix, which were given by the Hecke algebra decomposition matrix, indicate that of the four irreps of $h_c$-weight $3$, namely $105_c,$ $210_b,$ $420_a,$ and $84_a$, $\m(420_a)$ does not appear in $\el(21_a)$ while $-\m(210_b)$ does (the dot product of the Verma-vector for $\el(21_a)$ with those columns must be $0$). So $$-294=-210+c_{84}84+c_{105}105$$
where $c_{84},\;c_{105}$ are the multiplicities of $\m(84_a),\;\m(105_c)$ in $\el(21_a)$. Therefore $c_{84}=-1$ and $c_{105}=0$, and  
\begin{align*}
\el(21_a)&=\m(21_a)-\m(105_a')+\m(168_a)-\m(210_b)-\m(84_a)+\m(84_a')+\m(210_b')\\
&\qquad-\m(168_a')+\m(105_a)-\m(21_a')\\
\end{align*}
%Consequently, $[\m(21_a):\el(210_a)]=1$ and $[\m(21_a):\el(168_a)]=0$.

Now pulling at the thread of $\el(21_a)$ disentangles $\el(105_a')$ from the induced module in which it has appeared. Writing $\tau$ as shorthand for $\m(\tau)$, one has:
\begin{align*}
&\frac{1}{2}\left(\Ind\el(6_p)-\el(7_a')-\el(21_a)\right)=105_a'-210_a-168_a+315_a'+210_b+84_a+105_c'-2(210_b')\\
&\qquad -84_a'-315_a+70_a+2(168_a')+105_b'-105_a-35_b'+21_a'-21_b+7_a\\
\end{align*}
This coincides with the decomposition of $\el(105_a')$ up through $h_c$-weight $2$. Since\\ $\dim\Supp\frac{1}{2}\left(\Ind\el(6_p)-\el(7_a')-\el(21_a)\right)=1$ and there are no $\el(\tau)$ in this block such that $h_c(\tau)>2$ and $\dim\Supp\el(\tau)\leq1,$ there can be no further simples in the composition series of $\frac{1}{2}\left(\Ind\el(6_p)-\el(7_a')-\el(21_a)\right)$. Therefore $\el(105_a')=\frac{1}{2}\left(\Ind\el(6_p)-\el(7_a')-\el(21_a)\right)$. %And from this we learn that $[\m(105_a'):\el(84_a)]=0$ and therefore $[\m(21_a):\el(84_a)]=0$ as well.

There are now only four irreducible representations left to find in this block: $\el(35_b)$, $\el(210_a)$, $\el(168_a)$, and $\el(105_b)$. Once their Verma-decompositions are known, then so is the inverse of the decomposition matrix consisting of the multiplicities $([\el(\tau):\m(\sigma)])$, and then the decomposition matrix is known as well. %The difficulty right now is that we do not know if $\el(105_c)$ occurs in these representations, and because $\Res\el(105_c)=0$, it cannot be detected by restriction to $E_6$. So a bit of work will be required.

 Induce $\el(15_q)$ from $H_{\frac{1}{6}}(E_6)$:
\begin{align*}
\Ind\el(15_q)&=35_b+15_a'-168_a-70_a'+420_a+210_b'-105_c-84_a'-280_a-105_b'\\&\qquad+105_a+35_b'-7_a-1_a'\\
\end{align*}
Clearly $\Ind\el(15_q)$ contains $\el(35_b)$ and $\el(15_a')$ in its composition series, but using the decomposition numbers found so far to calculate $\el(35_b)$ up through Vermas of $h_c$-weight $1$:
$$\el(35_b)=\m(35_b)-\m(168_a)-\m(105_b)...$$
shows that $\Ind\el(15_q)$ also contains $2$ copies of $\el(105_b)$ (one copy will be produced by subtracting $\el(15_a')$, the other by subtracting $\el(35_b)$). So the dimensions of support of $\el(35_b)$ and $\el(105_b)$ are $1$ or $2$ since $\dim\Supp\Ind\el(15_q)=2$, and $\el(35_b),$ $\el(105_b)$ are not finite-dimensional. %If $\el(35_b)$ has dimension of support $2$ then its restriction to $H_{\frac{1}{6}}(E_6)$ must have dimension of support $1$, so must contain $\el(15_q)$, since that's the only $\el(\tau)$ for  $H_{\frac{1}{6}}(E_6)$ with $\dim\Supp\el(\tau)=1$. 

But

$$\Ind 15_q=35_b+15_a'+105_b$$

so $15_q$ does not appear in $\Res\el(35_b)$, and $1_p$ and $6_p$ cannot appear either since any $\tau$ with $1_p\subset\Res(\tau)$ has $h_c(\tau)<h_c(15_q)$, while the only $\tau$ with $6_p\subset\Res(\tau)$ and $h_c(\tau)>h_c(15_q)$ is $\tau=105_a'$, but $\m(105_a')$ does not occur in $\el(15_q)$. Moreover, $\m(15_q)$ appears only in $\el(15_q)$ and $\el(1_p)$. Therefore no irreducible representation of dimension of support $\leq1$ appears in $\Res\el(35_b)$. But $\Res$ lowers dimension of support by at least $1$, and $\dim\Supp(35_b)\leq2$, so it must be that $\Res\el(35_b)=0$. %We will use this to find the decomposition of $\el(35_b)$ up to a possible factor of $\el(105_c)$, but for now we bracket this question, make a note that $\dim\Supp\el(105_b)\leq2$, and get more information about $\el(105_b)$ by inducing a different representation of $H_{\frac{1}{6}}(E_6)$.

As a first pass at $\el(210_a)$ and $\el(168_a)$, their dimensions of support and restriction to $H_{\frac{1}{6}}(E_6)$ can be pinned down from some induced representations. Recall that $\el(20_p)\in\oh_{\frac{1}{6}}(E_6)$ has $2$-dimensional support and decomposition as found earlier.
%\begin{align*}
%&\el(20_p)=\m(20_p)-\m(30_p)-\m(24_p)-\m(60_p)+2\m(80_s)+\m(60_s)-\m(24_p')-2\m(60_p')\\
%&\qquad\qquad-\m(30_p')+\m(15_q')+2\m(20_p')-\m(6_p')-\m(1_p')\\
%\end{align*}
Inducing it up to $E_7$ then restricting back down gives:
$$\Res\circ\Ind(\el(20_p))=4\el(20_p)+\el(6_p)+\el(24_p)+\el(1_p)$$
Moreover, $\Ind(20_p)$, when its terms are put in the order of increasing $h_c$-weights, begins with the expression $\m(21_b')+\m(35_b)+\m(105_a')-\m(168_a)-2\m(105_b)....$. Since all Hom's between $\m(21_b')$, $\m(35_b')$, and $\m105_a')$ are zero, this means $\Ind(20_p)$ contains $\el(21_b')$, $\el(35_b)$, and $\el(105_a')$ each once in its composition series. %We don't yet know the full Verma decomposition for $\el(35_b)$, only the first few terms. But s
Subtracting off $\el(21_b')$ and $\el(105_a')$, $\Ind\el(20_p)-\el(21_b')-\el(105_a')=\m(35_b)+2\m(210_a)+\m(168_a)-\m(105_b)...$ where the remaining terms $\m(\tau)$ satisfy $h_c(\tau)\geq2$. Since $[\el(35_b):\m(210_a)]=0$, $2\el(210_a)$ belong to the composition series; since $[\el(35_b):\m(168_a)]=1$ but $[\el(210_a):\m(168_a)]=0$, $2\el(168_a)$ also belongs to the composition series. There are at most four irreducible representations in this block with support of dimension $3$: $\el(210_a)$, $\el(168_a)$, $\el(84_a)$, and $\el(105_c')$. A calculation shows that $\Res\el(84_a)=\Res\el(105_c')=\el(24_p)$. This implies that $\Res\el(168_a)$ and $\Res\el(210_a)$ must either (a) both equal $\el(20_p)$, or (b) one of them has restriction equal to $2\el(20_p)$ and the other to $0$.

To see that $\Res\el(168_a)=\Res\el(210_a)=\el(20_p)$, consider $\Ind^{E_7}_{E_6}\el(24_p)$ and $\Res_{E_6}^{E_7}\Ind_{E_6}^{E_7}\el(24_p)$:
\begin{align*}
\Ind\el(24_p)&=\m(168_a)-\m(315_a')+\m(105_c)-\m(210_b)+\m(84_a)+\m(105_c')+\m(84_a')\\
&+\m(70_a)-\m(210_a')+\m(168_a')-\m(15_a)-2\m(35_b')+\m(21_a')+2\m(7_a)+\m(1_a')\\
\end{align*}
Obviously $\el(168_a)$ appears in the decomposition series of this module. The columns of the decomposition matrix established so far show that $[\el(168_a):\m(84_a)]=-1$, so $2\el(84_a)$ appears in the composition series of $\Ind\el(24_p)$. Calculating that $\Res\circ\Ind\el(24_p)=4\el(24_p)+\el(20_p)$, we conclude that $\Ind\el(24_p)=\el(168_a)+\alpha\el(105_c)+2\el(84_a)+2\el(105_c')$
where $\alpha=0$ or $1$.Therefore  $\Res\el(168_a)=\el(20_p)$ and so $\Res\el(210_a)=\el(20_p)$ as well. Since $\dim\Supp\el(20_p)=2$, $\dim\Supp\el(168_a)\geq3$ and likewise $\dim\Supp\el(210_a)\geq3$. On the other hand, these dimensions of support are at most $3$ since $\el(210_a)$ and $\el(168_a)$ appear in $\Ind\el(20_p)$, a module of support of dimension $3$. Therefore $\dim\Supp\el(168_a)=\dim\Supp\el(210_a)=3$. 

All that's missing to calculate $\el(105_b)$ are two numbers: $[\el(105_b):\m(105_c)]$ and $[\el(105_b):\m(105_c')]$. All the rows of the decomposition matrix below this are complete now as well as all the other entries in the row of $\el(105_b)$. Since $\dim\Hom(\m(105_c),\m(315_a'))=1$ and $\dim\Hom(\m(105_c),\m(105_b))=0$ by (dim Hom) and the decomposition of $\el(105_c')$, $\alpha:=[\el(105_b):\m(105_c)]$ is either $0$ or $1$. In either case, a calculation shows that if $[\el(105_b):\m(105_c')]>0$ then $\dim\Supp\el(105_b)=3$, but if $[\el(105_b):\m(105_c')]=0$ then $\dim\Supp\el(105_b)=2$. But $\el(105_b)$ appears in the composition series of $\Ind\el(15_q')$ which is a module of $2$-dimensional support. Therefore $[\el(105_b):\m(105_c')]=0$ and there exist two possible decompositions for $\el(105_b)$, according to whether $\alpha=0$ or $1$.

The parabolic $A_5\times A_1$ can detect the correct decompositions of $\el(105_b),$ $\el(35_b)$, $\el(168_a)$, and $\el(210_a)$ from this point forward.
%Set $$\alpha:=[\m(105_b):\el(105_c)]=0\textrm{ or }1$$
%and write
Write
\begin{align*}
\el(105_b)&=\m(105_b)-\m(315_a')+(1-\alpha)\m(105_c)+\m(420_a)-\m(210_b)-\m(105_c')\\
&\qquad+(\alpha-1)\m(420_a')+\m(210_b')+(2-\alpha)\m(315_a)-\alpha\m(280_a)-\m(210_a')\\
&\qquad-\m(168_a')+(\alpha-2)\m(105_b')+(\alpha+1)\m(105_a)+(1-\alpha)\m(15_a)+\m(35_b')\\
&\qquad+\m(21_b)-(1+\alpha)\m(7_a)-\m(1_a')\\
\end{align*}
where $\alpha=[\m(105_b):\el(105_c)]=0\textrm{ or }1$. Consider $\Ind_{E_6}^{E_7}\el(15_q)=\el(15_a')+\el(35_b)+\alpha\el(105_c)+2\el(105_b)$. Restricting from $E_7$ to $A_5\times A_1$, one obtains $$\Res^{E_7}_{A_5\times A_1}\Ind_{E_6}^{E_7}\el(15_q)=0$$ However, $\Res^{E_7}_{A_5\times A_1}\el(105_c)\neq0$, so $\alpha=0$. This takes care of $\el(105_b)$, and immediately provides the Verma-decomposition of $\el(35_b)$ as well:
$$\el(35_b)=\Ind_{E_6}^{E_7}\el(15_q)-\el(15_a')-2\el(105_b)$$

Recall that $\Ind_{E_6}^{E_7}\el(24_p)=\el(168_a)+\alpha\el(105_c)+2\el(84_a)+2\el(105_c')$. Calculate that $\Res_{A_5\times A_1}^{E_7}\Ind_{E_6}^{E_7}\el(24_p)=0$. Again, since $\Res_{A_5\times A_1}^{E_7}\el(105_c)\neq0$ it follows that $\alpha=0$, and so $$\el(168_a)=\Ind_{E_6}^{E_7}\el(24_p)-2\el(84_a)-2\el(105_c')$$

Likewise, $\Res^{E_7}_{A_5\times A_1}\Ind_{E_6}^{E_7}\el(20_p)=0$, and since $\el(210_a)$ is a composition factor of $\Ind_{E_6}^{E_7}\el(20_p)$, this implies $\Res^{E_7}_{A_5\times A_1}\el(210_a)=0$. It follows that $[\m(210_a):\el(105_c)]=1$ as otherwise the restriction of $\el(210_a)$ to $A_5\times A_1$ will not be $0$. $\el(210_a)$ will now be determined once $[\m(210_a):\el(105_c')]$ is found. It was already argued that $\Res^{E_7}_{E_6}\el(210_a)=\el(20_p)$. This forces $[\m(210_a):\el(105_c')]=1$ by calculating the restrictions of the resulting expression for $\el(210_a)$ depending on this decomposition number.

Since Verma-decompositions of all irreducibles in the block have now been found, the inverse of the decomposition matrix is complete, and thus the decomposition matrix as well.

\end{proof}

There are also three blocks mapping to blocks of defect $1$ under the KZ functor, yielding two irreps with dimension $3$ support and one irrep with dimension $2$ support:
\begin{comment}
\begin{center}
\begin{picture}(270,50)
\put(0,30){\line(1,0){270}}
\put(0,30){\circle*{5}}
\put(30,30){\circle*{5}}
\put(120,30){\circle*{5}}
\put(150,30){\circle*{5}}
\put(240,30){\circle*{5}}
\put(270,30){\circle*{5}}
\put(-8,40){$-1$}
\put(28,40){$0$}
\put(118,40){$3$}
\put(148,40){$4$}
\put(238,40){$7$}
\put(268,40){$8$}
\put(-8,15){$56_a'$}
\put(21,15){$120_a$}
\put(110,15){$336_a'$}
\put(145,15){$336_a$}
\put(233,15){$120_a'$}
\put(267,15){$56_a$}
\end{picture}
\end{center}
\normalsize
\end{comment}

\begin{align*}
\el(27_a)&=\m(27_a)-\m(189_b')+\m(378_a')-\m(405_a')+\m(189_c)\\
\chi_{\el(27_a)}(t)&=t^{-2\frac{1}{3}}\frac{189t^4+315t^3+270t^2+108t+27}{(1-t)^3}\\
&\dim\Supp\el(27_a)=3\\
\\
\el(189_c')&=\m(189_c')-\m(405_a)+\m(378_a)-\m(189_b)+\m(27_a')\\
\chi_{\el(189_c')}(t)&=t^{1\frac{1}{3}}\frac{27t^4+108t^3+270t^2+351t+189}{(1-t)^3}\\
&\dim\Supp\el(189_c')=3\\
\\
\el(56_a')&=\m(56_a')-\m(120_a)+\m(336_a')-\m(336_a)+\m(120_a')-\m(56_a)\\
\chi_{\el(56_a')}(t)&=t^{-1}\frac{56t^4+160t^3+240t^2+160t+56}{(1-t)^2}\\
&\dim\Supp\el(56_a')=2\\
\end{align*}

\subsection{{\Large$\mbf{c=\frac{1}{5}}$}} All six nontrivial blocks at this parameter are defect $1$ blocks:
\begin{itemize}
\item$\dim\Supp\el(1_a)=3$
\item$\dim\Supp\el(56_a')=3$
\item$\dim\Supp\el(7_a')=3$
\item$\dim\Supp\el(27_a)=3$
\item$\dim\Supp\el(21_b')=3$
\item$\dim\Supp\el(21_a)=3$
\end{itemize}

\begin{comment}
\begin{align*}
\el(1_a)&=\m(1_a)-\m(84_a)+\m(216_a)-\m(189_b)+\m(56_a)\\
\chi_{\el(1_a)}(t)&=t^{-9.1}\frac{56t^{14} + 224t^{13} + 371t^{12} +\sum_{k=0}^{11}
{k+3\choose 3}t^{k}}{(1-t)^3}\\
\\
\el(56_a')&=\m(56_a')-\m(189_b')+\m(216_a')-\m(84_a')+\m(1_a')\\
\chi_{\el(56_a')}(t)&=t^{-1.9}\frac{56 + 224t + 371t^2 + \sum_{k=0}^{11}{14-k\choose 3}t^{k+3}}{(1-t)^3}\\
\\
\el(7_a')&=\m(7_a')-\m(378_a')+\m(512_a)-\m(168_a')+\m(27_a')\\
\chi_{\el(7_a')}(t)&=\\
\\
\el(27_a)&=\m(27_a)-\m(168_a)+\m(512_a')-\m(378_a)+\m(7_a)\\
\chi_{\el(27_a)}(t)&=\\
\\
\el(21_b')&=\m(21_b')-\m(189_c')+\m(336_a')-\m(189_a')+\m(21_a')\\
\chi_{\el(21_b')}(t)&=\\
\\
\el(21_a)&=\m(21_a)-\m(189_a)+\m(336_a)-\m(189_c)+\m(21_b)\\
\chi_{\el(21_a)}(t)&=\\
\end{align*}
\end{comment}

\subsection{{\Large$\mbf{c=\frac{1}{4}}$}}

There are four nontrivial blocks in $\oh_{\frac{1}{4}}(E_7)$: one pair of ``dual blocks" containing the spherical rep $\el(1a)$ and its dual $\el(1a')$ respectively, and a second pair of ``dual blocks" containing the irrep with lowest weight the standard rep $\el(7a')$ and its dual $\el(7a)$.  The blocks containing $\el(7a')$ and $\el(7a)$ may be recovered immediately from the corresponding blocks of $\mathcal{H}_q(E_7)$ by applying (dim Hom) and (RR). The blocks containing $\el(1a)$ and $\el(1a')$ require additional attention.
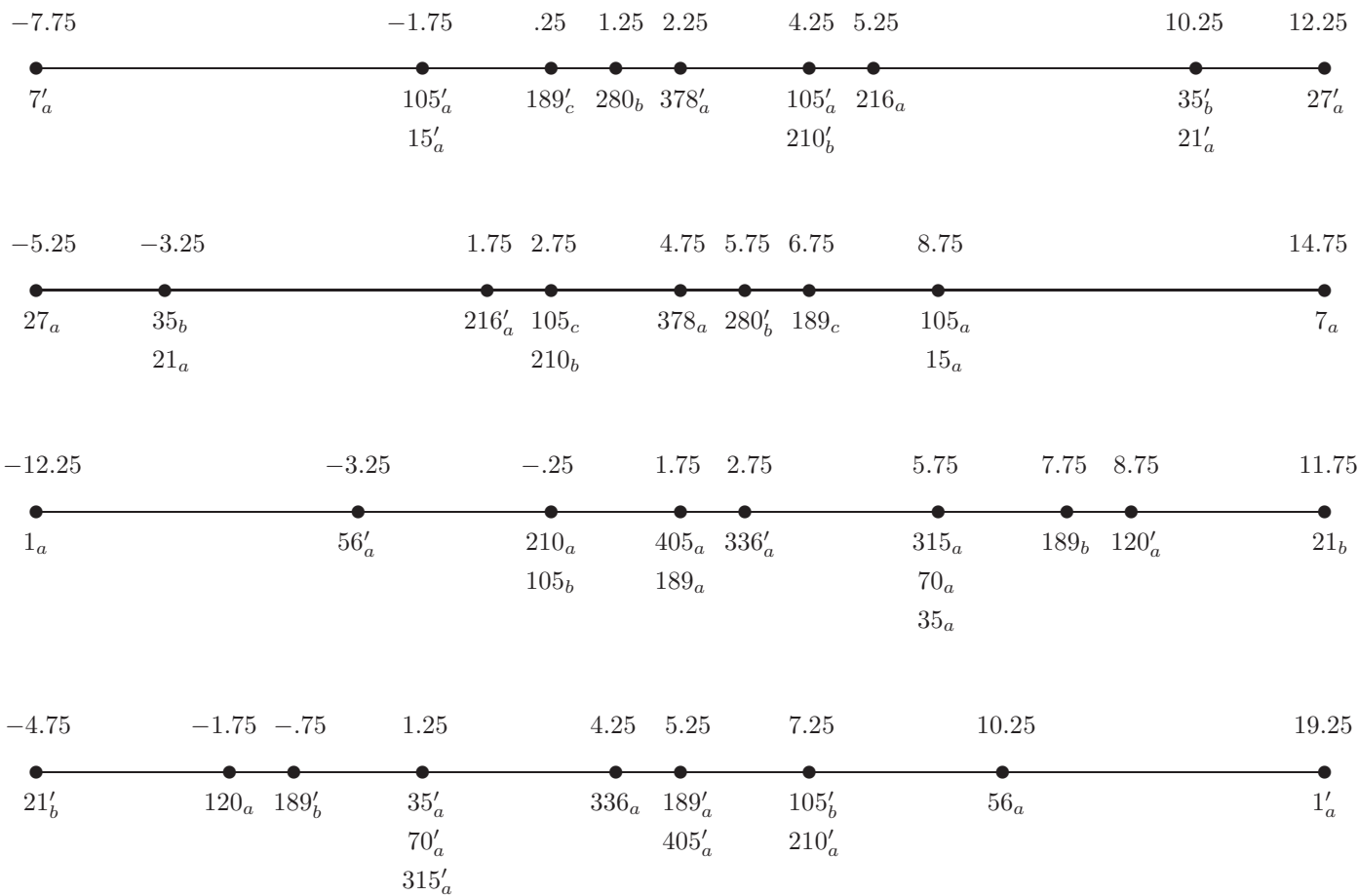
\small
\begin{center}
\begin{picture}(500,65)
\put(0,45){\line(1,0){500}}
\put(0,45){\circle*{5}}
\put(150,45){\circle*{5}}
\put(200,45){\circle*{5}}
\put(225,45){\circle*{5}}
\put(250,45){\circle*{5}}
\put(300,45){\circle*{5}}
\put(325,45){\circle*{5}}
\put(450,45){\circle*{5}}
\put(500,45){\circle*{5}}
\put(-10,60){$-7.75$}
\put(136,60){$-1.75$}
\put(193,60){$.25$}
\put(218,60){$1.25$}
\put(243,60){$2.25$}
\put(292,60){$4.25$}
\put(317,60){$5.25$}
\put(438,60){$10.25$}
\put(486,60){$12.25$}
\put(-3,30){$7_a'$}
\put(142,30){$105_a'$}
\put(144,15){$15_a'$}
\put(190,30){$189_c'$}
\put(217,30){$280_b$}
\put(242,30){$378_a'$}
\put(291,30){$105_a'$}
\put(291,15){$210_b'$}
\put(318,30){$216_a$}
\put(443,30){$35_b'$}
\put(443,15){$21_a'$}
\put(493,30){$27_a'$}
\end{picture}

\begin{picture}(500,85)
\put(0,45){\line(1,0){500}}
\put(0,45){\circle*{5}}
\put(50,45){\circle*{5}}
\put(175,45){\circle*{5}}
\put(200,45){\circle*{5}}
\put(250,45){\circle*{5}}
\put(275,45){\circle*{5}}
\put(300,45){\circle*{5}}
\put(350,45){\circle*{5}}
\put(500,45){\circle*{5}}
\put(-10,60){$-5.25$}
\put(40,60){$-3.25$}
\put(167,60){$1.75$}
\put(192,60){$2.75$}
\put(242,60){$4.75$}
\put(267,60){$5.75$}
\put(292,60){$6.75$}
\put(342,60){$8.75$}
\put(486,60){$14.75$}
\put(-5,30){$27_a$}
\put(45,30){$35_b$}
\put(45,15){$21_a$}
\put(166,30){$216_a'$}
\put(192,30){$105_c$}
\put(192,15){$210_b$}
\put(241,30){$378_a$}
\put(267,30){$280_b'$}
\put(293,30){$189_c$}
\put(343,30){$105_a$}
\put(345,15){$15_a$}
\put(496,30){$7_a$}
\end{picture}

\begin{picture}(500,120)
\put(0,80){\line(1,0){500}}
\put(0,80){\circle*{5}}
\put(125,80){\circle*{5}}
\put(200,80){\circle*{5}}
\put(250,80){\circle*{5}}
\put(275,80){\circle*{5}}
\put(350,80){\circle*{5}}
\put(400,80){\circle*{5}}
\put(425,80){\circle*{5}}
\put(500,80){\circle*{5}}
\put(-13,95){$-12.25$}
\put(112,95){$-3.25$}
\put(188,95){$-.25$}
\put(240,95){$1.75$}
\put(268,95){$2.75$}
\put(340,95){$5.75$}
\put(390,95){$7.75$}
\put(418,95){$8.75$}
\put(490,95){$11.75$}
\put(-5,65){$1_a$}
\put(117,65){$56_a'$}
\put(190,65){$210_a$}
\put(190,50){$105_b$}
\put(240,65){$405_a$}
\put(240,50){$189_a$}
\put(267,65){$336_a'$}
\put(340,65){$315_a$}
\put(342,50){$70_a$}
\put(342,35){$35_a$}
\put(390,65){$189_b$}
\put(417,65){$120_a'$}
\put(495,65){$21_b$}
\end{picture}

\begin{picture}(500,70)
\put(0,50){\line(1,0){500}}
\put(0,50){\circle*{5}}
\put(75,50){\circle*{5}}
\put(100,50){\circle*{5}}
\put(150,50){\circle*{5}}
\put(225,50){\circle*{5}}
\put(250,50){\circle*{5}}
\put(300,50){\circle*{5}}
\put(375,50){\circle*{5}}
\put(500,50){\circle*{5}}
\put(-12,65){$-4.75$}
\put(60,65){$-1.75$}
\put(92,65){$-.75$}
\put(142,65){$1.25$}
\put(215,65){$4.25$}
\put(244,65){$5.25$}
\put(292,65){$7.25$}
\put(365,65){$10.25$}
\put(488,65){$19.25$}
\put(-5,35){$21_b'$}
\put(65,35){$120_a$}
\put(92,35){$189_b'$}
\put(144,35){$35_a'$}
\put(144,20){$70_a'$}
\put(142,5){$315_a'$}
\put(215,35){$336_a$}
\put(243,35){$189_a'$}
\put(243,20){$405_a'$}
\put(292,35){$105_b'$}
\put(292,20){$210_a'$}
\put(369,35){$56_a$}
\put(495,35){$1_a'$}
\end{picture}

\end{center}
\normalsize

The decomposition matrices for the blocks containing $\el(7a')$ and $\el(7a)$ are as follows. All but the first five columns of each matrix are given by the Hecke algebra decomposition matrices, and for any $\tau$ labeling one of the missing columns, the decomposition of $\m(\tau')$ only involves $\el(\sigma)$'s which are not linked to each other but only to $\el(\tau')$, so that $[\m(\tau'):\el(\sigma)]=\dim\Hom(\m(\sigma),\m(\tau'))$ for such $\tau$. (dim Hom) and (RR) then determines the missing entries for the first five columns.

\[
\begin{blockarray}{cccccccccccccc}
&&7_a'&105_a'&15_a'&189_c'&280_b&378_a'&105_c'&210_b'&216_a&35_b'&21_a'&27_a'\\
\begin{block}{cc(cccccccccccc)}
\mathtt{(3)}&7_a'&1&\cdot&1&\cdot&1&\cdot&\cdot&1&\cdot&\cdot&\cdot&\cdot\\
\mathtt{(3)}&105_a'&\cdot&1&\cdot&1&1&1&\cdot&\cdot&\cdot&\cdot&\cdot&\cdot\\
\mathtt{(4)}&15_a'&\cdot&\cdot&1&\cdot&\cdot&\cdot&\cdot&1&\cdot&1&\cdot&\cdot\\
\mathtt{(4)}&189_c'&\cdot&\cdot&\cdot&1&\cdot&1&1&\cdot&\cdot&\cdot&\cdot&\cdot\\
\mathtt{(4)}&280_b&\cdot&\cdot&\cdot&\cdot&1&1&\cdot&1&1&\cdot&\cdot&\cdot\\
&378_a'&\cdot&\cdot&\cdot&\cdot&\cdot&1&1&\cdot&1&\cdot&1&\cdot\\
&105_c'&\cdot&\cdot&\cdot&\cdot&\cdot&\cdot&1&\cdot&\cdot&\cdot&1&\cdot\\
&210_b'&\cdot&\cdot&\cdot&\cdot&\cdot&\cdot&\cdot&1&1&1&\cdot&1\\
&216_a&\cdot&\cdot&\cdot&\cdot&\cdot&\cdot&\cdot&\cdot&1&\cdot&1&1\\
&35_b'&\cdot&\cdot&\cdot&\cdot&\cdot&\cdot&\cdot&\cdot&\cdot&1&\cdot&1\\
&21_a'&\cdot&\cdot&\cdot&\cdot&\cdot&\cdot&\cdot&\cdot&\cdot&\cdot&1&\cdot\\
&27_a'&\cdot&\cdot&\cdot&\cdot&\cdot&\cdot&\cdot&\cdot&\cdot&\cdot&\cdot&1\\
\end{block}\\
\begin{block}{cc(cccccccccccc)}
\mathtt{(3)}&7_a'&1&\cdot&-1&\cdot&-1&1&-1&1&-1&\cdot&1&\cdot\\
\mathtt{(3)}&105_a'&\cdot&1&\cdot&-1&-1&1&\cdot&1&-1&-1&\cdot&1\\
\mathtt{(4)}&15_a'&\cdot&\cdot&1&\cdot&\cdot&\cdot&\cdot&-1&1&\cdot&-1&\cdot\\
\mathtt{(4)}&189_c'&\cdot&\cdot&\cdot&1&\cdot&-1&\cdot&\cdot&1&\cdot&\cdot&-1\\
\mathtt{(4)}&280_b&\cdot&\cdot&\cdot&\cdot&1&-1&1&-1&1&1&-1&-1\\
&378_a'&\cdot&\cdot&\cdot&\cdot&\cdot&1&-1&\cdot&-1&\cdot&1&1\\
&105_c'&\cdot&\cdot&\cdot&\cdot&\cdot&\cdot&1&\cdot&\cdot&\cdot&-1&\cdot\\
&210_b'&\cdot&\cdot&\cdot&\cdot&\cdot&\cdot&\cdot&1&-1&-1&1&1\\
&216_a&\cdot&\cdot&\cdot&\cdot&\cdot&\cdot&\cdot&\cdot&1&\cdot&-1&-1\\
&35_b'&\cdot&\cdot&\cdot&\cdot&\cdot&\cdot&\cdot&\cdot&\cdot&1&\cdot&-1\\
&21_a'&\cdot&\cdot&\cdot&\cdot&\cdot&\cdot&\cdot&\cdot&\cdot&\cdot&1&\cdot\\
&27_a'&\cdot&\cdot&\cdot&\cdot&\cdot&\cdot&\cdot&\cdot&\cdot&\cdot&\cdot&1\\
\end{block}
\end{blockarray}
\]

Reading off the decompositions of simples with less than full support and calculating their characters, we find that there are three irreps with $4$-dimensional support and two irreps with $3$-dimensional support.

Likewise with the dual block:

\[
\begin{blockarray}{cccccccccccccc}
&&27_a&35_b&21_a&216_a'&105_c&210_b&378_a&280_b'&189_c&105_a&15_a&7_a\\
\begin{block}{cc(cccccccccccc)}
\mathtt{(3)}&27_a&1&1&\cdot&1&\cdot&1&\cdot&\cdot&\cdot&\cdot&\cdot&\cdot\\
\mathtt{(4)}&35_b&\cdot&1&\cdot&\cdot&\cdot&1&\cdot&\cdot&\cdot&\cdot&1&\cdot\\
\mathtt{(3)}&21_a&\cdot&\cdot&1&1&1&\cdot&1&\cdot&\cdot&\cdot&\cdot&\cdot\\
\mathtt{(4)}&216_a'&\cdot&\cdot&\cdot&1&\cdot&1&1&1&\cdot&\cdot&\cdot&\cdot\\
\mathtt{(4)}&105_c&\cdot&\cdot&\cdot&\cdot&1&\cdot&1&\cdot&1&\cdot&\cdot&\cdot\\
&210_b&\cdot&\cdot&\cdot&\cdot&\cdot&1&\cdot&1&\cdot&\cdot&1&1\\
&378_a&\cdot&\cdot&\cdot&\cdot&\cdot&\cdot&1&1&1&1&\cdot&\cdot\\
&280_b'&\cdot&\cdot&\cdot&\cdot&\cdot&\cdot&\cdot&1&\cdot&1&\cdot&1\\
&189_c&\cdot&\cdot&\cdot&\cdot&\cdot&\cdot&\cdot&\cdot&1&1&\cdot&\cdot\\
&105_a&\cdot&\cdot&\cdot&\cdot&\cdot&\cdot&\cdot&\cdot&\cdot&1&\cdot&\cdot\\
&15_a&\cdot&\cdot&\cdot&\cdot&\cdot&\cdot&\cdot&\cdot&\cdot&\cdot&1&1\\
&7_a&\cdot&\cdot&\cdot&\cdot&\cdot&\cdot&\cdot&\cdot&\cdot&\cdot&\cdot&1\\
\end{block}\\
\begin{block}{cc(cccccccccccc)}
\mathtt{(3)}&27_a&1&-1&\cdot&-1&\cdot&1&1&-1&-1&1&\cdot&\cdot\\
\mathtt{(4)}&35_b&\cdot&1&\cdot&\cdot&\cdot&-1&\cdot&1&\cdot&-1&\cdot&\cdot\\
\mathtt{(3)}&21_a&\cdot&\cdot&1&-1&-1&1&1&-1&\cdot&\cdot&-1&1\\
\mathtt{(4)}&216_a'&\cdot&\cdot&\cdot&1&\cdot&-1&-1&1&1&-1&1&-1\\
\mathtt{(4)}&105_c&\cdot&\cdot&\cdot&\cdot&1&\cdot&-1&1&\cdot&\cdot&\cdot&-1\\
&210_b&\cdot&\cdot&\cdot&\cdot&\cdot&1&\cdot&-1&\cdot&1&-1&1\\
&378_a&\cdot&\cdot&\cdot&\cdot&\cdot&\cdot&1&-1&-1&1&\cdot&1\\
&280_b'&\cdot&\cdot&\cdot&\cdot&\cdot&\cdot&\cdot&1&\cdot&-1&\cdot&-1\\
&189_c&\cdot&\cdot&\cdot&\cdot&\cdot&\cdot&\cdot&\cdot&1&-1&\cdot&\cdot\\
&105_a&\cdot&\cdot&\cdot&\cdot&\cdot&\cdot&\cdot&\cdot&\cdot&1&\cdot&\cdot\\
&15_a&\cdot&\cdot&\cdot&\cdot&\cdot&\cdot&\cdot&\cdot&\cdot&\cdot&1&-1\\
&7_a&\cdot&\cdot&\cdot&\cdot&\cdot&\cdot&\cdot&\cdot&\cdot&\cdot&\cdot&1\\
\end{block}
\end{blockarray}
\]

Calculating characters, we find there are three irreps with $4$-dimensional support and two irreps with $3$-dimensional support. Note that the characters of the irreps $\el(\tau)=\m(\tau)-...+\m(\sigma)$ of minimal support in this block are obtained from the characters of the irreps $\el(\sigma)$ of minimal support in the dual block by replacing $t^k$ with $t^{N-k}$ in the numerator, where $N$ is the degree of the polynomial in the numerator of $\chi_{\el(\sigma)}(t)$. Thus $\el(27_a)$ and $\el(105_a')$ form a dual pair, and likewise $\el(21_a)$ and $\el(7_a')$.

\begin{comment}
\begin{align*}
\el(27_a)&=\m(27_a)-\m(35_b)-\m(216_a')+\m(210_b)+\m(378_a)-\m(280_b')-\m(189_c)+\m(105_a)\\
\chi_{\el(27_a)}(t)&=t^{-5.25}\frac{105t^{10} + 420t^9 + 861t^8 + 1064t^7 + 1043t^6 + 812t^5 + 595t^4 + 400t^3 + 235t^2 + 108t + 27}{(1-t)^3}\\
&\dim\Supp\el(27_a)=3\\
\\
\el(35_b)&=\m(35_b)-\m(210_b)+\m(280_b')-\m(105_a)\\
\chi_{\el(35_b)}(t)&=t^{-3.25}\frac{105t^9 + 315t^8 + 630t^7 + 770t^6 + 735t^5 + 525t^4 + 350t^3 + 210t^2 + 105t + 35}{(1-t)^4}\\
&\dim\Supp\el(35_b)=4\\
\\
\el(21_a)&=\m(21_a)-\m(216_a')-\m(105_c)+\m(210_b)+\m(378_a)-\m(280_b')-\m(15_a)+\m(7_a)\\
\chi_{\el(21_a)}(t)&=t^{-3.25}\frac{7t^{14} + 28t^{13} + 70t^{12} + 140t^{11} + 245t^{10} + 392t^9 + 573t^8 + 780t^7 + 1005t^6 + 960t^5}{(1-t)^3}\\
&\qquad+t^{-3.25}\frac{21\sum_{k=0}^4{k+3\choose 3}t^k}{(1-t)^3}\\
&\dim\Supp\el(21_a)=3\\
\\
\el(216_a')&=\m(216_a')-\m(210_b)-\m(378_a)+\m(280_b')+\m(189_c)-\m(105_a)+\m(15_a)-\m(7_a)\\
\chi_{\el(216_a')}(t)&=t^{1.75}\frac{7t^{10} + 21t^9 + 42t^8 + 70t^7 + 105t^6 + 147t^5 + 286t^4 + 522t^3 + 666t^2 + 438t + 216}{(1-t)^4}\\
&\dim\Supp\el(216_a')=4\\
\\
\el(105_c)&=\m(105_c)-\m(378_a)+\m(280_b')-\m(7_a)\\
\chi_{\el(105_c)}(t)&=t^{2.75}\frac{7t^9 + 21t^8 + 42t^7 + 70t^6 + 105t^5 + 147t^4 + 196t^3 + 252t^2 + 315t + 105}{(1-t)^4}\\
&\dim\Supp\el(105_c)=4\\
\end{align*}
\end{comment}

Next we consider the block containing the spherical representation.

\[
\begin{blockarray}{ccccccccccccccc}
&&1_a&56_a'&210_a&105_b&405_a&189_a&336_a'&315_a&70_a&35_a&189_b&120_a'&21_b\\
\begin{block}{cc(ccccccccccccc)}
\mathtt{(3)}&1_a&1&\cdot&\cdot&1&\cdot&\cdot&\cdot&\cdot&1&\cdot&\cdot&\cdot&\cdot\\
\mathtt{(3)}&56_a'&\cdot&1&1&1&1&\cdot&\cdot&\cdot&\cdot&\cdot&\cdot&\cdot&\cdot\\
\mathtt{(4)}&210_a&\cdot&\cdot&1&\cdot&1&1&1&\cdot&\cdot&\cdot&\cdot&\cdot&\cdot\\
\mathtt{(4)}&105_b&\cdot&\cdot&\cdot&1&1&\cdot&\cdot&\cdot&1&\cdot&1&\cdot&\cdot\\
&405_a&\cdot&\cdot&\cdot&\cdot&1&\cdot&1&1&\cdot&\cdot&1&\cdot&\cdot\\
\mathtt{(3)}&189_a&\cdot&\cdot&\cdot&\cdot&\cdot&1&1&\cdot&\cdot&1&\cdot&\cdot&\cdot\\
&336_a'&\cdot&\cdot&\cdot&\cdot&\cdot&\cdot&1&1&\cdot&1&\cdot&\cdot&\cdot\\
&315_a&\cdot&\cdot&\cdot&\cdot&\cdot&\cdot&\cdot&1&\cdot&\cdot&1&1&\cdot\\
&70_a&\cdot&\cdot&\cdot&\cdot&\cdot&\cdot&\cdot&\cdot&1&\cdot&1&\cdot&1\\
&35_a&\cdot&\cdot&\cdot&\cdot&\cdot&\cdot&\cdot&\cdot&\cdot&1&\cdot&\cdot&\cdot\\
&189_b&\cdot&\cdot&\cdot&\cdot&\cdot&\cdot&\cdot&\cdot&\cdot&\cdot&1&1&1\\
&120_a'&\cdot&\cdot&\cdot&\cdot&\cdot&\cdot&\cdot&\cdot&\cdot&\cdot&\cdot&1&\cdot\\
&21_b&\cdot&\cdot&\cdot&\cdot&\cdot&\cdot&\cdot&\cdot&\cdot&\cdot&\cdot&\cdot&1\\
\end{block}\\
%\end{blockarray}
%\]
%\[
%\begin{blockarray}{ccccccccccccccc}
%&&1_a&56_a'&210_a&105_b&405_a&189_a&336_a'&315_a&70_a&35_a&189_b&120_a'&21_b\\
\begin{block}{cc(ccccccccccccc)}
\mathtt{(3)}&1_a&1&\cdot&\cdot&-1&1&\cdot&-1&\cdot&\cdot&1&\cdot&\cdot&\cdot\\
\mathtt{(3)}&56_a'&\cdot&1&-1&-1&1&1&-1&\cdot&1&\cdot&-1&1&\cdot\\
\mathtt{(4)}&210_a&\cdot&\cdot&1&\cdot&-1&-1&1&\cdot&\cdot&\cdot&1&-1&-1\\
\mathtt{(4)}&105_b&\cdot&\cdot&\cdot&1&-1&\cdot&1&\cdot&-1&-1&1&-1&\cdot\\
&405_a&\cdot&\cdot&\cdot&\cdot&1&\cdot&-1&\cdot&\cdot&1&-1&1&1\\
\mathtt{(3)}&189_a&\cdot&\cdot&\cdot&\cdot&\cdot&1&-1&1&\cdot&\cdot&-1&\cdot&1\\
&336_a'&\cdot&\cdot&\cdot&\cdot&\cdot&\cdot&1&-1&\cdot&-1&1&\cdot&-1\\
&315_a&\cdot&\cdot&\cdot&\cdot&\cdot&\cdot&\cdot&1&\cdot&\cdot&-1&\cdot&1\\
&70_a&\cdot&\cdot&\cdot&\cdot&\cdot&\cdot&\cdot&\cdot&1&\cdot&-1&1&\cdot\\
&35_a&\cdot&\cdot&\cdot&\cdot&\cdot&\cdot&\cdot&\cdot&\cdot&1&\cdot&\cdot&\cdot\\
&189_b&\cdot&\cdot&\cdot&\cdot&\cdot&\cdot&\cdot&\cdot&\cdot&\cdot&1&-1&-1\\
&120_a'&\cdot&\cdot&\cdot&\cdot&\cdot&\cdot&\cdot&\cdot&\cdot&\cdot&\cdot&1&\cdot\\
&21_b&\cdot&\cdot&\cdot&\cdot&\cdot&\cdot&\cdot&\cdot&\cdot&\cdot&\cdot&\cdot&1\\
\end{block}
\end{blockarray}
\]

\begin{comment}
\begin{align*}
\el(1_a)&=\m(1_a)-\m(105_b)+\m(405_a)-\m(336_a')+\m(35_a)\\
\chi_{\el(1_a)}(t)&=t^{-12.25}\frac{35t^{14} + 140t^{13} + 350t^{12} +\sum_{k=0}^{11}{k+3\choose3}t^k }{(1-t)^3}\\
&\dim\Supp\el(1_a)=3\\
\\
\el(56_a')&=\m(56_a')-\m(210_a)-\m(105_b)+\m(405_a)+\m(189_a)-\m(336_a')+\m(70_a)-\m(189_b)+\m(120_a')\\
\chi_{\el(56_a')}(t)&=t^{-3.25}\frac{120t^8 + 291t^7 + 444t^6 + 580t^5 + 700t^4 + 805t^3 + 560t^2 + 224t + 56}{(1-t)^3}\\
&\dim\Supp\el(56_a')=3\\
\\
\el(210_a)&=\m(210_a)-\m(405_a)-\m(189_a)+\m(336_a')+\m(189_b)-\m(120_a')-\m(21_b)\\
\chi_{\el(210_a)}(t)&=t^{-.25}\frac{21t^9 + 63t^8 + 126t^7 + 330t^6 + 486t^5 + 594t^4 + 654t^3 + 666t^2 + 630t + 210}{(1-t)^4}\\
&\dim\Supp\el(210_a)=4\\
\\
\el(105_b)&=\m(105_b)-\m(405_a)+\m(336_a')-\m(70_a)-\m(35_a)+\m(189_b)-\m(120_a')\\
\chi_{\el(105_b)}(t)&=t^{-.25}\frac{120t^6 + 171t^5 + 153t^4 + 171t^3 +225t^2 + 315t + 105}{(1-t)^4}\\
&\dim\Supp\el(105_b)=4\\
\\
\el(189_a)&=\m(189_a)-\m(336_a')+\m(315_a)-\m(189_b)+\m(21_b)\\
\chi_{\el(189_a)}(t)&=t^{1.75}\frac{21t^6 + 84t^5 + 210t^4 + 420t^3 + 546t^2 + 420t + 189}{(1-t)^3}\\
&\dim\Supp\el(189_a)=3\\
\end{align*}
\end{comment}

In this block, all columns except columns $1_a$, $56_a'$, $210_a$, $105_b$, and $189_a$ come from the Hecke algebra decomposition matrix. The first four columns then come for free by (dim Hom) and (RR). As for column $189_a$, all of its entries also come for free by (dim Hom) and (RR) except for $[\m(56_a'):\el(189_a)]$. To get this entry, refer to the decomposition matrix of the dual block just below where all entries are found from scratch, and then take $\el(56_a')=\el(120_a)^\vee$. This gives the Verma-decomposition of $\el(56_a')$ as the same as that of $\el(120_a)$ except for with the underlying lowest weights of the Vermas tensored by sign.

Finally, the block dual to that containing $\el(1a)$:

\[
\begin{blockarray}{ccccccccccccccc}
&&21_b'&120_a&189_b'&315_a'&70_a'&35_a'&336_a&405_a'&189_a'&210_a'&105_b'&56_a&1_a'\\
\begin{block}{cc(ccccccccccccc)}
\mathtt{(3)}&21_b'&1&\cdot&1&\cdot&1&\cdot&\cdot&\cdot&\cdot&\cdot&\cdot&\cdot&\cdot\\
\mathtt{(3)}&120_a&\cdot&1&1&1&\cdot&\cdot&1&\cdot&\cdot&\cdot&\cdot&\cdot&\cdot\\
\mathtt{(4)}&189_b'&\cdot&\cdot&1&1&1&\cdot&\cdot&\cdot&\cdot&\cdot&1&\cdot&\cdot\\
&315_a'&\cdot&\cdot&\cdot&1&\cdot&\cdot&1&1&\cdot&\cdot&1&\cdot&\cdot\\
&70_a'&\cdot&\cdot&\cdot&\cdot&1&\cdot&\cdot&\cdot&\cdot&\cdot&1&\cdot&1\\
\mathtt{(3)}&35_a'&\cdot&\cdot&\cdot&\cdot&\cdot&1&1&\cdot&1&\cdot&\cdot&\cdot&\cdot\\
\mathtt{(4)}&336_a&\cdot&\cdot&\cdot&\cdot&\cdot&\cdot&1&1&1&1&\cdot&\cdot&\cdot\\
&405_a'&\cdot&\cdot&\cdot&\cdot&\cdot&\cdot&\cdot&1&\cdot&1&1&1&\cdot\\
&189_a'&\cdot&\cdot&\cdot&\cdot&\cdot&\cdot&\cdot&\cdot&1&1&\cdot&\cdot&\cdot\\
&210_a'&\cdot&\cdot&\cdot&\cdot&\cdot&\cdot&\cdot&\cdot&\cdot&1&\cdot&1&\cdot\\
&105_b'&\cdot&\cdot&\cdot&\cdot&\cdot&\cdot&\cdot&\cdot&\cdot&\cdot&1&1&1\\
&56_a&\cdot&\cdot&\cdot&\cdot&\cdot&\cdot&\cdot&\cdot&\cdot&\cdot&\cdot&1&\cdot\\
&1_a'&\cdot&\cdot&\cdot&\cdot&\cdot&\cdot&\cdot&\cdot&\cdot&\cdot&\cdot&\cdot&1\\
\end{block}\\
%\end{blockarray}
%\]
%\[
%\begin{blockarray}{ccccccccccccccc}
%&&21_b'&120_a&189_b'&315_a'&70_a'&35_a'&336_a&405_a'&189_a'&210_a'&105_b'&56_a&1_a'\\
\begin{block}{cc(ccccccccccccc)}
\mathtt{(3)}&21_b'& 1  & \cdot & -1 &  1 &  \cdot & \cdot& -1 & \cdot & 1 & \cdot&  \cdot&  \cdot & \cdot\\
\mathtt{(3)}&120_a& \cdot&  1& -1&  \cdot&  1&  \cdot& -1&  1&  1& -1& -1&  1&  \cdot\\
\mathtt{(4)}&189_b'& \cdot & \cdot & 1 & -1&  -1&  \cdot&  1&  \cdot &-1&  \cdot&  1& -1 & \cdot\\
&315_a'& \cdot&  \cdot&  \cdot&  1&  \cdot&  \cdot& -1&  \cdot & 1 & \cdot &-1 & 1 & 1\\
&70_a'& \cdot & \cdot & \cdot & \cdot & 1 & \cdot & \cdot & \cdot & \cdot & \cdot &-1 & 1 & \cdot\\
\mathtt{(3)}&35_a'& \cdot & \cdot & \cdot & \cdot & \cdot & 1 &-1 & 1 & \cdot & \cdot &-1 & \cdot & 1\\
\mathtt{(4)}&336_a& \cdot & \cdot & \cdot & \cdot & \cdot & \cdot & 1 &-1& -1&  1&  1& -1& -1\\
&405_a'& \cdot & \cdot & \cdot & \cdot & \cdot & \cdot & \cdot  &1  &\cdot &-1& -1&  1&  1\\
&189_a' &\cdot&  \cdot&  \cdot&  \cdot&  \cdot & \cdot & \cdot & \cdot & 1 &-1&  \cdot&  1&  \cdot\\
&210_a'& \cdot&  \cdot & \cdot  &\cdot & \cdot&  \cdot&  \cdot & \cdot & \cdot & 1 & \cdot &-1 & \cdot\\
&105_b'& \cdot & \cdot & \cdot & \cdot & \cdot & \cdot & \cdot & \cdot & \cdot & \cdot & 1 &-1 &-1\\
&56_a& \cdot & \cdot & \cdot & \cdot&  \cdot & \cdot & \cdot & \cdot&  \cdot&  \cdot & \cdot & 1 & \cdot\\
&1_a'& \cdot & \cdot & \cdot & \cdot & \cdot & \cdot & \cdot&  \cdot&  \cdot & \cdot & \cdot & \cdot & 1\\
\end{block}
\end{blockarray}
\]

\begin{comment}

The decompositions and characters of irreps of less than full support in this block are as follows:

\begin{align*}
\el(21_b')&=\m(21_b')-\m(189_b')+\m(315_a')-\m(336_a)+\m(189_a')\\
\chi_{\el(21_b')}(t)&=t^{-4.75}\frac{189t^6 + 420t^5 + 546t^4 + 420t^3 + 210t^2 + 84t + 21}{(1-t)^3}\\
&\dim\Supp\el(21_b')=3\\
\\
\el(120_a)&=\m(120_a)-\m(189_b')+\m(70_a')-\m(336_a)+\m(405_a')+\m(189_a')-\m(210_a')-\m(105_b')+\m(56_a)\\
\chi_{\el(120_a)}(t)&=t^{-1.75}\frac{56t^8 + 224t^7 + 560t^6 + 805t^5 + 700t^4 + 580t^3 + 444t^2 + 291t + 120}{(1-t)^3}\\
&\dim\Supp\el(120_a)=3\\
\\
\el(189_b')&=\m(189_b')-\m(315_a')-\m(70_a')+\m(336_a)-\m(189_a')+\m(105_b')-\m(56_a)\\
\chi_{\el(189_b')}(t)&=t^{-.75}\frac{56t^8 + 168t^7 + 336t^6 + 455t^5 + 525t^4 + 735t^3 + 749t^2 + 567t + 189}{(1-t)^4}\\
&\dim\Supp\el(189_b')=4\\
\\
\el(35_a')&=\m(35_a')-\m(336_a)+\m(405_a')-\m(105_b')+\m(1_a')\\
\chi_{\el(35_a')}(t)&=t^{1.25}\frac{\left( \sum_{k=0}^{11}{k+3\choose3}t^{14-k}\right) + 350t^2 + 140t + 35}{(1-t)^3}\\
&\dim\Supp\el(35_a')=3\\
\\
\el(336_a)&=\m(336_a)-\m(405_a')-\m(189_a')+\m(210_a')+\m(105_b')-\m(56_a)-\m(1_a')\\
\chi_{\el(336_a)}(t)&=t^{4.25}\frac{ \left(\sum_{k=0}^8{k+2\choose 2}t^{12-k}\right)+ 111t^3 + 234t^2 +414t + 336}{(1-t)^4}\\
&\dim\Supp\el(336_a)=4\\
\end{align*}

\end{comment}

All the entries of the decomposition matrix of this block containing $\el(21_b')$ are directly determined by the corresponding decomposition matrix of the Hecke algebra (by (dim Hom) and (RR) in the case of entries in columns labeled by irreps killed by KZ functor), except for three entries: $[\m(21_b'):\el(336_a)],\;[\m(120_a):\el(336_a)],\;\textrm{and}\;[\m(189_b'):\el(336_a)]$. 

To determine $[\m(189_b'):\el(336_a)]$, use restriction to the rational Cherednik algebra of $E_6$ at $c=\frac{1}{4}$ to eliminate one of the two possible cases. Write the decomposition of $\el(189_b')$ into Vermas as:
\begin{align*}
\el(189_b')&=\m(189_b')-\m(315_a')-\m(70_a')+c\m(336_a)+(1-c)\m(405_a')-c\m(189_a')\\&\qquad+(c-1)\m(210_a')+c\m(105_b')-c\m(56_a)+(1-c)\m(1_a')\\
\end{align*}
Since $\dim\Hom(\m(336_a),\m(70_a))=0$ and $\dim\Hom(\m(336_a),\m(315_a'))=1$, $c$ can only be $0$ or $1$. Suppose $c=0$, and consider the restriction of the resulting expression $$M:=\m(189_b')-\m(315_a')-\m(70_a')+\m(405_a')-\m(210_a')+\m(1_a')$$ to $E_6$:
\begin{align*}
\Res^{E_7}_{E_6}M&=%\Res\m(189_b')-\Res\m(315_a')-\Res\m(70_a')+\Res\m(405_a')-\Res\m(210_a')+\Res\m(1_a')\\
%&=
\m(20_p)+\m(15_q)-\m(81_p)-\m(60_p)-\m(10_s)+\m(60_p')+\m(90_s)\\
&\qquad-\m(15_p')-\m(20_p')+\m(1_p')\\
&=\el(20_p)+\m(15_q)-\m(81_p)-\m(10_s)+\m(90_s)-\m(15_p')+\m(1_p')\\
\end{align*}
$\Res^{E_7}_{E_6}M$ contains $\el(20_p)$ and $\el(15_q)$ in its composition series; subtracting them leaves an expression whose lowest-$h_c$-weight term is $-\m(81_p)$. So  $\Res^{E_7}_{E_6}M$ is a virtual module but not a module, so $M$ can't be $\el(189_b')$, $c=1$, and $[\m(189_b'):\el(336_a)]=0$.

Projecting onto the block of $21_b'$ the module induced from the spherical rep from $H_{\frac{1}{4}}(E_6)$ determines the decomposition of $\el(21_b')$ into Vermas:
%\begin{align*}
%\Ind^{E_7}_{E_6}\el(1_p)&=\Ind\m(1_p)-\Ind\m(15_q)+\Ind\m(80_s)-\Ind\m(81_p')+\Ind\m(15_p')\\
%&=\m(1_a)+\m(27_a)+\m(7_a')+\m(21_b')-\m(35_b)-\m(105_b)-\m(280_b)-\m(189_b')\\&\qquad-\m(15_a')-\m(216_a')
%+\m(405_a)+\m(378_a)+\m(210_b)+\m(315_a')+\m(378_a')\\&\qquad+\m(210_b')-\m(189_c)-\m(336_a)-\m(216_a)-\m%(336_a')-\m(280_b')+\m(35_a)\\
%&\qquad+\m(105_a)+\m(21_a')+\m(189_a')\\
%\end{align*}
%The part belonging to the block of $\el(21_b')$ is:
$$\Ind^{E_7}_{E_6}|_{\mathcal{B}(21b')}=\m(21_b')-\m(189_b')+\m(315_a')-\m(336_a)+\m(189_a')$$
and as the first three terms coincide with the first three terms of $\el(21_b')$, the fourth must as well, as otherwise subtracting $\el(21_b')$ would leave an expression whose lowest $h_c$-weight Verma has a negative coefficient and this could not be the expression for a module. There can be no composition factors of $\Ind^{E_7}_{E_6}|_{\mathcal{B}(21b')}$ for terms bigger than $336_a$ in the partial order, since there are no more modules of less than full support there. So this is $\el(21_b')$ in full. It follows that $[\m(21_b'):\el(336_a)]=0$ .

The last remaining entry to determine is $[\m(120_a):\el(336_a)]$. Here it suffices to consider dimensions of supports. Set $\alpha=[\m(120_a):\el(336_a)]$ and write:
\begin{align*}
\el(120_a)&=\m(120_a)-\m(189_b')+\m(70_a')-\alpha\m(336_a)+\alpha\m(405_a')+\alpha\m(189_a')-\alpha\m(210_a')\\
&\qquad-\alpha\m(105_b')+\alpha\m(56_a)+(\alpha-1)\m(1_a')\\
\end{align*}
A calculation of the resulting characters for different possible values of $\alpha$ shows that $\dim\Supp\el(120_a)=4$ unless $\alpha=1$, in which case $\dim\Supp\el(120_a)=3$. Induction from the rational Cherednik algebra of a maximal parabolic raises the dimension of support of a module by $1$. $\el(6_p)\in\oh_{\frac{1}{4}}(E_6)$ has $2$-dimensional support; $\Ind_{E_6}^{E_7}\el(6_p)$ therefore has $3$-dimensional support, and any simple summand of $\Ind_{E_6}^{E_7}\el(6_p)$ in the Grothendieck group has dimension of support at most $3$. Writing down $\Ind_{E_6}^{E_7}\el(6_p)$ it is immediate that $\el(120_a)$ belongs to its composition series since $\m(120_a)$ appears in the expression and  $h_c(120_a)$ is minimal among the terms of  $\Ind_{E_6}^{E_7}\el(6_p)$ belonging to this block. It follows that $\dim\Supp\el(120_a)=3$ and $\alpha=1$.

\subsection{{\Large$\mbf{c=\frac{1}{3}}$}} 

\theorem The decomposition matrices of the principal block and its dual for $H_{\frac{1}{3}}(E_7)$ are as follows:

\[
\begin{blockarray}{cccccccccccccccccccccccc}
\begin{block}{cc(cccccccccccccccccccccc)}
\mathtt{(1)}&1_a&1&1&\cdot&1&\cdot&1&\cdot&\cdot&\cdot&\cdot&\cdot&\cdot&\cdot&\cdot&\cdot&\cdot&\cdot&\cdot&1&\cdot&\cdot&\cdot\\
\mathtt{(3)}&35_b&\cdot&1&\cdot&\cdot&\cdot&1&\cdot&1&1&\cdot&\cdot&1&\cdot&\cdot&\cdot&\cdot&\cdot&\cdot&1&\cdot&\cdot&\cdot\\
\mathtt{(1)}&21_a&\cdot&\cdot&1&1&1&\cdot&\cdot&\cdot&\cdot&\cdot&\cdot&\cdot&\cdot&\cdot&\cdot&1&\cdot&\cdot&\cdot&\cdot&\cdot&\cdot\\
\mathtt{(1)}&120_a&\cdot&\cdot&\cdot&1&1&1&1&1&\cdot&1&1&\cdot&1&\cdot&\cdot&1&\cdot&\cdot&1&\cdot&\cdot&\cdot\\
\mathtt{(5)}&210_a&\cdot&\cdot&\cdot&\cdot&1&\cdot&\cdot&1&1&1&\cdot&\cdot&1&1&\cdot&1&1&\cdot&\cdot&\cdot&\cdot&\cdot\\
\mathtt{(5)}&168_a&\cdot&\cdot&\cdot&\cdot&\cdot&1&\cdot&1&\cdot&1&1&1&1&\cdot&\cdot&\cdot&1&1&1&\cdot&\cdot&\cdot\\
\mathtt{(3)}&105_b&\cdot&\cdot&\cdot&\cdot&\cdot&\cdot&1&1&\cdot&\cdot&1&\cdot&1&\cdot&\cdot&\cdot&\cdot&\cdot&\cdot&\cdot&\cdot&1\\
\mathtt{(5)}&280_b&\cdot&\cdot&\cdot&\cdot&\cdot&\cdot&\cdot&1&1&\cdot&\cdot&1&1&1&\cdot&\cdot&1&1&1&\cdot&1&1\\
\mathtt{(3)}&105_c&\cdot&\cdot&\cdot&\cdot&\cdot&\cdot&\cdot&\cdot&1&\cdot&\cdot&\cdot&\cdot&1&\cdot&\cdot&\cdot&\cdot&\cdot&\cdot&1&\cdot\\
\mathtt{(5)}&420_a&\cdot&\cdot&\cdot&\cdot&\cdot&\cdot&\cdot&\cdot&\cdot&1&\cdot&\cdot&1&1&1&\cdot&1&1&\cdot&\cdot&\cdot&\cdot\\
\mathtt{(3)}&210_b&\cdot&\cdot&\cdot&\cdot&\cdot&\cdot&\cdot&\cdot&\cdot&\cdot&1&\cdot&1&\cdot&\cdot&1&\cdot&1&\cdot&1&\cdot&1\\
&84_a&\cdot&\cdot&\cdot&\cdot&\cdot&\cdot&\cdot&\cdot&\cdot&\cdot&\cdot&1&\cdot&\cdot&\cdot&\cdot&\cdot&1&1&\cdot&1&\cdot\\
&512_a'&\cdot&\cdot&\cdot&\cdot&\cdot&\cdot&\cdot&\cdot&\cdot&\cdot&\cdot&\cdot&1&1&1&1&\cdot&1&1&1&1&1\\
&336_a&\cdot&\cdot&\cdot&\cdot&\cdot&\cdot&\cdot&\cdot&\cdot&\cdot&\cdot&\cdot&\cdot&1&1&\cdot&1&1&\cdot&\cdot&1&\cdot\\
&280_a&\cdot&\cdot&\cdot&\cdot&\cdot&\cdot&\cdot&\cdot&\cdot&\cdot&\cdot&\cdot&\cdot&\cdot&1&\cdot&\cdot&1&\cdot&1&\cdot&\cdot\\
&70_a&\cdot&\cdot&\cdot&\cdot&\cdot&\cdot&\cdot&\cdot&\cdot&\cdot&\cdot&\cdot&\cdot&\cdot&\cdot&1&\cdot&\cdot&\cdot&1&\cdot&\cdot\\
\mathtt{(5)}&35_a&\cdot&\cdot&\cdot&\cdot&\cdot&\cdot&\cdot&\cdot&\cdot&\cdot&\cdot&\cdot&\cdot&\cdot&\cdot&\cdot&1&1&\cdot&\cdot&\cdot&\cdot\\
&105_a&\cdot&\cdot&\cdot&\cdot&\cdot&\cdot&\cdot&\cdot&\cdot&\cdot&\cdot&\cdot&\cdot&\cdot&\cdot&\cdot&\cdot&1&\cdot&1&1&1\\
&15_a&\cdot&\cdot&\cdot&\cdot&\cdot&\cdot&\cdot&\cdot&\cdot&\cdot&\cdot&\cdot&\cdot&\cdot&\cdot&\cdot&\cdot&\cdot&1&\cdot&1&\cdot\\
&56_a&\cdot&\cdot&\cdot&\cdot&\cdot&\cdot&\cdot&\cdot&\cdot&\cdot&\cdot&\cdot&\cdot&\cdot&\cdot&\cdot&\cdot&\cdot&\cdot&1&\cdot&1\\
&21_b&\cdot&\cdot&\cdot&\cdot&\cdot&\cdot&\cdot&\cdot&\cdot&\cdot&\cdot&\cdot&\cdot&\cdot&\cdot&\cdot&\cdot&\cdot&\cdot&\cdot&1&1\\
&7_a&\cdot&\cdot&\cdot&\cdot&\cdot&\cdot&\cdot&\cdot&\cdot&\cdot&\cdot&\cdot&\cdot&\cdot&\cdot&\cdot&\cdot&\cdot&\cdot&\cdot&\cdot&1\\
\end{block}
\end{blockarray}
\]
\\
\newgeometry{left=3cm,right=6cm,top=2cm,bottom=2cm}
\begin{turn}{270}
\begin{minipage}{\linewidth}
\small
\[
\begin{blockarray}{cccccccccccccccccccccccc}
&&1_a&35_b&21_a&120_a&210_a&168_a&105_b&280_b&105_c&420_a&210_b&84_a&512_a'&336_a&280_a&70_a&35_a&105_a&15_a&56_a&21_b&7_a\\
\begin{block}{cc(cccccccccccccccccccccc)}
\mathtt{(1)}&1_a&1&-1& \cdot&-1& 1& 1& 1&-1& 1&-1&-1& 1& 1&-1& 1& \cdot& 1&-1&-1& \cdot& 1& \cdot\\
\mathtt{(3)}&35_b&\cdot& 1& \cdot& \cdot& \cdot&-1& \cdot& \cdot&-1& 1& 1& \cdot&-1& 1&-1& \cdot&-1& 1& 1& \cdot&-1& \cdot\\
\mathtt{(1)}&21_a&\cdot& \cdot& 1&-1& \cdot& 1& 1&-1& 1& \cdot&-1& \cdot& 1&-1& \cdot& \cdot& 1& \cdot& \cdot& \cdot& \cdot& \cdot\\
\mathtt{(1)}&120_a&\cdot& \cdot& \cdot& 1&-1&-1&-1& 2&-1& 1& 1&-1&-2& 1& \cdot& 1&-2& 1& 1&-1&-1& 1\\
\mathtt{(5)}&210_a&\cdot& \cdot& \cdot& \cdot& 1& \cdot& \cdot&-1& \cdot&-1& \cdot& 1& 1& \cdot& \cdot&-2& 1&-1&-1& 2& 1&-2\\
\mathtt{(5)}&168_a&\cdot& \cdot& \cdot& \cdot& \cdot& 1& \cdot&-1& 1&-1&-1& \cdot& 2&-1& \cdot&-1& 2&-1&-2& 1& 2&-2\\
\mathtt{(3)}&105_b&\cdot& \cdot& \cdot& \cdot& \cdot& \cdot& 1&-1& 1& \cdot&-1& 1& 1&-1& \cdot& \cdot& 2&-1&-1& 1& 1&-1\\
\mathtt{(5)}&280_b&\cdot& \cdot& \cdot& \cdot& \cdot& \cdot& \cdot& 1&-1& \cdot& \cdot&-1&-1& 1& \cdot& 1&-2& 2& 1&-2&-2& 2\\
\mathtt{(3)}&105_c&\cdot& \cdot& \cdot& \cdot& \cdot& \cdot& \cdot& \cdot& 1& \cdot& \cdot& \cdot& \cdot&-1& 1& \cdot& 1&-1& \cdot& \cdot& 1& \cdot\\
\mathtt{(5)}&420_a&\cdot& \cdot& \cdot& \cdot& \cdot& \cdot& \cdot& \cdot& \cdot& 1& \cdot& \cdot&-1& \cdot& \cdot& 1&-1& 1& 1&-1&-1& 2\\
\mathtt{(3)}&210_b&\cdot& \cdot& \cdot& \cdot& \cdot& \cdot& \cdot& \cdot& \cdot& \cdot& 1& \cdot&-1& 1& \cdot& \cdot&-1& \cdot& 1& \cdot&-1& 1\\
&84_a&\cdot& \cdot& \cdot& \cdot& \cdot& \cdot& \cdot& \cdot& \cdot& \cdot& \cdot& 1& \cdot& \cdot& \cdot& \cdot& \cdot&-1&-1& 1& 1&-1\\
&512_a'&\cdot& \cdot& \cdot& \cdot& \cdot& \cdot& \cdot& \cdot& \cdot& \cdot& \cdot& \cdot& 1&-1& \cdot&-1& 1&-1&-1& 1& 2&-3\\
&336_a&\cdot& \cdot& \cdot& \cdot& \cdot& \cdot& \cdot& \cdot& \cdot& \cdot& \cdot& \cdot& \cdot& 1&-1& \cdot&-1& 1& \cdot& \cdot&-2& 1\\
&280_a&\cdot& \cdot& \cdot& \cdot& \cdot& \cdot& \cdot& \cdot& \cdot& \cdot& \cdot& \cdot& \cdot& \cdot& 1& \cdot& \cdot&-1& \cdot& \cdot& 1& \cdot\\
&70_a&\cdot& \cdot& \cdot& \cdot& \cdot& \cdot& \cdot& \cdot& \cdot& \cdot& \cdot& \cdot& \cdot& \cdot& \cdot& 1& \cdot& \cdot& \cdot&-1& \cdot& 1\\
\mathtt{(5)}&35_a&\cdot& \cdot& \cdot& \cdot& \cdot& \cdot& \cdot& \cdot& \cdot& \cdot& \cdot& \cdot& \cdot& \cdot& \cdot& \cdot& 1&-1& \cdot& 1& 1&-1\\
&105_a&\cdot& \cdot& \cdot& \cdot& \cdot& \cdot& \cdot& \cdot& \cdot& \cdot& \cdot& \cdot& \cdot& \cdot& \cdot& \cdot& \cdot& 1& \cdot&-1&-1& 1\\
&15_a&\cdot& \cdot& \cdot& \cdot& \cdot& \cdot& \cdot& \cdot& \cdot& \cdot& \cdot& \cdot& \cdot& \cdot& \cdot& \cdot& \cdot& \cdot& 1& \cdot&-1& 1\\
&56_a&\cdot& \cdot& \cdot& \cdot& \cdot& \cdot& \cdot& \cdot& \cdot& \cdot& \cdot& \cdot& \cdot& \cdot& \cdot& \cdot& \cdot& \cdot& \cdot& 1& \cdot&-1\\
&21_b&\cdot& \cdot& \cdot& \cdot& \cdot& \cdot& \cdot& \cdot& \cdot& \cdot& \cdot& \cdot& \cdot& \cdot& \cdot& \cdot& \cdot& \cdot& \cdot& \cdot& 1&-1\\
&7_a&\cdot& \cdot& \cdot& \cdot& \cdot& \cdot& \cdot& \cdot& \cdot& \cdot& \cdot& \cdot& \cdot& \cdot& \cdot& \cdot& \cdot& \cdot& \cdot& \cdot& \cdot& 1\\
\end{block}
\end{blockarray}
\]
\end{minipage}
\end{turn}

\restoregeometry
\normalsize

\[
\begin{blockarray}{cccccccccccccccccccccccc}
\begin{block}{cc(cccccccccccccccccccccc)}
\mathtt{(1)}&7_a'&1&1&1&1&\cdot&\cdot&\cdot&1&\cdot&\cdot&\cdot&\cdot&\cdot&\cdot&\cdot&\cdot&\cdot&1&\cdot&\cdot&\cdot&\cdot\\
\mathtt{(1)}&21_b'&\cdot&1&\cdot&1&1&\cdot&\cdot&\cdot&\cdot&1&1&\cdot&\cdot&1&1&\cdot&\cdot&1&\cdot&\cdot&\cdot&\cdot\\
\mathtt{(3)}&56_a'&\cdot&\cdot&1&1&\cdot&1&1&\cdot&\cdot&1&\cdot&\cdot&1&\cdot&\cdot&\cdot&\cdot&1&\cdot&\cdot&\cdot&\cdot\\
\mathtt{(5)}&105_a'&\cdot&\cdot&\cdot&1&\cdot&1&\cdot&1&1&1&1&\cdot&1&1&1&\cdot&1&1&\cdot&\cdot&\cdot&\cdot\\
\mathtt{(3)}&15_a'&\cdot&\cdot&\cdot&\cdot&1&\cdot&\cdot&\cdot&\cdot&1&\cdot&\cdot&\cdot&1&1&\cdot&\cdot&\cdot&\cdot&\cdot&\cdot&1\\
\mathtt{(5)}&280_a'&\cdot&\cdot&\cdot&\cdot&\cdot&1&\cdot&\cdot&1&1&\cdot&1&1&\cdot&\cdot&\cdot&1&\cdot&\cdot&\cdot&\cdot&\cdot\\
\mathtt{(3)}&70_a'&\cdot&\cdot&\cdot&\cdot&\cdot&\cdot&1&\cdot&\cdot&1&\cdot&\cdot&1&\cdot&\cdot&\cdot&\cdot&\cdot&\cdot&\cdot&1&\cdot\\
\mathtt{(1)}&35_a'&\cdot&\cdot&\cdot&\cdot&\cdot&\cdot&\cdot&1&1&\cdot&\cdot&\cdot&\cdot&\cdot&\cdot&\cdot&1&\cdot&\cdot&\cdot&\cdot&\cdot\\
\mathtt{(5)}&336_a'&\cdot&\cdot&\cdot&\cdot&\cdot&\cdot&\cdot&\cdot&1&1&1&1&\cdot&\cdot&1&1&1&\cdot&\cdot&\cdot&\cdot&\cdot\\
\mathtt{(5)}&512_a&\cdot&\cdot&\cdot&\cdot&\cdot&\cdot&\cdot&\cdot&\cdot&1&\cdot&1&1&\cdot&1&1&1&1&1&\cdot&1&1\\
\mathtt{(5)}&105_c'&\cdot&\cdot&\cdot&\cdot&\cdot&\cdot&\cdot&\cdot&\cdot&\cdot&1&\cdot&\cdot&\cdot&1&1&\cdot&\cdot&\cdot&\cdot&\cdot&\cdot\\
&420_a'&\cdot&\cdot&\cdot&\cdot&\cdot&\cdot&\cdot&\cdot&\cdot&\cdot&\cdot&1&\cdot&\cdot&\cdot&1&1&\cdot&1&\cdot&\cdot&\cdot\\
&210_b'&\cdot&\cdot&\cdot&\cdot&\cdot&\cdot&\cdot&\cdot&\cdot&\cdot&\cdot&\cdot&1&\cdot&\cdot&\cdot&1&1&1&\cdot&1&\cdot\\
\mathtt{(3)}&84_a'&\cdot&\cdot&\cdot&\cdot&\cdot&\cdot&\cdot&\cdot&\cdot&\cdot&\cdot&\cdot&\cdot&1&1&\cdot&1&\cdot&\cdot&1&\cdot&1\\
&280_b'&\cdot&\cdot&\cdot&\cdot&\cdot&\cdot&\cdot&\cdot&\cdot&\cdot&\cdot&\cdot&\cdot&\cdot&1&1&1&1&1&1&\cdot&1\\
&210_a'&\cdot&\cdot&\cdot&\cdot&\cdot&\cdot&\cdot&\cdot&\cdot&\cdot&\cdot&\cdot&\cdot&\cdot&\cdot&1&\cdot&\cdot&1&\cdot&1&\cdot\\
&168_a'&\cdot&\cdot&\cdot&\cdot&\cdot&\cdot&\cdot&\cdot&\cdot&\cdot&\cdot&\cdot&\cdot&\cdot&\cdot&\cdot&1&\cdot&1&1&\cdot&1\\
&105_b'&\cdot&\cdot&\cdot&\cdot&\cdot&\cdot&\cdot&\cdot&\cdot&\cdot&\cdot&\cdot&\cdot&\cdot&\cdot&\cdot&\cdot&1&1&\cdot&\cdot&\cdot\\
&120_a'&\cdot&\cdot&\cdot&\cdot&\cdot&\cdot&\cdot&\cdot&\cdot&\cdot&\cdot&\cdot&\cdot&\cdot&\cdot&\cdot&\cdot&\cdot&1&\cdot&1&1\\
&35_b'&\cdot&\cdot&\cdot&\cdot&\cdot&\cdot&\cdot&\cdot&\cdot&\cdot&\cdot&\cdot&\cdot&\cdot&\cdot&\cdot&\cdot&\cdot&\cdot&1&\cdot&1\\
&21_a'&\cdot&\cdot&\cdot&\cdot&\cdot&\cdot&\cdot&\cdot&\cdot&\cdot&\cdot&\cdot&\cdot&\cdot&\cdot&\cdot&\cdot&\cdot&\cdot&\cdot&1&\cdot\\
&1_a'&\cdot&\cdot&\cdot&\cdot&\cdot&\cdot&\cdot&\cdot&\cdot&\cdot&\cdot&\cdot&\cdot&\cdot&\cdot&\cdot&\cdot&\cdot&\cdot&\cdot&\cdot&1\\
\end{block}
\end{blockarray}
\]
\\

\newgeometry{left=3cm,right=6cm,top=2cm,bottom=2cm}
\begin{turn}{270}
\begin{minipage}{\linewidth}
\small
\[
\begin{blockarray}{cccccccccccccccccccccccc}
&&7_a'&21_b'&56_a'&105_a'&15_a'&280_a'&70_a'&35_a'&336_a'&512_a&105_c'&420_a'&210_b'&84_a'&280_b'&210_a'&168_a'&105_b'&120_a'&35_b'&21_a'&1_a'\\
\begin{block}{cc(cccccccccccccccccccccc)}
\mathtt{(1)}&7_a'&1&-1&-1&1&1&\cdot&1&-2&1&-2&-1&1&1&-1&2&-1&-1&-1&1&\cdot&\cdot&\cdot\\
\mathtt{(1)}&21_b'&\cdot&1&\cdot&-1&-1&1&\cdot&1&-1&1&1&-1&-1&1&-1&1&1&1&-1&-1&\cdot&1\\
\mathtt{(3)}&56_a'&\cdot&\cdot&1&-1&\cdot&\cdot&-1&1&\cdot&1&1&-1&\cdot&1&-2&1&1&1&-1&\cdot&\cdot&\cdot\\
\mathtt{(5)}&105_a'&\cdot&\cdot&\cdot&1&\cdot&-1&\cdot&-1&1&-1&-2&1&1&-1&2&-1&-2&-3&3&1&-2&-2\\
\mathtt{(3)}&15_a'&\cdot&\cdot&\cdot&\cdot&1&\cdot&\cdot&\cdot&\cdot&-1&\cdot&1&1&-1&1&-1&-1&-1&1&1&\cdot&-1\\
\mathtt{(5)}&280_a'&\cdot&\cdot&\cdot&\cdot&\cdot&1&\cdot&\cdot&-1&\cdot&1&\cdot&-1&\cdot&\cdot&\cdot&1&1&-1&-1&2&1\\
\mathtt{(3)}&70_a'&\cdot&\cdot&\cdot&\cdot&\cdot&\cdot&1&\cdot&\cdot&-1&\cdot&1&\cdot&\cdot&1&-1&-1&\cdot&1&\cdot&\cdot&\cdot\\
\mathtt{(1)}&35_a'&\cdot&\cdot&\cdot&\cdot&\cdot&\cdot&\cdot&1&-1&1&1&\cdot&-1&\cdot&-1&\cdot&1&1&-1&\cdot&1&\cdot\\
\mathtt{(5)}&336_a'&\cdot&\cdot&\cdot&\cdot&\cdot&\cdot&\cdot&\cdot&1&-1&-1&\cdot&1&\cdot&1&\cdot&-2&-1&2&1&-2&-1\\
\mathtt{(5)}&512_a&\cdot&\cdot&\cdot&\cdot&\cdot&\cdot&\cdot&\cdot&\cdot&1&\cdot&-1&-1&\cdot&-1&1&2&1&-2&-1&1&1\\
\mathtt{(5)}&105_c'&\cdot&\cdot&\cdot&\cdot&\cdot&\cdot&\cdot&\cdot&\cdot&\cdot&1&\cdot&\cdot&\cdot&-1&\cdot&1&1&-1&\cdot&1&1\\
&420_a'&\cdot&\cdot&\cdot&\cdot&\cdot&\cdot&\cdot&\cdot&\cdot&\cdot&\cdot&1&\cdot&\cdot&\cdot&-1&-1&\cdot&1&1&\cdot&-1\\
&210_b'&\cdot&\cdot&\cdot&\cdot&\cdot&\cdot&\cdot&\cdot&\cdot&\cdot&\cdot&\cdot&1&\cdot&\cdot&\cdot&-1&-1&1&1&-2&-1\\
\mathtt{(3)}&84_a'&\cdot&\cdot&\cdot&\cdot&\cdot&\cdot&\cdot&\cdot&\cdot&\cdot&\cdot&\cdot&\cdot&1&-1&1&\cdot&1&-1&\cdot&\cdot&1\\
&280_b'&\cdot&\cdot&\cdot&\cdot&\cdot&\cdot&\cdot&\cdot&\cdot&\cdot&\cdot&\cdot&\cdot&\cdot&1&-1&-1&-1&2&\cdot&-1&-2\\
&210_a'&\cdot&\cdot&\cdot&\cdot&\cdot&\cdot&\cdot&\cdot&\cdot&\cdot&\cdot&\cdot&\cdot&\cdot&\cdot&1&\cdot&\cdot&-1&\cdot&\cdot&1\\
&168_a'&\cdot&\cdot&\cdot&\cdot&\cdot&\cdot&\cdot&\cdot&\cdot&\cdot&\cdot&\cdot&\cdot&\cdot&\cdot&\cdot&1&\cdot&-1&-1&1&1\\
&105_b'&\cdot&\cdot&\cdot&\cdot&\cdot&\cdot&\cdot&\cdot&\cdot&\cdot&\cdot&\cdot&\cdot&\cdot&\cdot&\cdot&\cdot&1&-1&\cdot&1&1\\
&120_a'&\cdot&\cdot&\cdot&\cdot&\cdot&\cdot&\cdot&\cdot&\cdot&\cdot&\cdot&\cdot&\cdot&\cdot&\cdot&\cdot&\cdot&\cdot&1&\cdot&-1&-1\\
&35_b'&\cdot&\cdot&\cdot&\cdot&\cdot&\cdot&\cdot&\cdot&\cdot&\cdot&\cdot&\cdot&\cdot&\cdot&\cdot&\cdot&\cdot&\cdot&\cdot&1&\cdot&-1\\
&21_a'&\cdot&\cdot&\cdot&\cdot&\cdot&\cdot&\cdot&\cdot&\cdot&\cdot&\cdot&\cdot&\cdot&\cdot&\cdot&\cdot&\cdot&\cdot&\cdot&\cdot&1&\cdot\\
&1_a'&\cdot&\cdot&\cdot&\cdot&\cdot&\cdot&\cdot&\cdot&\cdot&\cdot&\cdot&\cdot&\cdot&\cdot&\cdot&\cdot&\cdot&\cdot&\cdot&\cdot&\cdot&1\\
\end{block}
\end{blockarray}
\]
\end{minipage}
\end{turn}

\restoregeometry
\normalsize

\begin{proof} The $h_c$-weight lines for these two blocks are:\\

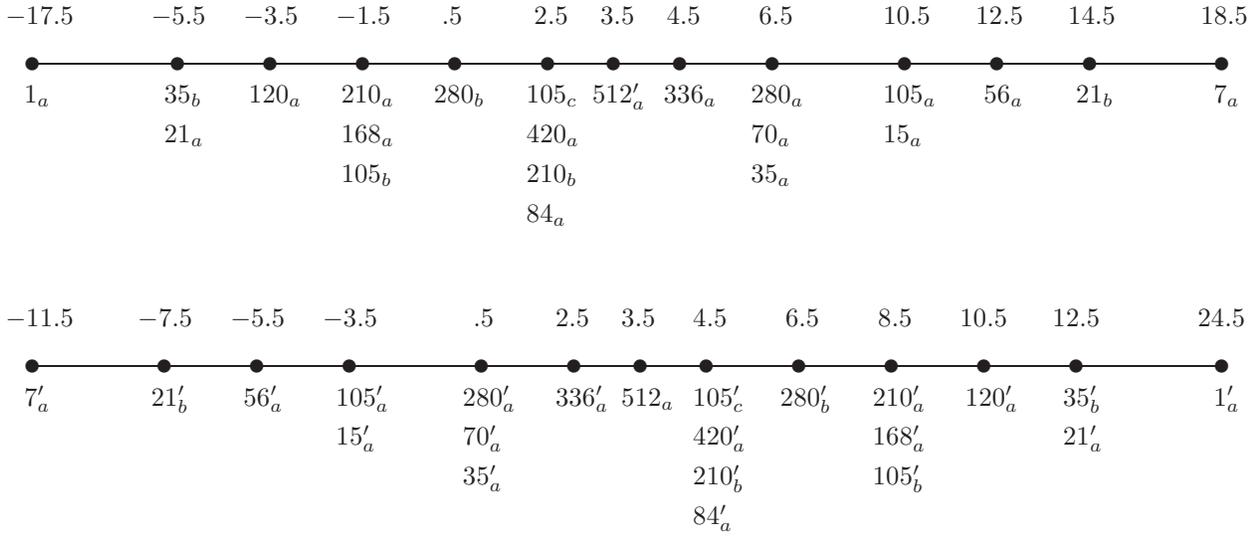
\small
\begin{center}
\begin{picture}(450,85)
\put(0,65){\line(1,0){450}}
\put(0,65){\circle*{5}}
\put(55,65){\circle*{5}}
\put(90,65){\circle*{5}}
\put(125,65){\circle*{5}}
\put(160,65){\circle*{5}}
\put(195,65){\circle*{5}}
\put(220,65){\circle*{5}}
\put(245,65){\circle*{5}}
\put(280,65){\circle*{5}}
\put(330,65){\circle*{5}}
\put(365,65){\circle*{5}}
\put(400,65){\circle*{5}}
\put(450,65){\circle*{5}}
\put(-10,80){$-17.5$}
\put(45,80){$-5.5$}
\put(80,80){$-3.5$}
\put(115,80){$-1.5$}
\put(155,80){$.5$}
\put(190,80){$2.5$}
\put(215,80){$3.5$}
\put(240,80){$4.5$}
\put(275,80){$6.5$}
\put(322,80){$10.5$}
\put(357,80){$12.5$}
\put(392,80){$14.5$}
\put(442,80){$18.5$}
\put(-3,50){$1_a$}
\put(50,50){$35_b$}
\put(50,35){$21_a$}
\put(82,50){$120_a$}
\put(117,50){$210_a$}
\put(117,35){$168_a$}
\put(117,20){$105_b$}
\put(152,50){$280_b$}
\put(187,50){$105_c$}
\put(187,35){$420_a$}
\put(187,20){$210_b$}
\put(187,5){$84_a$}
\put(212,50){$512_a'$}
\put(239,50){$336_a$}
\put(272,50){$280_a$}
\put(272,35){$70_a$}
\put(272,20){$35_a$}
\put(322,50){$105_a$}
\put(322,35){$15_a$}
\put(360,50){$56_a$}
\put(395,50){$21_b$}
\put(447,50){$7_a$}
\end{picture}

\vspace*{1cm}

\begin{picture}(450,85)
\put(0,65){\line(1,0){450}}
\put(0,65){\circle*{5}}
\put(50,65){\circle*{5}}
\put(85,65){\circle*{5}}
\put(120,65){\circle*{5}}
\put(170,65){\circle*{5}}
\put(205,65){\circle*{5}}
\put(230,65){\circle*{5}}
\put(255,65){\circle*{5}}
\put(290,65){\circle*{5}}
\put(325,65){\circle*{5}}
\put(360,65){\circle*{5}}
\put(395,65){\circle*{5}}
\put(450,65){\circle*{5}}
\put(-10,80){$-11.5$}
\put(40,80){$-7.5$}
\put(75,80){$-5.5$}
\put(110,80){$-3.5$}
\put(167,80){$.5$}
\put(198,80){$2.5$}
\put(223,80){$3.5$}
\put(250,80){$4.5$}
\put(285,80){$6.5$}
\put(320,80){$8.5$}
\put(351,80){$10.5$}
\put(386,80){$12.5$}
\put(441,80){$24.5$}
\put(-3,50){$7_a'$}
\put(45,50){$21_b'$}
\put(80,50){$56_a'$}
\put(115,50){$105_a'$}
\put(115,35){$15_a'$}
\put(163,50){$280_a'$}
\put(163,35){$70_a'$}
\put(163,20){$35_a'$}
\put(198,50){$336_a'$}
\put(223,50){$512_a$}
\put(250,50){$105_c'$}
\put(250,35){$420_a'$}
\put(250,20){$210_b'$}
\put(250,5){$84_a'$}
\put(283,50){$280_b'$}
\put(318,50){$210_a'$}
\put(318,35){$168_a'$}
\put(318,20){$105_b'$}
\put(353,50){$120_a'$}
\put(390,50){$35_b'$}
\put(390,35){$21_a'$}
\put(447,50){$1_a'$}
\end{picture}
\end{center}
\normalsize

The principal block contains $22$ irreps of which $12$ have less than full support. After the ten columns from the corresponding block of the Hecke algebra have been copied and (dim Hom), (RR) have been applied, the decomposition matrix of the principal block contains the following known entries:
\[
\begin{blockarray}{cccccccccccccccccccccccc}
\begin{block}{cc(cccccccccccccccccccccc)}
\bullet&1_a&1&1&\cdot&1&{\color{red}?}&{\color{red}?}&{\color{red}?}&{\color{red}?}&{\color{red}?}&{\color{red}?}&{\color{red}?}&\cdot&\cdot&\cdot&\cdot&\cdot&{\color{red}?}&\cdot&1&\cdot&\cdot&\cdot\\
\bullet&35_b&\cdot&1&\cdot&\cdot&\cdot&1&\cdot&{\color{red}?}&{\color{red}?}&{\color{red}?}&{\color{red}?}&1&\cdot&\cdot&\cdot&\cdot&{\color{red}?}&\cdot&1&\cdot&\cdot&\cdot\\
\bullet&21_a&\cdot&\cdot&1&1&{\color{red}?}&{\color{red}?}&{\color{red}?}&{\color{red}?}&{\color{red}?}&{\color{red}?}&{\color{red}?}&\cdot&\cdot&\cdot&\cdot&1&{\color{red}?}&\cdot&\cdot&\cdot&\cdot&\cdot\\
\bullet&120_a&\cdot&\cdot&\cdot&1&1&1&1&{\color{red}?}&{\color{red}?}&{\color{red}?}&{\color{red}?}&\cdot&1&\cdot&\cdot&1&{\color{red}?}&\cdot&1&\cdot&\cdot&\cdot\\
\bullet&210_a&\cdot&\cdot&\cdot&\cdot&1&\cdot&\cdot&1&{\color{red}?}&1&\cdot&\cdot&1&1&\cdot&1&{\color{red}?}&\cdot&\cdot&\cdot&\cdot&\cdot\\
\bullet&168_a&\cdot&\cdot&\cdot&\cdot&\cdot&1&\cdot&1&{\color{red}?}&1&1&1&1&\cdot&\cdot&\cdot&{\color{red}?}&1&1&\cdot&\cdot&\cdot\\
\bullet&105_b&\cdot&\cdot&\cdot&\cdot&\cdot&\cdot&1&1&{\color{red}?}&\cdot&1&\cdot&1&\cdot&\cdot&\cdot&{\color{red}?}&\cdot&\cdot&\cdot&\cdot&1\\
\bullet&280_b&\cdot&\cdot&\cdot&\cdot&\cdot&\cdot&\cdot&1&1&\cdot&\cdot&1&1&1&\cdot&\cdot&{\color{red}?}&1&1&\cdot&1&1\\
\bullet&105_c&\cdot&\cdot&\cdot&\cdot&\cdot&\cdot&\cdot&\cdot&1&\cdot&\cdot&\cdot&\cdot&1&\cdot&\cdot&{\color{red}?}&\cdot&\cdot&\cdot&1&\cdot\\
\bullet&420_a&\cdot&\cdot&\cdot&\cdot&\cdot&\cdot&\cdot&\cdot&\cdot&1&\cdot&\cdot&1&1&1&\cdot&{\color{red}?}&1&\cdot&\cdot&\cdot&\cdot\\
\bullet&210_b&\cdot&\cdot&\cdot&\cdot&\cdot&\cdot&\cdot&\cdot&\cdot&\cdot&1&\cdot&1&\cdot&\cdot&1&{\color{red}?}&1&\cdot&1&\cdot&1\\
&84_a&\cdot&\cdot&\cdot&\cdot&\cdot&\cdot&\cdot&\cdot&\cdot&\cdot&\cdot&1&\cdot&\cdot&\cdot&\cdot&{\color{red}?}&1&1&\cdot&1&\cdot\\
&512_a'&\cdot&\cdot&\cdot&\cdot&\cdot&\cdot&\cdot&\cdot&\cdot&\cdot&\cdot&\cdot&1&1&1&1&{\color{red}?}&1&1&1&1&1\\
&336_a&\cdot&\cdot&\cdot&\cdot&\cdot&\cdot&\cdot&\cdot&\cdot&\cdot&\cdot&\cdot&\cdot&1&1&\cdot&{\color{red}?}&1&\cdot&\cdot&1&\cdot\\
&280_a&\cdot&\cdot&\cdot&\cdot&\cdot&\cdot&\cdot&\cdot&\cdot&\cdot&\cdot&\cdot&\cdot&\cdot&1&\cdot&\cdot&1&\cdot&1&\cdot&\cdot\\
&70_a&\cdot&\cdot&\cdot&\cdot&\cdot&\cdot&\cdot&\cdot&\cdot&\cdot&\cdot&\cdot&\cdot&\cdot&\cdot&1&\cdot&\cdot&\cdot&1&\cdot&\cdot\\
\bullet&35_a&\cdot&\cdot&\cdot&\cdot&\cdot&\cdot&\cdot&\cdot&\cdot&\cdot&\cdot&\cdot&\cdot&\cdot&\cdot&\cdot&1&1&\cdot&\cdot&\cdot&\cdot\\
&105_a&\cdot&\cdot&\cdot&\cdot&\cdot&\cdot&\cdot&\cdot&\cdot&\cdot&\cdot&\cdot&\cdot&\cdot&\cdot&\cdot&\cdot&1&\cdot&1&1&1\\
&15_a&\cdot&\cdot&\cdot&\cdot&\cdot&\cdot&\cdot&\cdot&\cdot&\cdot&\cdot&\cdot&\cdot&\cdot&\cdot&\cdot&\cdot&\cdot&1&\cdot&1&\cdot\\
&56_a&\cdot&\cdot&\cdot&\cdot&\cdot&\cdot&\cdot&\cdot&\cdot&\cdot&\cdot&\cdot&\cdot&\cdot&\cdot&\cdot&\cdot&\cdot&\cdot&1&\cdot&1\\
&21_b&\cdot&\cdot&\cdot&\cdot&\cdot&\cdot&\cdot&\cdot&\cdot&\cdot&\cdot&\cdot&\cdot&\cdot&\cdot&\cdot&\cdot&\cdot&\cdot&\cdot&1&1\\
&7_a&\cdot&\cdot&\cdot&\cdot&\cdot&\cdot&\cdot&\cdot&\cdot&\cdot&\cdot&\cdot&\cdot&\cdot&\cdot&\cdot&\cdot&\cdot&\cdot&\cdot&\cdot&1\\
\end{block}
\end{blockarray}
\]

The only irrep of less than full support whose decomposition is immediately given is $\el(35_a)=\m(35_a)-\m(105_a)+\m(56_a)+\m(21_b)-\m(7_a)$; its dimension of support is $5$.

The rows with a single missing entry can be recovered from induction and restriction to $E_6$. 

$$\Ind_{E_6}^{E_7}\el(64_p')=\m(336_a)-\m(280_a)-\m(35_a)...$$ so $[\m(336_a):\el(35_a)]\geq1$. On the other hand, $\Res^{E_7}_{E_6}\el(35_a)=\el(20_s)$, and if $\alpha:=[\m(336_a):\el(35_a)]$ then $\Res^{E_7}_{E_6}\el(336_a)=\el(64_p)-(\alpha-1)\el(20_s)$, so $\alpha\leq1$, and therefore $[\m(336_a):\el(35_a)]=1$.

$$\Ind_{E_6}^{E_7}\el(60_p')=\el(512_a')-\m(336_a)+\m(70_a)+\m(35_a)...$$
and since the decomposition of $\el(512_a')$ begins $\m(512_a')-\m(336_a)-\m(70_a)$, $2\el(70_a)$ belongs to the composition series of $\Ind_{E_6}^{E_7}\el(60_p')$. Set $\alpha=[\m(512_a'):\el(35_a)]$. Then
\begin{align*}
\el(512_a)&=\Ind_{E_6}^{E_7}\el(60_p')-2\el(70_a)-\alpha\el(35_a)\\
&=\m(512_a)-\m(336_a)-\m(70_a)+\m(35_a)-\m(105_a)-\m(15_a)\\
&\qquad\qquad+\m(56_a)+2\m(21_b)-3\m(7_a)-\alpha\el(35_a)\\
\end{align*}
Restrict to $E_6$ and repeatedly subtract off $c_\tau\el(\tau)$ where $\tau$ is the leftmost nonzero place in the vector of the restricted module and $c_\tau>0$ the entry in this place. This algorithm either terminates with the $0$ vector or eventually produces a leftmost nonzero entry that is negative. In the latter case, it means the vector encodes a virtual module. This is what happens in the case of $\el(512_a')$: if $\alpha>0$ this algorithm shows that $\Res\el(512_a')$ is only a virtual module, whereas if $\alpha=0$ then the restriction is an actual module. Therefore $\alpha=0$.

Likewise for $\alpha=[\m(84_a):\el(35_a)]$, write the Verma-decomposition of $\el(84_a)$ with the variable $\alpha$ and compute the restriction to the rational Cherednik algebra of $E_6$, and apply the procedure just stated to see that $\Res\el(84_a)$ is a module if and only if $\alpha=0$. The same procedure applied to $\el(210_b)$ shows that $[\m(210_b):\el(35_a)]=0$.

For $\el(420_a)$ with $\alpha:=[\m(420_a):\el(35_a)]$, $\Res^{E_7}_{E_6}\el(420_a)$ is a module if and only if $\alpha\leq1$. On the other hand, $\Ind_{E_6}^{E_7}\el(80_s)-\el(210_b)-2\m(512_a')-\el(336_a)$ coincides with $\el(420_a)$ up through the place just to the left of $35_a$ and has a coefficient of $-1$ in the $35_a$ place, showing that $[\el(420_a):\m(35_a)]\leq-1$, from which it follows that $[\m(420_a):\el(35_a)]\geq1$. Therefore $[\m(420_a):\el(35_a)]=1$.

A restriction argument identical to those for $\el(84_a)$ and following shows that $[\m(105_c):\el(35_a)]=0$.

In the case of $\el(280_b)$, an induction argument similar to that for $\el(420_a)$ starts from $\Ind_{E_6}\el(60_s)$, whose leading term is $\m(280_b)$, and after peeling back layers of simple factors reveals that $\el(280_b)$ contains the term $c\m(35_a)$ with $c\leq-2$. This implies $[\m(280_b):\el(35_a)]\geq1$.  $\Res^{E_7}_{E_6}\el(420_a)$ is merely a virtual module if $[\m(280_b):\el(35_a)]>1$, so $[\m(280_b):\el(35_a)]=1$.

Using the entries recovered so far for column $35a$ together with those given by the Hecke algebra decomposition matrix, applying (dim Hom) and (RR) produces many of the decomposition numbers for the block dual to the principal block.

\[
\begin{blockarray}{cccccccccccccccccccccccc}
\begin{block}{cc(cccccccccccccccccccccc)}
\bullet&7_a'&1&1&1&{\color{red}?}&{\color{red}?}&{\color{red}?}&{\color{red}?}&{\color{red}?}&{\color{red}?}&{\color{red}?}&{\color{red}?}&\cdot&\cdot&{\color{red}?}&\cdot&\cdot&\cdot&1&\cdot&\cdot&\cdot&\cdot\\
\bullet&21_b'&\cdot&1&\cdot&1&1&{\color{red}?}&\cdot&{\color{red}?}&{\color{red}?}&{\color{red}?}&{\color{red}?}&\cdot&\cdot&{\color{red}?}&1&\cdot&\cdot&1&\cdot&\cdot&\cdot&\cdot\\
\bullet&56_a'&\cdot&\cdot&1&1&\cdot&{\color{red}?}&1&{\color{red}?}&{\color{red}?}&{\color{red}?}&{\color{red}?}&\cdot&1&{\color{red}?}&\cdot&\cdot&\cdot&1&\cdot&\cdot&\cdot&\cdot\\
\bullet&105_a'&\cdot&\cdot&\cdot&1&\cdot&1&\cdot&1&{\color{red}?}&{\color{red}?}&{\color{red}?}&\cdot&1&1&1&\cdot&1&1&\cdot&\cdot&\cdot&\cdot\\
\bullet&15_a'&\cdot&\cdot&\cdot&\cdot&1&\cdot&\cdot&\cdot&\cdot&1&\cdot&\cdot&\cdot&1&1&\cdot&\cdot&\cdot&\cdot&\cdot&\cdot&1\\
\bullet&280_a'&\cdot&\cdot&\cdot&\cdot&\cdot&1&\cdot&\cdot&1&{\color{red}?}&{\color{red}?}&1&1&\cdot&\cdot&\cdot&1&\cdot&\cdot&\cdot&\cdot&\cdot\\
\bullet&70_a'&\cdot&\cdot&\cdot&\cdot&\cdot&\cdot&1&\cdot&\cdot&1&\cdot&\cdot&1&\cdot&\cdot&\cdot&\cdot&\cdot&\cdot&\cdot&1&\cdot\\
\bullet&35_a'&\cdot&\cdot&\cdot&\cdot&\cdot&\cdot&\cdot&1&1&{\color{red}?}&{\color{red}?}&\cdot&\cdot&\cdot&\cdot&\cdot&1&\cdot&\cdot&\cdot&\cdot&\cdot\\
\bullet&336_a'&\cdot&\cdot&\cdot&\cdot&\cdot&\cdot&\cdot&\cdot&1&1&1&1&\cdot&\cdot&1&1&1&\cdot&\cdot&\cdot&\cdot&\cdot\\
\bullet&512_a&\cdot&\cdot&\cdot&\cdot&\cdot&\cdot&\cdot&\cdot&\cdot&1&\cdot&1&1&\cdot&1&1&1&1&1&\cdot&1&1\\
\bullet&105_c'&\cdot&\cdot&\cdot&\cdot&\cdot&\cdot&\cdot&\cdot&\cdot&\cdot&1&\cdot&\cdot&\cdot&1&1&\cdot&\cdot&\cdot&\cdot&\cdot&\cdot\\
&420_a'&\cdot&\cdot&\cdot&\cdot&\cdot&\cdot&\cdot&\cdot&\cdot&\cdot&\cdot&1&\cdot&\cdot&\cdot&1&1&\cdot&1&\cdot&\cdot&\cdot\\
&210_b'&\cdot&\cdot&\cdot&\cdot&\cdot&\cdot&\cdot&\cdot&\cdot&\cdot&\cdot&\cdot&1&\cdot&\cdot&\cdot&1&1&1&\cdot&1&\cdot\\
\bullet&84_a'&\cdot&\cdot&\cdot&\cdot&\cdot&\cdot&\cdot&\cdot&\cdot&\cdot&\cdot&\cdot&\cdot&1&1&\cdot&1&\cdot&\cdot&1&\cdot&1\\
&280_b'&\cdot&\cdot&\cdot&\cdot&\cdot&\cdot&\cdot&\cdot&\cdot&\cdot&\cdot&\cdot&\cdot&\cdot&1&1&1&1&1&1&\cdot&1\\
&210_a'&\cdot&\cdot&\cdot&\cdot&\cdot&\cdot&\cdot&\cdot&\cdot&\cdot&\cdot&\cdot&\cdot&\cdot&\cdot&1&\cdot&\cdot&1&\cdot&1&\cdot\\
&168_a'&\cdot&\cdot&\cdot&\cdot&\cdot&\cdot&\cdot&\cdot&\cdot&\cdot&\cdot&\cdot&\cdot&\cdot&\cdot&\cdot&1&\cdot&1&1&\cdot&1\\
&105_b'&\cdot&\cdot&\cdot&\cdot&\cdot&\cdot&\cdot&\cdot&\cdot&\cdot&\cdot&\cdot&\cdot&\cdot&\cdot&\cdot&\cdot&1&1&\cdot&\cdot&\cdot\\
&120_a'&\cdot&\cdot&\cdot&\cdot&\cdot&\cdot&\cdot&\cdot&\cdot&\cdot&\cdot&\cdot&\cdot&\cdot&\cdot&\cdot&\cdot&\cdot&1&\cdot&1&1\\
&35_b'&\cdot&\cdot&\cdot&\cdot&\cdot&\cdot&\cdot&\cdot&\cdot&\cdot&\cdot&\cdot&\cdot&\cdot&\cdot&\cdot&\cdot&\cdot&\cdot&1&\cdot&1\\
&21_a'&\cdot&\cdot&\cdot&\cdot&\cdot&\cdot&\cdot&\cdot&\cdot&\cdot&\cdot&\cdot&\cdot&\cdot&\cdot&\cdot&\cdot&\cdot&\cdot&\cdot&1&\cdot\\
&1_a'&\cdot&\cdot&\cdot&\cdot&\cdot&\cdot&\cdot&\cdot&\cdot&\cdot&\cdot&\cdot&\cdot&\cdot&\cdot&\cdot&\cdot&\cdot&\cdot&\cdot&\cdot&1\\
\end{block}
\end{blockarray}
\]

Playing with induction and restriction of the finite-dimensional irreps of $H_{\frac{1}{3}}(E_6)$ it turns out that $\Ind\el(1_p)=\el(21_b')$ (after projecting to the block in question here, and likewise for the next two --), $\Ind\el(6_p=\el(7_a')$, and $\Ind\el(15_p)=\el(35_a')$. Moreover, $\Ind\el(30_p)=\el(56_a')$ (after projecting the induced module to this block). This gives the Verma-decompositions of those simples, whose dimension of support is (obviously) $1$, for $\el(7_a')$, $\el(21_b')$, and $\el(35_a')$; and $3$, for $\el(56_a')$. That leaves $\el(280_a')$ and $\el(105_a')$ to find; once they are found, the decomposition matrix will be complete. 

Set $\alpha:=[\m(280_a'):\el(512_a)]$ and $\beta:=[\m(280_a'):\el(105_c')]$. Checking restriction to $E_6$ shows that $\alpha\leq2$ and $\beta\leq1$.  The induced rep $\el(64_p)$ from $E_6$ shows that $\alpha\geq1$. On the other hand, $\Res^{E_7}_{A_6}\el(280_a')$ is a virtual module if $\beta>0$ or if $\alpha>1$. Thus $\alpha=1$ and $\beta=0$.

There are three missing entries in row $105_a'$: $\alpha:=[\m(105_a'):\el(336_a')]$, $\beta:=[\m(105_a'):\el(512_a)]$, and $\gamma:=[\m(105_a'):\el(105_c')]$. Since $\el(56_a')$ is known, setting equal to $0$ the dot product of its Verma-vector with column $336_a'$ gives the equation $1+[\m(56_a'):\el(336_a')]-\alpha=0$, and since $[\m(56_a'):\el(336_a')]$ is nonnegative, $\alpha\geq1$. Identical arguments shows $\gamma\geq1$ and that $\beta=[\m(56_a'):\el(512_a)]=[\m(21_b'):\el(512_a)]$. Using the equality involving $\beta$, $\el(7_a')$ then gives that $\beta=[\m(7_a'):\el(512_a)]+1$, so $\beta\geq1$. Restricting $\Ind^{E_7}_{E_6}\el(20_p)-\el(21_b')-\el(56_a')-\el(35_a')$, which contains $\el(105_a')$ with multiplicity $2$, back down to $E_6$, one finds it contains $\el(60_p)=\Res\el(512_a)$ just once, $\el(64_p)=\Res\el(336_a')$ twice, and $\el(24_p)=\Res\el(105_c')$. This immediately implies that $\beta=\gamma\leq1$, so they are equal to $1$. Restricting the possible decomposition of $\el(105_a')$ to $A_6$ then shows that $\el(4,1^3)+\el(5,1^2)=\Res\el(336_a')$ will appear with negative multiplicity unless $\alpha\leq1$. So $\alpha=1$ also. This completes the decomposition matrix of the block dual to the principal block.

It remains to complete the decomposition matrix of the principal block. From duality, one has $\el(1_a)=\el(21_b')^\vee$, $\el(21_a)=\el(35_a)^\vee$, $\el(120_a)=\el(7_a')^\vee$. Moreover, $\el(35_b)=\Ind\el(15_q)$. So the decomposition matrix will be complete once the six entries in rows $210_a$, $168_a$, $105_b$ and columns $105_c$, $35_a$ are found.  Setting the dot product of $\el(120_a)$ and $\el(21_a)$ with columns $105_c$ and $35_a$ equal to $0$ gives some equations which imply that $[\m(210_a):\el(105_c)]$ and $[\m(210_a):\el(35_a)]$ are at least $1$. On the other hand, $\Res^{E_7}_{A_6}\el(210_a)=\el(6,1)+\el(5,1^2)$ if these decomposition numbers are exactly $1$, the sum of the two simples of $4$-dimensional support in the dual defect $1$ blocks of $H_{\frac{1}{3}}(A_6)$. And $\Res\el(105_c)$ and $\Res\el(35_a)$ land in a completely different block, the principal block of $H_{\frac{1}{3}}(A_6)$ -- so if $[\m(210_a):\el(105_c)]$ or $[\m(210_a):\el(35_a)]$ is greater than $1$ then the restriction would be a virtual module. Considering restriction to $H_{\frac{1}{3}}(E_6)$ shows that $[\m(168_a):\el(105_c)]=0$ and that $[\m(168_a):\el(35_a)]\leq1$. Looking at $\Ind_{E_6}^{E_7}\el(20_p)$, one sees that $[\m(168_a):\el(35_a)]\geq1$. As for $\el(105_b)$, its restriction to $H_{\frac{1}{3}}(E_6)$ is $\el(30_p)-[\m(105_b):\el(105_c)]\el(24_p')-[\m(105_b):\el(35_a)]\el(20_s)$, so both those decomposition numbers must be $0$. This concludes the calculation of the decomposition matrix of the principal block.

\end{proof}

There are also four blocks of defect $1$, each containing an irrep of $5$-dimensional support.

\begin{align*}
\el(27_a)&=\m(27_a)-\m(216_a)+\m(189_c)\\
&\dim\Supp\el(27_a)=5\\
\\
\el(189_c')&=\m(189_c')-\m(216_a')+\m(27_a')\\
&\dim\Supp\el(189_c')=5\\
\\
\el(189_b')&=\m(189_b')-\m(378_a')+\m(189_a')\\
&\dim\Supp\el(189_b')=5\\
\\
\el(189_a)&=\m(189_a)-\m(378_a)+\m(189_b')\\
&\dim\Supp\el(189_a)=5\\
\end{align*}

\vspace*{.2cm}

\huge
\begin{center}
\leafleft\leafleft\leafleft\\
\end{center}
\normalsize

\vspace*{.3cm}

\bibliography{decO1}
\bibliographystyle{plain}

\end{document}